\documentclass{amsart}
\usepackage{amsthm}
\usepackage{amsmath}
\usepackage{amssymb}
\usepackage[margin=0.5in]{geometry}
\usepackage{enumerate}
\usepackage{enumitem}
\usepackage{mathrsfs}
\usepackage{float}
\usepackage{graphicx}
\usepackage{mathtools}
\usepackage{hyperref}
\usepackage{caption}
\usepackage{subcaption}
\usepackage{faktor}
\usepackage[dvipsnames]{xcolor}
\usepackage{lineno}

\newtheorem{theorem}{Theorem}
\newtheorem{corollary}[theorem]{Corollary}
\newtheorem{proposition}[theorem]{Proposition}
\newtheorem{lemma}[theorem]{Lemma}
\newtheorem{definition}[theorem]{Definition}
\newtheorem{example}[theorem]{Example}
\newtheorem{remark}[theorem]{Remark}
\newtheorem{question}{Question}
\newtheorem{openquestion}{Open Question}
\numberwithin{theorem}{section}
\numberwithin{equation}{section}

\DeclareMathOperator{\sech}{sech}
\DeclareMathOperator{\Beta}{B}

\newcounter{relctr} 
\everydisplay\expandafter{\the\everydisplay\setcounter{relctr}{0}} 

\newcommand\labelrel[2]{%
  \begingroup
    \refstepcounter{relctr}%
    \stackrel{\textnormal{(\alph{relctr})}}{\mathstrut{#1}}%
    \originallabel{#2}%
  \endgroup
}
\AtBeginDocument{\let\originallabel\label} 

\makeatletter
\newcommand{\mylabel}[2]{#2\def\@currentlabel{#2}\label{#1}}
\makeatother

\allowdisplaybreaks

\begin{document}
	\title[Weighted $L^p$-Improving Estimates on constant curvature spaces]{Certain Weighted $L^p$-improving estimates for the totally-geodesic $k$-plane transform on simply connected spaces of constant curvature}
	\author[A Deshmukh]{Aniruddha Deshmukh}
	\author[A Kumar]{Ashisha Kumar}
	\date{\today}
	\begin{abstract}
		In this article we study the $L^p$-improving mapping properties of the totally-geodesic $k$-plane transform on simply connected spaces of constant curvature, namely, $\mathbb{R}^n$, $\mathbb{H}^n$ and $\mathbb{S}^n$. We begin our study by answering the question of the existence of the totally-geodesic $k$-plane transform on weighted $L^p$ spaces {with radial weights arising from the volume growth on these spaces}. These weights arise naturally from the geometry of these spaces. We then derive necessary and sufficient conditions for the $k$-plane transform of radial functions to be bounded on weighted Lebesgue spaces, with radial power weights. {Following an idea of Kurusa,} we also {derive} formulae for the totally-geodesic $k$-plane transform of general functions on the hyperbolic space and the sphere. Using this formula, and an elementary technique of Minkowski inequality, we prove weighted $L^p$-$L^p$ boundedness of the $k$-plane transform of general functions as well. Along with this, we also study the end-point behaviour of the transform, where the ``end-point" naturally arises due to either the existence conditions or the necessary conditions for boundedness.
	\end{abstract}
	\maketitle
	\section{Introduction}
		\label{IntroductionSection}
		The study of $k$-plane transform (and more specially, the Radon and $X$-Ray transforms) originated from the problem of reconstructing a function from its integral averages. For a ``nice" function $f: \mathbb{R}^n \rightarrow \mathbb{C}$, its $k$-plane transform is a function $R_kf: G \left( n, k \right) \rightarrow \mathbb{C}$, where $G \left( n, k \right)$ is the affine Grassmann manifold of all $k$-dimensional planes in $\mathbb{R}^n$, defined as
		\begin{equation}
			\label{KPlaneTransformDefinitionRn}
			R_kf \left( \xi \right) = \int\limits_{\xi} f \left( x \right) \mathrm{d}_{\xi}x.
		\end{equation}
		Here, $\mathrm{d}_{\xi}x$ is the $k$-dimensional Lebesgue measure on the $k$-plane $\xi$. {The reconstruction problem for the $k$-plane transform asks whether the function $f$ can be reconstructed solely from the values of $R_kf$ at points\footnote{``Points" in $G \left( n, k \right)$ are $k$-dimensional planes in $\mathbb{R}^n$.} in $G \left( n, k \right)$.} Much research has been done in this direction in the past century. We refer the reader to \cite{SmithSolmon} for a survey on this topic.
		
		{A question of interest in the study of such operators is about their mapping properties.} Particularly, one may ask that given a $p \geq 1$, what are the admissible values of $q \geq 1$ for which $R_k: L^p \left( \mathbb{R}^n \right) \rightarrow L^q \left( G \left( n, k \right) \right)$ is a bounded operator. {Moving a step ahead, one might also try to add some weights to the measures on $\mathbb{R}^n$ and $G \left( n, k \right)$, and ask the same question.} Among the various choices of weights one can use with the measures, the simple ones are those which depend only on the distance of the ``points" from the origin of $\mathbb{R}^n$. Such weights are referred to as radial weights. One of the early results in this direction was given by Quinto in \cite{QuintoNullSpace} about the $L^2$-boundedness with weights. Particularly, the author proved the following result for the Radon transform.
		\begin{theorem}[Quinto]
			\label{QuintoL2Boundedness}
			For $n \geq 3$, the Radon transform $R_{n - 1}: L^2 \left( \mathbb{R}^n, \| x \|^{n - 1} \mathrm{d}x \right) \rightarrow L^2 \left( G \left( n, n - 1 \right) \right)$ is a bounded operator.
		\end{theorem}
		Continuing the study of the $L^p$-$L^q$ mapping properties of the Radon transform, Oberlin and Stein in \cite{OberlinStein} gave $L^p$-improving boundedness for the Radon transform on the admissible range of $p$ (i.e., $1 \leq p < \frac{n}{n - 1}$). In fact, they prove the following mixed-norm estimate.
		\begin{theorem}[Oberlin and Stein \cite{OberlinStein}]
			\label{OSMixedNormEstimate}
			For $1 \leq p < \frac{n}{n - 1}$, $f \in L^p \left( \mathbb{R}^n \right)$, $q \leq p'$, and $\frac{1}{r} = \frac{n}{p} - n + 1$, the following mixed-norm inequality holds.
			\begin{equation}
				\label{OSMixedNormInequality}
				\| Rf \|_{L^q \left( L^r \right)} \leq C \| f \|_{L^p \left( \mathbb{R}^n \right)},
			\end{equation}
			where, $R$ is the Radon transform, $C$ is a constant independent of $f$, and for a function $\varphi$ on $\mathbb{S}^{n - 1} \times \mathbb{R}$, its mixed norm is given by
			$$\| \varphi \|_{L^q \left( L^r \right)} = \left[ \int\limits_{\mathbb{S}^{n - 1}} \left( \int\limits_{\mathbb{R}} \left| \varphi \left( \xi, t \right) \right|^r \mathrm{d}t \right)^{\frac{q}{r}} {\mathrm{d}\xi} \right]^{\frac{1}{q}}.$$
		\end{theorem}
		We notice that by taking $q = r$ in Theorem \ref{OSMixedNormEstimate}, we get an $L^p$-$L^q$ {norm estimate for} the Radon transform with the condition $p \leq \frac{n + 1}{n}$ and ${\frac{n}{p} - \frac{1}{q} = n - 1}$. Indeed, it was shown in \cite{OberlinStein} that the conditions mentioned in Theorem \ref{OSMixedNormEstimate} are also necessary. In a way, therefore, Oberlin and Stein obtained necessary and sufficient conditions for the Radon transform to be a bounded operator from $L^p \left( \mathbb{R}^n \right)$ to $L^q \left( G \left( n, n - 1 \right) \right)$.
		
		For the $k$-plane transform, analogous partial results were given by Drury (see \cite{Drury}). Particularly, in \cite{Drury}, the author proves the following result.
		\begin{theorem}[Drury \cite{Drury}]
			\label{DruryTheorem}
			Let $1 \leq p \leq \frac{n + 1}{k + 1}$, $1 \leq q \leq n + 1$ and $\frac{n}{p} - \frac{n - k}{q} = k$ with $n \leq 2k + 1$. Then, the $k$-plane transform $R_k: L^p \left( \mathbb{R}^n \right) \rightarrow L^q \left( G \left( n, k \right) \right)$ is a bounded operator.
		\end{theorem}
		We notice that Drury's technique does not give the $L^p$-$L^q$ boundedness for $k < \frac{n - 1}{2}$. In fact, this direction still seems to be open. As mentioned by {Duoandikoetxea et al. in} \cite{Duoandikoetxea}, improving the results for all values of $k$ seems to be a difficult task. Nonetheless, better understanding of boundedness can be obtained if one restricts their attention to radial functions. In \cite{Duoandikoetxea}, the authors {obtained} ``if and only if" conditions for the $L^p$-$L^q$ boundedness of the $k$-plane transform of radial functions. However, they consider the unweighted case.
		
		This line of study was further pursued by Kumar and Ray (see \cite{KumarRayWE}), where the authors {proved certain} weighted $L^p$-improving bounds of the $k$-plane transform of radial functions {on the Euclidean space}. We state their result for reference.
		\begin{theorem}[Kumar and Ray]
			\label{KumarRayWeightedTheorem}
			For $\alpha, \beta > k - n$, $1 \leq p < \frac{\alpha + n}{k}$, $k \geq 2$, and
			$$\frac{\alpha + n}{p} - \frac{\beta + n - k}{q} = k,$$
			the $k$-plane transform $R_k: L^p_{rad} \left( \mathbb{R}^n, \| x \|^{\alpha} \mathrm{d}x \right) \rightarrow L^q \left( G \left( n, k \right), d \left( 0, \xi \right)^{\beta} \mathrm{d}\xi \right)$ is a bounded operator. {For $k = 1$, the additional assumption that $\frac{1}{p} - \frac{1}{q} \leq \frac{1}{2}$, for $p > 1$, is required.}
		\end{theorem}
		It is well-known (see for instance, \cite{RubinKPlane}) that the $k$-plane transform $R_k$ is well-defined on the weighted Lebesgue space $L^p \left( \mathbb{R}^n, \| x \|^{\alpha} \mathrm{d}x \right)$ exactly when $1 \leq p < \frac{\alpha + n}{k}$. Thus, $p = \frac{\alpha + n}{k}$ is an ``end-point" for the existence of the $k$-plane transform. It is easily shown that there are radial functions in $L^{\frac{\alpha + n}{k}} \left( \mathbb{R}^n, \| x \|^{\alpha} \mathrm{d}x \right)$ for which the $k$-plane transform is infinite almost everywhere. Nonetheless, one may ask whether the $k$-plane transform is well-defined on a {smaller} subspace of $L^p \left( \mathbb{R}^n, \| x \|^{\alpha} \mathrm{d}x \right)$, when $p = \frac{\alpha + n}{k}$. {The answer to this question can be obtained from the following Lorentz norm inequality.}
		$$\| R_kf \|_{{L^{\infty} \left( G \left( n, k \right) \right)}} \leq C \| f \|_{L^{p, 1} \left( \mathbb{R}^n, \| x \|^{\alpha} \mathrm{d}x \right)}.$$
		For the unweighted case, this inequality occurred in \cite{OberlinStein} for the Radon transform. As a generalization, Duoandikoetxea et al. in \cite{Duoandikoetxea}, proved the estimate for the $k$-plane transform of radial functions. In \cite{KumarRayWE}, analogous results {were proved} by considering radial power weights on the domain. We state their result for reference, and remark that the result of Oberlin and Stein (for radial functions) and that of Duoandikoetxea et al. can be seen by considering $\alpha = 0$.
		\begin{theorem}[Kumar and Ray]
			\label{KumarRayWeightedEndPointRn}
			Let $\alpha > k - n$, $p = \frac{\alpha + n}{k}$. For $k = 1$, assume additionally that $\alpha + n \geq 2$. Then, for all radial functions $f \in L^{p, 1} \left( \mathbb{R}^n, \| x \|^{\alpha} \mathrm{d}x \right)$, we have
			$$\| R_kf \|_{\infty} \leq C \| f \|_{L^{p, 1} \left( \mathbb{R}^n, \| x \|^{\alpha} \mathrm{d}x \right)}.$$
		\end{theorem}
		Along with the study of these transforms on the Euclidean space $\mathbb{R}^n$, {many} researchers turned their gaze towards {the} study {of such transforms} on non-Euclidean spaces. We mention a few contributions to non-Euclidean spaces that have motivated us. First, we remark that among the non-Euclidean spaces, those with a Riemannian metric and a constant sectional curvature are easy to deal with. A result from Riemannian geometry, in fact, classifies all simply connected constant curvature spaces upto isometries as either the Euclidean space $\mathbb{R}^n$, the hyperbolic space $\mathbb{H}^n$, or the sphere $\mathbb{S}^n$. These spaces, respectively have curvatures $0$, $-1$ and $+1$. Helgason has been a pioneer in the study of the $k$-plane transform on these spaces. We would like to mention his works \cite{HelgasonCCS}, \cite{HelgasonSupport}, and \cite{HelgasonDuality}, where he has dealt with the inversion problem and provided support theorems for the totally-geodesic $k$-plane transform. Range characterizations of the totally{-}geodesic Radon transforms on constant curvature spaces were obtained by Berenstein et al. in \cite{BerensteinKurusaCCS}.
		
		{On the other hand,} the mapping properties of the totally-geodesic $k$-plane transform on the spaces of constant curvature, {has not} been {studied much}. A partial result, similar to that of Quinto (Theorem \ref{QuintoL2Boundedness}) was proved by Kurusa in \cite{KurusaHyperbolicSpace} for the hyperbolic space and in \cite{KurusaSphere} for the sphere. {In these papers, the author also obtained a formula for the Radon transform of functions on $\mathbb{H}^n$ and $\mathbb{S}^n$.}
		
		{In the study of constant curvature spaces with} curvature $c \in \left\lbrace -1, 0, 1 \right\rbrace$, the following function becomes important for various reasons.
		\begin{equation}
			\label{CurvatureFunction}
			s_c \left( t \right) = \begin{cases}
										t, & c = 0. \\
										\sinh t, & c = -1. \\
										\sin t, & c = +1.
									\end{cases}
		\end{equation}
		One sees its use in studying Jacobi fields on constant curvature spaces, and it plays an important role in the ``polar decomposition" of the Riemannian measure. {This also gives us the volume-growth of geodesic balls in terms of $s_c$.} We, {therefore}, call {the} function {$s_c$} as {the} ``\textit{{volume-growth} function}". Based on the above motivation, and the knowledge of the function $s_c$, we pose the following questions.
		\begin{question}
			\label{MainQuestion1}
			\normalfont
			Let $X$ be a space of constant curvature and $\Xi_k \left( X \right)$ be the collection of all totally-geodesic submanifolds in $X$. Let $s_c: \left[ 0, \infty \right) \rightarrow \mathbb{R}$ be the ``curvature function" given in Equation \eqref{CurvatureFunction}. What are the values of $\alpha_1, \alpha_2 \in \mathbb{R}$ and $p \geq 1$ for which the totally-geodesic $k$-plane transform exists on $L^p \left( X, s_c^{\alpha_1} \left( d \left( 0, x \right) \right) \left( s_c' \right)^{\alpha_2} \left( d \left( 0, x \right) \right) \mathrm{d}x \right)$?
		\end{question}
		\begin{question}
			\label{MainQuestion2}
			\normalfont
			In the setting of Question \ref{MainQuestion1}, let $p \geq 1$ and $\alpha_1, \alpha_2 \in \mathbb{R}$ be such that the totally-geodesic $k$-plane transform exists on $L^p \left( X, s_c^{\alpha_1} \left( d \left( 0, x \right) \right) \left( s_c' \right)^{\alpha_2} \left( d \left( 0, x \right) \right) \mathrm{d}x \right)$. {What are the values of $q \geq p$ and $\beta_1, \beta_2 \in \mathbb{R}$ such that the map $R_k: L^p \left( X, s_c^{\alpha_1} \left( d \left( 0, x \right) \right) \left( s_c' \right)^{\alpha_2} \left( d \left( 0, x \right) \right) \mathrm{d}x \right) \rightarrow L^q \left( {\Xi_k \left( X \right)}, s_c^{\alpha_1} \left( d \left( 0, \xi \right) \right) \left( s_c' \right)^{\alpha_2} \left( d \left( 0, \xi \right) \right) \mathrm{d}\xi \right)$ is bounded?}
		\end{question}
		The article is developed in answering Questions \ref{MainQuestion1} and \ref{MainQuestion2}. Indeed, we revisit these questions again in various sections, where the meaning becomes clear by fixing $X$ to be one of the ``model" Riemannian spaces ($\mathbb{R}^n, \mathbb{H}^n$, or $\mathbb{S}^n$).
		
		{As discussed before}, answering the {Question \ref{MainQuestion2}} posed above for general $L^p$ functions, {is difficult} even in the Euclidean case. Thus, as a first step, we try to understand the boundedness for radial functions. {In the Euclidean case, this theory is well-established (see for instance, \cite{Duoandikoetxea} and \cite{KumarRayWE}). In our work, we focus on $\mathbb{H}^n$ and $\mathbb{S}^n$.} {Here}, we obtain necessary and sufficient conditions for certain norm inequalities to hold. In this process, we observe that certain values of $p$, or {powers of} the weights, are natural ``end-points". We then prove the corresponding end-point estimates {on Lorentz spaces}.
		
		Going a step forward from radial functions, we ask the following question. The question is motivated from a similar work of Rubin (\cite{RubinKPlane}), where the author proves weighted $L^p$-$L^p$ estimates for the $k$-plane transform on the Euclidean space.
		\begin{question}
			\label{MainQuestion3}
			\normalfont
			Let $X$ be a simply connected space of constant curvature, and $s_c$ be the function given in Equation \eqref{CurvatureFunction}. What are the admissible values of $\alpha_1, \alpha_2, \beta_1, \beta_2 \in \mathbb{R}$ such that the $k$-plane transform $R_k: L^p \left( X, s_c^{\alpha_1} \left( d \left( 0, x \right) \right) \left( s_c' \right)^{\alpha_2} \left( d \left( 0, x \right) \right) \mathrm{d}x \right) \rightarrow L^p \left( \Xi_k, s_c^{\alpha_1} \left( d \left( 0, \xi \right) \right) \left( s_c' \right)^{\alpha_2} \left( d \left( 0, \xi \right) \right) \mathrm{d}\xi \right)$ is bounded?
		\end{question}
		The article is organized as follows: in Section \ref{PreliminariesSection}, we discuss the prerequisites needed for our results. {In this section, we recall a few results from Riemannian geometry}, some known results about the $k$-plane transform and its dual, Lorentz spaces and a weighted interpolation result. For the convenience of the reader, we keep the notations intact from the source when stating well-known results. However, in certain instances, for the needs of this article, we have stated some results not in a way originally stated, but in a way convenient to our use. Next, in Section \ref{ExistenceSection}, we discuss the conditions of existence of the totally-geodesic $k$-plane transform on weighted $L^p$ spaces for the spaces of constant curvature{,} $\mathbb{H}^n$ and $\mathbb{S}^n$. In Section \ref{RadialFunctionSection}, we {discuss} the $L^p$-improving bounds for the $k$-plane transform of radial functions on $\mathbb{H}^n$ and $\mathbb{S}^n$. Here, we give both necessary and sufficient conditions for boundedness. The necessary conditions are obtained by mimicking the dilations of the Euclidean space. We observe in the proofs of these results that some of the values of $p$ or the weights become ``end-points"{, where the boundedness is not obtained}. In Section \ref{EndPointSection}, we {take up the study of boundedness of the totally-geodesic $k$-plane transform at these end-points}. Particularly, we improve the results of \cite{KumarRay} concerning the end-point estimates. Next, we move to the {weighted} $L^p$-$L^p$ boundedness of the $k$-plane transform on $\mathbb{H}^n$ and $\mathbb{S}^n$ for general functions {in Section \ref{FormulaSection}}. {In fact, similar formula for the Radon transform was given by} Kurusa (see \cite{KurusaHyperbolicSpace} and \cite{KurusaSphere}) . With the formulae at hand, we {proceed} to {understand} the $L^p$-$L^p$ boundedness in Section \ref{LPLPEstimatesSection}. Lastly, in Section \ref{ConclusionSection}, we give some comments and observations about the results presented and conclude our article. At the end of this article, we provide three appendices which include material that we have often used throughout our work. In Appendix \ref{EvaluationAppendix}, we have evaluated explicitly, the $k$-plane transform and the dual of various (radial) functions. These expressions are helpful at various places in obtaining necessary conditions and/or counter-examples. Next, in Appendix \ref{HypergoemetricFunctionSubsection}, we recall certain definitions and properties of the Gauss' hypergeometric function that have been helpful in our proofs. Finally, Appendix \ref{FIBoundednessSection} deals with the weighted boundedness of Riemann-Liouville Fractional Integrals. These results are used in the study of the $X$-ray transform in Section \ref{RadialFunctionSection}. In Appendix \ref{FIBoundednessSection}, we have also given certain counterexamples about the weighted $L^1$-$L^2$ and $L^2$-$L^{\infty}$ boundedness of fractional integrals. Again, these serve us in justifying the conditions of boundedness of the $X$-ray transform.
		
		In this article, we deal with various type of inequalities. {We use the letter $C$ to denote a constant that might change at every step, but only depends on the dimension of the space ($n$), that of the submanifolds ($k$), the exponents $p, r \geq 1$, and the values of $\alpha_1, \alpha_2, \beta_1, \beta_2$ mentioned in Questions \ref{MainQuestion1}--\ref{MainQuestion3}. Also, we use the phrase ``$k$-plane transform" to mean ``totally-geodesic $k$-plane transform".}

	\section{Preliminaries}
		\label{PreliminariesSection}
		{In this section, we look at a few preliminary results required throughout this article.} For the convenience of the reader, we have divided this section into {various} subsections, each concerning a particular topic of interest to us. 
		\subsection{Differential Geometry and Spaces of Constant Curvature}
		\label{DGCCSSubsection}
		We begin with a few facts about manifolds with constant sectional curvature. For a description on sectional curvature (and other curvature tensors), we refer the reader to \cite{LeeRM}. One of the major results in this direction is the classification of all simply connected spaces of constant curvature. The result can be found in \cite{LeeRM}.
		\begin{theorem}
			\label{ClassificationCCS}
			Let $X$ be a simply connected Riemannian manifold with constant curvature $c$. Then, $X$ is, upto an isometry, one of the three spaces: $\mathbb{R}^n$, $\mathbb{H}^n_R$, or $\mathbb{S}^n_R$, in the cases $c = 0$, $c = - \frac{1}{R^2} < 0$, and $c = \frac{1}{R^2} > 0$, respectively. Here, $\mathbb{H}^n_R$ and $\mathbb{S}^n_R$, respectively, denote the hyperboloid and the sphere of radius $R$.
		\end{theorem}
		{In the course of this article, we denote by $\mathbb{H}^n$ and $\mathbb{S}^n$, the hyperbolic space and the sphere with radius $1$.} We observe that by a simple rescaling, the spaces $\mathbb{H}^n_R$ and $\mathbb{H}^n$ are diffeomorphic, and so are $\mathbb{S}^n_R$ and $\mathbb{S}^n$. Therefore, the measures on the spaces $\mathbb{H}^n_R$ and $\mathbb{S}^n_R$ are a constant multiple of those on $\mathbb{H}^n$ and $\mathbb{S}^n$, respectively. Since we deal with norm estimates, we can restrict our attention to the spaces $\mathbb{R}^n$, $\mathbb{H}^n$, and $\mathbb{S}^n$, which have curvatures $0, -1$, and $1$, respectively. 
		
		{When it comes to describing the model Riemannian spaces, the real hyperbolic space, $\mathbb{H}^n$, presents itself with various interesting results. Particularly, it} is known that there are many ways to describe the hyperbolic space $\mathbb{H}^n$, such as the upper sheet of hyperboloid, the Poincar\'{e} ball model, the Klein disc model, and the upper half space model{,} to mention a few. For a description on this topic, we refer the reader to \cite{HyperbolicGeometryNotesCannon}, where the equivalence between all these descriptions is proved. Among these {descriptions}, we use the Poincar\'{e} ball model {in our considerations of} the hyperbolic space $\mathbb{H}^n$. Our motivation for this consideration lies in Kurusa's work (\cite{KurusaHyperbolicSpace}), where the author has obtained a formula for the Radon transform on $\mathbb{H}^n$. {We wish to do the same for the $k$-plane transform in Section \ref{FormulaSection}.}
		
		In this article we use totally-geodesic submanifolds as a generalization of planes in the Euclidean space. We {now} give a formal definition of totally-geodesic submanifolds in a Riemannian manifold.
		\begin{definition}[Totally-Geodesic Submanifold]
			\label{TGSDefinition}
			Let $\left( M, g \right)$ be a Riemannian manifold. An embedded submanifold $S \subseteq M$ is totally-geodesic if the geodesics in $S$ are also geodesics in $M$.
		\end{definition}
		\begin{remark}
			\normalfont
			It is often stated that an embedded submanifold $S$ is totally-geodesic if every geodesic in $M$ that is tangent to $S$ at a point lies entirely in $S$. However, this is equivalent to Definition \ref{TGSDefinition} (see \cite{LeeRM}).
		\end{remark}
		Next, we give a comprehensive result {that characterizes} the totally-geodesic submanifolds in the spaces of constant curvature. We have written a compiled version of results available {at various instances in} the literature. For details, we refer the reader to \cite{LeeRM}, \cite{HelgasonGGA} and \cite{Berenstein}.
		\begin{theorem}
			\label{TGSClassification}
			Let $X$ be a simply connected Riemannian manifold of constant curvature.
			\begin{enumerate}
				\item If $X = \mathbb{R}^n$, the $k$-dimensional totally-geodesic submanifolds are precisely the $k$-dimensional planes in $\mathbb{R}^n$.
				\item If $X = \mathbb{H}^n$, given by the Poincar\'{e} ball model, the $k$-dimensional totally-geodesic submanifolds are precisely the $k$-dimensional {spherical caps} of {those} spheres {that are} orthogonal to the unit sphere of $\mathbb{R}^n$.
				\item If $X = \mathbb{S}^n$, then the $k$-dimensional totally-geodesic submanifolds are precisely the intersections of $\mathbb{S}^n$ with $\left( k + 1 \right)$-dimensional subspaces of $\mathbb{R}^{n + 1}$.
			\end{enumerate}
		\end{theorem}
		Using Theorem \ref{TGSClassification}, we give a parametrization of the set of all $k$-dimensional totally-geodesic {submanifolds} in the space of constant curvature $X$. {Although we believe that this result must be known to experts, we provide a quick proof for the benefit of the reader.} We denote by $\Xi_k \left( X \right)$, the collection of all $k$-dimensional totally-geodesic submanifolds of $X$. First, we consider the Euclidean space $\mathbb{R}^n$, in which case to get any $k$-dimensional plane ${\eta}$, we require a $k$-dimensional subspace ${\xi}$ and a ``foot of perpendicular" ${u}$ {from $0$ to $\eta$}. That is, we have $\Xi_k \left( \mathbb{R}^n \right) \simeq \left\lbrace \left( \xi, u \right) | \xi \in G_{n, k}, u \in \xi^{\perp} \right\rbrace$. Here, $G_{n, k}$ is the Grassmannian of $k$-dimensional subspaces of $\mathbb{R}^n$.
		
		The idea of getting a ``foot of perpendicular" is also valid when we consider  totally-geodesic submanifolds of a Riemannian manifold. Particularly, we have the following result.
		\begin{theorem}[\cite{LeeRM}]
			\label{FootOfPerpendicular}
			Let $\left( M, g \right)$ be a Riemannian manifold {and} $S_1, S_2 \subseteq M$ be totally-geodesic submanifolds. Then, there exists point $x_1 \in S_1$ and $x_2 \in S_2$ such that $d \left( x_1, x_2 \right) = d \left( S_1, S_2 \right)$. Moreover, the shortest geodesic joining $x_1$ and $x_2$ is orthogonal to both $S_1$ and $S_2$.
		\end{theorem}
		As a consequence of Theorem \ref{FootOfPerpendicular}, we have the following corollary.
		\begin{corollary}
			\label{FOPCor}
			Let $\left( M, g \right)$ be a Riemannian manifold, $S \subseteq M$ be a totally-geodesic submanifold and $p \in M$. Let $x \in S$ be a point closest to $p$, i.e., $d \left( p, x \right) = d \left( p, S \right)$. Then, the shortest geodesic joining $p$ and $x$ is orthogonal to $S$.
		\end{corollary}
		We name the point $x$ (which is the point in $S$ closest to a given point $p$) the {``}foot of perpendicular{"} to $S$ {from} $p$. We employ this {observation} in the paramterization of $\Xi_k \left( X \right)$, where $X = \mathbb{H}^n$ and $X = \mathbb{S}^n$. {To get a parametrization of $\Xi_k \left( X \right)$, we also use the following result.}
		\begin{theorem}[\cite{Gudmundsson}]
			\label{UniqueTGS}
			Let $\left( M, g \right)$ be a Riemannian Manifold, $p \in M$ and $V \subseteq T_pM$ be a subspace. Then, locally, there is at most one totally-geodesic submanifold $S$ through $p$ such that $T_pS = V$.
		\end{theorem}
		We now give the parametrization of the collection $\Xi_k \left( X \right)$, of all $k$-dimensional totally-geodesic submanifolds in $X$.
		\begin{proposition}
			\label{TGSParametrization}
			Let $X$ be a simply connected Riemannian manifold of constant curvature. Then, we have $\Xi_k \left( X \right) \simeq \left\lbrace \left( \xi, u \right) | \xi \in G_{n, k}, u \in \xi^{\perp} \cap B \left( 0, \frac{\delta \left( X \right)}{2} \right) \right\rbrace$, where $G_{n, k}$ is the Grassmannian of $k$-dimensional subspaces of $\mathbb{R}^n$, and $\delta \left( X \right)$ is the diameter of $X$. 
		\end{proposition}
		\begin{proof}
			When $X = \mathbb{R}^n$, the parametrization is clear, as discussed before Theorem \ref{FootOfPerpendicular}.
			
			We start with the case of $\mathbb{H}^n$. We observe that given a $k$-dimensional totally-geodesic submanifold {$\eta$}, we get a point $u {\in \eta}$ that is closest to the origin, and the tangent space $\xi$ at this point is $k$-dimensional and perpendicular to the line joining $0$ and $u$. Conversely, given a pair $\left( \xi, u \right)$ where $\xi \in G_{n, k}$ and $u \in \xi^{\perp}$, consider the geodesic through $0$ in the direction $\frac{u}{\| u \|}$ (whenever $u \neq 0$), and let the vector $u \in \xi^{\perp}$ be identified with the point on this geodesic at a {hyperbolic} distance $\| u \|$ from $0$. Let $\varphi: \mathbb{H}^n \rightarrow \mathbb{H}^n$ be an isometry that takes $0$ to $u$ {without any} rotational component. With Theorem \ref{UniqueTGS} at hand, we know that there is (locally) a unique totally-geodesic submanifold through $u$ with the tangent space at $u$ given by $d\varphi_0 \left( \xi \right)$. However, $\xi \cap \mathbb{H}^n$ is a totally-geodesic subamnifold through $0$ with $T_0 \left( \xi \cap \mathbb{H}^n \right) = \xi$. Hence, $\varphi \left( {\xi \cap \mathbb{H}^n} \right)$ is a totally-geodesic submanifold through $\varphi \left( 0 \right) = u$ with $T_u \left( \varphi \left( {\xi \cap \mathbb{H}^n} \right) \right) = d\varphi_0 \left( \xi \right)$. Thus, for every pair $\left( \xi, u \right)$ with $\xi \in G_{n, k}$ and $u \in \xi^{\perp}$, we have a totally-geodesic submanifold, ${\eta}$, through $u$ {with $T_u\eta \simeq \xi$}. {In the sequel, we identify $d\varphi_0 \left( \xi \right)$ by $\xi$, and make no difference in writing the two.} Also, we observe that $d \left( 0, \varphi \left( {\xi \cap \mathbb{H}^n} \right) \right) = d \left( 0, u \right)$. This gives us the desired parametrization (see Figure \ref{TGSParameterizationHn}).
		
			Similarly, we have the paramterization of $\Xi_k \left( \mathbb{S}^n \right) \simeq \left\lbrace \left( \xi, u \right) | \xi \in G_{n, k}, u \in \xi^{\perp} \cap B \left( 0, \frac{\pi}{2} \right) \right\rbrace$, where $B \left( 0, \frac{\pi}{2} \right)$ is the ball of radius $\frac{\pi}{2}$ centered at $0$. This ball suffices because if we move a distance bigger than $\frac{\pi}{2}$ in the direction $u$, the totally-geodesic submanifold obtained would have the closest point to $e_{n{+1}}$ in the direction $-u$, thereby making repetitions. We also exclude the equatorial totally-geodesic submanifolds (i.e., those submanifolds at a distance equal to $\frac{\pi}{2}$ from $e_{n{+1}}$) because for {such totally-geodesic submanifolds}, we do not have a unique point closest to $e_{n{+1}}$ (see Figure \ref{TGSSphere}).
			\begin{figure}[ht!]
				\centering
				\begin{subfigure}[t]{0.45\textwidth}
					\centering
					\includegraphics[scale=0.525]{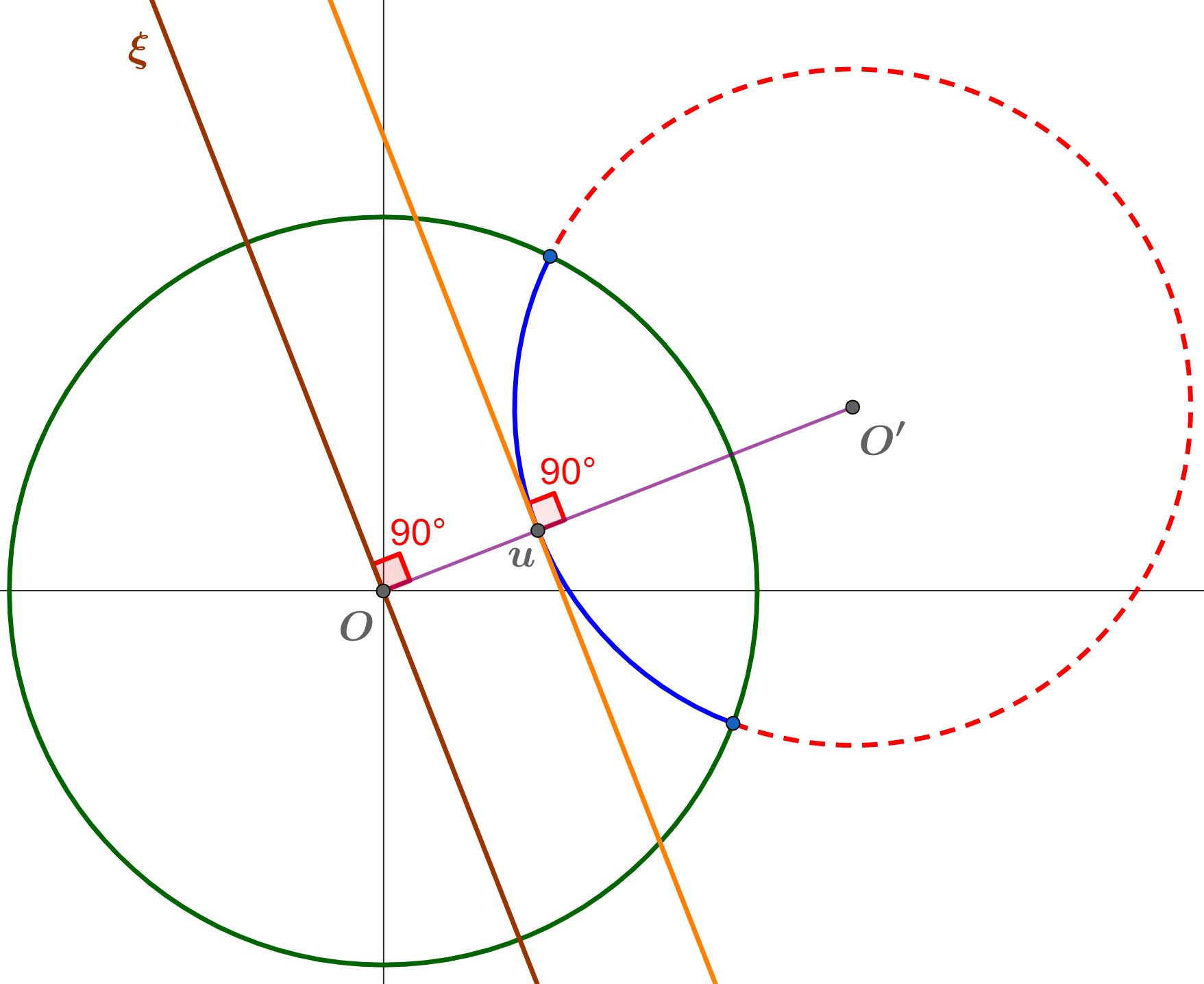}
					\caption{Totally-geodesic subamnifolds of $\mathbb{H}^n$.}
					\label{TGSParameterizationHn}
				\end{subfigure}
				~~~~~~~~
				\begin{subfigure}[t]{0.45\textwidth}
					\centering
					\includegraphics[scale=0.35]{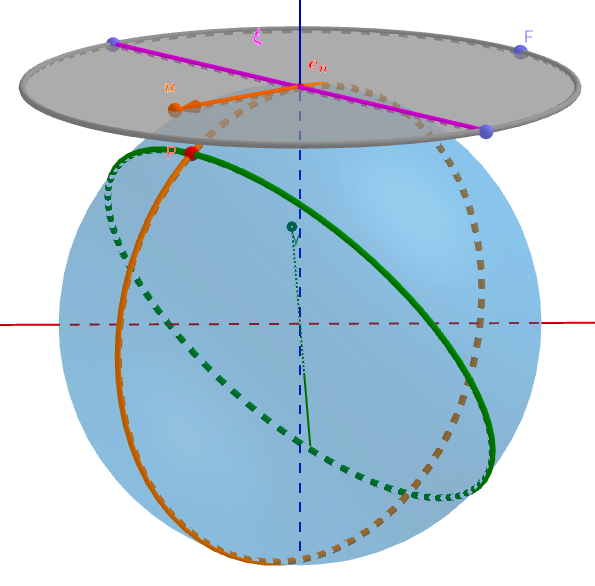}
					\caption{Totally-geodesic submanifolds on the sphere.}
					\label{TGSSphere}
				\end{subfigure}
				\caption{Parameterizing the totally-geodesic submanifolds on spaces of constant curvature}
				\label{RightTriangleFigure}
			\end{figure}
		\end{proof}
		In this article, we wish to look at integrals of functions defined on $\Xi_k \left( X \right)$. {For the same, we require a measure on $\Xi_k \left( X \right)$.} We postpone the discussion of the measure to Subsection \ref{MeasureSubsection}, and instead proceed to give a few more results related to differential geometry.
		
		An important technique that we use throughout the article is the polar decomposition of geodesic balls on Riemannian manifolds. The motivation lies in the fact that any element of $\mathbb{R}^n \setminus \left\lbrace 0 \right\rbrace$ can be described uniquely by an element of $\mathbb{S}^{n - 1} \times \left( 0, \infty \right)$. The same can be said about geodesic balls (diffeomorphic images of the Euclidean balls {in} the tangent space through the exponential map). However, since the exponential map {deforms} the Euclidean structure of the tangent space, the measure changes! On spaces of constant curvature, this change is given by what {the} ``{volume-growth} function", $s_c$ (see Equation \eqref{CurvatureFunction}).
		\begin{theorem}[\cite{LeeRM}]
			\label{PolarDecomposition}
			Let $\left( M, g \right)$ be an $n$-dimensional Riemannian manifold of constant sectional curvature $c$. Let $U \subseteq M$ be a geodesic ball (open or closed) around a point $p \in M$ of radius $b$. If $f: U \rightarrow \mathbb{R}$ is any integrable function, then, we have
			\begin{equation}
				\label{PolarDecompositionEquation}
				\int\limits_{U} f \ \mathrm{d}V_g = \int\limits_{\mathbb{S}^{n - 1}} \int\limits_{0}^{b} \left( f \circ \varphi \right) \left( r, \omega \right) \left( s_c \left( r \right) \right)^{n - 1} \ \mathrm{d}r {\mathrm{d}\omega}.
			\end{equation}
			Here, $\mathrm{d}V_g$ is the volume form induced by the Riemannian metric $g$, and $\varphi: \left( 0, b \right) \times \mathbb{S}^{n - 1} \rightarrow U \setminus \left\lbrace p \right\rbrace$ is the normal coordinate chart map given by
			\begin{equation}
				\label{NormalCoordinates}
				\varphi \left( r, \omega \right) = \exp_p \left( r \omega \right).
			\end{equation}
		\end{theorem}
		\begin{remark}
			\normalfont
			It is well-known that the hyperbolic space and the punctured sphere can be {parametrized} by geodesic normal coordinates. Particularly, the exponential map of $\mathbb{H}^n$ becomes a diffeomorphism on $T_0 \mathbb{H}^n$, and that of $\mathbb{S}^n$ becomes a diffeomorphism on $B \left( 0, \pi \right) \subseteq T_{e_{n + 1}} \mathbb{S}^n$. Therefore, we can think of polar decomposition on the entirety of $\mathbb{H}^n$ and $\mathbb{S}^n$. {For this, one would require an ``origin" of the spaces, which we denote by $0$.} Indeed, in the Poincar\'{e} ball model of $\mathbb{H}^n$, this is the origin of $\mathbb{R}^n$. However, in the case of the sphere, the symbol `$0$' means $e_{n + 1} \in \mathbb{S}^n$.
		\end{remark}
		With the understanding of polar decomposition on $\mathbb{H}^n$ and $\mathbb{S}^n$, one can think of looking at functions that depend only on the distance from the origin. These functions, when looked at in $\mathbb{R}^n$, are called \textbf{radial} functions. {We keep the terminology same, and define the following.} Some authors prefer to call such functions ``\textit{zonal}" (see for instance, \cite{Rubin}).
		\begin{definition}[Radial Functions on $X$]~
			\label{RadialFunctions}
			\begin{enumerate}
				\item A function $f: \mathbb{R}^n \rightarrow \mathbb{C}$ if radial if there is a function $\tilde{f}: \left[ 0, \infty \right) \rightarrow \mathbb{C}$ such that $f \left( x \right) = \tilde{f} \left( d \left( 0, x \right) \right)$.
				\item A function $f: \mathbb{H}^n \rightarrow \mathbb{C}$ is radial if there is a function $\tilde{f}: \left[ 1, \infty \right) \rightarrow \mathbb{C}$ such that $f \left( x \right) = \tilde{f} \left( \cosh d \left( 0, x \right) \right)$.
				\item A function $f: \mathbb{S}^n \rightarrow \mathbb{C}$ is radial if there is a function $\tilde{f}: \left[ -1, 1 \right] \rightarrow \mathbb{C}$ such that $f \left( x \right) = \tilde{f} \left( \cos d \left( 0, x \right) \right)$.
			\end{enumerate}
		\end{definition}
		\begin{remark}
			\normalfont
			The choices of using $\cosh d \left( 0, x \right)$ and $\cos d \left( 0, x \right)$ for the hyperbolic space and the sphere, respectively, are for computational convenience. However, they do have a geometric meaning to them. If we view $\mathbb{H}^n$ as the upper sheet of the two sheeted hyperboloid $\left\lbrace \left( x_1, \cdots, x_{n + 1} \right) \in \mathbb{R}^{n{+1}} \Big| \sum\limits_{i = 1}^{n} x_i^2 - x_{n + 1}^2 = - 1 \right\rbrace$, then the collection of all points on $\mathbb{H}^n$ at a fixed distance $d$ from $e_{n + 1}$ (the ``origin" of $\mathbb{H}^n$) are on the plane $z = \cosh d$ (see Figure \ref{ZonalFunctionsHn}). 

			Similarly, the points on the sphere that are at a fixed distance $d$ from $e_{n + 1}$ are on the plane $z = \cos d$ (see Figure \ref{ZonalFunctionsSn}). We see that this consideration helps us a lot when it comes to integration of radial functions. It also helps that the functions $\cosh$ and $\cos$ are complimentary functions to $\sinh$ and $\sin$, respectively, which are the ``curvature functions" of $\mathbb{H}^n$ and $\mathbb{S}^n$. The sections of $\mathbb{H}^n$ and $\mathbb{S}^n$ by the planes $z = \cosh d$ and $z = \cos d$, respectively, are called ``\textit{zones}".

			\begin{figure}[ht!]
				\centering
				\begin{subfigure}[t]{0.45\textwidth}
					\centering
					\includegraphics[scale=0.225]{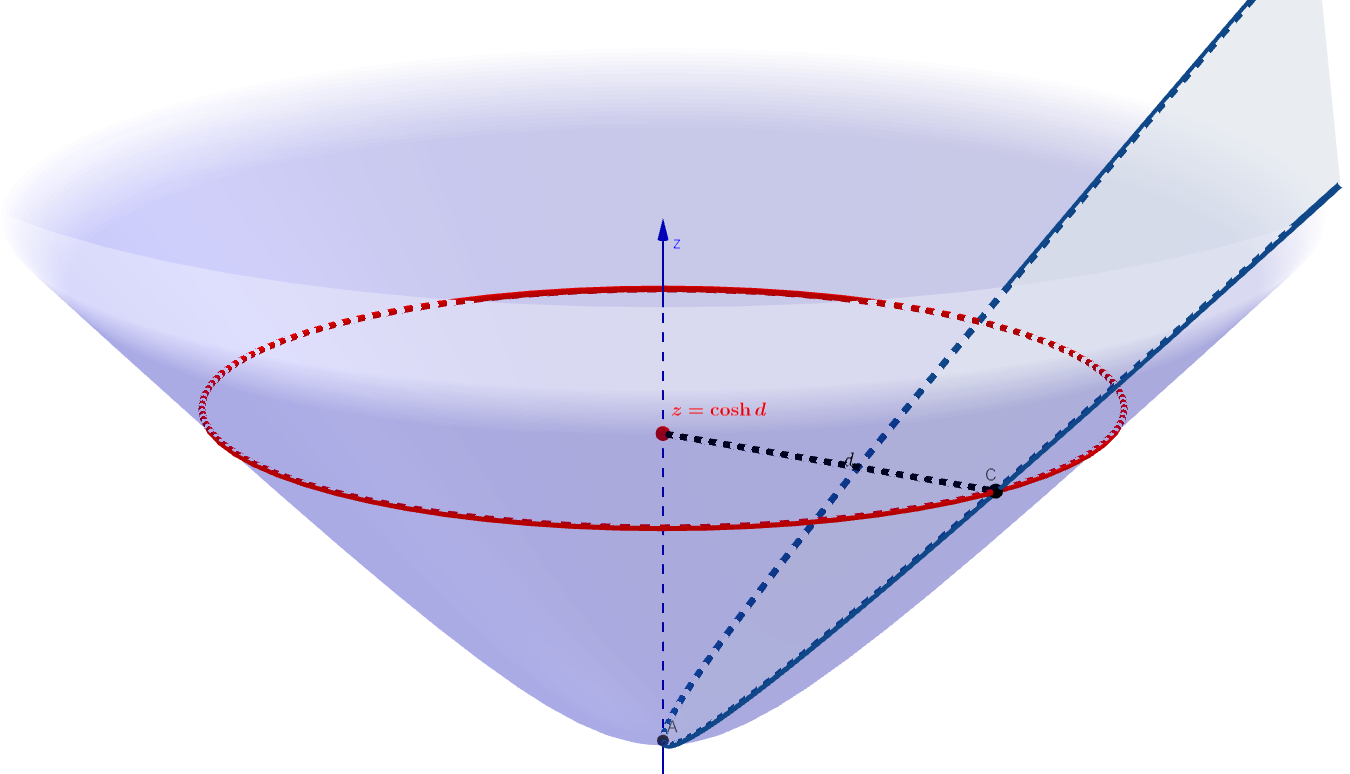}
					\caption{Zones of $\mathbb{H}^n$.}
					\label{ZonalFunctionsHn}
				\end{subfigure}
				~~~~~~~~
				\begin{subfigure}[t]{0.45\textwidth}
					\centering
					\includegraphics[scale=0.325]{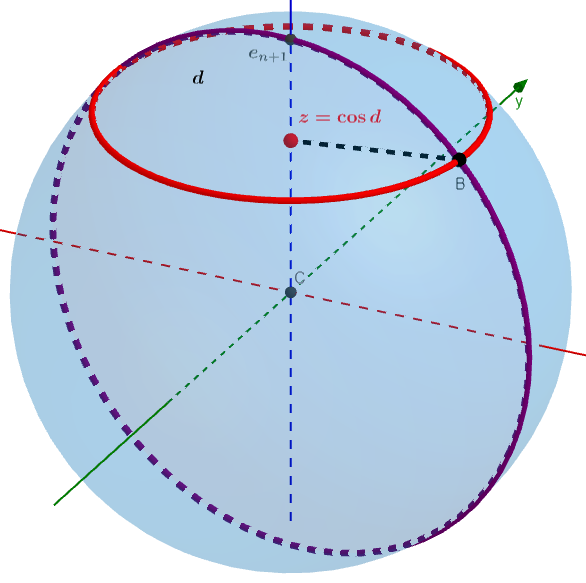}
					\caption{Zones of $\mathbb{S}^n$.}
					\label{ZonalFunctionsSn}
				\end{subfigure}
				\caption{Collection of points at a fixed distance from the origin in a space of consatnt curvature.}
				\label{RightTriangleFigure2}
			\end{figure}
		\end{remark}
		Before we move on to the next topic of discussion, we also give the definition of radial functions on $\Xi_k \left( X \right)$.
		\begin{definition}[Radial Functions on $\Xi_k \left( X \right)$]~
			\label{RadialFunctionsXiK}
			\begin{enumerate}
				\item A function $\varphi: \Xi_k \left( \mathbb{R}^n \right) \rightarrow \mathbb{C}$ if radial if there is a function $\tilde{\varphi}: \left[ 0, \infty \right) \rightarrow \mathbb{C}$ such that $\varphi \left( \xi \right) = \tilde{\varphi} \left( d \left( 0, \xi \right) \right)$.
				\item A function $\varphi: \Xi_k \left( \mathbb{H}^n \right) \rightarrow \mathbb{C}$ is radial if there is a function $\tilde{\varphi}: \left[ 0, \infty \right) \rightarrow \mathbb{C}$ such that $\varphi \left( \xi \right) = \tilde{\varphi} \left( \sinh d \left( 0, {\xi} \right) \right)$.
				\item A function $\varphi: \Xi_k \left( \mathbb{S}^n \right) \rightarrow \mathbb{C}$ is radial if there is a function $\tilde{\varphi}: \left[ {0}, 1 \right] \rightarrow \mathbb{C}$ such that $\varphi \left( \xi \right) = \tilde{\varphi} \left( \sin d \left( 0, \xi \right) \right)$.
			\end{enumerate}
		\end{definition}
		\subsection{Measure on  \texorpdfstring{$\Xi_k \left( X \right)$}{Affine Grassmann}}
		\label{MeasureSubsection}
		We require a measure on $\Xi_k \left( X \right)$, that allows us to talk about the Lebesgue spaces. The descriptions of these measures are well known. We only wish to write these descriptions in a form suitable to us. It is well known that when $\Xi_k \left( \mathbb{H}^n \right)$ is identified with $\faktor{SO_0 \left( n, 1 \right)}{S \left( O \left( n - k \right) \times O \left( k + 1 \right) \right)}$, the measure on $\Xi_k \left( \mathbb{H}^n \right)$ is given by (see, for instance, \cite{BerensteinRubin})
		\begin{equation}
			\label{MeasureXiKBerenstein}
			\int\limits_{\Xi_k \left( \mathbb{H}^n \right)} \varphi \left( \eta \right) \mathrm{d}\eta = \left| \mathbb{S}^{n - k - 1} \right| \int\limits_{0}^{\infty} \int\limits_{SO(n)} \varphi \left( \gamma g_{\theta}^{-1} \cdot \mathbb{R}^k \right) \cosh^{k} \theta \sinh^{n - k - 1} \theta \ \mathrm{d}\gamma \mathrm{d}\theta,
		\end{equation}
		where $g_{\theta} = \left[ \begin{matrix}
										\cosh \theta & 0 & \sinh \theta \\
										0 & I_{n - 1} & 0 \\
										\sinh \theta & 0 & \cosh \theta
									\end{matrix} \right]$ is a hyperbolic rotation\footnote{The action is effectively the translation of $e_{n + 1}$ to a point at a hyperbolic distance $\theta$ in the $\left\lbrace e_1, e_{n + 1} \right\rbrace$-plane.}. The right-hand-side of Equation \eqref{MeasureXiKBerenstein} can be understood in the following way. First, to obtain any $k$-dimensional totally-geodesic submanifold $\eta$ of $\mathbb{H}^n$ at a distance $\theta$ from the origin, one can start with a submanifold through the origin, ``translate" it by a distance $\theta$, and then perform an appropriate rotation to obtain $\eta$. Essentially, this mimics ``polar decomposition", from where we get the powers of $\cosh$ and $\sinh$. In our parametrization of {Proposition} \ref{TGSParametrization}, we have that any $k$-dimensional totally-geodesic submanifold can be given by the pair $\left( g \cdot \mathbb{R}^k, h g \cdot e_{k + 1} \right)$, where $g \in SO(n)$ and $h > 0$ is the hyperbolic distance of the submanifold $\left( g \cdot \mathbb{R}^k, h g \cdot e_{k + 1} \right)$ from the hyperbolic origin. Using this parametrization {in} Equation \eqref{MeasureXiKBerenstein}, we get
		\begin{equation}
			\label{MeasureXiK}
			\int\limits_{\Xi_k \left( \mathbb{H}^n \right)} \varphi \left( \eta \right) \mathrm{d}\eta = \left| \mathbb{S}^{n - k - 1} \right| \int\limits_{SO(n)} \int\limits_{0}^{\infty} \varphi \left( g \cdot \mathbb{R}^{k}, h g \cdot e_{k + 1} \right) \cosh^{k} h \sinh^{n - k - 1} h \ \mathrm{d}h \mathrm{d}g.
		\end{equation}
		Similarly, for suitable functions $\varphi$ on $\Xi_k \left( \mathbb{S}^n \right)$, we have (see, for instance, \cite{KumarRay})
		\begin{equation}
			\label{MeasureXiKSNKR}
			\int\limits_{\Xi_k \left( \mathbb{S}^n \right)} \varphi \left( \eta \right) \mathrm{d}\eta = C_{n, k} \int\limits_{SO(n)} \int\limits_{0}^{\frac{\pi}{2}} \varphi \left( \gamma g_{\theta}^{-1} \xi_{0} \right) \sin^{n - k - 1} \theta \cos^{k} \theta \ {\mathrm{d}\theta \mathrm{d}\gamma},
		\end{equation}
		where $g_{\theta} = \left[ \begin{matrix}
										\cos \theta & 0 & - \sin \theta \\
										0 & I_{n - 1} & 0 \\
										\sin \theta & 0 & \cos \theta
									\end{matrix} \right]$ is a spherical ``translation" that takes $e_{n + 1}$ to a point in $\left\lbrace e_1, e_{n + 1} \right\rbrace$-plane at a (spherical) distance of $\theta$ from $e_{n + 1}$, and $\xi_0$ is a $k$-dimensional totally-geodesic submanifold passing through $e_{n + 1}$. We may very well assume $\xi_0$ to have the tangent space $\mathbb{R}^k$. Therefore, in our parametrization of $\Xi_k \left( \mathbb{S}^n \right)$ given in {Proposition} \ref{TGSParametrization}, we have
		\begin{equation}
			\label{MeasureXiKSN}
			\int\limits_{\Xi_k \left( \mathbb{S}^n \right)} \varphi \left( \eta \right) \mathrm{d}\eta = C_{n, k} \int\limits_{SO(n)} \int\limits_{0}^{\frac{\pi}{2}} \varphi \left( g \cdot \mathbb{R}^k, hg \cdot e_{k + 1} \right) \sin^{n - k - 1} h \cos^{k} h \ \mathrm{d}h \mathrm{d}g.
		\end{equation}
		With the introduction of measure on $\Xi_k \left( X \right)$, we can define the Lebesgue spaces. In this article we are interested in the weighted Lebesgue spaces. For the same, we consider radial weights in terms of the powers of ``curvature function" and its complimentary function. Particularly, we define the following weighted Lebesgue spaces.
		\begin{equation}
			\label{LPX}
			L^p_{\alpha_1, \alpha_2} \left( X \right) := \left\lbrace f : X \rightarrow \mathbb{C} \Bigg| {f \text{ is measurable and }} \int\limits_{X} \left| f \left( x \right) \right|^p s_c^{\alpha_1} d \left( 0, x \right) \left| s_c' \right|^{\alpha_2} d \left( 0, x \right) \mathrm{d}x < + \infty \right\rbrace,
		\end{equation}
		and
		\begin{equation}
			\label{LPXiK}
			L^p_{\beta_1, \beta_2} \left( \Xi_k \left( X \right) \right) := \left\lbrace \varphi: \Xi_k \left( X \right) \rightarrow \mathbb{C} \Bigg| {\varphi \text{ is measurable and }} \int\limits_{\Xi_k \left( X \right)} \left| \varphi \left( \xi \right) \right|^p s_c^{\beta_1} d \left( 0, \xi \right) \left( s_c' \right)^{\beta_2} d \left( 0, \xi \right) \mathrm{d}x < + \infty \right\rbrace.
		\end{equation}
		For a function $f: X \rightarrow \mathbb{C}$, its $L^p_{\alpha_1, \alpha_2} \left( X \right)${-}norm is given by
		\begin{equation}
			\label{LPXNorm}
			\| f \|_{L^p_{\alpha_1, \alpha_2} \left( X \right)} = \left[ \int\limits_{X} \left| f \left( x \right) \right|^p s_c^{\alpha_1} d \left( 0, x \right) \left| s_c' \right|^{\alpha_2} d \left( 0, x \right) \mathrm{d}x \right]^{\frac{1}{p}}{.}
		\end{equation}
		Similarly, for a function $\varphi: \Xi_k \left( X \right) \rightarrow \mathbb{C}$, its $L^p_{\beta_1, \beta_2} \left( \Xi_k \left( X \right) \right)${-}norm is given by
		\begin{equation}
			\label{LPXikNorm}
			\| \varphi \|_{L^p_{\beta_1, \beta_2} \left( \Xi_k \left( X \right) \right)} = \left[ \int\limits_{\Xi_k \left( X \right)} \left| \varphi \left( \xi \right) \right|^p s_c^{\beta_1} d \left( 0, \xi \right) \left( s_c' \right)^{\beta_2} d \left( 0, \xi \right) \mathrm{d}{\xi} \right]^{\frac{1}{p}}.
		\end{equation}
		\begin{remark}
			\normalfont
			For the Euclidean space $\mathbb{R}^n$, we have $c = 0$ and hence $s_c' \left( t \right) = 1$, for any $t \in \mathbb{R}$. Consequently, the weighted Lebesgue spaces have only one weight ($\| x \|^{\alpha}$ and $\left( d \left( 0, \xi \right) \right)^{\beta}$).
		\end{remark}
		Here, $1 \leq p < \infty$, $\alpha_1, \alpha_2, \beta_1, \beta_2 \in \mathbb{R}$, and the ``curvature function" is given by Equation \eqref{CurvatureFunction}. We recall that one of the aims of this article is answering Question \ref{MainQuestion2}. {Since} answering the question for general (non-radial) functions is difficult, we give necessary and sufficient conditions for the boundedness of the $k$-plane transform acting on radial functions. For the same, we now state certain {formulae} (and their variations) for the weighted norms of radial functions.
		\begin{theorem}~
			\label{LpNormsRadialX}
			\begin{enumerate}
				\item For a radial function $f: \mathbb{R}^n \rightarrow \mathbb{C}$, from polar decomposition, we get
				\begin{equation}
					\label{LpNormRadialRn1}
					\| f \|_{L^p_{\alpha} \left( \mathbb{R}^n \right)} = C \left[ \int\limits_{0}^{\infty} \left| \tilde{f} \left( t \right) \right|^p t^{\alpha + n - 1} \mathrm{d}t \right]^{\frac{1}{p}}.
				\end{equation}
				\item For a radial function $f: \mathbb{H}^n \rightarrow \mathbb{C}$, using the polar decomposition, we get
				\begin{equation}
					\label{LpNormRadialHn1}
					\| f \|_{L^p_{\alpha_1, \alpha_2} \left( \mathbb{H}^n \right)} = C \left[ \int\limits_{0}^{\infty} \left| \tilde{f} \left( \cosh t \right) \right|^p \sinh^{\alpha_1 + n - 1} t \cosh^{\alpha_2} t \ \mathrm{d}t \right]^{\frac{1}{p}}.
				\end{equation}
				Moreover, by replacing $\cosh t$ with $t$, we get
				\begin{equation}
					\label{LpNormRadialHn2}
					\| f \|_{L^p_{\alpha_1, \alpha_2} \left( \mathbb{H}^n \right)} = C \left[ \int\limits_{1}^{\infty} \left| \tilde{f} \left( t \right) \right|^p \left( t^2 - 1 \right)^{\frac{\alpha_1 + n}{2} - 1} t^{\alpha_2} \mathrm{d}t \right]^{\frac{1}{p}}.
				\end{equation}
				Further, by replacing $\sinh^2 t$ with $t$ in Equation \eqref{LpNormRadialHn1}, we get
				\begin{equation}
					\label{LpNormRadialHn3}
					\| f \|_{L^p_{\alpha_1, \alpha_2} \left( \mathbb{H}^n \right)} = C \left[ \int\limits_{0}^{\infty} \left| \tilde{f} \left( \sqrt{t + 1} \right) \right|^p t^{{\frac{\alpha_1 + n}{2} - 1}} \left( t + 1 \right)^{{\frac{\alpha_2 - 1}{2}}} \frac{\mathrm{d}t}{t} \right]^{\frac{1}{p}}.
				\end{equation}
				\item {Using the polar decomposition, for even radial functions $f: \mathbb{S}^n \rightarrow \mathbb{C}$, we have,}
				\begin{equation}
					\label{LpNormRadialSn1}
					\| f \|_{L^p_{\alpha_1, \alpha_2} \left( \mathbb{S}^n \right)} = C \left[ \int\limits_{0}^{\frac{\pi}{2}} \left| \tilde{f} \left( \cos t \right) \right|^p \sin^{\alpha_1 + n - 1} t \cos^{\alpha_2} t \ \mathrm{d}t \right]^{\frac{1}{p}}.
				\end{equation}
				Moreover, by replacing $\cos^{{2}} t$ with $t$, we get
				\begin{equation}
					\label{LpNormRadialSn2}
					\| f \|_{L^p_{\alpha_1, \alpha_2} \left( \mathbb{S}^n \right)} = C \left[ \int\limits_{0}^{1} \left| \tilde{f} \left( t \right) \right|^p \left( 1 - t \right)^{\frac{\alpha_1 + n}{2} - 1} t^{\frac{\alpha_2 - 1}{2}} \mathrm{d}t \right]^{\frac{1}{p}}.
				\end{equation}
				Further, by replacing $\sin^2 t$ with $t$ in Equation \eqref{LpNormRadialSn1}, we get
				\begin{equation}
					\label{LpNormRadialSn3}
					\| f \|_{L^p_{\alpha_1, \alpha_2} \left( \mathbb{S}^n \right)} = C \left[ \int\limits_{0}^{1} \left| \tilde{f} \left( \sqrt{1 - t} \right) \right|^p \left( 1 - t \right)^{\frac{\alpha_{{2}} - 1}{2}} t^{\frac{\alpha_1 + n}{2} - 1} \mathrm{d}t \right]^{\frac{1}{p}} {.}
				\end{equation}
				Lastly, by replacing $\tan^2 t$ by $t$ in Equation \eqref{LpNormRadialSn1}, we get
				\begin{equation}
					\label{LpNormRadialSn4}
					\| f \|_{L^p_{\alpha_1, \alpha_2} \left( \mathbb{S}^n \right)} = C \left[ \int\limits_{0}^{1} \left| \tilde{f} \left( \frac{1}{\sqrt{1 + t}} \right) \right|^p \frac{t^{{\frac{\alpha_1 + n}{2}}}}{\left( 1 + t \right)^{{\frac{\alpha_1 + \alpha_2 + n + 1}{2}}}} \frac{\mathrm{d}t}{t} \right]^{\frac{1}{p}}.
				\end{equation}
			\end{enumerate}
		\end{theorem}
		We again mention that Equations \eqref{LpNormRadialRn1}--\eqref{LpNormRadialSn4} find their use at various places throughout this article.		
		\subsection{The \texorpdfstring{$k$}{k}-plane Transform and its Dual}
		\label{RTDualRTSubsection}
		We now briefly talk about the $k$-plane transform and its dual on the spaces of constant curvature. We begin with basic definitions about these maps, and then state a formula for the totally-geodesic $k$-plane transform of radial functions. As discussed earlier, the $k$-plane transform takes in functions on the space $X$ and give{s} out a function on the ``affine Grassmannian", $\Xi_k \left( X \right)$. It does so by integrating the {restriction of the} given function on the given $k$-dimensional totally-geodesic submanifold. Precisely, we have the following definition.
		\begin{definition}[$k$-plane transform]
			\label{KPlaneTransformDefinition}
			Let $X$ be a space of constant curvature and $f: X \rightarrow \mathbb{C}$ be a ``nice" function\footnote{By ``nice", we mean a function for which all integrals in the definition make sense.}. The $k$-plane transform of $f$ is a function $R_kf: \Xi_k \left( X \right) \rightarrow \mathbb{C}$, given by
			\begin{equation}
				\label{KPlaneTransformEquation}
				R_kf \left( \xi \right) := \int\limits_{\xi} f \left( x \right) \mathrm{d}_{\xi}x,
			\end{equation}
			where, $\mathrm{d}_{\xi}x$ is the induced measure on $\xi \subseteq X$ from the Riemannian metric on $X$.
		\end{definition}
		For details about induced measures on embedded (Riemannian) submanifolds of a Riemannian manifold, we refer the reader to \cite{LeeRM}. Indeed, when $X = \mathbb{R}^n$, we get back to the definition as given in Equation \eqref{KPlaneTransformDefinitionRn}.
		
		We also have the dual $k$-plane transform, which integrates a function $\varphi: \Xi_k \left( X \right) \rightarrow \mathbb{C}$ on the set of all the $k$-dimensional totally-geodesic submanifolds that contain a particular point.
		\begin{definition}[Dual $k$-plane transform]
			\label{DualKPlaneTransformDefinition}
			Let $X$ be a space of constant curvature and $\varphi: \Xi_k \left( X \right) \rightarrow \mathbb{C}$ be a ``nice" function. The dual $k$-plane transform of $\varphi$ is a function $R_k^*\varphi: X \rightarrow \mathbb{C}$, given by
			\begin{equation}
				\label{DualKPlaneTransformEquation}
				R_k^*\varphi \left( x \right) := \int\limits_{\left\lbrace \xi \in \Xi_k \left( X \right) | x \in \xi \right\rbrace} \varphi \left( \xi \right) \mathrm{d}_S\xi.
			\end{equation}
			Here, the measure $\mathrm{d}_S\xi$ is on the set $S = \left\lbrace \xi \in \Xi_k \left( X \right) | x \in \xi \right\rbrace$.
		\end{definition}
		Measure $\mathrm{d}_S \xi$ on the set $S = \left\lbrace \xi \in \Xi_k \left( X \right) | x \in \xi \right\rbrace$ is often obtained by considering $S$ as a quotient of the isometry group of $X$. For details, we refer the reader to \cite{HelgasonCCS}. We have the following duality relation that is of help later.
		\begin{theorem}[\cite{Rubin}]
			\label{DualityKPlaneTransform}
			Let $f: X \rightarrow \mathbb{C}$ and $\varphi: \Xi_k \left( X \right) \rightarrow \mathbb{C}$ be ``nice" functions. Then, the following duality holds whenever the integrals on both sides exist in the Lebesgue sense.
			\begin{equation}
				\label{DualityEquation}
				\int\limits_{\Xi_k \left( X \right)} R_kf \left( \xi \right) \varphi \left( \xi \right) \mathrm{d}\xi = \int\limits_{X} f \left( x \right) R_k^*\varphi \left( x \right) \mathrm{d}x.
			\end{equation}
		\end{theorem}
		Of importance to us are the formulae of the $k$-plane transform and its dual for radial functions. They can be easily obtained from the polar decomposition of the submanifold $\xi$. For details, we refer the reader to \cite{BerensteinRubin}, \cite{RubinInversionSphere}, and \cite{HelgasonRTBook}. We only state the formulae here.
		\begin{proposition}
			\label{KPlaneRadialFormulaeRn}
			For a radial function $f: \mathbb{R}^n \rightarrow \mathbb{C}$, the $k$-plane transform is given by
			\begin{equation}
				\label{KPlaneTransformRadialRnEquation}
				R_kf \left( \xi \right) = \left| \mathbb{S}^{k - 1} \right| \int\limits_{d \left( 0, \xi \right)}^{\infty} \tilde{f} \left( t \right) \left( t^2 - \left( d \left( 0, \xi \right) \right)^2 \right)^{\frac{k}{2} - 1} t \ \mathrm{d}t.
			\end{equation}
			Equation \eqref{KPlaneTransformRadialRnEquation} can also be written as
			\begin{equation}
				\label{KPlaneTransformRadialRnEquation2}
				R_kf \left( \xi \right) = \left| \mathbb{S}^{k - 1} \right| \int\limits_{d \left( 0, \xi \right)}^{\infty} \tilde{f} \left( t \right) \left( 1 - \frac{\left( d \left( 0, \xi \right) \right)^2}{t^2} \right)^{\frac{k}{2} - 1} t^{k - 1} \ \mathrm{d}t.
			\end{equation}
		\end{proposition}
		{We now state a few formulae for the $k$-plane transform on the real hyperbolic space $\mathbb{H}^n$. These formulae play important roles at various places throughout the article.}
		\begin{proposition}
			\label{KPlaneRadialFormulaeHn}
			For a radial function $f: \mathbb{H}^n \rightarrow \mathbb{C}$, the $k$-plane transform is given by
			\begin{equation}
				\label{KPlaneTransformRadialHnEquation}
				R_kf \left( \xi \right) = \frac{| \mathbb{S}^{k - 1} |}{{\cosh}^{k - 1} d \left( 0, \xi \right)} \int\limits_{d \left( 0, \xi \right)}^{\infty} \tilde{f} \left( \cosh t \right) \left( \cosh^2 t - \cosh^2 d \left( 0, \xi \right) \right)^{\frac{k}{2} - 1} \sinh t \ \mathrm{d}t.
			\end{equation}
			Equation \eqref{KPlaneTransformRadialHnEquation} can also be written as
			\begin{equation}
				\label{KPlaneTransformRadialHnEquation3}
				R_kf \left( \xi \right) = \frac{\left| \mathbb{S}^{k - 1} \right|}{\cosh^{k - 1} d \left( 0, \xi \right)} \int\limits_{d \left( 0, \xi \right)}^{\infty} \tilde{f} \left( \cosh t \right) \left( 1 - \frac{\sinh^2 d \left( 0, \xi \right)}{\sinh^2 t} \right)^{\frac{k}{2} - 1} \sinh^{k - 1} t \ \mathrm{d}t.
			\end{equation}
			Further, using the relations between hyperbolic trigonometric functions, we have,
			\begin{equation}
				\label{KPlaneTransformRadialHnEquation2}
				R_kf \left( \xi \right) = \frac{\left| \mathbb{S}^{k - 1} \right|}{\cosh d \left( 0, \xi \right)} \int\limits_{d \left( 0, \xi \right)}^{\infty} \tilde{f} \left( \cosh t \right) \left( 1 - \frac{\tanh^2 d \left( 0, \xi \right)}{\tanh^2 t} \right)^{\frac{k}{2} - 1} \sinh^{k - 1} t \ \mathrm{d}t.
			\end{equation}
		\end{proposition}
		Lastly, we state the formulae for the $k$-plane transform of radial functions on the sphere $\mathbb{S}^n$.
		\begin{proposition}
			\label{KPlaneRadialFormulaeSn}
			For a radial (even) function $f: \mathbb{S}^n \rightarrow \mathbb{C}$, the $k$-plane transform is given by
			\begin{equation}
				\label{KPlaneTransformRadialSnEquation}
				R_kf \left( \xi \right) = \frac{2 | \mathbb{S}^{k - 1} |}{\cos^{k - 1} d \left( 0, \xi \right)} \int\limits_{d \left( 0, \xi \right)}^{\frac{\pi}{2}} \tilde{f} \left( \cos t \right) \left( \cos^2 d \left( 0, \xi \right) - \cos^2 t \right)^{\frac{k}{2} - 1} \sin t \ \mathrm{d}t.
			\end{equation}
			Using the relations between trigonometric functions, Equation \eqref{KPlaneTransformRadialSnEquation} can be written as
			\begin{equation}
				\label{KPlaneTransformRadialSnEquation2}
				R_kf \left( \xi \right) = \frac{2 \left| \mathbb{S}^{k - 1} \right|}{\cos d \left( 0, \xi \right)} \int\limits_{d \left( 0, \xi \right)}^{{\frac{\pi}{2}}} \tilde{f} \left( \cos t \right) \left( 1 - \frac{\tan^2 d \left( 0, \xi \right)}{\tan^2 t} \right)^{\frac{k}{2} - 1} \sin^{k - 1} t \ \mathrm{d}t.
			\end{equation}
		\end{proposition}
		Next, we see a few formulae for the dual $k$-plane transform of radial functions on $\Xi_k \left( X \right)$, where $X = \mathbb{R}^n$, $\mathbb{H}^n$, or $\mathbb{S}^n$. We begin with the Euclidean case. {These formulae can be found in \cite{Rubin}, \cite{BerensteinRubin}, and \cite{RubinInversionSphere}.}
		\begin{proposition}
			\label{DualKPlaneRadialFormulaeX}
			{We have the following formulae for the dual $k$-plane transform of radial functions on $\Xi_k \left( X \right)$, where $X$ is a space of constant curvature.
			\begin{enumerate}
				\item \label{DualKPlaneRadialFormulaeRn}
				For a radial function $\varphi: \Xi_k \left( \mathbb{R}^n \right) \rightarrow \mathbb{C}$, its dual $k$-plane transform is given by
				\begin{equation}
					\label{DualKPlaneTransformRadialRnEquation}
					R_k^*\varphi \left( x \right) = \frac{\left| \mathbb{S}^{k - 1} \right| \left| \mathbb{S}^{n - k - 1} \right|}{\left| \mathbb{S}^{n - 1} \right|} \frac{1}{\| x \|^{k - 1}} \int\limits_{0}^{\| x \|} \tilde{\varphi} \left( s \right) \left( \| x \|^2 - s^2 \right)^{\frac{k}{2} - 1} \mathrm{d}s.
				\end{equation}
				\item \label{DualKPlaneRadialFormulaeHn}
				For a radial function $\varphi: \Xi_k \left( \mathbb{H}^n \right) \rightarrow \mathbb{C}$, the dual $k$-plane transform of $\varphi$ is given by
				\begin{equation}
					\label{DualKPlaneTransformRadialHnEquation}
					R_k^*\varphi \left( x \right) = \frac{| \mathbb{S}^{k - 1} | | \mathbb{S}^{n - k - 1} |}{| \mathbb{S}^{n - 1} |} \frac{1}{\sinh^{n - 2} d \left( 0, \xi \right)} \int\limits_{0}^{\sinh d \left( 0, \xi \right)} \tilde{\varphi} \left( s \right) \left( \sin^2 d \left( 0, \xi \right) - s^2 \right)^{\frac{k}{2} - 1} s^{n - k - 1} \mathrm{d}s.
				\end{equation}
				\item \label{DualKPlaneRadialFormulaeSn}
				For a radial function $\varphi: \Xi_k \left( \mathbb{S}^n \right) \rightarrow \mathbb{C}$, the dual $k$-plane transform of $\varphi$ is given by
				\begin{equation}
					\label{DualKPlaneTransformRadialSnEquation}
					R_k^*\varphi \left( \xi \right) = \frac{| \mathbb{S}^{k - 1} | | \mathbb{S}^{n - k - 1} |}{| \mathbb{S}^{n - 1} |} \frac{1}{\sin^{n - 2} d \left( 0, x \right)} \int\limits_{0}^{\sin d \left( 0, x \right)} \tilde{\varphi} \left( s \right) \left( \sin^2 d \left( 0, x \right) - s^2 \right)^{\frac{k}{2} - 1} s^{n - k - 1} \mathrm{d}s.
				\end{equation}
			\end{enumerate}}	
		\end{proposition}
		We are interested in estimating the $L^p_{\beta_1, \beta_2} \left( \Xi_k \left( X \right) \right)$-norms of the functions $R_kf$. By making suitable changes of variables in the formulae for the $k$-plane transform of radial functions and Equation \eqref{LPXikNorm}, we get the following well-known result. We only state it here for reference.
		\begin{theorem}~
			\label{NormRkRadialExpression}
			\begin{enumerate}
				\item On the Euclidean space, we have for radial functions,
				\begin{equation}
					\label{NormRkRn}
					\| R_kf \|_{L^p_{\beta} \left( \mathbb{R}^n \right)} = C \left[ \int\limits_{0}^{\infty} \left| \int\limits_{s}^{\infty} \tilde{f} \left( t \right) \left( t^2 - s^2 \right)^{\frac{k}{2} - 1} t \ \mathrm{d}t \right|^p s^{\beta + n - k - 1} \mathrm{d}s \right]^{\frac{1}{p}}.
				\end{equation}
				\item On the Hyperbolic space, {by Equation \eqref{KPlaneTransformRadialHnEquation}, with a change of variables, $\cosh t$ replaced by $t$} and $\cosh d \left( 0, \xi \right) = s$, we get for radial functions
				\begin{equation}
					\label{NormRkHn}
					\| R_kf \|_{L^p_{\beta_1, \beta_2} \left( \mathbb{H}^n \right)} = C \left[ \int\limits_{1}^{\infty} \left| \int\limits_{s}^{\infty} \frac{\tilde{f} \left( t \right) \left( t^2 - s^2 \right)^{\frac{k}{2} - 1}}{s^{k - 1}} \ \mathrm{d}t \right|^p s^{\beta_2 + k - 1} \left( s^2 - 1 \right)^{\frac{\beta_1 + n - k}{2} - 1} \mathrm{d}s \right]^{\frac{1}{p}}.
				\end{equation}
				\item On the sphere, {by Equation \eqref{KPlaneTransformRadialSnEquation}, with a change of variables, $\cos t$ replaced by $t$} and $\cos d \left( 0, \xi \right) = s$, we get
				\begin{equation}
					\label{NormRkSn}
					\| R_kf \|_{L^p_{\beta_1, \beta_2}} \left( \mathbb{S}^n \right) = C \left[ \int\limits_{0}^{1} \left| \int\limits_{0}^{s} \frac{\tilde{f} \left( t \right) \left( s^2 - t^2 \right)^{\frac{k}{2} - 1}}{s^{k - 1}} \ \mathrm{d}t \right|^p s^{\beta_2 + k - 1} \left( 1 - s^2 \right)^{\frac{\beta_1 + n - k}{2} - 1} \mathrm{d}s \right]^{\frac{1}{p}}.
				\end{equation}
			\end{enumerate}
		\end{theorem}
		\subsection{Lorentz Spaces and a weighted interpolation result}
		\label{FunctionSpacesSubsection}
		In this section, we mention a few definitions and results from function spaces that we require throughout this article. We have already seen the weighted Lebesgue spaces on the spaces of constant curvature in Subsection \ref{MeasureSubsection}.
		
		Apart from the weighted Lebesgue spaces, we also work on Lorentz spaces for obtaining end-point behaviour of the $k$-plane transform. We give the definition of {such} spaces here for an abstract measure space.
		\begin{definition}[Lorentz space $L^{p, q}$ (\cite{GrafakosCFA})]
			\label{LorentzSpaceDefinition}
			Let $\left( X, \mu \right)$ be a measure space and $p, q \geq 1$. The Lorentz space $L^{p, q} \left( X, \mu \right)$ is the space of all measurable functions $f$ on $X$ for which the Lorentz norm,
			\begin{equation}
				\label{LorentzNorm}
				\| f \|_{L^{p, q} \left( X, \mu \right)} = \begin{cases}
																p^{\frac{1}{q}} \left[ \int\limits_{0}^{\infty} \left( \lambda_f \left( s \right)^{\frac{1}{p}} s \right)^q \frac{\mathrm{d}s}{s} \right]^{\frac{1}{q}}, & q < + \infty, \\
																\sup\limits_{s > 0} \left\lbrace s \lambda_f \left( s \right)^{\frac{1}{p}} \right\rbrace, & q = + \infty,
															\end{cases}
			\end{equation}
			is finite. Here, $\lambda_f$ is the distribution function of $f$, given by
			\begin{equation}
				\label{DistributionFunction}
				\lambda_f \left( s \right) = \mu \left\lbrace x \in X | \left| f \left( x \right) \right| > s \right\rbrace.
			\end{equation}
		\end{definition}
		\begin{remark}
			\normalfont
			It is known (see for instance, \cite{GrafakosCFA}) that $\| \chi_E \|_{L^{p, 1} \left( X, \mu \right)} = \mu \left( E \right)^{\frac{1}{p}}$, for measurable sets $E \subseteq X$.
		\end{remark}
		The following theorem about boundedness of operators on Lorentz spaces is crucial.
		\begin{theorem}[\cite{SteinWeissFA}]
			\label{LorentzSpaceBoundedness}
			Let $T$ be a linear operator that maps finite linear combinations of characteristic functions {of measurable sets} $E \subseteq X$ {with} finite measure into a vector space $B$ {of functions on $X$}, with an order preserving norm, {$\| \cdot \|$}. {By an ``\textit{order preserving norm}" we mean that whenever $f, g \in B$ are such that $\left| f \left( x \right) \right| \leq \left| g \left( x \right) \right|$, for almost every $x \in X$, we have $\| f \| \leq \| g \|$.} If, for any measurable set with finite measure $E \subseteq X$, we have
			$$\| T \chi_E \| \leq C \| \chi_E \|_{L^{p, 1} \left( X, \mu \right)},$$
			where, the constant $C$ is independent of the choice of the set $E$, then, for every $f \in L^{p, 1} \left( X, \mu \right)$, we have
			$$\| Tf \| \leq C \| f \|_{L^{p, 1} \left( X, \mu \right)}.$$
		\end{theorem}
		We now give an interpolation result due to Stein and Weiss (\cite{SteinWeissWeightedInterpolation}). It can also be found in \cite{Bennet}, stated in a slightly different manner. We state it here in a way convenient to us.
		\begin{theorem}[Stein and Weiss]
			\label{RealInterpolationWeighted}
			Let $\left( X, \mu \right)$, and $\left( Y, \nu \right)$ be measure spaces. Let $\alpha_1, \alpha_2$ be positive functions on $X$, and $\beta_1, \beta_2$ be positive functions on $Y$. Let $1 \leq p_i, q_i \leq \infty$ for $i = 0, 1$, and suppose that $T$ is a operator such that for $i = 0, 1$, we have
			$$\| \left( Tf \right) \beta_i \|_{q_i} \leq C_i \| f \alpha_i \|_{p_i}.$$
			Then, for $\theta \in \left( 0, 1 \right)$, $\frac{1}{p} = \frac{1 - \theta}{p_0} + \frac{\theta}{p_1}$, $\frac{1}{q} = \frac{1 - \theta}{q_0} + \frac{\theta}{q_1}$, and $\alpha = \alpha_0^{1 - \theta} \alpha_1^{\theta}$ and $\beta = \beta_0^{1 - \theta} \beta_1^{\theta}$ such that $p < + \infty$, we have
			$$\| \left( Tf \right) \beta \|_{q} \leq C \| f \alpha \|_p.$$
		\end{theorem}
	\section{Existence of \texorpdfstring{$k$}{k}-plane transform}
		\label{ExistenceSection}
		In this section, we study the existence of the $k$-plane transforms on the spaces of constant curvature. Essentially, we wish to answer Question \ref{MainQuestion1} posed in Section \ref{IntroductionSection}. The existence of $k$-plane transforms for (weighted) $L^p${-}functions on $\mathbb{R}^n$ is well studied (see for instance, \cite{RubinKPlane}). We look at the weighted Lebesgue spaces on $\mathbb{H}^n$ and $\mathbb{S}^n$, given by Equation \eqref{LPX}, and ask for which $p \geq 1$ and choice of ${\alpha_1, \alpha_2 \in \mathbb{R}}$ does the $k$-plane transform exist.
		\subsection{The Hyperbolic Space \texorpdfstring{$\mathbb{H}^n$}{}}
			We start by getting the existence of the $k$-plane transform on the weighted Lebesgue spaces on the real hyperbolic space. Particularly, we ask the following question.
			\begin{question}
				\label{ExistenceQuestionHn}
				\normalfont
				For which values of $p \geq 1$, and $\alpha_1, \alpha_2 \in \mathbb{R}$, does the $k$-plane transform {exist} for functions in $L^p_{\alpha_1, \alpha_2} \left( \mathbb{H}^n \right)$?
			\end{question}
			In obtaining the desired existence result, we employ the duality given in Theorem \ref{DualityKPlaneTransform}.
		\begin{theorem}
			\label{ExistenceMixed}
			For $\alpha_1 + \alpha_2 > k - n$, and $1 \leq p < \frac{\alpha_1 + \alpha_2 + n - 1}{k - 1}$, the $k$-plane transform exists for all functions $f \in {L^p_{\alpha_1, \alpha_2} \left( \mathbb{H}^n \right)}$. Also for functions in ${L^1_{\alpha_1, \alpha_2} \left( \mathbb{H}^n \right)}$ with $\alpha_1 + \alpha_2 = k - n$, the $k$-plane transform is well-defined. {Moreover, when the conditions stated above are not satisfied, there are radial functions $f \in L^p_{\alpha_1, \alpha_2} \left( \mathbb{H}^n \right)$ for which $R_kf \equiv + \infty$.}
		\end{theorem}
		\begin{proof}
			First, let us fix $\alpha_1 + \alpha_2 > k - n$ and $1 < p < \frac{\alpha_1 + \alpha_2 + n - 1}{k - 1}$. We choose 
			$$\gamma_1 > \max \left\lbrace k - n, - \frac{\alpha_1 + n}{p'} \right\rbrace \text{ and } k - n - \alpha_1 - \alpha_2 < \gamma_1 + \gamma_2 < \frac{1 - n - \alpha_1 - \alpha_2}{p'}.$$
			Such a choice of $\gamma_1$ and $\gamma_2$ is possible because under the assumptions on $\alpha_1, \alpha_2$, and $p$, we get an unbounded region of $\mathbb{R}^2$ with the constraints mentioned above (see Figure \ref{GammaChoice}). \\
		\begin{figure}[ht!]
			\centering
			\includegraphics[scale=0.75]{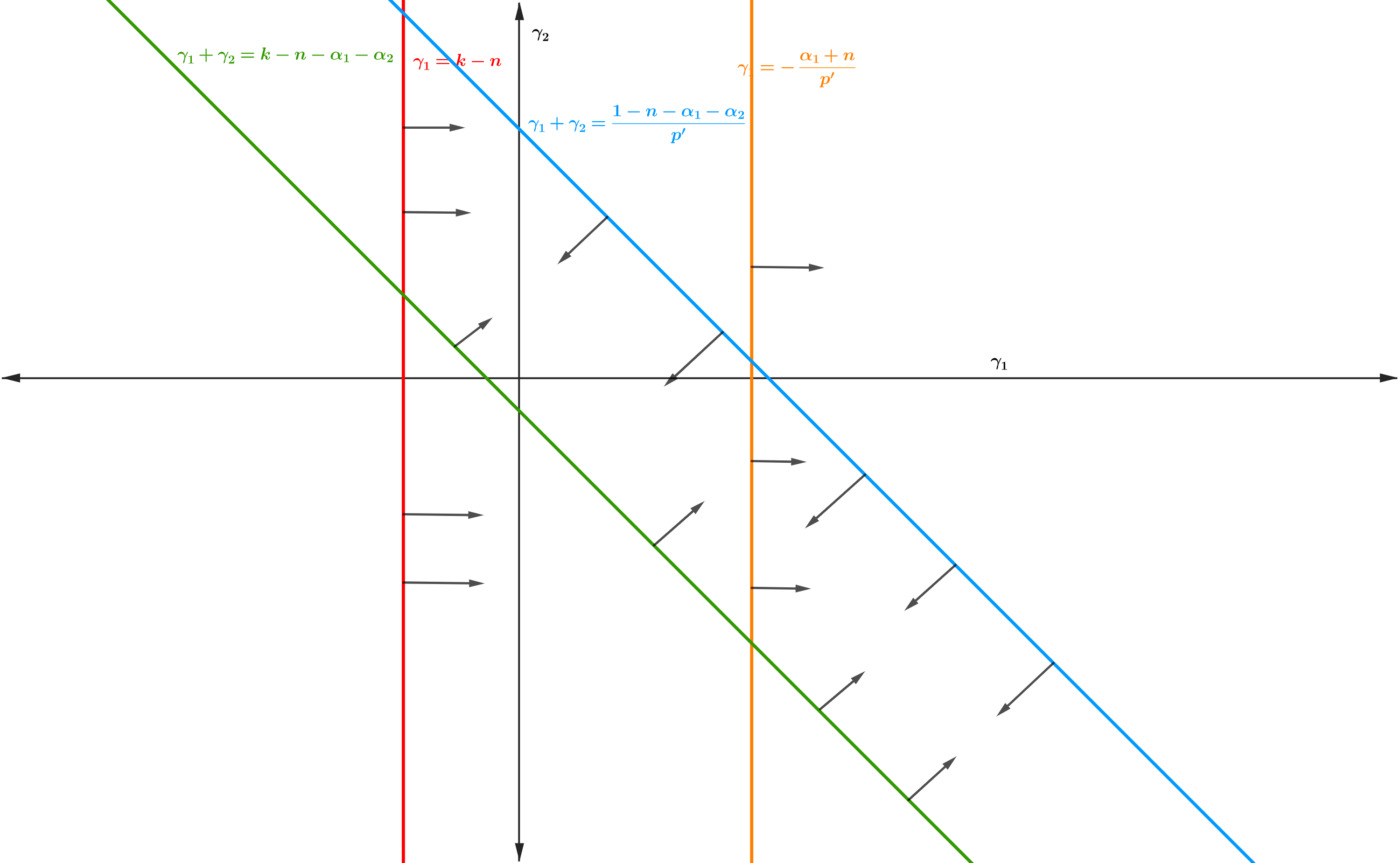}
			\caption{Region in $\mathbb{R}^2$ for choices of $\gamma_1$ and $\gamma_2$.}
			\label{GammaChoice}
		\end{figure}
		
		Now, using the duality of the $k$-plane transform given in Equation \eqref{DualityEquation} and {using} Equation \eqref{DualMixedEquation}, we get for any $f \in L^p_{\alpha_1, \alpha_2} \left( \mathbb{H}^n \right)$,
		\begin{align*}
			&\left| \int\limits_{\Xi_k \left( \mathbb{H}^n \right)} R_kf \left( \xi \right) \sinh^{\alpha_1 + \gamma_1} d \left( 0, \xi \right) \cosh^{\alpha_2 + \gamma_2} d \left( 0, \xi \right) \mathrm{d}\xi \right| \\
			&= C \left| \int\limits_{\mathbb{H}^n} f \left( x \right) \sinh^{\alpha_1 + \gamma_1} d \left( 0, x \right) \cosh^{\alpha_2 + \gamma_2} d \left( 0, x \right) {}_2F_1 \left( - \frac{\alpha_2 + \gamma_2}{2}, \frac{k}{2}; \frac{\alpha_1 + \gamma_1 + n}{2}; \tanh^2 d \left( 0, x \right) \right) \mathrm{d}x  \right|
		\end{align*}
		{We write $\alpha_1 = \frac{\alpha_1}{p} + \frac{\alpha_1}{p'}$ and $\alpha_2 = \frac{\alpha_2}{p} + \frac{\alpha_2}{p'}$, and use H\"{o}lder inequality with exponents $p$ and $p'$ along with polar decomposition on $\mathbb{H}^n$ to obtain,}
		\begin{align*}
			&{\left| \int\limits_{\Xi_k \left( \mathbb{H}^n \right)} R_kf \left( \xi \right) \sinh^{\alpha_1 + \gamma_1} d \left( 0, \xi \right) \cosh^{\alpha_2 + \gamma_2} d \left( 0, \xi \right) \mathrm{d}\xi \right|} \\
			&\leq C \| f \|_{{L^p_{\alpha_1, \alpha_2} \left( \mathbb{H}^n \right)}} \left( \int\limits_{0}^{\infty} \sinh^{\alpha_1 + p'\gamma_1 + n - 1} t \cosh^{\alpha_2 + p'\gamma_2} t \ {\Bigg|} {}_2F_1 \left( - \frac{\alpha_2 + \gamma_2}{2}, \frac{k}{2}; \frac{\alpha_1 + \gamma_1 + n}{2}; \tanh^2 t \right) {\Bigg|}^{{p'}} \ \mathrm{d}t \right)^{\frac{1}{p'}},
		\end{align*}
		By the choice of $\gamma_1$ and $\gamma_2$, we have from Equation \eqref{Behaviour2F1Case1} that ${}_2F_1 \left( - \frac{\alpha_2 + \gamma_2}{2}, \frac{k}{2}; \frac{\alpha_1 + \gamma_1 + n}{2}; \tanh^2 t \right)$ is bounded, and the integral is convergent. This gives the existence {of} $R_kf$ under the {given} assumptions.
		
		To check the existence for $p = 1$, let us choose $\gamma_1 > 0$ such that $\alpha_1 + \gamma_1 > k - n$ and $\gamma_2 = - \gamma_1$. Again, from Equation \eqref{DualityEquation} and Equation \eqref{DualMixedEquation}, we get for any $f \in L^1_{\alpha_1, \alpha_2} \left( \mathbb{H}^n \right)$,
		\begin{align*}
			&\left| \int\limits_{\Xi_k \left( \mathbb{H}^n \right)} R_kf \left( \xi \right) \sinh^{\alpha_1 + \gamma_1} d \left( 0, \xi \right) \cosh^{\alpha_2 - \gamma_1} d \left( 0, \xi \right) \mathrm{d}\xi \right| \\
			&= C \left| \int\limits_{\mathbb{H}^n} f \left( x \right) \sinh^{\alpha_1 + \gamma_1} d \left( 0, x \right) \cosh^{\alpha_2 - \gamma_1} d \left( 0, x \right) {}_2F_1 \left( - \frac{\alpha_2 - \gamma_1}{2}, \frac{k}{2}; \frac{\alpha_1 + \gamma_1 + n}{2}; \tanh^2 d \left( 0, x \right) \right) \mathrm{d}x  \right| \\
			&\leq C \int\limits_{\mathbb{H}^n} \left| f \left( x \right) \right| \sinh^{\alpha_1} d \left( 0, x \right) \cosh^{\alpha_2} d \left( 0, x \right) \tanh^{\gamma_1} d \left( 0, x \right) \left| {}_2F_1 \left( - \frac{\alpha_2 - \gamma_1}{2}, \frac{k}{2}; \frac{\alpha_1 + \gamma_1 + n}{2}; \tanh^2 d \left( 0, x \right) \right) \right| \mathrm{d}x.
		\end{align*}
		It is clear that with the choice of $\gamma_1$ and $\gamma_2$, we have from Equation \eqref{Behaviour2F1Case1} that ${}_2F_1 \left( - \frac{\alpha_2 - \gamma_1}{2}, \frac{k}{2}; \frac{\alpha_1 + \gamma_1 + n}{2}; \tanh^2 d \left( 0, x \right) \right)$ is bounded. Also, $\tanh^{\gamma_1} d \left( 0, x \right) \leq 1$. Therefore, we get
		\begin{align*}
			&\left| \int\limits_{\Xi_k \left( \mathbb{H}^n \right)} R_kf \left( \xi \right) \sinh^{\alpha_1 + \gamma_1} d \left( 0, \xi \right) \cosh^{\alpha_2 - \gamma_1} d \left( 0, \xi \right) \mathrm{d}\xi \right| \leq C \int\limits_{\mathbb{H}^n} \left| f \left( x \right) \right| \sinh^{\alpha_1} d \left( 0, x \right) \cosh^{\alpha_2} d \left( 0, x \right) \mathrm{d}x = C \| f \|_{{L^1_{\alpha_1, \alpha_2} \left( \mathbb{H}^n \right)}} < + \infty.
		\end{align*}
		Hence, we have shown the existence of the $k$-plane transform for functions in $L^p_{\alpha_1, \alpha_2} \left( \mathbb{H}^n \right)$ whenever $\alpha_1 + \alpha_2 > k - n$ and $1 \leq p < \frac{\alpha_1 + \alpha_2 + n - 1}{k - 1}$.
		
		It only remains to check the existence in the case when $\alpha_1 + \alpha_2 = k - n$ and $p = 1$. To see this case, let us choose $\gamma_1 > 0$ such that $\gamma_1 + \alpha_1 > k - n$ and $\gamma_2 < - \gamma_1$. To make things precise, let us fix $\gamma_2 = - \gamma_1 - 1$. {For any $f \in L^1_{\alpha_1, \alpha_2} \left( \mathbb{H}^n \right)$, we have by using the duality relation given in Equation \eqref{DualityEquation}, and Equation \eqref{DualMixedEquation},}
		\begin{align*}
			&\left| \int\limits_{\Xi_k \left( \mathbb{H}^n \right)} R_kf \left( \xi \right) \sinh^{\alpha_1 + \gamma_1} d \left( 0, \xi \right) \cosh^{\alpha_2 - \gamma_1 - 1} d \left( 0, \xi \right) \mathrm{d}\xi \right| \\
			&= C \left| \int\limits_{\mathbb{H}^n} f \left( x \right) \sinh^{\alpha_1 + \gamma_1} d \left( 0, x \right) \cosh^{\alpha_2 - \gamma_1 - 1} d \left( 0, x \right) {}_2F_1 \left( - \frac{\alpha_2 - \gamma_1 - 1}{2}, \frac{k}{2}; \frac{\alpha_1 + \gamma_1 + n}{2}; \tanh^2 d \left( 0, x \right) \right) \mathrm{d}x  \right|.
		\end{align*}
		{For our} choice of parameters, we have from Equation \eqref{Behaviour2F1Case3} that $\lim\limits_{d \left( 0, x \right) \rightarrow \infty} {\frac{{}_2F_1 \left( - \frac{\alpha_2 - \gamma_1 - 1}{2}, \frac{k}{2}; \frac{\alpha_1 + \gamma_1 + n}{2}; \tanh^2 d \left( 0, x \right) \right)}{\cosh d \left( 0, x \right)} = C}$. {Therefore, for sufficiently large $R > 0$, using the fact that hypergeometric function is continuous, hence bounded on compact sets, and that $\gamma_1 > 0$, we get}
		\begin{align*}
			&\left| \int\limits_{\Xi_k \left( \mathbb{H}^n \right)} R_kf \left( \xi \right) \sinh^{\alpha_1 + \gamma_1} d \left( 0, \xi \right) \cosh^{\alpha_2 - \gamma_1 - 1} d \left( 0, \xi \right) \mathrm{d}\xi \right| \\
			&\leq C \int\limits_{B \left( 0, R \right)} \left| f \left( x \right) \right| \sinh^{\alpha_1} d \left( 0, x \right) \cosh^{\alpha_2} d \left( 0, x \right) \tanh^{\gamma_1} d \left( 0, x \right) \mathrm{d}x \\
			&+ C \int\limits_{\mathbb{H}^n \setminus B \left( 0, R \right)} \left| f \left( x \right) \right| \sinh^{\alpha_1} d \left( 0, x \right) \cosh^{\alpha_2} d \left( 0, x \right) \tanh^{\gamma_1} d \left( 0, x \right) \cosh^{-1} d \left( 0, x \right) \cosh d \left( 0, x \right) \mathrm{d}x \\
			&\leq C \int\limits_{\mathbb{H}^n} \left| f \left( x \right) \right| \sinh^{\alpha_1} d \left( 0, x \right) \cosh^{\alpha_2} d \left( 0, x \right) \mathrm{d}x = C \| f \|_{{L^1_{\alpha_1, \alpha_2} \left( \mathbb{H}^n \right)}} < + \infty.
		\end{align*}
		Hence, $R_kf$ is finite for almost every $k$-dimensional totally-geodesic submanifold in $\mathbb{H}^n$.
		\end{proof}
		We now see through a few examples that the constraints on $\alpha_1, \alpha_2$ and $p$ imposed in Theorem \ref{ExistenceMixed} are optimal. That is, beyond these restrictions, we do not have the existence of the $k$-plane transform on weighted Lebesgue spaces, $L^p_{\alpha_1, \alpha_2} \left( \mathbb{H}^n \right)$.
		\begin{example}
			\label{CounterExampleHn1}
			\normalfont
			In this example, we see that the condition $p < \frac{\alpha_1 + \alpha_2 + n - 1}{k - 1}$ imposed in Theorem \ref{ExistenceMixed} is optimal{, when $\alpha_1 + \alpha_2 > k - n$}. {Let assume that $p \geq \frac{\alpha_1 + \alpha_2 + n - 1}{k - 1} > 1$, and consider the function}
		$${f} \left(  x \right) = \frac{\sinh^{\frac{2 - n - \alpha_1}{p}} d \left( 0, x \right) \cosh^{- \frac{\alpha_2}{p}} d \left( 0, x \right)}{\left( 1 + \cosh d \left( 0, x \right) \right)^{\frac{1}{p}} \ln \left( 1 + \cosh d \left( 0, x \right) \right)}.$$
		Using the polar decomposition {on} $\mathbb{H}^n$, we have,
		\begin{align*}
			{\| f \|_{L^p_{\alpha_1, \alpha_2} \left( \mathbb{H}^n \right)}^{p}} &= C \int\limits_{0}^{\infty} \frac{\sinh t}{\left( 1 + \cosh t \right) \ln^p \left( 1 + \cosh t \right)} \ \mathrm{d}t.
		\end{align*}
		Upon substituting $\ln \left( 1 + \cosh t \right) = u$, we get
		\begin{align*}
			{\| f \|_{L^p_{\alpha_1, \alpha_2} \left( \mathbb{H}^n \right)}^p} &= C \int\limits_{\ln 2}^{\infty} \frac{1}{u^p} \ \mathrm{d}u < + \infty,
		\end{align*}
		for every $p > 1$. That is, whenever $p > 1$, we have ${f} \in L^p_{\alpha_1, \alpha_2} \left( \mathbb{H}^n \right)$. For the $k$-plane transform, we have from Equation \eqref{KPlaneTransformRadialHnEquation2},
		\begin{align*}
			R_kf \left( \xi \right) &= \frac{C}{\cosh  d \left( 0, \xi \right)} \int\limits_{d \left( 0, \xi \right)}^{\infty} \frac{\sinh^{\frac{2 - n - \alpha_1}{p}} t \cosh^{- \frac{\alpha_2}{p}} t}{\left( 1 + \cosh t \right)^{\frac{1}{p}} \ln \left( 1 + \cosh t \right)} \left[ 1 - \frac{\tanh^2 d \left( 0, \xi \right)}{\tanh^2 t} \right]^{\frac{k}{2} - 1} \sinh^{k - 1} t \ \mathrm{d}t.
		\end{align*}
		Substitute $\ln \left( 1 + \cosh t \right) = u$. Then, we have,
		\begin{align*}
			R_kf \left( \xi \right) &= \frac{C}{\cosh d \left( 0, \xi \right)} \int\limits_{\ln \left( 1 + \cosh d \left( 0, \xi \right) \right)}^{\infty} \frac{\left( \left( e^u - 1 \right)^2 - 1 \right)^{\frac{2 - n - \alpha_1}{2p} + \frac{k - 2}{2}} \left( e^u - 1 \right)^{- \frac{\alpha_2}{p}}}{e^{\left( \frac{1}{p} - 1 \right) u} u} \left[ 1 - {\frac{\left( e^u - 1 \right)^2}{\left( e^u - 1 \right)^2 - 1}} \tanh^2 d \left( 0, \xi \right) \right]^{\frac{k}{2} - 1} \mathrm{d}{u}.
		\end{align*}
		The total power of $e^u$ towards infinity is $k - 1 - \frac{\alpha_1 + \alpha_2 + n - 1}{p} \geq 0$, for $p \geq \frac{\alpha_1 + \alpha_2 + n - 1}{k - 1}$. Therefore, the integral does not converge and $R_kf \equiv + \infty$ in this case.
		\end{example}
		Example \ref{CounterExampleHn1} suggests that for given $\alpha_1, \alpha_2 \in \mathbb{R}$ with $\alpha_1 + \alpha_2 > k - n$, the restriction on $p$ (i.e., $p < \frac{\alpha_1 + \alpha_2 + n - 1}{k - 1}$) is optimal. We now see whether the restriction on $\alpha_1 + \alpha_2$ is also optimal.
		\begin{example}
			\label{CounterExampleHn2}
			\normalfont
			Let us now fix $\alpha_1 + \alpha_2 \leq k - n$, and consider the function,
			$${f} \left( x \right) = \frac{\sinh^{- \frac{\alpha_1}{p}} d \left( 0, x \right) \cosh^{\frac{1- \alpha_2}{p}} d \left( 0, x \right) \left( 2 + \sinh d \left( 0, x \right) \right)^{- \frac{n}{p}}}{\ln^{1 + \frac{1}{p}} \left( 2 + \sinh d \left( 0, x \right) \right)}.$$
			Then, we have using the polar decomposition on $\mathbb{H}^n$,
			\begin{align*}
				\| {f} \|_{L^p_{\alpha_1, \alpha_2} \left( \mathbb{H}^n \right)}^p &= C \int\limits_{0}^{\infty} \frac{\sinh^{n - 1} t \cosh t}{\left( 2 + \sinh t \right)^n \ln^{p + 1} \left( 2 + \sinh t \right)} \ \mathrm{d}t.
			\end{align*}
			By substituting $\ln \left( 2 + \sinh t \right) = u$, we get
			\begin{align*}
				\| {f} \|_{L^p_{\alpha_1, \alpha_2} \left( \mathbb{H}^n \right)}^p &= C \int\limits_{\ln 2}^{\infty} \frac{\left( e^u - 2 \right)^{n - 1}}{e^{\left( n - 1 \right) u} u^{p + 1}} \mathrm{d}u < + \infty,
			\end{align*}
			for every $p \geq 1$, {and hence} ${f} \in L^p_{\alpha_1, \alpha_2} \left( \mathbb{H}^n \right)$. For the $k$-plane transform, we have from Equation \eqref{KPlaneTransformRadialHnEquation2},
			\begin{align*}
				R_kf \left( \xi \right) &= \frac{C}{\cosh d \left( 0, \xi \right)} \int\limits_{d \left( 0, \xi \right)} \frac{\sinh^{- \frac{\alpha_1}{p}} t \cosh^{\frac{1 - \alpha_2}{p}} t}{\left( 2 + \sinh t \right)^{\frac{n}{p}} \ln^{1 + \frac{1}{p}} \left( 2 + \sinh t \right)} \left[ 1 - \frac{\tanh^2 d \left( 0, \xi \right)}{\tanh^2 t} \right]^{\frac{k}{2} - 1} \sinh^{k - 1} t \ \mathrm{d}t.
			\end{align*}
			By substituting $\ln \left( 2 + \sinh t \right) = u$, we get
			\begin{align*}
				R_kf \left( \xi \right) &= \frac{C}{\cosh d \left( 0, \xi \right)} \int\limits_{\ln \left( 2 + \sinh d \left( 0, \xi \right) \right)}^{\infty} \frac{\left( e^u - 2 \right)^{- \frac{\alpha_1}{p} + k - 1} \left( \left( e^u - 2 \right)^2 + 1 \right)^{\frac{1 - \alpha_2}{2p} - \frac{1}{2}}}{e^{\left( \frac{n}{p} - 1 \right) u} u^{1 + \frac{1}{p}}} \left[ 1 - \left( \frac{\left( e^u - 2 \right)^2 + 1}{\left( e^u - 2 \right)^2} \right) \tanh^2 d \left( 0, \xi \right) \right]^{\frac{k}{2} - 1} \mathrm{d}t.
			\end{align*}
			The total power of $e^u$ towards infinity is $k - 1 - \frac{\alpha_1 + \alpha_2 + n - 1}{p} > 0$ for $\alpha_1 + \alpha_2 < k - n$ and $p \geq 1$ or for $\alpha_1 + \alpha_2 {=} k - n$ and $p > 1$. Therefore, in these cases, we have $R_kf \equiv + \infty$.
		\end{example}
		We notice that {the proof of} Theorem \ref{ExistenceMixed} together with Examples \ref{CounterExampleHn1} and \ref{CounterExampleHn2}, completely answer{s} Question \ref{ExistenceQuestionHn}.
		\begin{remark}
			\normalfont
			Theorem \ref{ExistenceMixed} and Example \ref{CounterExampleHn1} give a natural ``end-point" of existence for the $k$-plane transform. Particularly, we know that for $p = \frac{\alpha_1 + \alpha_2 + n - 1}{k - 1}$, the $k$-plane transform is not well-defined on $L^p_{\alpha_1, \alpha_2} \left( \mathbb{H}^n \right)$. In Section \ref{EndPointSection}, we take up this end-point and prove a Lorentz norm estimate for radial functions.
		\end{remark}
		\subsection{The Sphere \texorpdfstring{$\mathbb{S}^n$}{}}
		We now look at the existence conditions for the totally-geodesic $k$-plane transform on $\mathbb{S}^n$. Here, we consider the weighted Lebesgue spaces $L^p_{\alpha_1, \alpha_2} \left( \mathbb{S}^n \right)$, as given in Equation \eqref{LPX}. Particularly, we ask the following question.
		\begin{question}
			\label{ExistenceQuestionSn}
			For what values of $p \geq 1$, and $\alpha_1, \alpha_2 \in \mathbb{R}$, does the $k$-plane transform {exist} for function{s} in $L^p_{\alpha_1, \alpha_2} \left( \mathbb{S}^n \right)$?
		\end{question}
		To get the required conditions, we use the technique described in Theorem \ref{ExistenceMixed}. We, then have the following result.
		\begin{theorem}
			\label{ExistenceSN}
			For any $\alpha_1, \alpha_2 \in \mathbb{R}$, and $p > 1 + \alpha_2$, the $k$-plane transform is well-defined on the weighted Lebesgue space $L^p_{\alpha_1, \alpha_2} \left( \mathbb{S}^n \right)$. {Moreover}, {f}or $\alpha_2 {=} 0$ and $p {=} 1$, the $k$-plane transform is well-defined on the weighted Lebesgue space $L^p_{\alpha_1, \alpha_2} \left( \mathbb{S}^n \right)$. {Further, if the hypotheses stated above are not satisfied, then there are radial functions $f \in L^p_{\alpha_1, \alpha_2} \left( \mathbb{S}^n \right)$ for which $R_kf \equiv + \infty$.}
		\end{theorem}
		\begin{proof}
			Let $\alpha_1, \alpha_2 \in \mathbb{R}$ be given. We prove this result in two cases.
			
			\textbf{Case I:} We first assume that $p > 1 + \alpha_2$ {and $p \neq 1$}. Let us fix $\gamma_1 > \max \left\lbrace k - n - \alpha_1, - \frac{\alpha_1 + n}{p'} \right\rbrace$ and $\gamma_2 \in \mathbb{R}$ such that $k + \alpha_2 + \gamma_2 \notin 2 \mathbb{Z}$ and $- \frac{1 + \alpha_2}{p'} < \gamma_2 < - \alpha_2$. Then, from the duality of the $k$-plane transform (Equation \eqref{DualityEquation}), and Equation \eqref{DualMixedSphereEquation}, and an application of H\"{o}lder's inequality, we have
			\begin{align}
				&\left| \int\limits_{\Xi_k} R_kf \left( \xi \right) \sin^{\alpha_1 + \gamma_1} d \left( 0, \xi \right) \cos^{\alpha_2 + \gamma_2} d \left( 0, \xi \right) \mathrm{d}\xi \right| \nonumber \\
				&\leq C \int\limits_{\mathbb{S}^n} \left| f \left( x \right) \right| \sin^{\alpha_1 + \gamma_1} d \left( 0, x \right) \left| \cos^{\alpha_2 + \gamma_2} d \left( 0, x \right) \right| \left| {}_2F_1 \left( - \frac{\alpha_2 + \gamma_2}{2}, \frac{k}{2}; \frac{\alpha_1 + \gamma_1 + n}{2}; - \tan^2 d \left( 0, x \right) \right) \right| \mathrm{d}x \nonumber \\
				\label{DualityHolderSn}
				&\leq C \| f \|_{{L^p_{\alpha_1, \alpha_2} \left( \mathbb{S}^n \right)}} I_{\alpha_1, \alpha_2, \gamma_1, \gamma_2}^{\frac{1}{p'}},
			\end{align}
			where,
			\begin{align}
				\label{ISn}
				I_{\alpha_1, \alpha_2, \gamma_1, \gamma_2} &= \int\limits_{\mathbb{S}^n} \sin^{\alpha_1 + p' \gamma_1} d \left( 0, x \right) \left| \cos^{\alpha_2 + p' \gamma_2} d \left( 0, x \right) \right| \left| {}_2F_1 \left( - \frac{\alpha_2 + \gamma_2}{2}, \frac{k}{2}; \frac{\alpha_1 + \gamma_1 + n}{2}; - \tan^2 d \left( 0, x \right) \right) \right|^{p'} \mathrm{d}x \\
				&= C \int\limits_{0}^{\frac{\pi}{2}} \sin^{\alpha_1 + p'\gamma_1 + n - 1} t \cos^{\alpha_2 + p'\gamma_2}  t \left| {}_2F_1 \left( - \frac{\alpha_2 + \gamma_2}{2}, \frac{k}{2}; \frac{\alpha_1 + \gamma_1 + n}{2}; - \tan^2 t \right) \right|^{p'} \mathrm{d}t \nonumber \\
				&= C \left( I_1 + I_2 \right). \nonumber
			\end{align}
			We have used the polar decomposition for $\mathbb{S}^n$ in the penultimate step{. Here,}
			$$I_1 = \int\limits_{0}^{{\frac{\pi}{4}}} \sin^{\alpha_1 + p'\gamma_1 + n - 1} t \cos^{\alpha_2 + p'\gamma_2} t \left| {}_2F_1 \left( - \frac{\alpha_2 + \gamma_2}{2}, \frac{k}{2}; \frac{\alpha_1 + \gamma_1 + n}{2}; - \tan^2 {t} \right) \right|^{p'} \mathrm{d}t,$$ 
			and 
			$$I_2 = \int\limits_{{\frac{\pi}{4}}}^{\frac{\pi}{2}} \sin^{\alpha_1 + p'\gamma_1 + n - 1} t \cos^{\alpha_2 + p'\gamma_2} t \left| {}_2F_1 \left( - \frac{\alpha_2 + \gamma_2}{2}, \frac{k}{2}; \frac{\alpha_1 + \gamma_1 + n}{2}; - \tan^2 {t} \right) \right|^{p'} \mathrm{d}t.$$
			The integral $I_1$ is finite due to the choice of $\gamma_1${, and the fact that the cosine and hypergeometric terms are bounded when $t \in \left( 0, \frac{\pi}{4} \right)$}. For the integral $I_2$, we need to {consider} the behaviour of the hypergeometric function. From Equation \eqref{Transformation2F12}, we get
			\begin{align*}
				&{}_2F_1 \left( - \frac{\alpha_2 + \gamma_2}{2}, \frac{k}{2}; \frac{{\alpha_1 +} \gamma_1 + n}{2}; - \tan^2 t \right) \\
				&= \frac{\pi}{\sin \left( k + \alpha_2 + \gamma_2 \right) \frac{\pi}{2}} \left[ \frac{\tan^{\alpha_2 + \gamma_2} t}{\Gamma \left( \frac{k}{2} \right) \Gamma \left( \frac{\alpha_1 + \alpha_2 + \gamma_1 + \gamma_2 + n}{2} \right)} {}_2F_1 \left( - \frac{\alpha_2 + \gamma_2}{2}, {1} - \frac{\alpha_1 + \alpha_2 + \gamma_1 + \gamma_2 + n}{2}; {1} - \frac{\alpha_2 + \gamma_2 + k}{2}; - \cot^2 t \right) \right. \\
				&\left. - \frac{\tan^{-k} t}{\Gamma \left( - \frac{\alpha_2 + \gamma_2}{2} \right) \Gamma \left( \frac{\alpha_1 + \gamma_1 + n - k}{2} \right)} {}_2F_1 \left( \frac{k}{2}, {1 - \frac{\alpha_1 + \gamma_1 + n - k}{2}}; {1 + } \frac{\alpha_2 + \gamma_2 + k}{2}; - \cot^2 t \right) \right].
			\end{align*}
			Since, $\gamma_2 < - \alpha_2$, we have that the first term in this expression is bounded towards $\frac{\pi}{2}$. On the other hand, the second term is anyway bounded, so that overall, the hypergeometric function is bounded. The power of $\cos$ in the integral $I_2$ is integrable since $\gamma_2 > - \frac{1 + \alpha_2}{p'}$. This completes the proof for existence in this case.
			
			\textbf{Case II:} We now consider $\alpha_2 = 0$ and $p = 1$. First, we notice that by taking $\gamma_2 = 0$ in {Equation \eqref{DualMixedSphereEquation} of} Theorem \ref{DualMixedSphere}, we get
			\begin{equation}
				\label{DualSineSphere}
				R_k^* \left( \sin^{\gamma_1} d \left( 0, \cdot \right) \right) \left( x \right) = C \sin^{\gamma_1} d \left( 0, x \right),
			\end{equation}
			provided $\gamma_1 > k - n$. Therefore, for a given $\alpha_1 \in \mathbb{R}$, let us choose $\gamma_1 \geq 0$ such that $\alpha_1 + \gamma_1 > k - n$. Then, we have, for functions $f \in L^1_{\alpha_1, 0} \left( \mathbb{S}^n \right)$, using Equations \eqref{DualityEquation} and \eqref{DualSineSphere},
			\begin{align*}
				\left| {\int\limits_{\Xi_k \left( \mathbb{S}^n \right)}} R_kf \left( \xi \right) \sin^{\alpha_1 + \gamma_1} d \left( 0, \xi \right) \mathrm{d}\xi \right| &\leq C \int\limits_{\mathbb{S}^n} \left| f \left( x \right) \right| \sin^{\alpha_1 + \gamma_1} d \left( 0, x \right) \mathrm{d}x \leq C \| f \|_{L^1_{\alpha_1, 0} \left( \mathbb{S}^n \right)} < + \infty.
			\end{align*}
			{It only remains to be checked that when $\alpha_2 < 0$ and $p = 1 > 1 + \alpha_2$, the $k$-plane transform of functions in $L^1_{\alpha_1, \alpha_2} \left( \mathbb{S}^n \right)$ exists. However, we notice that for $\alpha_2 < 0$, we have for any $t \in \left( 0, \frac{\pi}{2} \right)$, $\cos^{\alpha_2} t \geq 1$. Consequently,
			$$\| f \|_{L^1_{\alpha_1, 0} \left( \mathbb{S}^n \right)} = \int\limits_{\mathbb{S}^n} \left| f \left( x \right) \right| \sin^{\alpha_1} d \left( 0, x \right) \mathrm{d}x \leq \int\limits_{\mathbb{S}^n} \left| f \left( x \right) \right| \sin^{\alpha_1} d \left( 0, x \right) \cos^{\alpha_2} d \left( 0, x \right) \mathrm{d}x = \| f \|_{L^1_{\alpha_1, \alpha_2} \left( \mathbb{S}^n \right)},$$
			for any $\alpha_1 \in \mathbb{R}$. That is, for $\alpha_2 < 0$, we have $L^1_{\alpha_1, \alpha_2} \left( \mathbb{S}^n \right) \subseteq L^1_{\alpha_1, 0} \left( \mathbb{S}^n \right)$. In Case II, we have seen that for functions in $L^1_{\alpha_1, 0} \left( \mathbb{S}^n \right)$, the $k$-plane transform is well-defined. Therefore, it is also defined on $L^1_{\alpha_1, \alpha_2} \left( \mathbb{S}^n \right)$, when $\alpha_2 < 0$. This completes the proof!}
		\end{proof}
		Next, we show that for $\alpha_2 > 0$ and $1 \leq p \leq 1 + \alpha_2$, there are radial functions in $L^p_{\alpha_1, \alpha_2} \left( \mathbb{S}^n \right)$ for which $R_kf \equiv + \infty$.
		\begin{example}
			\label{CounterExampleSn}
			\normalfont
			Let us fix $\alpha_2 > 0$ and $1 \leq p < 1 + \alpha_2$. We define a function ${f}: \mathbb{S}^n \rightarrow \mathbb{C}$, as
			$${f} \left( x \right) = {\frac{\sin^{- \frac{\alpha_1}{p}} d \left( 0, x \right)}{\cos d \left( 0, x \right)}}.$$
			Then, {by using the polar decomposition on $\mathbb{S}^n$,} we have,
			\begin{align*}
				\| {f} \|_{L^{p}_{\alpha_1, \alpha_2} \left( \mathbb{S}^n \right)}^{{p}}&= C \int\limits_{0}^{\frac{\pi}{2}} \cos^{\alpha_2 - p} t \sin^{n - 1} t \ \mathrm{d}t.
			\end{align*}
			By substituting $\cos^2 t = u$, we get
			\begin{align*}
				\| {f} \|_{L^{p}_{\alpha_1, \alpha_2} \left( \mathbb{S}^n \right)}^{{p}} &= C \int\limits_{0}^{1} u^{\frac{\alpha_2 - p - 1}{2}} \left( 1 - u \right)^{\frac{n}{2} - 1} \mathrm{d}t < + \infty,
			\end{align*}
			since $p < 1 + \alpha_2$. That is, ${f} \in L^p_{\alpha_1, \alpha_2} \left( \mathbb{S}^n \right)$. {Further,} for the $k$-plane transform {of $f$ at $\xi \in \Xi_k \left( \mathbb{S}^n \right)$}, we have from Equation \eqref{KPlaneTransformRadialSnEquation},
			\begin{align*}
				R_kf \left( \xi \right) &= \frac{C}{\cos^{k - 1} d \left( 0, \xi \right)} \int\limits_{d \left( 0, \xi \right)}^{\frac{\pi}{2}} {\frac{\sin^{- \frac{\alpha_1}{p}} t}{\cos t}} \left( \cos^2 d \left( 0, \xi \right) - \cos^{2} t \right)^{\frac{k}{2} - 1} \sin t \ \mathrm{d}t = +\infty,
			\end{align*}
			since ${\frac{1}{\cos t}}$ is not integrable near $\frac{\pi}{2}$.
		\end{example}
		Example \ref{CounterExampleSn} gives us that the restriction $p \geq 1 + \alpha_2$ is optimal, for any $\alpha_2 {>} 0$. However, Theorem \ref{ExistenceSN} has a stronger restriction on $p$, when $\alpha_2 > 0$. The next example shows that for $\alpha_2 > 0$ and $p = 1 + \alpha_2$, there are radial functions for which the $k$-plane transform is not well-defined. Particularly, the restriction $p > 1 + \alpha_2$ is optimal.
		\begin{example}
			\label{CounterExampleSn2}
			\normalfont
			Let us fix $\alpha_2 > 0$ and $p = 1 + \alpha_2$. We consider the function ${f} {:} \mathbb{S}^n \rightarrow \mathbb{C}$, defined as
			$${f} \left( x \right) = {\frac{\sin^{- \frac{\alpha_1}{p}} d \left( 0, x \right)}{\cos d \left( 0, x \right)}} \left( \ln \frac{1}{\cos d \left( 0, x \right)} \right)^{- \gamma} \chi_{\left( \frac{\pi}{4}, \frac{\pi}{2} \right)} \left( d \left( 0, x \right) \right),$$
			where, $\frac{1}{p} < \gamma \leq 1$. We notice that such a choice of $\gamma$ is possible since $\alpha_2 > 0$.
			We have,
			\begin{align*}
				\| {f} \|_{L^p_{\alpha_1, \alpha_2} \left( \mathbb{S}^n \right)}^p = C \int\limits_{\frac{\pi}{4}}^{\frac{\pi}{2}} \cos^{\alpha_2 - p} t \sin^{n - 1} t \left( \ln \frac{1}{\cos t} \right)^{- \gamma p} \mathrm{d}t = C \int\limits_{\frac{\pi}{4}}^{\frac{\pi}{2}} \frac{\sin^{n - 1} t}{\cos t} \left( \ln \frac{1}{\cos t} \right)^{- \gamma {p}} \mathrm{d}t.
			\end{align*}
			By substituting $\cos^2 t = e^{-v}$, we get
			\begin{align*}
				\| {f} \|_{L^p_{\alpha_1, \alpha_2} \left( \mathbb{S}^n \right)}^p = C \int\limits_{\ln 2}^{\infty} \left( 1 - e^{-v} \right)^{\frac{n}{2} - 1} v^{- p \gamma} \mathrm{d}v < + \infty,
			\end{align*}
			since $\gamma p > 1$. On the other hand, from Equation \eqref{KPlaneTransformRadialSnEquation}, we get for any $\xi \in \Xi_k \left( \mathbb{S}^n \right)$,
			\begin{align*}
				R_kf \left( \xi \right) = \frac{C}{\cos^{k - 1} d \left( 0, \xi \right)} \int\limits_{\max \left\lbrace \frac{\pi}{4}, d \left( 0, \xi \right) \right\rbrace}^{\frac{\pi}{2}} \sin^{1 - \frac{\alpha_1}{p}} t \left( \cos^2 d \left( 0, \xi \right) - \cos^2 t \right)^{\frac{k}{2} - 1} \frac{1}{\cos t} \left( \ln \frac{1}{\cos t} \right)^{- \gamma} \mathrm{d}t.
			\end{align*}
			Substituting $\ln \frac{1}{\cos t} = u$, we get
			\begin{align*}
				R_kf \left( \xi \right) &= \frac{C}{\cos^{k - 1} d \left( 0, \xi \right)} \int\limits_{\max \left\lbrace \ln 2, \ln \frac{1}{\cos d \left( 0, \xi \right)} \right\rbrace}^{\infty} \left( 1 - e^{-2u} \right)^{- \frac{\alpha_1}{2p}} \left( \cos^2 d \left( 0, \xi \right) - e^{-2u} \right)^{\frac{k}{2} - 1} u^{- \gamma} \mathrm{d}u = + \infty,
			\end{align*}
			since $\gamma \leq 1$.
		\end{example}
		We notice that {the proof of} Theorem \ref{ExistenceSN} together with Examples \ref{CounterExampleSn} and \ref{CounterExampleSn2}, completely answer{s} Question \ref{ExistenceQuestionSn}.
		\begin{remark}
			\normalfont
			We would like to mention that the existence conditions on the sphere give a natural end-point, namely $p = 1 + \alpha_2$, when $\alpha_2 > 0$. We consider this end-point in Section \ref{EndPointSection}, and obtain boundedness of the $k$-plane transform (for radial functions) on Lorentz space, $L^{p, 1}_{\alpha_1, \alpha_2} \left( \mathbb{S}^n \right)$.
		\end{remark}
	\section{ \texorpdfstring{$L^p$}{Lp}-improving estimates for radial functions}
		\label{RadialFunctionSection}
		{In this section, we find necessary and sufficient conditions for the $L^p$-improving boundedness of the $k$-plane transform on weighted Lebesgue spaces of radial functions.} This section is divided into three subsections, each dealing with the $k$-plane transform on the Euclidean space, Hyperbolic space, and the {S}phere, respectively.
		\subsection{The Euclidean space \texorpdfstring{$\mathbb{R}^n$}{}}
			\label{LPLInfinityEuclideanSection}
			The {$L^p$-improving} boundedness of the Euclidean $k$-plane transform has been a subject of study for a long time, and much theory is developed around it. Some of the known results are mentioned in Section \ref{IntroductionSection}. The motivation of our study in this section lies in the works of Kumar and Ray (see \cite{KumarRayWE}), where they proved weighted $L^p$-improving estimates for the $k$-plane transform of radial functions on the Euclidean space. In their paper, the authors give necessary and sufficient conditions for the boundedness. We move a step further and ask the following question:
			\begin{question}
				\label{EuclideanLPLInfinityQuestion}
				{Given $\alpha > k - n$ and $1 \leq p < \frac{\alpha + n}{k}$} {or $\alpha = k - n$ and $p = 1$}, {w}hat are the necessary and sufficient conditions on $\beta \in \mathbb{R}$ such that the following inequality holds for all radial functions?
				\begin{equation}
					\label{RequiredEuclideanInequalityLPLInfinity}
					\| \left( d \left( 0, \cdot \right) \right)^{\beta} R_kf \|_{\infty} \leq C \| f \|_{{L^p_{\alpha} \left( \mathbb{R}^n \right)}}.
				\end{equation}
			\end{question}
			We know (see for instance \cite{RubinKPlane}) that the $k$-plane transform exists on ${L^{p}_{\alpha} \left( \mathbb{R}^n \right)}$, {exactly when} $\alpha > k - n$ and $1 \leq p < \frac{\alpha + n}{k}$ {or $\alpha = k - n$ and $p = 1$}.
			
			First, let us obtain necessary conditions for Inequality \eqref{RequiredEuclideanInequalityLPLInfinity} to hold.
			\begin{theorem}
				\label{NecessaryLPLInfinityRn}
				Let $\alpha > k - n$ and $1 \leq p < \frac{\alpha + n}{k}$ {or $\alpha = k - n$ and $p = 1$}. Then, Inequality \eqref{RequiredEuclideanInequalityLPLInfinity} holds only if $\beta = \frac{\alpha + n}{p} - k$.
			\end{theorem}
			\begin{proof}
				For a fixed $\lambda > 0$, let us consider the function $f_{\lambda} = \chi_{B \left( 0, \lambda \right)}$. Then, we have, using the polar decomposition of $\mathbb{R}^n$,
			\begin{align*}
				\| f \|_{L^p_{\alpha} \left( \mathbb{R}^n \right)} &= C \left[ \int\limits_{0}^{\lambda} t^{\alpha + n - 1} \ \mathrm{d}t \right]^{\frac{1}{p}} = C \lambda^{\frac{\alpha + n}{p}}.
			\end{align*}
			The $k$-plane transform of $f_{\lambda}$ is given in Equation \eqref{LPMassBallRnEquation}. By choosing $d \left( 0, \xi \right) = \frac{1}{\sqrt{2}} \lambda$, it is easy to see that
			$$\| \left( d \left( 0, \cdot \right) \right)^{\beta} R_kf \|_{\infty} \geq C \lambda^{\beta + k}.$$
			Therefore, a necessary condition for Inequality \eqref{RequiredEuclideanInequalityLPLInfinity} to hold is
			$$\lambda^{\beta + k} \leq C \lambda^{\frac{\alpha + n}{p}},$$
			for all $\lambda > 0$. This is possible only if $\beta = \frac{\alpha + n}{p} - k$.
			\end{proof}
			We now see that the necessary condition of Theorem \ref{NecessaryLPLInfinityRn} is also sufficient.
			\begin{theorem}
				\label{SufficientLPLInifityRn}
				Let $\alpha > k - n$ and $1 \leq p < \frac{\alpha + n}{k}$.
				\begin{enumerate}
					\item[\mylabel{WeightLPLRRnA}{(A)}] For $k \geq 2$, Inequality \eqref{RequiredEuclideanInequalityLPLInfinity} holds if and only if $\beta = \frac{\alpha + n}{p} - k$.
					\item[\mylabel{WeightLPLRRnB}{(B)}] {For $k = 1$, Inequality \eqref{RequiredEuclideanInequalityLPLInfinity} holds if and only if $\beta = \frac{\alpha + n}{p} - k$ and $p > 2$.}
				\end{enumerate}
			\end{theorem}
			\begin{proof}
				The necessity of the condition $\beta = \frac{\alpha + n}{p} - k$ was proved in Theorem \ref{NecessaryLPLInfinityRn}. We now check the sufficiency.
				\begin{enumerate}
					\item[\mylabel{WeightLPLRRnAProof}{(A)}] First, we consider the case $k \geq 2$ and $p > 1$, which is fairly simple. Here, we have, from Equation \eqref{KPlaneTransformRadialRnEquation},
					\begin{align*}
						\left| \left( d \left( 0, \xi \right) \right)^{\frac{\alpha + n}{p} - k} R_kf \left( \xi \right) \right| &\leq C \left( d \left( 0, \xi \right) \right)^{\frac{\alpha + n}{p} - k} \int\limits_{d \left( 0, \xi \right)}^{\infty} \left| \tilde{f} \left( t \right) \right| \left( t^2 - \left( d \left( 0, \xi \right) \right)^2 \right)^{\frac{k}{2} - 1} t \ \mathrm{d}t \\
						&= C \left( d \left( 0, \xi \right) \right)^{\frac{\alpha + n}{p} - k} \int\limits_{d \left( 0, \xi \right)}^{\infty} \left| \tilde{f} \left( t \right) \right| t^{\frac{\alpha + n - 1}{p}} \left( 1 - \left( \frac{d \left( 0, \xi \right)}{t} \right)^2 \right)^{\frac{k}{2} - 1} t^{k - \frac{\alpha + n}{p} - \frac{1}{p'}} \mathrm{d}t.
					\end{align*}
					Using H\"{o}lder's inequality {and Equation \eqref{LpNormRadialRn1}}, we get
					\begin{align*}
						\left| \left( d \left( 0, \xi \right) \right)^{\frac{\alpha + n}{p} - k} R_kf \left( \xi \right) \right| &\leq C \left( d \left( 0, \xi \right) \right)^{\frac{\alpha + n}{p} - k} \| f \|_{L^p_{\alpha} \left( \mathbb{R}^n \right)} I^{\frac{1}{p'}},
					\end{align*}
					where,
					\begin{align*}
						I &= \int\limits_{d \left( 0, \xi \right)}^{\infty} t^{p' \left( k - \frac{\alpha + n}{p} \right) - 1} \left( 1 - \left( \frac{d \left( 0, \xi \right)}{t} \right)^2 \right)^{p' \left( \frac{k}{2} - 1 \right)} \mathrm{d}t.
					\end{align*}
					By substituting $\frac{d \left( 0, \xi \right)}{t} = u$, we obtain
					\begin{equation}
						\label{IntermediateIEuclideanLPLInfinity}
						I = \left( d \left( 0, \xi \right) \right)^{p' \left( k - \frac{\alpha + n}{p} \right)} \int\limits_{0}^{1} u^{p' \left( \frac{\alpha + n}{p} - k \right) - 1} \left( 1 - u^2 \right)^{p' \left( \frac{k}{2} - 1 \right)} \mathrm{d}u = C \left( d \left( 0, \xi \right) \right)^{p' \left( k - \frac{\alpha + n}{p} \right)}.
					\end{equation}
					In the above equation, $C = \frac{\Gamma \left( \frac{p'}{2} \left( \frac{\alpha + n}{p} - k \right) \right) \Gamma \left( p' \left( \frac{k}{2} - 1 \right) + 1 \right)}{\Gamma \left( p' \left( \frac{\alpha + n}{2p} - 1 \right) + 1 \right)}$. Hence, for $k \geq 2$ and $p > 1$, we have Inequality \eqref{RequiredEuclideanInequalityLPLInfinity}.
					Now, let us consider $p = 1$. From Equation \eqref{KPlaneTransformRadialRnEquation2} and the fact that $\alpha {\geq} k - n$, we have the following.
					\begin{align*}
						\left| \left( d \left( 0, \xi \right) \right)^{{\alpha + n} - 1} R_kf \left( \xi \right) \right| &\leq C \left( d \left( 0, \xi \right) \right)^{\alpha + n - k} \int\limits_{d \left( 0, \xi \right)}^{\infty} \left| \tilde{f} \left( t \right) \right| t^{k - 1} \left( 1 - \left( \frac{d \left( 0, \xi \right)}{t} \right)^2 \right)^{\frac{k}{2} - 1} \mathrm{d}t \\
						&\leq C \int\limits_{d \left( 0, \xi \right)}^{\infty} \left| \tilde{f} \left( t \right) \right| t^{\alpha + n - 1} \mathrm{d}t \leq C \| f \|_{L^1_{\alpha} \left( \mathbb{R}^n \right)}.
					\end{align*}
					\item[\mylabel{WeightLPLRRnBProof}{(B)}] For the case $k = 1$ with $p > 2$, the proof is verbatim to {that of \ref{WeightLPLRRnAProof}}, {where the condition that $p > 2$ is used to make the integral $I$, given in Equation \eqref{IntermediateIEuclideanLPLInfinity} finite. The necessity of $p > 2$ when $k = 1$ is given in Example \ref{EuclideanAnnulusCounterExample} and Remark \ref{RemarkL2LInfinityRn}.}
				\end{enumerate}
			\end{proof}
			{Now we show that} for $k = 1$ and $p \leq 2$, we cannot expect Inequality \eqref{RequiredEuclideanInequalityLPLInfinity} to hold. We begin with the case $p < 2$, where we produce a counter-example.
			{\begin{example}
				\label{EuclideanAnnulusCounterExample}
				\normalfont
				Let $k = 1$ and $p < 2$. For $0 < a < b < \infty$, let us consider the characteristic function of an annulus $A := \left\lbrace x \in \mathbb{R}^n | a < \| x \| < b \right\rbrace$. Then, we have,
				\begin{equation}
					\label{AnnulusMeasure}
					\| \chi_A \|_{L^{p}_{\alpha} \left( \mathbb{R}^n \right)} = \left[ \int\limits_{a}^{b} r^{\alpha + n - 1} \mathrm{d}r \right]^{\frac{1}{p}} = C \left( b^{\alpha + n} - a^{\alpha + n} \right)^{\frac{1}{p}}.
				\end{equation}
				As for the X-ray transform of $\chi_A$, we have using Equation \eqref{KPlaneTransformRadialRnEquation},
				\begin{equation}
					\label{XRayAnnulus}
					\begin{aligned}
						R_1{\chi_A} \left( \xi \right) = C \int\limits_{\max \left\lbrace a, d \left( 0, \xi \right) \right\rbrace}^{b} \left( t^2 - \left( d \left( 0, \xi \right) \right)^2 \right)^{- \frac{1}{2}} t \ \mathrm{d}t &= C \left[ \left( b^2 - \left( d \left( 0, \xi \right) \right)^2 \right)^{\frac{1}{2}} - \left( \max \left\lbrace a^2, \left( d \left( 0, \xi \right) \right)^2 \right\rbrace - \left( d \left( 0, \xi \right) \right)^2 \right)^{\frac{1}{2}} \right].
					\end{aligned}
				\end{equation}
				Using $d \left( 0, \xi \right) = a$ in Equation \eqref{XRayAnnulus}, we easily see that
				\begin{equation}
					\label{XRayAnnulusEstimate}
					\| \left( d \left( 0, \cdot \right) \right)^{\beta} R_1\chi_A \|_{L^{\infty} \left( \mathbb{R}^n \right)} \geq C a^{\beta} \left( b^2 - a^2 \right)^{\frac{1}{2}}.
				\end{equation}
				Now, testing Inequality \eqref{RequiredEuclideanInequalityLPLInfinity} against $\chi_A$, we get from Equations \eqref{AnnulusMeasure} and \eqref{XRayAnnulusEstimate},
				\begin{equation}
					\label{RequiredXRayAnnulusEstimateRn}
					a^{\beta} \left( b^2 - a^2 \right)^{\frac{1}{2}} \leq C \left( b^{\alpha + n} - a^{\alpha + n} \right)^{\frac{1}{p}},
				\end{equation}
				for $0 < a < b < \infty$. For a fixed $a > 0$, let us consider the function $H: \left( a, \infty \right) \rightarrow \mathbb{R}$, defined as
				$$H \left( b \right) = \frac{a^{p \beta} \left( b^2 - a^2 \right)^{\frac{p}{2}}}{b^{\alpha + n} - a^{\alpha + n}}.$$
				Inequality \eqref{RequiredXRayAnnulusEstimateRn} is equivalent to saying that $H$ is bounded. However, we see that since $p < 2$, we have,
				$$\lim\limits_{b \rightarrow a} H \left( b \right) = \lim\limits_{b \rightarrow a} \frac{a^{p \beta} \left( b^2 - a^2 \right)^{\frac{p}{2}}}{b^{\alpha + n} - a^{\alpha + n}} = \lim\limits_{b \rightarrow a} \frac{a^{p \beta} p b \left( b^2 - a^2 \right)^{\frac{p}{2} - 1}}{\left( \alpha + n \right) b^{\alpha + n - 1}} = + \infty.$$
				Therefore, Inequality \eqref{RequiredEuclideanInequalityLPLInfinity} cannot be expected for $p < 2$ and $k = 1$.
			\end{example}}
			\begin{remark}
				\label{RemarkL2LInfinityRn}
				\normalfont
				Example \ref{EuclideanAnnulusCounterExample} shows that for Inequality \eqref{RequiredEuclideanInequalityLPLInfinity} to hold {with} $k = 1$, we must have $p \geq 2$.	Let us now see that we cannot expect Inequality \eqref{RequiredEuclideanInequalityLPLInfinity} for $p = 2$ when $k = 1$. To see this, we first observe from Equation \eqref{KPlaneTransformRadialRnEquation}, by replacing $t^2$ by $t$ {and Equation \eqref{InfiniteUpperRLFractionalIntegral}}, that
				$$R_1f \left( \xi \right) = I^{\frac{1}{2}}_{-}\varphi \left( \left( d \left( 0, \xi \right) \right)^2 \right),$$
				where,
				$$\varphi \left( t \right) = \tilde{f} \left( \sqrt{t} \right).$$
				Also, it is easy to see that
				$$\| f \|_{L^p_{\alpha} \left( \mathbb{R}^n \right)} = C \left[ \int\limits_{0}^{\infty} \left| \tilde{f} \left( t \right) \right|^p t^{\alpha + n - 1} \mathrm{d}t \right]^{\frac{1}{p}} = C \left[ \int\limits_{0}^{\infty} \left| \tilde{f} \left( \sqrt{t} \right) \right|^p t^{\frac{\alpha + n}{2} - 1} \mathrm{d}t \right]^{\frac{1}{p}} = C \| \varphi \|_{L^p_{\frac{\alpha + n}{2} - 1} \left( 0, \infty \right)}.$$
				Hence, proving Inequality \eqref{RequiredEuclideanInequalityLPLInfinity} for $p = 2$ and $k = 1$ is equivalent to proving
				\begin{equation}
					\label{RequiredFIInequalityEuclidean}
					\| x^{\frac{\alpha + n}{4} - \frac{1}{2}} I^{\frac{1}{2}}_{-}\varphi \|_{L^{\infty} \left( 0, \infty \right)} \leq C \| \varphi \|_{L^2_{\frac{\alpha + n}{2} - 1} \left( 0, \infty \right)}.
				\end{equation}
				{However, Inequality \eqref{RequiredFIInequalityEuclidean} is not valid as seen in Proposition \ref{L2LInfinityIHalfPlusNotPossible}.}
			\end{remark}
			The analysis done here completely answers Question \ref{EuclideanLPLInfinityQuestion}. We now move on to consider $L^p$-improving mapping properties for the $k$-plane transform of radial functions on the hyperbolic space.
		\subsection{The Hyperbolic Space \texorpdfstring{$\mathbb{H}^n$}{}}
			In this section, we get necessary and sufficient conditions for the totally-geodesic $k$-plane transform on $\mathbb{H}^n$ to be {an $L^p$-improving} bounded operator {on} weighted Lebesgue spaces of radial functions {on $\mathbb{H}^n$}. That is, we ask the following question.
			\begin{question}
				\label{LPLRHnQuestion}
				What are the admissible values of $\alpha_1, \alpha_2, \beta_1, \beta_2 \in \mathbb{R}$, and ${r \geq p} \geq 1$ such that the following inequality holds for all radial functions in $L^p_{\alpha_1, \alpha_2} \left( \mathbb{H}^n \right)$?
				\begin{equation}
					\label{LPLRRequiredHn}
					\| R_kf \|_{L^r_{\beta_1, \beta_2} \left( \Xi_k \left( \mathbb{H}^n \right) \right)} \leq C \| f \|_{L^p_{\alpha_1, \alpha_2} \left( \mathbb{H}^n \right)}.
				\end{equation}
			\end{question}
			Further, we also ask the following question.
			\begin{question}
				\label{LPLInfinityHnQuestion}
				What are the admissible values of $\alpha_1, \alpha_2, \gamma_1, \gamma_2 \in \mathbb{R}$, and $p \geq 1$ such that the following inequality holds for all radial functions $f \in L^p_{\alpha_1, \alpha_2} \left( \mathbb{H}^n \right)$?
				\begin{equation}
					\label{RequiredLPLInfinityHn}
					\| \sinh^{\gamma_1} d \left( 0, \cdot \right) \cosh^{\gamma_2} d \left( 0, \cdot \right) R_kf \|_{L^{\infty} \left( \Xi_k \left( \mathbb{H}^n \right) \right)} \leq C \| f \|_{L^p_{\alpha_1, \alpha_2} \left( \mathbb{H}^n \right)}.
				\end{equation}
			\end{question}
			{First, we recall from Theorem \ref{ExistenceMixed} that we must have $\alpha_1 + \alpha_2 > k - n$ and $1 \leq p < \frac{\alpha_1 + \alpha_2 + n - 1}{k - 1}$ or $\alpha_1 + \alpha_2 = k - n$ and $p = 1$ for the existence of the $k$-plane transform. {Without} these conditions, Inequalities \eqref{LPLRRequiredHn} and \eqref{RequiredLPLInfinityHn} do not make sense.}
			
			We begin with answering Question \ref{LPLRHnQuestion}. First, we obtain necessary conditions for Inequality \eqref{LPLRRequiredHn} to hold. To do so, we test the desired {in}equality against certain functions that mimic dilations.
			\begin{theorem}
				\label{MixedWeightNecessary}
				{Let} $\alpha_1 + \alpha_2 > k - n$ {and} $1 \leq p < \frac{\alpha_1 + \alpha_2 + n - 1}{k - 1}$ {or $\alpha_1 + \alpha_2 = k - n$ and $p = 1$. Then}, Inequality \eqref{LPLRRequiredHn} holds only if
				\begin{equation}
					\label{NecessaryMixedWeightInequality}
					{\beta_1 > k - n \text{ and }} \frac{\beta_2 + k - 1}{r} - \frac{\alpha_2 - 1}{p} \leq \frac{\alpha_1 + n}{p} - \frac{\beta_1 + n - k}{r} \leq k.
				\end{equation}
			\end{theorem}
			\begin{proof}
				We prove the necessary condition stated in the statement of this theorem in two cases.
				
				\textbf{Case I:} Let us first assume that $\alpha_1 > - n$. Here, we consider, for $\lambda > 0$, the function $f_{\lambda} \left( x \right) = \chi_{B \left( 0, \lambda \right)} \left( x \right) \cosh d \left( 0, x \right)$. Then, we have, from Equation {\eqref{LpNormRadialHn1}},
				\begin{equation}
					\label{LPMassHn1}
					\| f_{\lambda} \|_{L^p_{\alpha_1, \alpha_2} \left( \mathbb{H}^n \right)} = C \left[ \int\limits_{0}^{\lambda} \sinh^{\alpha_1 + n - 1} t \cosh^{\alpha_2 + p} t \ \mathrm{d}t \right]^{\frac{1}{p}} \leq C \begin{cases}
				\lambda^{\frac{\alpha_1 + n}{p}}, & \lambda \text{ is small}. \\
				e^{\left( \frac{\alpha_1 + \alpha_2 + n - 1}{p} + 1 \right) \lambda}, & \lambda \text{ is large}.
			\end{cases}
				\end{equation}
				On the other hand, from Equation {\eqref{NormRkHn}}, we have
				\begin{align*}
					\| R_kf_{\lambda} \|_{L^r_{\beta_1, \beta_2} \left( \Xi_k \left( \mathbb{H}^n \right) \right)} & = C \left[ \int\limits_{1}^{\infty} \left| \int\limits_{s}^{\infty} \frac{\chi_{\left( 1, \cosh \lambda \right)} \left( t \right) t \left( t^2 - s^2 \right)^{\frac{k}{2} - 1}}{s^{k - 1}} \ \mathrm{d}t \right|^r s^{\beta_2 + k} \left( s^2 - 1 \right)^{\frac{\beta_1 + n - k}{2} - 1} \mathrm{d}s \right]^{\frac{1}{r}} \\
					&= C \left[ \int\limits_{1}^{\cosh \lambda} \left| \int\limits_{s}^{\cosh \lambda} t \left( t^2 - s^2 \right)^{\frac{k}{2} - 1} \mathrm{d}t \right|^r s^{\beta_2 + k - r \left( k - 1 \right)} \left( s^2 - 1 \right)^{\frac{\beta_1 + n - k}{2} - 1} \mathrm{d}s \right]^{\frac{1}{r}} \\
					&= C \left[ \int\limits_{1}^{\cosh \lambda} \left( \cosh^2 \lambda - s^2 \right)^{\frac{rk}{2}} s^{\beta_2 + k - r \left( k - 1 \right)} \left( s^2 - 1 \right)^{\frac{\beta_1 + n - k}{2} - 1} \mathrm{d}s \right]^{\frac{1}{r}}.
				\end{align*}
				By substituting $s^2 = \left( 1 - x \tanh^2 \lambda \right) \cosh^2 \lambda$, and simplifying, we get
				\begin{align}
					\label{Beta1Case1Hn}
					&\| R_kf_{\lambda} \|_{L^r_{\beta_1, \beta_2} \left( \Xi_k \left( \mathbb{H}^n \right) \right)} \nonumber \\
					&= C \cosh^{\frac{\beta_1 + \beta_2 + n - 1}{r} + 1} \lambda \tanh^{k + \frac{\beta_1 + n - k}{r}} \lambda \left[ \int\limits_{0}^{1} \frac{x^{\frac{rk}{2}} \left( 1 - x \right)^{\frac{\beta_1 + n - k}{2} - 1}}{\left( 1 - x \tanh^2 \lambda \right)^{\frac{\left( r - 1 \right) \left( k - 1 \right) - \beta_2}{2}}} \ \mathrm{d}x \right]^{\frac{1}{r}} \\
					&= C \cosh^{\frac{\beta_2 + k - 1}{r} + 1 - k} \lambda \sinh^{k + \frac{\beta_1 + n - k}{r}} \lambda \left[ {}_2F_1 \left( \frac{\left( r - 1 \right) \left( k - 1 \right) - \beta_2}{2}, {1 +} \frac{rk}{2}; {1 +} \frac{\beta_1 + n - k + rk}{2}; \tanh^2 \lambda \right) \right]^{\frac{1}{r}}. \nonumber
				\end{align}
				The last equality follows from the integral representation of the hypergeometric function (Equation \eqref{IntegralForm2F1}). {Indeed, Equation \eqref{Beta1Case1Hn} gives the necessity of $\beta_1 > k - n$, for otherwise the hypergeometric integral does not converge and $R_kf_{\lambda} \notin L^r_{\beta_1, \beta_2} \left( \Xi_k \left( \mathbb{H}^n \right) \right)$.} Depending on the behaviour of the hypergeometric function in the above expression, we have from Theorem \ref{Behaviour2F1},
				\begin{equation}
					\label{LRNormRKFMixedWeight}
					\| R_kf_{\lambda} \|_{L^r_{\beta_1, \beta_2} \left( \Xi_k \left( \mathbb{H}^n \right) \right)} \geq C \begin{cases}
																			\lambda^{k + \frac{\beta_1 + n - k}{r}}, & \lambda \text{ is small}. \\
																			e^{\left( \frac{\beta_1 + \beta_2 + n - 1}{r} + 1 \right) \lambda}, & \lambda \text{ is large and } \beta_1 + \beta_2 \geq r \left( k - 1 \right) - \left( n - 1 \right). \\
																			e^{k \lambda}, & \lambda \text{ is large and } \beta_1 + \beta_2 < r \left( k - 1 \right) - \left( n - 1 \right).
																		\end{cases}
				\end{equation}
				We want Inequality \eqref{LPLRRequiredHn} to hold for all $\lambda > 0$. Particularly, for small $\lambda > 0$, we get from Inequalities \eqref{LPMassHn1} and \eqref{LRNormRKFMixedWeight},
				$$\lambda^{k + \frac{\beta_1 + n - k}{r}} \leq C \ \lambda^{\frac{\alpha_1 + n}{p}}.$$
				As $\lambda \rightarrow 0$, the above inequality forces
				$$k + \frac{\beta_1 + n - k}{r} \geq \frac{\alpha_1 + n}{p}.$$
				
				For large $\lambda > 0$, we have two cases. First, if $\beta_1 + \beta_2 \geq r \left( k - 1 \right) - \left( n - 1 \right)$, then, from Inequalities \eqref{LPLRRequiredHn}, \eqref{LPMassHn1}, and \eqref{LRNormRKFMixedWeight}, we get
				$$e^{\left( \frac{\beta_1 + \beta_2 + n - 1}{r} + 1 \right) \lambda} \leq C \ e^{\left( \frac{\alpha_1 + \alpha_2 + n - 1}{p} + 1 \right) \lambda}.$$
				As $\lambda \rightarrow \infty$, the above inequality forces
				$$\frac{\beta_1 + \beta_2 + n - 1}{r} \leq \frac{\alpha_1 + \alpha_2 + n - 1}{p}.$$
				This is {equivalent to $\frac{\beta_2 + k - 1}{r} - \frac{\alpha_2 - 1}{p} \leq \frac{\alpha_1 + n}{p} - \frac{\beta_1 + n - k}{r}$}. {In the case when $\beta_1 + \beta_2 < r \left( k - 1 \right) - \left( n - 1 \right)$, we first observe that from the existence condition of Theorem \ref{ExistenceMixed},}
				$$\frac{\beta_1 + \beta_2 + n - 1}{r} < k - 1 < \frac{\alpha_1 + \alpha_2 + n - 1}{p}.$$
				{That is, in any case, we must have $\frac{\beta_1 + \beta_2 + n - 1}{r} \leq \frac{\alpha_1 + \alpha_2 + n - 1}{p}$, for the validity of Inequality \eqref{LPLRRequiredHn}.}
				
				Also, from Inequalities \eqref{LPLRRequiredHn}, \eqref{LPMassHn1}, and \eqref{LRNormRKFMixedWeight}, we get
				$$e^{k \lambda} \leq C e^{\left( \frac{\alpha_1 + \alpha_2 + n - 1}{p} + 1 \right) \lambda}.$$
				As $\lambda \rightarrow \infty$, the above inequality forces $p \leq \frac{\alpha_1 + \alpha_2 + n - 1}{k - 1}$, which is redundant! This gives the necessity of the conditions of Equation \eqref{NecessaryMixedWeightInequality} for $\alpha_1 > -n$.
				
				\textbf{Case II:} Let us now assume $\alpha_1 \leq -n$, and consider the function $f_{\lambda} \left( x \right) = \chi_{B \left( 0, \lambda \right)} \left( x \right) \sinh^{- \frac{\alpha_1}{p}} d \left( 0, x \right) \cosh d \left( 0, x \right)$. Then, we have, using Equation {\eqref{LpNormRadialHn1}},
				\begin{equation}
					\label{LPMassHn2}
					\| f_{\lambda} \|_{L^p_{\alpha_1, \alpha_2} \left( \mathbb{H}^n \right)} = C \left[ \int\limits_{0}^{{\lambda}} \sinh^{n - 1} t \cosh^{\alpha_2 + p} t \ \mathrm{d}t \right]^{\frac{1}{p}} \leq C \begin{cases}
	\lambda^{\frac{n}{p}}, & \lambda \text{ is small}. \\
	e^{\left( \frac{\alpha_2 + n - 1}{p} + 1 \right) \lambda}, & \lambda \text{ is large}.
\end{cases}
				\end{equation}
				{We now use Equation \eqref{NormRkHn} to estimate the $L^r_{\beta_1, \beta_2}$-norm of $R_kf_{\lambda}$.}
				\begin{align*}
					\| R_kf \|_{L^r_{\beta_1, \beta_2} \left( \Xi_k \left( \mathbb{H}^n \right) \right)} &= C \left[ \int\limits_{1}^{\infty} \left| \int\limits_{s}^{\infty} \frac{\chi_{\left( 1, \cosh \lambda \right)} \left( t \right) t \left( t^2 - 1 \right)^{- \frac{\alpha_1}{{2}p}} \left( t^2 - s^2 \right)^{\frac{k}{2} - 1}}{s^{k- 1}} \mathrm{d}t \right|^r s^{\beta_2 + k} \left( s^2 - 1 \right)^{\frac{\beta_1 + n - k}{2} - 1} \mathrm{d}s \right]^{\frac{1}{r}} \\
					&= C \left[ \int\limits_{1}^{\cosh \lambda} \left| \int\limits_{s}^{\cosh {\lambda}} t \left( t^2 - 1 \right)^{- \frac{\alpha_1}{{2}p}} \left( t^2 - s^2 \right)^{\frac{k}{2} - 1} \mathrm{d}t \right|^r s^{\beta_2 + k - r \left( k - 1 \right)} \left( s^2 - 1 \right)^{\frac{\beta_1 + n - k}{2} - 1} \mathrm{d}s \right]^{\frac{1}{r}}.
				\end{align*}
				{In the inner integral, let} us substitute $t^2 = x \left( \cosh^2 \lambda - s^2 \right) + s^2$. Then, after simplification, we get
				\begin{align*}
					&\| R_kf \|_{L^r_{\beta_1, \beta_2} \left( \Xi_k \left( \mathbb{H}^n \right) \right)} \\
					&= C \left[ \int\limits_{1}^{\cosh \lambda} \left| \int\limits_{0}^{1} \left( 1 + \left( \frac{\cosh^2 \lambda - s^2}{s^2 - 1} \right) x \right)^{- \frac{\alpha_1}{{2}p}} \left( s^2 - 1 \right)^{- \frac{\alpha_1}{{2}p}} \left( \cosh^2 \lambda - s^2 \right)^{\frac{k}{2}} x^{\frac{k}{2} - 1} \mathrm{d}x \right|^r s^{\beta_2 + k - r \left( k - 1 \right)} \left( s^2 - 1 \right)^{\frac{\beta_1 + n - k}{2} - 1} \mathrm{d}s \right]^{\frac{1}{r}} \\
					&= C \left[ \int\limits_{1}^{\cosh \lambda} \left| {}_2F_1 \left( \frac{\alpha_1}{{2}p}, \frac{k}{2}; 1 + \frac{k}{2}; \frac{s^2 - \cosh^2 \lambda}{s^2 - 1} \right) \right|^r \left( \cosh^2 \lambda - s^2 \right)^{\frac{{r}k}{2}} s^{\beta_2 + k - r \left( k - 1 \right)} \left( s^2 - 1 \right)^{\frac{\beta_1 + n - k}{2} - 1 - \frac{{r}\alpha_1}{{2}p}} \mathrm{d}s \right]^{\frac{1}{r}}.
				\end{align*}
				In the last step, we have used the integral form of the hypergeometric function given in Equation \eqref{IntegralForm2F1}. Now, we use the transformation of the hypergeometric function given in Equation \eqref{Transformation2F1} \eqref{2F1Transformation1} to obtain
				\begin{align*}
					&\| R_kf \|_{L^r_{\beta_1, \beta_2} \left( \Xi_k \left( \mathbb{H}^n \right) \right)} \\
					&= C \sinh^{- \frac{\alpha_1}{p}} {\lambda} \left[ \int\limits_{1}^{\cosh \lambda} \left| {}_2F_1 \left( \frac{\alpha_1}{{2}p}, 1; 1 + \frac{k}{2}; \frac{\cosh^2 \lambda - s^2}{\sinh^2 \lambda} \right) \right|^r \left( \cosh^2 \lambda - s^2 \right)^{\frac{rk}{2}} s^{\beta_2 + k - r \left( k - 1 \right)} \left( s^2 - 1 \right)^{\frac{\beta_1 + n - k}{2} - 1} \mathrm{d}s \right]^{\frac{1}{r}}.
				\end{align*}
				{Since $\alpha_1 \leq -n < 0$, from Remark \ref{HypergemetricFunctionBoundedBelow}, we see} that the hypergeometric function in the above integral is bounded below {by a positive constant}. Hence, we have,
				\begin{align*}
					\| R_kf \|_{L^r_{\beta_1, \beta_2} \left( \Xi_k \left( \mathbb{H}^n \right) \right)} &\geq C \sinh^{- \frac{\alpha_1}{p}} \lambda \left[ \int\limits_{1}^{\cosh \lambda} \left( \cosh^2 \lambda - s^2 \right)^{\frac{rk}{2}} s^{\beta_2 + k - r \left( k - 1 \right)} \left( s^2 - 1 \right)^{\frac{\beta_1 + n - k}{2} - 1} \mathrm{d}s \right]^{\frac{1}{r}}.
				\end{align*}
				Now, we substitute $s^2 = 1 + z \sinh^2 \lambda$ and simplify to get
				\begin{align}
					\label{Beta1Case2Hn}
					\| R_kf \|_{L^r_{\beta_1, \beta_2} \left( \Xi_k \left( \mathbb{H}^n \right) \right)} &\geq C \sinh^{\frac{\beta_1 + n - k}{r} - \frac{\alpha_1}{p} + k} \lambda \left[ \int\limits_{0}^{1} \frac{\left( 1 - z \right)^{\frac{rk}{2}} z^{\frac{\beta_1 + n - k}{2} - 1}}{\left( 1 + z \sinh^2 \lambda \right)^{\frac{\left( r - 1 \right) \left( k - 1 \right) - \beta_2}{2}}} \ \mathrm{d}z \right]^{\frac{1}{r}} \\
					&= C \sinh^{\frac{\beta_1 + n - k}{r} - \frac{\alpha_1}{p} + k} \lambda \left[ {}_2F_1 \left( \frac{\left( r - 1 \right) \left( k - 1 \right) {- \beta_2}}{2}, \frac{\beta_1 + n - k}{2}; 1 + {\frac{\beta_1 + n - k + rk}{2}}; - \sinh^2 \lambda \right) \right]^{\frac{1}{r}}. \nonumber
				\end{align}
				In the last equality, we have employed the integral form of the hypergeometric function given in Equation \eqref{IntegralForm2F1}. {Again, from Equation \eqref{Beta1Case2Hn}, it is evident that $\beta_1 > k - n$ is necessary.} Now, using the transformation of Equation \eqref{Transformation2F1} {\eqref{2F1Transformation1}}, we get
				\begin{align*}
					\| R_kf \|_{L^r_{\beta_1, \beta_2} \left( \Xi_k \left( \mathbb{H}^n \right) \right)} &{\geq} C \sinh^{\frac{\beta_1 + n - k}{r} - \frac{\alpha_1}{p} + k} \lambda \ \cosh^{\beta_2 - \left( r - 1 \right) \left( k - 1 \right)} \lambda \times \\
					&\left[ {}_2F_1 \left( \frac{\left( r - 1 \right) \left( k - 1 \right) - \beta_2}{2}, 1 + \frac{rk}{2}; {1 +} \frac{\beta_1 + n - k + rk}{2}; \tanh^2 \lambda \right) \right]^{\frac{1}{r}}.
				\end{align*}
				Finally, using Theorem \ref{Behaviour2F1}, we have
				\begin{equation}
					\label{LRNormMixed2}
					\| R_kf_{\lambda} \|_{L^r_{\beta_1, \beta_2} \left( \Xi_k \left( \mathbb{H}^n \right) \right)} \geq  C \begin{cases}
																				\lambda^{\frac{\beta_1 + n - k}{r} - \frac{\alpha_1}{p} + k}, & \lambda \text{ is samll}. \\
																			e^{\left( \frac{\beta_1 + \beta_2 + n - 1}{r} - \frac{\alpha_1}{p} + 1 \right) \lambda}, & \lambda \text{ is large and } \beta_1 + \beta_2 \geq r \left( k - 1 \right) - \left( n - 1 \right). \\
																			e^{\left( k - \frac{\alpha_1}{p} \right) \lambda}, & \lambda \text{ is large and } \beta_1 + \beta_2 < r \left( k - 1 \right) - \left( n - 1 \right).
																		\end{cases}
				\end{equation}
				We desire the inequality \eqref{LPLRRequiredHn}. From Inequalities \eqref{LPMassHn2} and \eqref{LRNormMixed2}, using $\lambda \rightarrow 0$ and $\lambda \rightarrow \infty$, following the arguments of Case I, we get the necessary conditions stated in Equation \eqref{NecessaryMixedWeightInequality}.
			\end{proof}
			\begin{remark}
				\normalfont
				{We again emphasize that the restrictions on $\alpha_1, \alpha_2$ and $p$, stated in Theorem \ref{NecessaryMixed} are required for the existence of the $k$-plane transform on the weighted Lebesgue space $L^p_{\alpha_1, \alpha_2} \left( \mathbb{H}^n \right)$.}
			\end{remark}
			We now see that conditions of Theorem \ref{MixedWeightNecessary} are also sufficient {for $L^p$-improving norm estimate of Equation \eqref{LPLRRequiredHn}}.
			\begin{theorem}
				\label{MixedWeightLPLRFullRadial}
				Let $\alpha_1 + \alpha_2 > k - n$ {and} $1 \leq p < \frac{\alpha_1 + \alpha_2 + n - 1}{k - 1}$ {or $\alpha_1 + \alpha_2 = k - n$ and $p = 1$}, {and} ${r \geq p}$.
				\begin{enumerate}
					\item[\mylabel{MixedWeightLPLRHnA}{(A)}] For $k \geq 2$, Inequality \eqref{LPLRRequiredHn} holds for any radial function $f \in L^p_{\alpha_1, \alpha_2} \left( \mathbb{H}^n \right)$ {if and only if $\beta_1, \beta_2 \in \mathbb{R}$ satisfy the conditions of Equation \eqref{NecessaryMixedWeightInequality}}.
					\item[\mylabel{MixedWeightLPLRHnB}{(B)}] For $k = 1$, we assume in addition that $\frac{1}{p} - \frac{1}{r} \leq \frac{1}{2}$ for $p > 1$ or that $r < 2$ for $p = 1$. Then, Inequality \eqref{LPLRRequiredHn} holds for any radial function $f \in L^p_{\alpha_1, \alpha_2} \left( \mathbb{H}^n \right)$ {if and only if $\beta_1, \beta_2 \in \mathbb{R}$ satisfy the conditions of Equation \eqref{NecessaryMixedWeightInequality}}.
				\end{enumerate}
			\end{theorem}
			Before we begin the proof, we present its outline {by fixing $\alpha_1 + \alpha_2 > k - n$, $1 \leq p < \frac{\alpha_1 + \alpha_2 + n - 1}{k - 1}$, and $r \geq p$.} First, we observe that the necessary condition stated in Equation \eqref{NecessaryMixedWeightInequality} gives two separate information about the weights $\beta_1$ and $\beta_2$. It is clear that $\beta_2$ has an upper bound (but no lower bound), while $\beta_1$ has both an upper bound and a lower bound. Such a condition is natural since for any $\beta_2' \leq \beta_2$, and a fixed $\beta_1 \in \mathbb{R}$, we have for any $\varphi: \Xi_k \left( \mathbb{H}^n \right) \rightarrow \mathbb{C}$,
			\begin{equation}
				\label{NormComparisonHn}
				\| \varphi \|_{L^r_{\beta_1, \beta_2'} \left( \Xi_k \left( \mathbb{H}^n \right) \right)} \leq \| \varphi \|_{L^r_{\beta_1, \beta_2} \left( \Xi_k \left( \mathbb{H}^n \right) \right)}.
			\end{equation}
			Due to the above inequality, we proceed in the following manner. First, we prove {Inequality \eqref{LPLRRequiredHn}} of the $k$-plane transform {by fixing $\beta_1^{\left( 0 \right)}, \beta_2^{\left( 0 \right)} \in \mathbb{R}$ by the following equations}
			$$\frac{\beta_2^{{\left( 0 \right)}} + k - 1}{r} - \frac{\alpha_2 - 1}{p} = k = \frac{\alpha_1 + n}{p} - \frac{\beta_1^{{\left( 0 \right)}} + n - k}{r}.$$
			{With the norm estimate obtained for the above choices of $\beta_1^{\left( 0 \right)}$ and $\beta_2^{\left( 0 \right)}$, we see with the help of} Inequality \eqref{NormComparisonHn}, {that we have} the boundedness for $\frac{\alpha_1 + n}{p} - \frac{\beta_1^{{\left( 0 \right)}} + n - k}{r} = k$ and {any $\beta_2 \in \mathbb{R}$ satisfying} $\frac{\beta_2 + k - 1}{r} - \frac{\alpha_2 - 1}{p} \leq k {= \frac{\beta_2^{{\left( 0 \right)}} + k - 1}{r} - \frac{\alpha_2 - 1}{p}}$. 
			
			{Now, we fix $\beta_2 \in \mathbb{R}$ such that $\frac{\beta_2 + k - 1}{r} - \frac{\alpha_2 - 1}{p} \leq k$, and choose $\beta_1^{\left( 1 \right)}  \in \mathbb{R}$ by the following equation.}
			$${\frac{\beta_2 + k - 1}{r} + \frac{\alpha_2 - 1}{p} = \frac{\alpha_1 + n}{p} - \frac{\beta_1^{\left( 1 \right)} + n - k}{r}}.$$
			That is, for a fixed $\alpha_1, \alpha_2 \in \mathbb{R}$ with $\alpha_1 + \alpha_2 > k - n$, $1 \leq p < \frac{\alpha_1 + \alpha_2 + n - 1}{k - 1}$, $r \geq p$, and {$\beta_2 \in \mathbb{R}$ satisfying} $\frac{\beta_2 + k - 1}{r} - \frac{\alpha_2 - 1}{p} \leq k$, we get two ``{end-points}", {$\beta_1^{\left( 0 \right)}$ and $\beta_1^{\left( 1 \right)}$}, for the weight $\beta_1$. {We may now use the weighted interpolation stated in} Theorem \ref{RealInterpolationWeighted} to complete the proof.
			
			{The proof for the case when $\alpha_1 + \alpha_2 = k - n$ and $p = 1$ is fairly simple and is given separately.}
			\begin{proof}
				{We prove the result in two cases.}
				
				{\textbf{Case I:} First, we consider $\alpha_1 + \alpha_2 > k - n$ and $1 \leq p < \frac{\alpha_1 + \alpha_2 + n - 1}{k - 1}$.}
				\begin{enumerate}
					\item[\mylabel{MixedWeightLPLRHnAProof}{(A)}] {Let us start by proving} the result for $k \geq 2$.
					
					\textbf{Step I:} We fix {$\beta_1^{\left( 0 \right)}, \beta_2^{\left( 0 \right)} \in \mathbb{R}$ such that} $\frac{\beta_2^{{\left( 0 \right)}} + k - 1}{r} - \frac{\alpha_2 - 1}{p} = k = \frac{\alpha_1 + n}{p} - \frac{\beta_1^{{\left( 0 \right)}} + n - k}{r}$, {with $\beta_1^{\left( 0 \right)} > k - n$}.
				
					By substituting $s^2 = u + 1$ and $t^2 = v + 1$ in Equation \eqref{NormRkHn}, and simplifying, we get,
					\begin{align}
						\label{NormRkHnConvolution}
						\| R_kf \|_{L^r_{\beta_1^{{\left( 0 \right)}}, \beta_2^{{\left( 0 \right)}}} \left( \Xi_k \left( \mathbb{H}^n \right) \right)} &= C \left[ \int\limits_{0}^{\infty} \left| \int\limits_{u}^{\infty} \tilde{f} \left( \sqrt{v + 1} \right) v^{\frac{k}{2}} \left( 1 - \frac{u}{v} \right)^{\frac{k}{2} - 1} \left( u + 1 \right)^{\frac{\beta_2^{{\left( 0 \right)}} + k - 1}{2r} - \frac{k - 1}{2}} u^{\frac{\beta_1^{{\left( 0 \right)}} + n - k}{2r}} \left( v + 1 \right)^{- \frac{1}{2}} \ \frac{\mathrm{d}v}{v} \right|^r \frac{\mathrm{d}u}{u} \right]^{\frac{1}{r}} \\
						&= C \left[ \int\limits_{0}^{\infty} \left| \int\limits_{u}^{\infty} \tilde{f} \left( \sqrt{v + 1} \right) v^{\frac{\alpha_1 + n}{2p}} \left( v + 1 \right)^{\frac{\alpha_2 - 1}{2p}} \left( \frac{u + 1}{v + 1} \right)^{\frac{1}{2} + \frac{\alpha_2 - 1}{2p}} \left( \frac{u}{v} \right)^{\frac{\beta_1^{{\left( 0 \right)}} + n - k}{2r}} \left( 1 - \frac{u}{v} \right)^{\frac{k}{2} - 1} \frac{\mathrm{d}v}{v} \right|^r \frac{\mathrm{d}u}{u} \right]^{\frac{1}{r}}. \nonumber
					\end{align}
					We break the analysis {by considering two situations about the values of $p$}.
				
					{\textbf{Situation I:}} Let us consider $1 + \frac{\alpha_2 - 1}{p} \geq 0$.
					Owing to the fact that $u \leq v$, we get
					\begin{align*}
						\| R_kf \|_{L^r_{\beta_1^{{\left( 0 \right)}}, \beta_2^{{\left( 0 \right)}}} \left( \Xi_k \left( \mathbb{H}^n \right) \right)} &\leq C \left[ \int\limits_{0}^{\infty} \left| \int\limits_{0}^{\infty} \tilde{f} \left( \sqrt{v + 1} \right) v^{\frac{\alpha_1 + n}{2p}} \left( v + 1 \right)^{\frac{\alpha_2 - 1}{2p}} \left( \frac{u}{v} \right)^{\frac{\beta_1^{{\left( 0 \right)}} + n - k}{2r}} \left( 1 - \frac{u}{v} \right)^{\frac{k}{2} - 1} \chi_{\left( 0, 1 \right)} \left( \frac{u}{v} \right) \frac{\mathrm{d}v}{v} \right|^r \frac{\mathrm{d}u}{u} \right]^{\frac{1}{r}} \\
						&= C \| F * G \|_{L^r \left( \left( 0, \infty \right), \frac{\mathrm{d}t}{t} \right)},
					\end{align*}
					where,
					\begin{equation}
						\label{ConvolutionFHn1}
						F \left( t \right) = \tilde{f} \left( \sqrt{t + 1} \right) t^{\frac{\alpha_1 + n}{2p}} \left( t + 1 \right)^{\frac{\alpha_2 - 1}{2p}},
					\end{equation}
					and
					\begin{equation}
						\label{ConvoltionGHn1}
						G \left( t \right) = t^{\frac{\beta_1^{{\left( 0 \right)}} + n - k}{2r}} \left( 1 - t \right)^{\frac{k}{2} - 1} \chi_{\left( 0, 1 \right)} \left( t \right).
					\end{equation}
					By Young's convolution inequality, we have
					$$\| R_kf \|_{L^r_{\beta_1^{{\left( 0 \right)}}, \beta_2^{{\left( 0 \right)}}} \left( \Xi_k \left( \mathbb{H}^n \right) \right)} \leq C \| F \|_{L^p \left( \left( 0, \infty \right), \frac{\mathrm{d}t}{t} \right)} \| G \|_{L^q \left( \left( 0, \infty \right), \frac{\mathrm{d}t}{t} \right)},$$
					where, $1 + \frac{1}{r} = \frac{1}{p} + \frac{1}{q}$. From Equation \eqref{LpNormRadialHn3}, we know that $\| F \|_{L^p \left( \left( 0, \infty \right), \frac{\mathrm{d}t}{t} \right)} = \| f \|_{L^p_{\alpha_1, \alpha_2} \left( \mathbb{H}^n \right)}$. Also, {since $k \geq 2$ and $\beta_1^{{\left( 0 \right)}} > k - n$, we have}
					\begin{align*}
						\| G \|_{L^q \left( \left( 0, \infty \right), \frac{\mathrm{d}t}{t} \right)} &= \left[ \int\limits_{0}^{{1}} t^{q \left( \frac{\beta_1^{{\left( 0 \right)}} + n - k}{2r} \right) - 1} \left( 1 - t \right)^{q \left( \frac{k}{2} - 1 \right)} \mathrm{d}t \right]^{\frac{1}{q}} < + \infty,
					\end{align*}
					Hence, we {get}
					$$\| R_kf \|_{L^r_{\beta_1^{{\left( 0 \right)}}, \beta_2^{{\left( 0 \right)}}} \left( \Xi_k \left( \mathbb{H}^n \right) \right)} \leq C \| f \|_{L^p_{\alpha_1, \alpha_2} \left( \mathbb{H}^n \right)}.$$
				
					{\textbf{Situation II:}} Now, let us consider $1 + \frac{\alpha_2 - 1}{p} < 0$. From Equation \eqref{NormRkHnConvolution} and the fact $\frac{\alpha_1 + n}{p} - \frac{\beta_1^{{\left( 0 \right)}} + n - k}{r} = k$, we get,
					\begin{align*}
						\| R_kf \|_{L^r_{\beta_1^{{\left( 0 \right)}}, \beta_2^{{\left( 0 \right)}}} \left( \Xi_k \left( \mathbb{H}^n \right) \right)} &= C \left[ \int\limits_{0}^{\infty} \left| \int\limits_{u}^{\infty} \tilde{f} \left( \sqrt{v + 1} \right) v^{\frac{\alpha_1 + n}{2p}} \left( v + 1 \right)^{\frac{\alpha_2 - 1}{2p}} \left( \frac{u + 1}{v + 1} \right)^{\frac{1}{2} + \frac{\alpha_2 - 1}{2p}} \left( \frac{u}{v} \right)^{\frac{\alpha_1 + n}{2p} - \frac{k}{2}} \left( 1 - \frac{u}{v} \right)^{\frac{k}{2} - 1} \frac{\mathrm{d}v}{v} \right|^r \frac{\mathrm{d}u}{u} \right]^{\frac{1}{r}} \\
						&\leq C \left[  \int\limits_{0}^{\infty} \left| \int\limits_{0}^{\infty} \tilde{f} \left( \sqrt{v + 1} \right) v^{\frac{\alpha_1 + n}{2p}} \left( v + 1 \right)^{\frac{\alpha_2 - 1}{2p}} \left( \frac{u}{v} \right)^{\frac{\alpha_1 + \alpha_2 + n - 1}{2p} - \frac{k - 1}{2}} \left( 1 - \frac{u}{v} \right)^{\frac{k}{2} - 1} \chi_{\left( 0, 1 \right)} \left( \frac{u}{v} \right) \frac{\mathrm{d}v}{v} \right|^r \frac{\mathrm{d}u}{u} \right]^{\frac{1}{r}},
					\end{align*}
					where the last inequality is due to the fact that $u \leq v$ and $1 + \frac{\alpha_2 - 1}{p} < 0$. From the above inequality, it is clear that
					$$\| R_kf \|_{L^r_{\beta_1^{{\left( 0 \right)}}, \beta_2^{{\left( 0 \right)}}} \left( \Xi_k \left( \mathbb{H}^n \right) \right)} \leq C\| F * G \|_{L^r \left( \left( 0, \infty \right), \frac{\mathrm{d}t}{t} \right)},$$
					where $F$ is as given in Equation \eqref{ConvolutionFHn1}, and
					\begin{equation}
						\label{ConvolutionGHn2}
						G \left( t \right) = t^{\frac{\alpha_1 + \alpha_2 + n - 1}{2p} - \frac{k - 1}{2}} \left( 1 - t \right)^{\frac{k}{2} - 1} \chi_{\left( 0, 1 \right)} \left( t \right).
					\end{equation}
					Now, {since $p < \frac{\alpha_1 + \alpha_2 + n - 1}{k - 1}$ and $k \geq 2$, we have}
					$$\| G \|_{L^q \left( \left( 0, \infty \right), \frac{\mathrm{d}t}{t} \right)} = \left[ \int\limits_{0}^{1} t^{q \left( \frac{\alpha_1 + \alpha_2 + n - 1}{2p} - \frac{k - 1}{2} \right) - 1} \left( 1 - t \right)^{q \left( \frac{k}{2} - 1 \right)} \mathrm{d}t \right]^{\frac{1}{q}} < + \infty{.}$$
					{Hence,} from Equation \eqref{LpNormRadialHn3}, {we conclude that}
					$$\| R_kf \|_{L^r_{\beta_1^{{\left( 0 \right)}}, \beta_2^{{\left( 0 \right)}}} \left( \Xi_k \left( \mathbb{H}^n \right) \right)} \leq C \| f \|_{L^p_{\alpha_1, \alpha_2} \left( \mathbb{H}^n \right)}.$$
					
					{As stated in the outline of the proof, we now have for any $\beta_2 \in \mathbb{R}$ with $\frac{\beta_2 + k - 1}{r} - \frac{\alpha_2 - 1}{p} \leq k$,
					$$\| R_kf \|_{L^r_{\beta_1^{{\left( 0 \right)}}, \beta_2} \left( \Xi_k \left( \mathbb{H}^n \right) \right)} \leq C \| f \|_{L^p_{\alpha_1, \alpha_2} \left( \mathbb{H}^n \right)}.$$}
					\textbf{Step II:} {Let us} now consider {$\beta_2 \in \mathbb{R}$ with $\frac{\beta_2 + k - 1}{r} - \frac{\alpha_2 - 1}{p} \leq k$ and $\beta_1^{\left( 1 \right)} \in \mathbb{R}$ given by} $\frac{\beta_2 + k - 1}{r} - \frac{\alpha_2 - 1}{p} = \frac{\alpha_1 + n}{p} - \frac{\beta_1^{{\left( 1 \right)}} + n - k}{r}$. We have from Equation \eqref{NormRkHnConvolution},
					\begin{align}
						\| R_kf \|_{L^r_{\beta_1^{{\left( 1 \right)}}, \beta_2} \left( \Xi_k \left( \mathbb{H}^n \right) \right)} &= C \left[ \int\limits_{0}^{\infty} \left| \int\limits_{u}^{\infty} \tilde{f} \left( \sqrt{v + 1} \right) v^{\frac{k}{2}} \left( 1 - \frac{u}{v} \right)^{\frac{k}{2} - 1} \left( u + 1 \right)^{\frac{\beta_2 + k - 1}{2r} - \frac{k - 1}{2}} u^{\frac{\beta_1^{{\left( 1 \right)}} + n - k}{2r}} \left( v + 1 \right)^{- \frac{1}{2}} \frac{\mathrm{d}v}{v} \right|^r \frac{\mathrm{d}u}{u} \right]^{\frac{1}{r}} \nonumber \\
						&= C \left[ \int\limits_{0}^{\infty} \left| \int\limits_{u}^{\infty} \tilde{f} \left( \sqrt{v + 1} \right) v^{\frac{\alpha_1 + n}{2p}} \left( v + 1 \right)^{\frac{\alpha_2 - 1}{2p}} \left( \frac{v}{v + 1} \right)^{\frac{k}{2} - \frac{\alpha_1 + n}{2p} - \frac{\beta_1^{{\left( 1 \right)}} + n - k}{2r}} \times \right. \right. \nonumber \\
						&\left. \left. \left( v + 1 \right)^{\frac{k}{2} - \frac{\alpha_1 + n}{2p} + \frac{\beta_1^{{\left( 1 \right)}} + n - k}{2r} - \frac{1}{2} - \frac{\alpha_2 - 1}{2p}} \left( u + 1 \right)^{\frac{\beta_2 + k - 1}{2r} - \frac{k - 1}{2}} \left( \frac{u}{v} \right)^{\frac{\beta_1^{{\left( 1 \right)}} + n - k}{2r}} \left( 1 - \frac{u}{v} \right)^{\frac{k}{2} - 1} \frac{\mathrm{d}v}{v} \right|^r \frac{\mathrm{d}u}{u} \right]^{\frac{1}{r}} \nonumber \\
						\label{Step2Inequality}
						&\leq C \left[ \int\limits_{0}^{\infty} \left| \int\limits_{u}^{\infty} \tilde{f} \left( \sqrt{v + 1} \right) v^{\frac{\alpha_1 + n}{2p}} \left( v + 1 \right)^{\frac{\alpha_2 - 1}{2p}} \left( \frac{u + 1}{v + 1} \right)^{\frac{\beta_2 + k - 1}{2r} - \frac{k - 1}{2}} \left( \frac{u}{v} \right)^{\frac{\beta_1^{{\left( 1 \right)}} + n - k}{2r}} \left( 1 - \frac{u}{v} \right)^{\frac{k}{2} - 1} \frac{\mathrm{d}v}{v} \right|^r \frac{\mathrm{d}u}{u} \right]^{\frac{1}{r}},
					\end{align}
					where, in the last inequality, we have used the relation $\frac{\alpha_1 + n}{p} - \frac{\beta_1^{{\left( 1 \right)}} + n - k}{r} = \frac{\beta_2 + k - 1}{r} - \frac{\alpha_2 - 1}{p} \leq k$. We again consider two {situations}. 
					
					{\textbf{Situation I:}} Let $\frac{\beta_2 + k - 1}{r} - k + 1 \geq 0$. Then, by virtue of $u \leq v$, we have from Inequality \eqref{Step2Inequality},
					\begin{align*}
						\| R_kf \|_{L^r_{\beta_1^{{\left( 1 \right)}}, \beta_2} \left( \Xi_k \left( \mathbb{H}^n \right) \right)} &\leq C \left[ \int\limits_{0}^{\infty} \left| \int\limits_{0}^{\infty} \tilde{f} \left( \sqrt{v + 1} \right) v^{\frac{\alpha_1 + n}{2p}} \left( v + 1 \right)^{\frac{\alpha_2 - 1}{2p}} \left( \frac{u}{v} \right)^{\frac{\beta_1^{{\left( 1 \right)}} + n - k}{2r}} \left( 1 - \frac{u}{v} \right)^{\frac{k}{2} - 1} \chi_{\left( 0, 1 \right)} \left( \frac{u}{v} \right) \frac{\mathrm{d}v}{v} \right|^r \frac{\mathrm{d}u}{u} \right]^{\frac{1}{r}} \\
						&= C \| F * G \|_{L^r \left( \left( 0, \infty \right), \frac{\mathrm{d}t}{t} \right)},
					\end{align*}
					where, $F$ and $G$ are given by Equations \eqref{ConvolutionFHn1} and \eqref{ConvoltionGHn1}, respectively, {with $\beta_1^{\left( 0 \right)}$ replaced by $\beta_1^{{\left( 1 \right)}}$}. From the arguments provided in Step I, {Situation} I, we have Inequality \eqref{LPLRRequiredHn} {for all radial functions in $L^p_{\alpha_1, \alpha_2} \left( \mathbb{H}^n \right)$}.
					
					{\textbf{Situation II:}} Next, let us consider $\frac{\beta_2 + k - 1}{r} - k + 1 < 0$. Then, we have from Inequality \eqref{Step2Inequality} and $u \leq v$,
					\begin{align*}
						\| R_kf \|_{L^r_{\beta_1^{{\left( 1 \right)}}, \beta_2} \left( \Xi_k \left( \mathbb{H}^n \right) \right)} &\leq C \left[ \int\limits_{0}^{\infty} \left| \int\limits_{0}^{\infty} \tilde{f} \left( \sqrt{v + 1} \right) v^{\frac{\alpha_1 + n}{2p}} \left( v + 1 \right)^{\frac{\alpha_2 - 1}{2p}} \left( \frac{u}{v} \right)^{\frac{\alpha_1 + \alpha_2 + n - 1}{2p} - \frac{k - 1}{2}} \left( 1 - \frac{u}{v} \right)^{\frac{k}{2} - 1} \chi_{\left( 0, 1 \right)} \left( \frac{u}{v} \right) \frac{\mathrm{d}v}{v} \right|^r \frac{\mathrm{d}u}{u} \right]^{\frac{1}{r}}.
					\end{align*}
					Here, we have utilized the fact that $\frac{\alpha_1 + n}{p} - \frac{\beta_1^{{\left( 1 \right)}} + n - k}{r} = \frac{\beta_2 + k - 1}{r} - \frac{\alpha_2 - 1}{p}$. Thus, we again have an inequality of the form
					$$\| R_kf \|_{{L^r_{\beta_1^{\left( 1 \right)}, \beta_2} \left( \Xi_k \left( \mathbb{H}^n \right) \right)}} \leq C \| F * G \|_{L^r \left( \left( 0, \infty \right), \frac{\mathrm{d}t}{t} \right)},$$
					where, $F$ and $G$ are given by Equations \eqref{ConvolutionFHn1} and \eqref{ConvolutionGHn2}, respectively. Following the arguments of Step I, {Situation} II, we now have the {desired} boundedness result.
					
					\textbf{Step III:} We now interpolate between the ``end-points" of $\beta_1$ to obtain the result. Let $1 \leq p < \frac{\alpha_1 + \alpha_2 + n - 1}{k - 1}$, ${p} \leq r < \infty$, $\frac{\beta_2 + k - 1}{r} - \frac{\alpha_2 - 1}{p} \leq k$ be fixed. Let $\beta_1^{\left( 0 \right)}, \beta_1^{\left( 1 \right)} \in \mathbb{R}$ be such that $\frac{\alpha_1 + n}{p} - \frac{\beta_1^{\left( 0 \right) + n - k}}{r} = k$ and $\frac{\alpha_1 + n}{p} - \frac{\beta_1^{\left( 1 \right)} + n - k}{r} = \frac{\beta_2 + k - 1}{r} - \frac{\alpha_2 - 1}{p}$. From Step I, we have an inequality of the form
					\begin{equation}
						\label{EndPoint1HnLPLR}
						\Bigg\| \sinh^{\frac{\beta_1^{\left( 0 \right)}}{r}} d \left( 0, \cdot \right) \cosh^{\frac{\beta_2}{r}} d \left( 0, \cdot \right) R_kf \Bigg\|_{L^r \left( \Xi_k \left( \mathbb{H}^n \right) \right)} \leq C \Bigg\| \sinh^{\frac{\alpha_1}{p}} d \left( 0, \cdot \right) \cosh^{\frac{\alpha_2}{p}} d \left( 0, \cdot \right) f \Bigg\|_{L^p \left( \mathbb{H}^n \right)}.
					\end{equation}
					From Step II, we have,
					\begin{equation}
						\label{EndPoint2HnLPLR}
						\Bigg\| \sinh^{\frac{\beta_1^{\left( 1 \right)}}{r}} d \left( 0, \cdot \right) \cosh^{\frac{\beta_2}{r}} d \left( 0, \cdot \right) R_kf \Bigg\|_{L^r \left( \Xi_k \left( \mathbb{H}^n \right) \right)} \leq C \Bigg\| \sinh^{\frac{\alpha_1}{p}} d \left( 0, \cdot \right) \cosh^{\frac{\alpha_2}{p}} d \left( 0, \cdot \right) f \Bigg\|_{L^p \left( \mathbb{H}^n \right)}.
					\end{equation}
					{Hence, } from Theorem \ref{RealInterpolationWeighted}, we have for every $\theta \in \left( 0, 1 \right)$,
					\begin{equation}
						\label{InterpolatedResult}
						\Bigg\| \sinh^{\frac{\beta_1^{\left( \theta \right)}}{r}} d \left( 0, \cdot \right) \cosh^{\frac{\beta_2}{r}} d \left( 0, \cdot \right) R_kf \Bigg\|_{L^r \left( \Xi_k \left( \mathbb{H}^n \right) \right)} \leq C \Bigg\| \sinh^{\frac{\alpha_1}{p}} d \left( 0, \cdot \right) \cosh^{\frac{\alpha_2}{p}} d \left( 0, \cdot \right) f \Bigg\|_{L^p \left( \mathbb{H}^n \right)},
					\end{equation}
					where,
					$$\beta_1^{\left( \theta \right)} = \left( 1 - \theta \right) \beta_1^{\left( 0 \right)} + \theta \beta_1^{\left( 1 \right)}.$$
					It is readily seen that we have
					$$\frac{\beta_2 + k - 1}{r} - \frac{\alpha_2 - 1}{p} \leq \frac{\alpha_1 + n}{p} - \frac{\beta_1^{\left( \theta \right)} + n - k}{r} \leq k,$$
					{and for any $\beta_1 \in \mathbb{R}$ that satisfies the conditions of Equation \eqref{NecessaryMixedWeightInequality}, there is some $\theta \in \left[ 0, 1 \right]$ such that $\beta_1^{\left( \theta \right)} = \beta_1$.}
					
					This completes the proof of \ref{MixedWeightLPLRHnAProof}.
					\item[\mylabel{MixedWeightLPLRHnBProof}{(B)}] We now prove the result for $k = 1$. First, we notice that for $\frac{1}{p} - \frac{1}{r} < \frac{1}{2}$, the proof is verbatim to that of \ref{MixedWeightLPLRHnAProof}. The condition $\frac{1}{p} - \frac{1}{r} < \frac{1}{2}$ is used to show that the functions $G$ of Equations \eqref{ConvoltionGHn1} and \eqref{ConvolutionGHn2} are in $L^q \left( \left( 0, \infty \right), \frac{\mathrm{d}t}{t} \right)$. Therefore, it only remains to see the case when $\frac{1}{p} - \frac{1}{r} = \frac{1}{2}$. As in the assumption of the statement of the result, in this case, we consider $p > 1$. To this end, we wish to employ Theorem \ref{LPLQBoundednessFI}.
					
					We start with Equation \eqref{NormRkHnConvolution}. We have,
					\begin{equation}
						\label{R1FIEquation}
						\| R_kf \|_{L^r_{\beta_1, \beta_2} \left( \Xi_k \left( \mathbb{H}^n \right) \right)} = C \left[ \int\limits_{0}^{\infty} \left| \int\limits_{u}^{\infty} \frac{\tilde{f} \left( \sqrt{v + 1} \right) \left( v + 1 \right)^{- \frac{1}{2}}}{\left( v - u \right)^{\frac{1}{2}}} \left( u + 1 \right)^{\frac{\beta_2}{2r}} u^{\frac{\beta_2 + n - 1}{2r} - \frac{1}{r}} \mathrm{d}v \right|^r \mathrm{d}u \right]^{\frac{1}{r}} = C \| \rho_- I_-^{\frac{1}{2}} \varphi \|_{L^r \left( \Omega \right)},
					\end{equation}
					where, in the notation of Theorem \ref{LPLQBoundednessFI}, $\Omega = \overline{\mathbb{R}_+}$, and
					$$\rho_- \left( x \right) = \left( 1 + x \right)^{\frac{\beta_2}{2r}} x^{\frac{\beta_1 + n - 1}{2r} - \frac{1}{r}}, \text{ and } \varphi \left( x \right) = \tilde{f} \left( \sqrt{x + 1} \right) \left( x + 1 \right)^{- \frac{1}{2}}.$$
					Now, let us consider
					$$\rho \left( x \right) = \left( 1 + x \right)^{\frac{\beta_2}{2r}} x^{\frac{\beta_1 + n - 1}{2r} - \frac{1}{r}}.$$
					Therefore, by Theorem \ref{LPLQBoundednessFI}, we have,
					\begin{equation}
						\label{NormR1HnFI}
						\| R_kf \|_{L^r_{\beta_1, \beta_2} \left( \Xi_k \left( \mathbb{H}^n \right) \right)} = C \| \rho_- I_-^{\frac{1}{2}} \varphi \|_{L^r \left( \Omega \right)} \leq C \| \rho \varphi \|_{L^p \left( \Omega \right)} = C \left[ \int\limits_{0}^{\infty} \left| \tilde{f} \left( \sqrt{x + 1} \right) \right|^p \left( x + 1 \right)^{p \left( \frac{\beta_2}{2r} - \frac{1}{2} \right)} x^{p \left( \frac{\beta_1 + n - 1}{2r} - \frac{1}{r} \right)} \mathrm{d}x \right]^{\frac{1}{p}}.
					\end{equation}
					Now, let us consider the two ``end-points" for the weight $\beta_1$ as done in the proof of (1).
					
					\textbf{Step I:} Let us {consider} {$\beta_1^{\left( 0 \right)}, \beta_2^{\left( 0 \right)} \in \mathbb{R}$ given by} $\frac{\beta_2^{{\left( 0 \right)}} + k - 1}{r} - \frac{\alpha_2 - 1}{p} = k = \frac{\alpha_1 + n}{p} - \frac{\beta_1^{{\left( 0 \right)}} + n - k}{r}$, {with $k = 1$}. Using this relation together with $\frac{1}{p} - \frac{1}{r} = \frac{1}{2}$ in Equation \eqref{NormR1HnFI}, we get
					\begin{align*}
						\| R_kf \|_{L^r_{\beta_1^{{\left( 0 \right)}}, \beta_2^{{\left( 0 \right)}}} \left( \Xi_k \left( \mathbb{H}^n \right) \right)} &\leq C \left[ \int\limits_{0}^{\infty} \left| \tilde{f} \left( \sqrt{x + 1} \right) \right|^p \left( x + 1 \right)^{\frac{\alpha_2 - 1}{2}} x^{\frac{\alpha_1 + n}{2} - 1} \mathrm{d}x \right]^{\frac{1}{p}}.
					\end{align*}
					From Equation \eqref{LpNormRadialHn3}, we conclude that
					$$\| R_kf \|_{L^r_{\beta_1^{{\left( 0 \right)}}, \beta_2^{{\left( 0 \right)}}} \left( \Xi_k \left( \mathbb{H}^n \right) \right)} \leq C \| f \|_{L^p_{\alpha_1, \alpha_2} \left( \mathbb{H}^n \right)}.$$
					
					{Hence, from Inequality \eqref{NormComparisonHn}, we have for any $\beta_2 \in \mathbb{R}$ with $\frac{\beta_2 + k - 1}{r} - \frac{\alpha_2 - 1}{p} \leq k$,
					$$\| R_kf \|_{L^r_{\beta_1^{{\left( 0 \right)}}, \beta_2} \left( \Xi_k \left( \mathbb{H}^n \right) \right)} \leq C \| f \|_{L^p_{\alpha_1, \alpha_2} \left( \mathbb{H}^n \right)},$$
					for all radial functions $f \in L^p_{\alpha_1, \alpha_2} \left( \mathbb{H}^n \right)$.}
					
					\textbf{Step II:} Now, we {consider} {$\beta_2 \in \mathbb{R}$ with $\frac{\beta_2 + k - 1}{r} - \frac{\alpha_2 - 1}{p} \leq k$ and $\beta_1^{\left( 1 \right)} \in \mathbb{R}$ given by} $\frac{\beta_2 + k - 1}{r} - \frac{\alpha_2 - 1}{p} = \frac{\alpha_1 + n}{p} - \frac{\beta_1 + n - k}{r}$, {with $k = 1$}. {Using this assumption} in Equation \eqref{NormR1HnFI}, we reach at
					\begin{align*}
						\| R_kf \|_{L^r_{\beta_1^{{\left( 1 \right)}}, \beta_2} \left( \Xi_k \left( \mathbb{H}^n \right) \right)} &\leq C \left[ \int\limits_{0}^{\infty} \left| \tilde{f} \left( \sqrt{x + 1} \right) \right|^p \left( x + 1 \right)^{\frac{\alpha_2 - 1}{2}} x^{\frac{\alpha_1 + n}{2} - 1} \left( \frac{x}{x + 1} \right)^{p \left( \frac{\alpha_2 - 1}{2p} - \frac{\beta_2}{2r} + \frac{1}{2} \right)} \mathrm{d}x \right]^{\frac{1}{p}} \\
						&\leq C \left[ \int\limits_{0}^{\infty} \left| \tilde{f} \left( \sqrt{x + 1} \right) \right|^p \left( x + 1 \right)^{\frac{\alpha_2 - 1}{2}} x^{\frac{\alpha_1 + n}{2} - 1} \mathrm{d}x \right]^{\frac{1}{p}} = C \| f \|_{L^p_{\alpha_1, \alpha_2} \left( \mathbb{H}^n \right)}.
					\end{align*}
					The last equality follows from Equation \eqref{LpNormRadialHn3}. The proof of the result is now complete by using the interpolation argument mentioned in Step III of \ref{MixedWeightLPLRHnAProof}.
				\end{enumerate}
				
				{\textbf{Case II}: We now take up the case when $\alpha_1 + \alpha_2 = k - n$ and $p = 1$. With the given assumptions, let $f \in L^1_{\alpha_1, \alpha_2} \left( \mathbb{H}^n \right)$ be a radial function. Then, from Equation \eqref{NormRkHn}, we get
				\begin{align*}
					\| R_kf \|_{L^r_{\beta_1, \beta_2} \left( \mathbb{H}^n \right)} &= C \left[ \int\limits_{1}^{\infty} \left| \int\limits_{s}^{\infty} \frac{\tilde{f} \left( t \right) \left( t^2 - s^2 \right)^{\frac{k}{2} - 1}}{s^{k - 1}} \mathrm{d}t \right|^r s^{\beta_2 + k - 1} \left( s^2 - 1 \right)^{\frac{\beta_1 + n - k}{2} - 1} \mathrm{d}s \right]^{\frac{1}{r}}.
				\end{align*}
				By an application of Minkowski's integral inequality, we get
				\begin{align*}
					\| R_kf \|_{L^r_{\beta_1, \beta_2} \left( \mathbb{H}^n \right)} &\leq C \int\limits_{1}^{\infty} \left| \tilde{f} \left( t \right) \right| \left[ \int\limits_{1}^{t} \left( t^2 - s^2 \right)^{r \left( \frac{k}{2} - 1 \right)} s^{\beta_2 - \left( r - 1 \right) \left( k - 1 \right)} \left( s^2 - 1 \right)^{\frac{\beta_1 + n - k}{2} - 1} \mathrm{d}s \right]^{\frac{1}{r}} \mathrm{d}t.
				\end{align*}
				In the inner integral, we substitute $s^2 - 1 = \left( t^2 - 1 \right) v$, and obtain
				\begin{align}
					\label{IntermediateEquationL1LRHn}
					\| R_kf \|_{L^r_{\beta_1, \beta_2} \left( \mathbb{H}^n \right)} &\leq C \int\limits_{1}^{\infty} \left| \tilde{f} \left( t \right) \right| \left( t^2 - 1 \right)^{\frac{\beta_1 + n - k}{2r} + \frac{k}{2} - 1} \left[ \int\limits_{0}^{1} \frac{\left( 1 - v \right)^{r \left( \frac{k}{2} - 1 \right)} v^{\frac{\beta_1 + n - k}{2} - 1}}{\left( 1 - \left( 1 - t^2 \right) v \right)^{\frac{\left( r - 1 \right) \left( k - 1 \right) + 1 - \beta_2}{2}}} \mathrm{d}v \right]^{\frac{1}{r}} \mathrm{d}t.
				\end{align}
				It is clear that for $k = 1$, we must have $r < 2$, for otherwise, the inner integral in Equation \eqref{IntermediateEquationL1LRHn} is not convergent. Now, using the integral form of the hypergeometric function given in Equation \eqref{IntegralForm2F1}, we get
				\begin{align*}
					&\| R_kf \|_{L^r_{\beta_1, \beta_2} \left( \mathbb{H}^n \right)} \\
					&\leq C \int\limits_{1}^{\infty} \left| \tilde{f} \left( t \right) \right| \left( t^2 - 1 \right)^{\frac{\beta_1 + n - k}{2r} + \frac{k}{2} - 1} \left[ {}_2F_1 \left( \frac{\left( r - 1 \right) \left( k - 1 \right) + 1 - \beta_2}{2}, \frac{\beta_1 + n - k}{2}; \frac{\beta_1 + n - k}{2} + r \left( \frac{k}{2} - 1 \right) + 1; 1 - t^2 \right) \right]^{\frac{1}{r}} \mathrm{d}t.
				\end{align*}
				Next, we use the transformation of the hypergeometric function given in Equation \eqref{Transformation2F1} \eqref{2F1Transformation1} and get
				\begin{align*}
					\| R_kf \|_{L^r_{\beta_1, \beta_2} \left( \mathbb{H}^n \right)} &\leq C \int\limits_{1}^{\infty} \left| \tilde{f} \left( t \right) \right| \left( t^2 - 1 \right)^{\frac{\beta_1 + n - k}{2r} + \frac{k}{2} - 1} t^{\frac{\beta_2 - \left( r - 1 \right) \left( k - 1 \right) - 1}{r}} \times \\
					&\left[ {}_2F_1 \left( \frac{\left( r - 1 \right) \left( k - 1 \right) + 1 - \beta_2}{2}, 1 + r \left( \frac{k}{2} - 1 \right); \frac{\beta_1 + n - k}{2} + r \left( \frac{k}{2} - 1 \right) + 1; 1 - \frac{1}{t^2} \right) \right]^{\frac{1}{r}} \mathrm{d}t.
				\end{align*}
				Now, by substituting $t = \cosh x$, we get
				\begin{align*}
					\| R_kf \|_{L^r_{\beta_1, \beta_2} \left( \mathbb{H}^n \right)} &\leq C \int\limits_{0}^{\infty} \left| \tilde{f} \left( \cosh x \right) \right| \sinh^{\frac{\beta_1 + n - k}{r} + k - 1} x \cosh^{\frac{\beta_2 - \left( r - 1 \right) \left( k - 1 \right) - 1}{r}} x \times \\
					&\left[ {}_2F_1 \left( \frac{\left( r - 1 \right) \left( k - 1 \right) + 1 - \beta_2}{2}, 1 + r \left( \frac{k}{2} - 1 \right); \frac{\beta_1 + n - k}{2} + r \left( \frac{k}{2} - 1 \right) + 1; \tanh^2 x \right) \right]^{\frac{1}{r}} \mathrm{d}x.
				\end{align*}
				From our assumption, we have $\frac{\beta_1 + \beta_2 + n - 1}{r} \leq \alpha_1 + \alpha_2 + n - 1	= k - 1$. Therefore,
				$$\frac{\beta_1 + n - k}{2} + \frac{\beta_2 - \left( r - 1 \right) \left( k - 1 \right) - 1}{2} = \frac{\beta_1 + \beta_2 + n - rk + r - 2}{2} < - \frac{\left( r - 1 \right) \left( k - 1 \right) + 1}{2} < 0.$$
				Therefore, from Equation \eqref{Behaviour2F1Case3} of Theorem \ref{Behaviour2F1}, we conclude that {there is a $C > 0$}
				$$\lim\limits_{x \rightarrow \infty} \frac{{}_2F_1 \left( \frac{\left( r - 1 \right) \left( k - 1 \right) + 1 - \beta_2}{2}, 1 + r \left( \frac{k}{2} - 1 \right); \frac{\beta_1 + n - k}{2} + r \left( \frac{k}{2} - 1 \right) + 1; \tanh^2 x \right)}{\cosh^{\left( r - 1 \right) \left( k - 1 \right) + 1 - \beta_1 - \beta_2 - n + k} x} = C.$$
				We also know that hypergeometric function is continuous and hence bounded on any compact set. Therefore, for sufficiently large $R > 0$, we have,
				\begin{align*}
					\| R_kf \|_{L^r_{\beta_1, \beta_2} \left( \mathbb{H}^n \right)} &\leq C \left[ \int\limits_{0}^{R} \left| \tilde{f} \left( \cosh x \right) \right| \sinh^{\alpha_1 + n - 1} x \cosh^{\alpha_2} x \sinh^{\frac{\beta_1 + n - k}{r} + k - n - \alpha_1} x \cosh^{ \frac{\beta_2 - \left( r - 1 \right) \left( k - 1 \right) - 1}{r} - \alpha_2} x \ \mathrm{d}x \right. \\
					&+\left. \int\limits_{R}^{\infty} \left| \tilde{f} \left( \cosh x \right) \right| \sinh^{\alpha_1 + n - 1} x \cosh^{\alpha_2} x \sinh^{\frac{\beta_1 + n - k}{r} + k - n - \alpha_2} x \cosh^{- \frac{\beta_1 + n - k}{r} - \alpha_2} x \ \mathrm{d}x \right] \\
					&= C \left[ \int\limits_{0}^{R} \left| \tilde{f} \left( \cosh x \right) \right| \sinh^{\alpha_1 + n - 1} x \cosh^{\alpha_2} x \sinh^{\frac{\beta_1 + n - k}{r} + k - n - \alpha_1} x \cosh^{ \frac{\beta_2 - \left( r - 1 \right) \left( k - 1 \right) - 1}{r} - \alpha_2} x \ \mathrm{d}x \right. \\
					&+\left. \int\limits_{R}^{\infty} \left| \tilde{f} \left( \cosh x \right) \right| \sinh^{\alpha_1 + n - 1} x \cosh^{\alpha_2} x \tanh^{\frac{\beta_1 + n - k}{r} + k - n - \alpha_2} x \cosh^{k - n - \alpha_1 - \alpha_2} x \ \mathrm{d}x \right].
				\end{align*}
				Now, we notice that under the assumptions of this result, we have $\frac{\beta_1 + n - k}{r} + k - n - \alpha_1 \geq 0$. Therefore, for $x \in \left( 0, R \right)$, we have $\sinh^{\frac{\beta_1 + n - k}{r} + k - n - \alpha_1} x \leq C$ {and $\cosh^{\frac{\beta_2 - \left( r - 1 \right) \left( k - 1 \right) - 1}{r} - \alpha_2} x \leq C$}. Also, when $x \in \left( 0, R \right)$, any power of $\cosh$ is bounded above. Further, we have for $x \in \left( R, \infty \right)$, $\tanh^{\frac{\beta_1 + n - k}{r} + k - n - \alpha_1} x \leq 1$. Hence, we get
				$$\| R_kf \|_{L^r_{\beta_1, \beta_2} \left( \mathbb{H}^n \right)} \leq C \| f \|_{L^1_{\alpha_1, \alpha_2} \left( \mathbb{H}^n \right)}.$$
				This completes the proof!}
			\end{proof}
			\begin{remark}
				\normalfont
				It was proved in \cite{KumarRayWE} that the condition $\frac{1}{p} - \frac{1}{r} \leq \frac{1}{2}$ is necessary for the boundedness of the $X$-Ray transform on the Euclidean case. However, a similar analysis cannot be done here for all admissible values of $\beta_1, \beta_2$. Instead, when $\frac{\beta_2 + k - 1}{r} - \frac{\alpha_2 - 1}{p} = k = \frac{\alpha_1 + n}{p} - \frac{\beta_1 + n - k}{r}$, we can proceed as in \cite{KumarRayWE}, and check that the condition $\frac{1}{p} - \frac{1}{r} \leq \frac{1}{2}$ is necessary for the boundedness of the $X$-ray transform. {However, in other cases, the necessity of the condition $\frac{1}{p} - \frac{1}{r} \leq \frac{1}{2}$ is unknown to us.}
				
				We now {argue} that we cannot expect the boundedness of the {$X$-ray} transform for the case $p = 1$ and $r = 2$. {This is, indeed, easily seen from Inequality \eqref{NormR1HnFI} and Example \ref{L1L2NotPossibleIHalfPlus}.}
			\end{remark}
			The analysis done above almost completely answers Question \ref{LPLRHnQuestion}. We now move on to answer Question \ref{LPLInfinityHnQuestion}. To answer this question, let us once again begin by first obtaining the necessary conditions.
			\begin{theorem}
				\label{NecessaryLPLInfinityHn}
				Let $\alpha_1 + \alpha_2 > k - n$ {and} $1 \leq p < \frac{\alpha_1 + \alpha_2 + n - 1}{k - 1}$ {or $\alpha_1 + \alpha_2 = k - n$ and $p = 1$}. Then, Inequality \eqref{RequiredLPLInfinityHn} holds for all radial functions in $L^p_{\alpha_1, \alpha_2} \left( \mathbb{H}^n \right)$ only if
				\begin{equation}
					\label{LPLInfinityHnNecessaryConditions}
					{\gamma_1 \geq \max \left\lbrace 0, \frac{\alpha_1 + n}{p} - k \right\rbrace \text{ and } \gamma_1 + \gamma_2 \leq \frac{\alpha_1 + \alpha_2 + n - 1}{p}.}
				\end{equation}
			\end{theorem}
			\begin{proof}
				We prove this result in two cases.
				
				\textbf{Case I:} First, let us assume that $\alpha_1 > -n$. In this case, we test Inequality \eqref{RequiredLPLInfinityHn} against the function {$f_{\lambda}: \mathbb{H}^n \rightarrow \mathbb{C}$, given by} $f_{\lambda} \left( x \right) = \chi_{B \left( 0, \lambda \right)} \left( x \right) \cosh d \left( 0, x \right)${, for $\lambda > 0$}. We have seen in Equation \eqref{LPMassHn1} that $f_{\lambda} \in L^p_{\alpha_1, \alpha_2} \left( \mathbb{H}^n \right)$. Using Equation \eqref{KPlaneBallCoshHnEquation}, we have,
				\begin{equation}
					\label{LPLInfinityHnNecessaryCaseIEquationRkf}
					\sinh^{\gamma_1} d \left( 0, \xi \right) \cosh^{\gamma_2} d \left( 0, \xi \right) R_kf_{\lambda} \left( \xi \right) = C \sinh^{\gamma_1} d \left( 0, \xi \right) \cosh^{\gamma_2 - k + 1} d \left( 0, \xi \right) \sinh^k \lambda \left( 1 - \frac{\sinh^2 d \left( 0, \xi \right)}{\sinh^2 \lambda} \right)^{\frac{k}{2}} {\chi_{\left( 0, \lambda \right)} \left( d \left( 0, \xi \right) \right)}.
				\end{equation}
				It is now clear that if $\gamma_1 < 0$, then $\| \sinh^{\gamma_1} d \left( 0, \cdot \right) \cosh^{\gamma_2} d \left( 0, \cdot \right) R_kf_{\lambda} \|_{L^{\infty} \left( \Xi_k \left( \mathbb{H}^n \right) \right)} = +\infty$.
				This gives that for Inequality \eqref{RequiredLPLInfinityHn}, $\gamma_1 \geq 0$ is {necessary}.
				
				Also, by choosing $d \left( 0, \xi \right)$ such that $\sinh^2 d \left( 0, \xi \right) = \frac{1}{2} \sinh^2 \lambda$ in Equation {\eqref{LPLInfinityHnNecessaryCaseIEquationRkf}}, we get
				\begin{equation}
					\label{RKEstimateInfinityHn}
					\| \sinh^{\gamma_1} d \left( 0, \cdot \right) \cosh^{\gamma_2} d \left( 0, \cdot \right) R_kf_{\lambda} \|_{L^{\infty} \left( \Xi_k \left( \mathbb{H}^n \right) \right)} \geq C \sinh^{\gamma_1 + k} \lambda \left( 1 + \cosh^2 \lambda \right)^{\frac{\gamma_2 - k + 1}{2}}.
				\end{equation}
				Now, using Inequalities \eqref{RKEstimateInfinityHn} and \eqref{LPMassHn1} in Inequality \eqref{RequiredLPLInfinityHn}, we get
				$$\sinh^{\gamma_1 + k} \lambda \left( 1 + \cosh^2 \lambda \right)^{\frac{\gamma_2 - k + 1}{2}} \leq C \begin{cases}
												\lambda^{\frac{\alpha_1 + n}{p}}, & \lambda \text{ is small}. \\
												e^{\left( \frac{\alpha_1 + \alpha_2 + n - 1}{p} + 1 \right) \lambda}, & \lambda \text{ is large}.
											\end{cases}$$
				Now, as $\lambda \rightarrow 0$, the above inequality is reduced to
				$$\lambda^{\gamma_1 + k} \leq C \lambda^{\frac{\alpha_1 + n}{p}},$$
				which forces $\gamma_1 \geq \frac{\alpha_1 + n}{p} - k$.
				
				On the other hand, as $\lambda \rightarrow \infty$, we get
				$$e^{\left( \gamma_1 + \gamma_2 + 1 \right) \lambda} \leq C e^{\left( \frac{\alpha_1 + \alpha_2 + n - 1}{p} + 1 \right) \lambda},$$
				which forces $\gamma_1 + \gamma_2 \leq \frac{\alpha_1 + \alpha_2 + n - 1}{p}$.
				
				\textbf{Case II:} Let us now consider $\alpha_1 \leq -n$ {and} the function $f_{\lambda} \left( x \right) = \chi_{B \left( 0, \lambda \right)} \left( x \right) \sinh^{- \frac{\alpha_1}{p}} d \left( 0, x \right) \cosh d \left( 0, x \right)$. We have seen in Equation \eqref{LPMassHn2} that $f_{\lambda} \in L^p_{\alpha_1, \alpha_2} \left( \mathbb{H}^n \right)$. From Equation \eqref{KPlaneBallSinhCoshHnEquation}, we have
				\begin{align*}
					&\sinh^{\gamma_1} d \left( 0, \xi \right) \cosh^{\gamma_2} d \left( 0, \xi \right) R_kf_{\lambda} \left( \xi \right) \\
					&= C \sinh^{\gamma_1} d \left( 0, \xi \right) \cosh^{\gamma_2 - k + 1} d \left( 0, \xi \right) \sinh^{k - \frac{\alpha_1}{p}} \lambda \left( 1 - \frac{\sinh^2 d \left( 0, \xi \right)}{\sinh^2 \lambda} \right)^{\frac{k}{2}} {}_2F_1 \left( \frac{\alpha_1}{2p}, 1; 1 - \frac{k}{2}; 1 - \frac{\sinh^2 d \left( 0, \xi \right)}{\sinh^2 \lambda} \right).
				\end{align*}
				Again, it is clear that $\gamma_1 \geq 0$ is necessary, for otherwise we have $\| \sinh^{\gamma_1} d \left( 0, \cdot \right) \cosh^{\gamma_2} d \left( 0, \cdot \right) R_kf_{\lambda} \|_{L^{\infty} \left( \Xi_k \left( \mathbb{H}^n \right) \right)} = +\infty$.
				
				Also, by choosing $d \left( 0, \xi \right)$ such that $\sinh^2 d \left( 0, \xi \right) = \frac{1}{2} \sinh^2 \lambda$, we get
				\begin{equation}
					\label{RkEstimateInfinityHn2}
					\| \sinh^{\gamma_1} d \left( 0, \cdot \right) \cosh^{\gamma_2} d \left( 0, \cdot \right) R_kf_{\lambda} \|_{L^{\infty} \left( \Xi_k \left( \mathbb{H}^n \right) \right)} \geq C \sinh^{\gamma_1 + k - \frac{\alpha_1}{p}} {\lambda} \left( 1 + \cosh^2 \lambda \right)^{\frac{\gamma_2 - k + 1}{2}}.
				\end{equation}
				By using Inequalities \eqref{RkEstimateInfinityHn2} and \eqref{LPMassHn2} in Inequality \eqref{RequiredLPLInfinityHn}, we get
				$$\sinh^{\gamma_1 + k - \frac{\alpha_1}{p}} {\lambda} \left( 1 + \cosh^2 \lambda \right)^{\frac{\gamma_2 - k + 1}{2}} \leq C \begin{cases}
																	\lambda^{\frac{n}{p}}, & \lambda \text{ is small}. \\
																	e^{\left( \frac{\alpha_2 + n - 1}{p} + 1 \right) \lambda}, & \lambda \text{ is large}.
																\end{cases}$$
				The necessary conditions are now obtained by taking $\lambda \rightarrow 0$ and $\lambda \rightarrow \infty$, as in Case I.
			\end{proof}
			{Now we} check the sufficiency of the conditions stated in Theorem \ref{NecessaryLPLInfinityHn}.
			\begin{theorem}
				\label{SufficientLPLInfinityHn}
				Let $\alpha_1 + \alpha_2 > k - n$ {and} $1 \leq p < \frac{\alpha_1 + \alpha_2 + n - 1}{k - 1}$ {or $\alpha_1 + \alpha_2 = k - n$ and $p = 1$}. {Let $\gamma_1, \gamma_2 \in \mathbb{R}$. We} assume that if {$p = \frac{\alpha_1 + n}{k}$, then $\gamma_1 > 0$}.
				\begin{enumerate}
					\item[\mylabel{LPLInfinityHnWeightedA}{(A)}] For $k \geq 2$, Inequality \eqref{RequiredLPLInfinityHn} holds for all radial functions in $L^p_{\alpha_1, \alpha_2} \left( \mathbb{H}^n \right)$ {if and only if $\gamma_1$ and $\gamma_2$ satisfy the conditions of Equation \eqref{LPLInfinityHnNecessaryConditions}}.
					\item[\mylabel{LPLInfinityHnWeightedB}{(B)}] For $k = 1$, Inequality \eqref{RequiredLPLInfinityHn} holds for all radial functions in $L^p_{\alpha_1, \alpha_2} \left( \mathbb{H}^n \right)$ {if and only if $\gamma_1$ and $\gamma_2$ satisfy the conditions of Equation \eqref{LPLInfinityHnNecessaryConditions} and} $p > 2$.
				\end{enumerate}
			\end{theorem}
			\begin{proof}~
				\begin{enumerate}
					\item[\mylabel{LPLInfinityHnWeightedAProof}{(A)}] Let us prove the result for $k \geq 2$ { and $p \neq 1$}. From Equation \eqref{KPlaneTransformRadialHnEquation2}, we have by an application of H\"{o}lder's inequality,
					\begin{align*}
						\left| \sinh^{\gamma_1} d \left( 0, \xi \right) \cosh^{\gamma_2} d \left( 0, \xi \right) R_kf \left( \xi \right) \right| &\leq C \sinh^{\gamma_1} d \left( 0, \xi \right) \cosh^{\gamma_2 - 1} d \left( 0, \xi \right) \int\limits_{d \left( 0, \xi \right)}^{\infty} \left| \tilde{f} \left( \cosh t \right) \right| \left( 1 - \frac{\tanh^2 {d \left( 0, \xi \right)}}{\tanh^2 {t}} \right)^{\frac{k}{2} - 1} \sinh^{k - 1} t \ \mathrm{d}t \\
						&\leq C \sinh^{\gamma_1} d \left( 0, \xi \right) \cosh^{\gamma_2 - 1} d \left( 0, \xi \right) \| f \|_{L^p_{\alpha_1, \alpha_2} \left( \mathbb{H}^n \right)} I^{\frac{1}{p'}},
					\end{align*}
					where,
					\begin{align*}
						I &= \int\limits_{d \left( 0, \xi \right)}^{\infty} \sinh^{p' \left( k - \frac{\alpha_1 + n}{p} \right) - 1} t \cosh^{- \frac{p' \alpha_2}{p}} {t} \left( 1 - \frac{\tanh^2 {d \left( 0, \xi \right)}}{\tanh^2 {t}} \right)^{p' \left( \frac{k}{2} - 1 \right)} \mathrm{d}t \\
						&= \cosh^{-2p{'} \left( \frac{k}{2} - 1 \right)} d \left( 0, \xi \right) \int\limits_{d \left( 0, \xi \right)}^{\infty} \sinh^{p' \left( 2 - \frac{\alpha_1 + n}{p} \right) - 1} t \cosh^{- \frac{p' \alpha_2}{p}} {t} \left( \cosh^2 t - \cosh^2 d \left( 0, \xi \right) \right)^{p' {(} \frac{k}{2} - 1 {)}} \mathrm{d}t.
					\end{align*}
					We now substitute $\cosh^2 t = \frac{\cosh^2 d \left( 0, \xi \right)}{y}$ and simplify to obtain
					\begin{align*}
						I &= C \cosh^{p' \left( 2 - \frac{\alpha_1 + \alpha_2 + n}{p} \right) {- 1}} d \left( 0, \xi \right) \int\limits_{0}^{1} \frac{y^{p' \left( \frac{\alpha_1 + \alpha_2 + n}{2p} - \frac{k}{2} \right) - \frac{1}{2}} \left( 1 - y \right)^{p' \left( \frac{k}{2} - 1 \right)}}{\left( 1 - y \sech^2 d \left( 0, \xi \right) \right)^{1 - p' \left( 1 - \frac{\alpha_1 + n}{2p} \right)}} \mathrm{d}y \\
						&= C \cosh^{p' \left( 2 - \frac{\alpha_1 + \alpha_2 + n}{p} \right) - 1} d \left( 0, \xi \right) \\
						&{}_2F_1 \left( 1 - p' \left( 1 - \frac{\alpha_1 + n}{2p} \right), \frac{{p'}}{2} \left( \frac{\alpha_1 + \alpha_2 + n - 1}{p} - k + 1 \right); p' \left( \frac{k}{2} - 1 \right) + \frac{{p'}}{2} \left( \frac{\alpha_1 + \alpha_2 + n - 1}{p} - k \right) {+ 1}; \sech^2 d \left( 0, \xi \right) \right).
					\end{align*}
					In the last step, we have used Equation \eqref{IntegralForm2F1}. That is, we have,
					\begin{align}
						\label{LPLInfinityEstimateHnSecondLast}
						&\left| \sinh^{\gamma_1} d \left( 0, \xi \right) \cosh^{\gamma_2} d \left( 0, \xi \right) R_kf \left( \xi \right) \right| \\
						&\leq C \sinh^{\gamma_1} d \left( 0, \xi \right) \cosh^{\gamma_2 - \frac{\alpha_1 + \alpha_2 + n - 1}{p}} d \left( 0, \xi \right) \| f \|_{L^p_{\alpha_1, \alpha_2} \left( \mathbb{H}^n \right)} \times \nonumber \\
						&\left| {}_2F_1 \left( 1 - p' \left( 1 - \frac{\alpha_1 + n}{2p} \right), \frac{1}{2} \left( \frac{\alpha_1 + \alpha_2 + n - 1}{p} - k + 1 \right); p' \left( \frac{k}{2} - 1 \right) + \frac{3}{2} + \frac{1}{2} \left( \frac{\alpha_1 + \alpha_2 + n - 1}{p} - k \right); \sech^2 d \left( 0, \xi \right) \right) \right|^{\frac{1}{p'}}. \nonumber
					\end{align}
					The behaviour of the hypergeometric function depends on the value of $p$, as given by the conditions of Theorem \ref{Behaviour2F1}.
					
					\textbf{Case I:} Let us assume $p > \frac{\alpha_1 + n}{k}$. In this case, we have from Equation \eqref{Behaviour2F1Case1} that the hypergeometric function is bounded. Therefore, for any $\xi \in \Xi_k \left( \mathbb{H}^n \right)$, we have,
					$$\left| \sinh^{\gamma_1} d \left( 0, \xi \right) \cosh^{\gamma_2} d \left( 0, \xi \right) R_kf \left( \xi \right) \right| \leq C \sinh^{\gamma_1} d \left( 0, \xi \right) \cosh^{\gamma_2 - \frac{\alpha_1 + \alpha_2 + n - 1}{p}} d \left( 0, \xi \right) \| f \|_{L^p_{\alpha_1, \alpha_2} \left( \mathbb{H}^n \right)}.$$
					When $d \left( 0, \xi \right)$ is close to $0$, the powers of hyperbolic sine and cosine are bounded (since $\gamma_1 \geq 0$ and $\cosh d \left( 0, \xi \right)$ is close to $1$). On the other hand, when $d \left( 0, \xi \right)$ is away from zero, both the hyperbolic trigonometric functions behave like exponential and produce $e^{\left( \gamma_1 + \gamma_2 - \frac{\alpha_1 + \alpha_2 + n - 1}{p} \right) d \left( 0, \xi \right)}$. However, from our assumptions, we have $\gamma_1 + \gamma_2 \leq \frac{\alpha_1 + \alpha_2 + n - 1}{p}$, making this term bounded. Hence, Inequality \eqref{RequiredLPLInfinityHn} holds in this case.
					
					\textbf{Case II:} For $p = \frac{\alpha_1 + n}{k}$, from Equation \eqref{Behaviour2F1Case2}, we know that the hypergeometric function on the right hand side of Inequality \eqref{LPLInfinityEstimateHnSecondLast} behaves as $- \ln \tanh^2 d \left( 0, \xi \right)$ as $d \left( 0, \xi \right) \rightarrow 0$. Also, from our assumptions, we must have $\gamma_1 > 0$. Therefore, when $d \left( 0, \xi \right)$ is near $0$, the term $\sinh^{\gamma_1} d \left( 0, \xi \right) \left( - \ln \tanh^2 d \left( 0, \xi \right) \right)$ is bounded. When $d \left( 0, \xi \right)$ is far from $0$, the terms on the right side of Inequality \eqref{LPLInfinityEstimateHnSecondLast} {that are powers of $\sinh$ and $\cosh$ are} bounded anyway due to the argument of Case I. {The hypergeometric function is bounded because $\sech d \left( 0, \xi \right)$ is near $0$ when $d \left( 0, \xi \right)$ is away from zero.} Hence, Inequality \eqref{RequiredLPLInfinityHn} holds.
					
					\textbf{Case III:} Lastly, let us take the case when $p < \frac{\alpha_1 + n}{k}$. Here, from Equation \eqref{Behaviour2F1Case3}, we know that the hypergeometric function behaves like $\tanh^{p' \left( k - \frac{\alpha_1 + n}{p} \right)} d \left( 0, \xi \right)$, for small values of $d \left( 0, \xi \right)$. When $d \left( 0, \xi \right)$ is away from $0$, the hypergometric function makes no change in the analysis as done in Case I. However, when $d \left( 0, \xi \right)$ is close to $0$, we have
					$$\left| \sinh^{\gamma_1} d \left( 0, \xi \right) \cosh^{\gamma_2} d \left( 0, \xi \right) R_kf \left( \xi \right) \right| \leq C \sinh^{\gamma_1 + k - \frac{\alpha_1 + n}{p}} d \left( 0, \xi \right) \cosh^{\gamma_2 - \frac{\alpha_1 + \alpha_2 + n - 1}{p} - k + \frac{\alpha_1 + n}{p}} d \left( 0, \xi \right) \| f \|_{L^p_{\alpha_1, \alpha_2} \left( \mathbb{H}^n \right)}.$$
					Again, the hyperbolic cosine is bounded. The power of hyperbolic sine is non-negative, since $\gamma_1 \geq \frac{\alpha_1 + n}{p} - k$. Hence, when $d \left( 0, \xi \right)$ is near $0$, the hyperbolic trigonometric functions on the right side of the inequality are bounded and we get Inequality \eqref{RequiredLPLInfinityHn}.
					
					{To complete the proof of \ref{LPLInfinityHnWeightedAProof} we now show the result for $p = 1$. First, we observe that by using Equation \eqref{KPlaneTransformRadialHnEquation3} and the fact that $k \geq 2$, we have,
					\begin{equation}
						\label{L1LInfinityHnInitialEstimate}
						\begin{aligned}
							&\left| \sinh^{\gamma_1} d \left( 0, \xi \right) \cosh^{\gamma_2} d \left( 0, \xi \right) R_kf \left( \xi \right) \right| \\
							&\leq C \sinh^{\gamma_1} d \left( 0, \xi \right) \cosh^{\gamma_2 - k + 1} d \left( 0, \xi \right) \int\limits_{d \left( 0, \xi \right)}^{\infty} \left| \tilde{f} \left( \cosh t \right) \right| \sinh^{\alpha_1 + n - 1} t \cosh^{\alpha_2} t \sinh^{k - n - \alpha_1} t \cosh^{- \alpha_2} t \ \mathrm{d}t.
						\end{aligned}
					\end{equation}
					We break the analysis in two cases.}
					
					{\textbf{Case I}: First, let us consider $d \left( 0, \xi \right) \geq 1$. Then, we have that all the hyperbolic trigonometric functions in the above inequality behave like exponentials. Therefore,
					$$\left| \sinh^{\gamma_1} d \left( 0, \xi \right) \cosh^{\gamma_2} d \left( 0, \xi \right) R_kf \left( \xi \right) \right| \leq C e^{\left( \gamma_1 + \gamma_1 - k + 1 \right) d \left( 0, \xi \right)} \int\limits_{d \left( 0, \xi \right)}^{\infty} \left| \tilde{f} \left( \cosh t \right) \right| \sinh^{\alpha_1 + n - 1} t \cosh^{\alpha_2} t \ e^{\left( k - n - \alpha_1 - \alpha_2 \right) t} \mathrm{d}t.$$
					We know that $\alpha_1 + \alpha_2 \geq k - n$, so that $e^{\left( k - n - \alpha_1 - \alpha_2 \right) t} \leq e^{\left( k - n - \alpha_1 - \alpha_2 \right) d \left( 0, \xi \right)}$. Hence, we get
					$$\left| \sinh^{\gamma_1} d \left( 0, \xi \right) \cosh^{\gamma_2} d \left( 0, \xi \right) R_kf \left( \xi \right) \right| \leq C e^{\left( \gamma_1 + \gamma_2 - \alpha_1 - \alpha_2 - n + 1 \right) d \left( 0, \xi \right)} \int\limits_{d \left( 0, \xi \right)}^{\infty} \left| \tilde{f} \left( \cosh t \right) \right| \sinh^{\alpha_1 + n - 1} t \cosh^{\alpha_2} t \ \mathrm{d}t.$$
					Further, we know from the necessary conditions of Equation \eqref{LPLInfinityHnNecessaryConditions} that $\gamma_1 + \gamma_2 \leq \alpha_1 + \alpha_2 + n - 1$. Thus, we have,
					$$\left| \sinh^{\gamma_1} d \left( 0, \xi \right) \cosh^{\gamma_2} d \left( 0, \xi \right) R_kf \left( \xi \right) \right| \leq C \| f \|_{L^1_{\alpha_1, \alpha_2} \left( \mathbb{H}^n \right)}.$$}
					
					{\textbf{Case II}: Let us now consider $d \left( 0, \xi \right) < 1$. From Inequality \eqref{L1LInfinityHnInitialEstimate}, we have,
					\begin{align*}
						&\left| \sinh^{\gamma_1} d \left( 0, \xi \right) \cosh^{\gamma_2} d \left( 0, \xi \right) R_kf \left( \xi \right) \right| \\
							&\leq C \sinh^{\gamma_1} d \left( 0, \xi \right) \cosh^{\gamma_2 - k + 1} d \left( 0, \xi \right) \int\limits_{d \left( 0, \xi \right)}^{1} \left| \tilde{f} \left( \cosh t \right) \right| \sinh^{\alpha_1 + n - 1} t \cosh^{\alpha_2} t \sinh^{k - n - \alpha_1} t \cosh^{- \alpha_2} t \ \mathrm{d}t \\
							&+ C \sinh^{\gamma_1} d \left( 0, \xi \right) \cosh^{\gamma_2 - k + 1} d \left( 0, \xi \right) \int\limits_{1}^{\infty} \left| \tilde{f} \left( \cosh t \right) \right| \sinh^{\alpha_1 + n - 1} t \cosh^{\alpha_2} t \sinh^{k - n - \alpha_1} t \cosh^{- \alpha_2} t \ \mathrm{d}t
					\end{align*}
					The second term in the above inequality is bounded above by a multiple $\| f \|_{L^1_{\alpha_1, \alpha_2} \left( \mathbb{H}^n \right)}$ due to Case I and the fact that both the hyperbolic trigonometric functions are bounded when $d \left( 0, \xi \right) < 1$. Therefore, we only need to check the boundedness of the first term. To do so, we observe that $\cosh^{\gamma_2 - k + 1} d \left( 0, \xi \right) \leq C$ and $\cosh^{- \alpha_2} t \leq C$, when $d \left( 0, \xi \right) \leq t \leq 1$. Let us now consider the following situations.}
					
					{\textbf{Case II(A)}: First, let us assume that $\alpha_1 > k - n$. Then, for $t \geq d \left( 0, \xi \right)$, we have, $\sinh^{k - n - \alpha_1} t \leq \sinh^{k - n - \alpha_1} d \left( 0, \xi \right)$. Consequently,
					\begin{align*}
						\sinh^{\gamma_1} d \left( 0, \xi \right) \int\limits_{d \left( 0, \xi \right)} \left| \tilde{f} \left( \cosh t \right) \right| \sinh^{\alpha_1 + n - 1} t \cosh^{\alpha_2} t \sinh^{k - n - \alpha_1} t \ \mathrm{d}t &\leq C \sinh^{\gamma_1 + k - n - \alpha_1} d \left( 0, \xi \right) \| f \|_{L^1_{\alpha_1, \alpha_2} \left( \mathbb{H}^n \right)} \\
						&\leq C \| f \|_{L^1_{\alpha_1, \alpha_2} \left( \mathbb{H}^n \right)},
					\end{align*}
					since $\gamma_1 \geq \alpha_1 + n - k$.}
					
					\textbf{Case II(B)}: Now, we assume that $\alpha_1 \leq k - n$. Then, we have for $t \leq 1$, $\sinh^{k - n - \alpha_1} t \leq C$. {Since $\gamma_1 \geq 0$, we have},
					\begin{align*}
						\sinh^{\gamma_1} d \left( 0, \xi \right) \int\limits_{d \left( 0, \xi \right)} \left| \tilde{f} \left( \cosh t \right) \right| \sinh^{\alpha_1 + n - 1} t \cosh^{\alpha_2} t \sinh^{k - n - \alpha_1} t \ \mathrm{d}t &\leq C \sinh^{\gamma_1} d \left( 0, \xi \right) \| f \|_{L^1_{\alpha_1, \alpha_2} \left( \mathbb{H}^n \right)} \\
						&\leq C \| f \|_{L^1_{\alpha_1, \alpha_2} \left( \mathbb{H}^n \right)},
					\end{align*}
					
					{From Cases I and II, we conclude that for any $\xi \in \Xi_k \left( \mathbb{H}^n \right)$, we have,
					$$\left| \sinh^{\gamma_1} d \left( 0, \xi \right) \cosh^{\gamma_2} d \left( 0, \xi \right) R_kf \left( \xi \right) \right| \leq C \| f \|_{L^1_{\alpha_1, \alpha_2} \left( \mathbb{H}^n \right)}.$$}
					
					This completes the proof for \ref{LPLInfinityHnWeightedA}.
					\item[\mylabel{LPLInfinityHnWeightedBProof}{(B)}] The proof for the case $k = 1$ and $p > 2$ is verbatim to that of \ref{LPLInfinityHnWeightedAProof}. The condition $p > 2$ is required for the convergence of the hypergeometric integral in $I$. {The necessity of the condition $p > 2$ is shown in Example \ref{CounterExampleLpLinfinityHn} and Remark \ref{L2LInfinityHnNotPossible}.}
				\end{enumerate}
			\end{proof}
			We now {show} that Theorem \ref{SufficientLPLInfinityHn} cannot hold when $k = 1$ and $p < 2$. The example is motivated from its Euclidean counterpart mentioned in {Example \ref{EuclideanAnnulusCounterExample}}.
			\begin{example}
				\label{CounterExampleLpLinfinityHn}
				\normalfont
				Let $1 < a < b < \infty$ and $f_{a, b} \left( x \right) = \chi_{\left( a, b \right)} \left( \cosh d \left( 0, x \right) \right)$. Clearly, $f_{a, b}$ is the characteristic function of the hyperbolic annulus of inner radius $\cosh^{-1} a$ and outer radius $\cosh^{-1} b$. We first observe that a necessary condition for Inequality \eqref{RequiredLPLInfinityHn} to hold is the following Lorentz norm inequality.
				\begin{equation}
					\label{RequiredLorentzNormLPLInfinityHn}
					\| \sinh^{\gamma_1} d \left( 0, \cdot \right) \cosh^{\gamma_2} d \left( 0, \cdot \right) R_1f_{a, b} \|_{L^{\infty} \left( \Xi_1 \left( \mathbb{H}^n \right) \right)} \leq C \| f_{a, b} \|_{L^{p, 1}_{\alpha_1, \alpha_2} \left( \mathbb{H}^n \right)}.
				\end{equation}
				We test Inequality \eqref{RequiredLorentzNormLPLInfinityHn} against the function $f_{a, b}$. Later, we shrink the annulus by taking $b \rightarrow a$, and obtain a necessary condition for Inequality \eqref{RequiredLorentzNormLPLInfinityHn} to hold.
				
				First, we see that by using $\cosh t = x$,
				\begin{equation}
					\label{LPMassWeightedAnnulusHn}
					\| f_{a, b} \|_{L^{p, 1}_{\alpha_1, \alpha_2} \left( \mathbb{H}^n \right)} = C \left[ \int\limits_{a}^{b} \left( x^2 - 1 \right)^{\frac{\alpha_1 + n}{2} - 1} x^{\alpha_2} \mathrm{d}x \right]^{\frac{1}{p}}.
				\end{equation}
				For the X-Ray transform of $f_{a, b}$, we have from Equation \eqref{KPlaneTransformRadialHnEquation},
				\begin{align*}
					R_1f_{a, b} \left( \xi \right) = C \int\limits_{d \left( 0, \xi \right)}^{\infty} \chi_{\left( a, b \right)} \left( \cosh t \right) \left( \cosh^2 t - \cosh^2 d \left( 0, \xi \right) \right)^{-\frac{1}{2}} \sinh t \ \mathrm{d}t.
				\end{align*}
				By substituting $\cosh t = x \cosh d \left( 0, \xi \right)$ and simplifying, we get
				\begin{align*}
					R_1f_{a, b} \left( \xi \right) &= C \int\limits_{1}^{\infty} \chi_{\left( a, b \right)} \left( x \cosh d \left( 0, \xi \right) \right) \left( x^2 - 1 \right)^{- \frac{1}{2}} \mathrm{d}x \\
					&= C \begin{cases}
								0, & \cosh d \left( 0, \xi \right) \geq b. \\
								\int\limits_{1}^{\frac{b}{\cosh d \left( 0, \xi \right)}} \left( x^2 - 1 \right)^{- \frac{1}{2}}, & a \leq \cosh d \left( 0, \xi \right) < b. \\
								\int\limits_{\frac{a}{\cosh d \left( 0, \xi \right)}}^{\frac{b}{\cosh d \left( 0, \xi \right)}} \left( x^2 - 1 \right)^{- \frac{1}{2}} \mathrm{d}x, & \cosh d \left( 0, {\xi} \right) < a.
							\end{cases}
				\end{align*}
				Therefore, by choosing $\cosh d \left( 0, \xi \right) = a$, we get
				\begin{equation}
					\label{R1EstimateLPLInfinityNecessaryHn}
					\| \sinh^{\gamma_1} d \left( 0, \cdot \right) \cosh^{\gamma_2} d \left( 0, \cdot \right) R_1f_{a, b} \|_{L^{\infty} \left( \Xi_k \left( \mathbb{H}^n \right) \right)} \geq \left( a^2 - 1 \right)^{\frac{\gamma_1}{2}} a^{\gamma_2} \int\limits_{1}^{\frac{b}{a}} \left( x^2 - 1 \right)^{- \frac{1}{2}} \mathrm{d}x.
				\end{equation}
				Plugging Inequality \eqref{R1EstimateLPLInfinityNecessaryHn} and Equation \eqref{LPMassWeightedAnnulusHn} in Inequality \eqref{RequiredLorentzNormLPLInfinityHn}, we get
				$$\left( a^2 - 1 \right)^{\frac{p\gamma_1}{2}} a^{p\gamma_2} \left[ \int\limits_{1}^{\frac{b}{a}} \left( x^2 - 1 \right)^{- \frac{1}{2}} \mathrm{d}x \right]^{p} \leq C \int\limits_{a}^{b} \left( x^2 - 1 \right)^{\frac{\alpha_1 + n}{2} - 1} x^{\alpha_2} \mathrm{d}x.$$
				For a fixed $a > 1$, let us consider the function
				$$H \left( b \right) = \frac{\left( \int\limits_{1}^{\frac{b}{a}} \left( x^2 - 1 \right)^{- \frac{1}{2}} \mathrm{d}x \right)^{p}}{\int\limits_{a}^{b} \left( x^2 - 1 \right)^{\frac{\alpha_1 + n}{2} - 1} x^{\alpha_2} \mathrm{d}x},$$
				defined for $b > a$. It is clear that {Inequality \eqref{RequiredLorentzNormLPLInfinityHn}} fails if $\lim\limits_{b \rightarrow a} H \left( b \right)$ is infinite. We use L'H\^{o}spital rule to {evaluate} the limit.
				\begin{align*}
					\lim\limits_{b \rightarrow a} H \left( b \right) &= \lim\limits_{b \rightarrow a} \frac{\left( \int\limits_{1}^{\frac{b}{a}} \left( x^2 - 1 \right)^{- \frac{1}{2}} \mathrm{d}x \right)^{p}}{\int\limits_{a}^{b} \left( x^2 - 1 \right)^{\frac{\alpha_1 + n}{2} - 1} x^{\alpha_2} \mathrm{d}x} \\
					&= \lim\limits_{b \rightarrow a} \frac{p \left( \int\limits_{1}^{\frac{b}{a}} \left( x^2 - 1 \right)^{- \frac{1}{2}} \mathrm{d}x \right)^{p - 1}}{\sqrt{b^2 - a^2} \left( b^2 - 1 \right)^{\frac{\alpha_1 + n}{2} - 1} b^{\alpha_2}} \\
					&= \lim\limits_{b \rightarrow a} \frac{p \left( p - 1 \right) \left( \int\limits_{1}^{\frac{b}{a}} \left( x^2 - 1 \right)^{- \frac{1}{2}} \mathrm{d}x \right)^{p - 2}}{\left( b^2 - 1 \right)^{\frac{\alpha_1 + n}{2} - 1} b^{\alpha_2 + 1} + \left( \alpha_1 + n - 2 \right) \left( b^2 - a^2 \right) \left( b^2 - 1 \right)^{\frac{\alpha_1 + n }{2} - 1} b^{\alpha_2 + 1} + \alpha_2 \left( b^2 - a^2 \right) \left( b^2 - {1} \right)^{\frac{\alpha_1 + n}{2} - 1} b^{\alpha_2 - 1}},
				\end{align*}
				which is clearly infinite if $p < 2$. That is, Inequality \eqref{RequiredLPLInfinityHn} cannot be expected unless $p \geq 2$.
			\end{example}
			\begin{remark}
				\label{L2LInfinityHnNotPossible}
				\normalfont
				{We now see that Inequality \eqref{RequiredLPLInfinityHn} does not hold for all radial functions $f \in L^2_{\alpha_1, \alpha_2} \left( \mathbb{H}^n \right)$. Let us consider its left-hand-side. We have for a radial function $f: \mathbb{H}^n \rightarrow \mathbb{C}$, from Equation \eqref{KPlaneTransformRadialHnEquation3},
				$$\sinh^{\gamma_1} d \left( 0, \xi \right) \cosh^{\gamma_2} d \left( 0, \xi \right) R_1f \left( \xi \right) = C \sinh^{\gamma_1} d \left( 0, \xi \right) \cosh^{\gamma_2} d \left( 0, \xi \right) \int\limits_{d \left( 0, \xi \right)}^{\infty} \tilde{f} \left( \cosh t \right) \left( \sinh^2 t - \sinh^2 d \left( 0, \xi \right) \right)^{- \frac{1}{2}} \sinh t \ \mathrm{d}t.$$
				By substituting $\sinh^2 t = u$, we obtain,
				\begin{align*}
					\sinh^{\gamma_1} d \left( 0, \xi \right) \cosh^{\gamma_2} d \left( 0, \xi \right) R_1f \left( \xi \right) &= C \sinh^{\gamma_1} d \left( 0, \xi \right) \cosh^{\gamma_2} d \left( 0, \xi \right) \int\limits_{\sinh^2 d \left( 0, \xi \right)}^{\infty} \frac{\tilde{f} \left( \sqrt{1 + u} \right) \left( 1 + u \right)^{- \frac{1}{2}}}{\left( u - \sinh^2 d \left( 0, \xi \right) \right)^{\frac{1}{2}}} \ \mathrm{d}u \\
					&= C \sinh^{\gamma_1} d \left( 0, \xi \right) \cosh^{\gamma_2} d \left( 0, \xi \right) \left( I^{\frac{1}{2}}_- \varphi \right) \left( \sinh^2 d \left( 0, \xi \right) \right),
				\end{align*}
				where,
				\begin{equation}
					\label{PhiL2LinfinityHn}
					\varphi \left( t \right) = \frac{\tilde{f} \left( \sqrt{1 + u} \right)}{\left( 1 + u \right)^{\frac{1}{2}}}.
				\end{equation}
				Further, by writing $\sinh^2 d \left( 0, \xi \right) = x$, we have that
				\begin{equation}
					\label{LInfinityNormR1fWeighted}
					\| \sinh^{\gamma_1} d \left( 0, \cdot \right) \cosh^{\gamma_2} d \left( 0, \cdot \right) R_1f \|_{L^{\infty} \left( \Xi_k \left( \mathbb{H}^n \right) \right)} = C \| x^{\frac{\gamma_1}{2}} \left( 1 + x \right)^{\frac{\gamma_2}{2}} I^{\frac{1}{2}}_- \varphi \|_{L^{\infty} \left( 0, \infty \right)},
				\end{equation}
				where, $\varphi$ is given in Equation \eqref{PhiL2LinfinityHn}. On the other hand, we also have,
				$$\| f \|_{L^2_{\alpha_1, \alpha_2} \left( \mathbb{H}^n \right)} = C \left[ \int\limits_{0}^{\infty} \left| \tilde{f} \left( \cosh t \right) \right|^2 \sinh^{\alpha_1 + n - 1} t \cosh^{\alpha_2} t \ \mathrm{d}t \right]^{\frac{1}{2}}.$$
				Again, by substituting $\sinh^2 t = u$, we get,
				\begin{equation}
					\label{L2NormfWeighted}
					\begin{aligned}
						\| f \|_{L^2_{\alpha_1, \alpha_2} \left( \mathbb{H}^n \right)} &= C \left[ \int\limits_{0}^{\infty} \left| \tilde{f} \left( \sqrt{1 + u} \right) \right|^2 u^{\frac{\alpha_1 + n}{2} - 1} \left( 1 + u \right)^{\frac{\alpha_2 - 1}{2}} \mathrm{d}u \right]^{\frac{1}{2}} \\
						&= C \left[ \int\limits_{0}^{\infty} \left| \frac{\tilde{f} \left( \sqrt{1 + u} \right)}{\left( 1 + u \right)^{\frac{1}{2}}} \right|^2 u^{\frac{\alpha_1 + n}{2} - 1} \left( 1 + u \right)^{\frac{\alpha_2 + 1}{2}} \mathrm{d}u \right]^{\frac{1}{2}} = \| \varphi \|_{L^2_{\frac{\alpha_1 + n}{2} - 1, \frac{\alpha_2 + 1}{2}} \left( 0, \infty \right)}.
					\end{aligned}
				\end{equation}
				From Equations \eqref{LInfinityNormR1fWeighted} and \eqref{L2NormfWeighted}, we easily see that Inequality \eqref{RequiredLPLInfinityHn} is equivalent to the following inequality for all $\varphi \in L^2_{\frac{\alpha_1 + n}{2} - 1, \frac{\alpha_2 + 1}{2}} \left( 0, \infty \right)$.
				\begin{equation}
					\label{InequalityForFIHn}
					\| x^{\frac{\gamma_1}{2}} \left( 1 + x \right)^{\frac{\gamma_2}{2}} I^{\frac{1}{2}}_- \varphi \|_{L^{\infty} \left( 0, \infty \right)} \leq C \| \varphi \|_{L^2_{\frac{\alpha_1 + n}{2} - 1, \frac{\alpha_2 + 1}{2}} \left( 0, \infty \right)}.
				\end{equation}
				However, Inequality \eqref{InequalityForFIHn} is not valid as seen in Proposition \ref{L2LInfinityIHalfPlusNotPossible}.}
			\end{remark}
		\subsection{The Sphere \texorpdfstring{$\mathbb{S}^n$}{}}
			In this section, we ask when is the $k$-plane transform bounded from $L^p_{\alpha_1, \alpha_2} \left( \mathbb{S}^n \right)$ to $L^r_{\beta_1, \beta_2} \left( \Xi_k \left( \mathbb{S}^n \right) \right)$. We have remarked in Section \ref{PreliminariesSection} that it is enough to study the $k$-plane transform of even functions. Consequently, for the functions considered in {our} analysis, we assume that the expressions involved are {defined on $\mathbb{S}^n_+ := \left\lbrace x \in \mathbb{S}^n | \langle x, e_{n + 1} \rangle > 0 \right\rbrace$ and are} extended evenly on $\mathbb{S}^n$. Analogous to the hyperbolic case, we try to answer the following questions.
			\begin{question}
				\label{LPLRSnQuestion}
				What are the admissible values of $\alpha_1, \alpha_2, \beta_1, \beta_2 \in \mathbb{R}$, and ${r \geq p} \geq 1$ such that the following inequality holds for all radial functions in $L^{p}_{\alpha_1, \alpha_2} \left( \mathbb{S}^n \right)$?
				\begin{equation}
					\label{RequiredLPLRSn}
					\| R_kf \|_{L^r_{\beta_1, \beta_2} \left( \Xi_k \left( \mathbb{S}^n \right) \right)} \leq C \| f \|_{L^p_{\alpha_1, \alpha_2} \left( \mathbb{S}^n \right)}.
				\end{equation}
			\end{question}
			\begin{question}
				\label{LPLInfinitySnQuestion}
				What are the admissible values of $\alpha_1, \alpha_2, \gamma_1, \gamma_2 \in \mathbb{R}$, and $p \geq 1$ such that the following inequality holds for all radial functions in $L^{p}_{\alpha_1, \alpha_2} \left( \mathbb{S}^n \right)$?
				\begin{equation}
					\label{RequiredLPLInfinitySn}
					\| \sin^{\gamma_1} d \left( 0, \cdot \right) \cos^{\gamma_2} d \left( 0, \cdot \right) R_kf \|_{L^{\infty} \left( \Xi_k \left( \mathbb{S}^n \right) \right)} \leq C \| f \|_{L^p_{\alpha_1, \alpha_2} \left( \mathbb{S}^n \right)}.
				\end{equation}
			\end{question}
			{We recall that $p > 1 + \alpha_2$ or $p = 1$ when $\alpha_2 = 0$ is required for the existence of the $k$-plane transform on $L^p_{\alpha_1, \alpha_2} \left( \mathbb{S}^n \right)$. That is, outside of these assumptions, Inequalities \eqref{RequiredLPLRSn} and \eqref{RequiredLPLInfinitySn} do not make sense. In the sequel, therefore, we assume these conditions on $p$ and $\alpha_2$.}
			
			{Now,} we begin answering Question \ref{LPLRSnQuestion}. As before, we start by getting the necessary conditions for Inequality \eqref{RequiredLPLRSn} to hold.
			\begin{theorem}
				\label{NecessaryMixed}
				Let $\alpha_1, \alpha_2, {\beta_1}, \beta_2 \in \mathbb{R}$, {and} $p \geq {\max \left\lbrace 1, 1 + \alpha_2 \right\rbrace}$, {with $p > 1 + \alpha_2$ if $\alpha_2 > 0$}. {Let}  $r \geq {p}$. Then, Inequality \eqref{RequiredLPLRSn} holds only if
				\begin{equation}
					\label{LPLRSNNecessaryConditions}
					{\beta_1 > k - n, \ \frac{\alpha_1 + n}{p} - \frac{\beta_1 + n - k}{r} \leq k \text{ and } \frac{\beta_2 + k + 1}{r} \geq \frac{\alpha_2 + 1}{p} \text{ with strict inequality if } p = 1 + \alpha_2.}
				\end{equation}
			\end{theorem}
			\begin{proof}
				We prove this result by testing the desired inequality against certain functions that mimic dilations. However, since in the sphere we lack a notion of ``infinity", we need to consider two such functions: one concentrated at the origin of the sphere, and the other concentrated at the equator.
				
				\textbf{{Case} I:} First, we consider functions concentrated at the origin. Particularly, we test Inequality \eqref{RequiredLPLRSn} against geodesic balls centered at $0$, with a prescribed density. For our purposes, we consider $\lambda$ ``very small". We do so in two cases.
				
				\textbf{Case I{(A)}:} First, let us assume $\alpha_1 > - n$. {Here}, we consider the function $f_{{\lambda}}: \mathbb{S}^n_+ \rightarrow \mathbb{C}$, defined as $f_{\lambda} \left( x \right) = \chi_{B \left( 0, \lambda \right)} \left( x \right) \cos d \left( 0, x \right)$.
				From Equation \eqref{LpNormRadialSn3}, we have
				\begin{align*}
					\| f_{\lambda} \|_{L^p_{\alpha_1, \alpha_2} \left( \mathbb{S}^n \right)} &= C \left[ \int\limits_{0}^{\sin^2 \lambda} \left( 1 - {t} \right)^{\frac{p + \alpha_2 - 1}{2}} {t}^{\frac{\alpha_1 + n}{2} - 1} \mathrm{d}{t} \right]^{\frac{1}{p}}.
				\end{align*}
				By using ${t} = x \sin^2 \lambda$, we get
				\begin{align}
					\| f_{\lambda} \|_{L^p_{\alpha_1, \alpha_2} \left( \mathbb{S}^n \right)} &= C \left[ \int\limits_{0}^{1} \left( 1 - x \sin^2 \lambda \right)^{\frac{\alpha_2 + p - 1}{2}} \sin^{\alpha_1 + n} \lambda \ x^{\frac{\alpha_1 + n}{2} - 1} \mathrm{d}x \right]^{\frac{1}{p}} \nonumber \\
					\label{LPMassSphereOrigin1}
					&= C \sin^{\frac{\alpha_1 + n}{p}} \lambda \left[ {}_2F_1 \left( \frac{1 - p - \alpha_2}{2}, \frac{\alpha_1 + n}{2}; 1 + \frac{\alpha_1 + n}{2}; \sin^2 \lambda \right) \right]^{\frac{1}{p}},
				\end{align}
				where, in the last equality, we have used the integral form of the hypergeometric function (Equation \eqref{IntegralForm2F1}). On the other hand, using Equation \eqref{NormRkSn}, we get
				\begin{align*}
					\| R_kf_{\lambda} \|_{L^r_{\beta_1, \beta_2} \left( \Xi_k \left( \mathbb{S}^n \right) \right)} &= C \left[ \int\limits_{0}^{1} \left| \int\limits_{0}^{s} \frac{\chi_{\left( \cos \lambda, 1 \right)} \left( t \right) t \left( s^2 - t^2 \right)^{\frac{k}{2} - 1}}{s^{k - 1}} \mathrm{d}t \right|^r s^{\beta_2 + k} \left( 1 - s^2 \right)^{\frac{\beta_1 + n - k}{2} - 1} \mathrm{d}s \right]^{\frac{1}{r}} \\
					&= C \left[ \int\limits_{\cos \lambda}^{1} \left| \int\limits_{\cos \lambda}^{s} t \left( s^2 - t^2 \right)^{\frac{k}{2} - 1} \mathrm{d}t \right|^r s^{\beta_2 + k - r \left( k - 1 \right)} \left( 1 - s^2 \right)^{\frac{\beta_1 + n - k}{2} - 1} \mathrm{d}s \right]^{\frac{1}{r}} \\
					&= C \left[ \int\limits_{\cos \lambda}^{1} \left( s^2 - \cos^2 \lambda \right)^{\frac{rk}{2}} s^{\beta_2 + k - r \left( k - 1 \right)} \left( 1 - s^2 \right)^{\frac{\beta_1 + n - k}{2} - 1} \mathrm{d}s \right]^{\frac{1}{r}}.
				\end{align*}
				Now, we substitute $s^2 = \left( 1 + z \tan^2 \lambda \right) \cos^2 \lambda$ and use Equation \eqref{IntegralForm2F1} to get
				\begin{align}
					\label{Beta1SnCondition1}
					\| R_kf_{\lambda} \|_{L^r_{\beta_1, \beta_2} \left( \Xi_k \left( \mathbb{S}^n \right) \right)} &= C \cos^{\frac{\beta_1 + \beta_2 + n - 1}{r} + 1} \lambda \tan^{k + \frac{\beta_1 + n - k}{r}} \lambda \left[ \int\limits_{0}^{1} \frac{z^{\frac{rk}{2}} \left( 1 - z \right)^{\frac{\beta_1 + n - k}{2} - 1}}{\left( 1 + z \tan^2 \lambda \right)^{\frac{\left( r - 1 \right) \left( k - 1 \right) - \beta_2}{2}}} \mathrm{d}z \right]^{\frac{1}{r}} \\
					&= C \cos^{\frac{\beta_2 + k - 1}{r} - k + 1} \lambda \sin^{\frac{\beta_1 + n - k}{r} + k} \lambda \times \nonumber \\
					&\left[ {}_2F_1 \left( \frac{\left( r - 1 \right) \left( k - 1 \right) - \beta_2}{2}, 1 + \frac{rk}{2}; 1 + \frac{\beta_1 + n - k + rk}{2}; - \tan^2 \lambda \right) \right]^{\frac{1}{r}} \nonumber
				\end{align}
				{Equation \eqref{Beta1SnCondition1} gives the necessity of $\beta_1 > k - n$. Indeed, if $\beta_1 \leq k - n$, then the integral involved in Equation \eqref{Beta1SnCondition1} is not convergent and hence $R_kf_{\lambda} \notin L^r_{\beta_1, \beta_2} \left( \Xi_k \left( \mathbb{S}^n \right) \right)$.}
				
				{Now}, using the transformation of Equation \eqref{Transformation2F1} {\eqref{2F1Transformation1}}, we get
				\begin{align}
					\label{RKFNormMixedZero1}
					\| R_kf_{\lambda} \|_{L^r_{\beta_1, \beta_2} \left( \Xi_k \left( \mathbb{S}^n \right) \right)} &= C \sin^{\frac{\beta_1 + n - k}{r} + k} \lambda \left[ {}_2F_1 \left( \frac{\left( r - 1 \right) \left( k - 1 \right) - \beta_2}{2}, \frac{\beta_1 + n - k}{2}; 1 + \frac{\beta_1 + n - k + rk}{2}; \sin^2 \lambda \right) \right]^{\frac{1}{r}}.
				\end{align}
				Plugging Equations \eqref{RKFNormMixedZero1} and \eqref{LPMassSphereOrigin1} in Inequality \eqref{RequiredLPLRSn}, we get {for small $\lambda > 0$},
				\begin{align*}
					&\sin^{\frac{\beta_1 + n - k}{r} + k} \lambda \left[ {}_2F_1 \left( \frac{\left( r - 1 \right) \left( k - 1 \right) - \beta_2}{2}, \frac{\beta_1 + n - k}{2}; 1 + \frac{\beta_1 + n - k + rk}{2}; \sin^2 \lambda \right) \right]^{\frac{1}{r}} \\
					&\leq C \sin^{\frac{\alpha_1 + n}{p}} \lambda \left[ {}_2F_1 \left( \frac{1 - p - \alpha_2}{2}, \frac{\alpha_1 + n}{2}; 1 + \frac{\alpha_1 + n}{2}; \sin^2 \lambda \right) \right]^{\frac{1}{p}}.
				\end{align*}
				As $\lambda \rightarrow 0$, the above inequality gives us $\frac{\beta_1 + n - k}{r} + k \geq \frac{\alpha_1 + n}{p}$.
				
				\textbf{Case {I(B)}:} Now, let us consider $\alpha_1 \leq -n$. Here, we {consider} the function {$f_{\lambda}: \mathbb{S}^n_+ \rightarrow \mathbb{C}$, defined as} $f_{\lambda} \left( x \right) = \chi_{B \left( 0, \lambda \right)} \left( x \right) \cos d \left( 0, x \right) \sin^{- \frac{{\alpha_1}}{p}} d \left( 0, x \right)$. Then, we have from Equation \eqref{LpNormRadialSn3},
				\begin{align*}
					\| f_{\lambda} \|_{L^p_{\alpha_1, \alpha_2} \left( \mathbb{S}^n \right)} &= C \left[ \int\limits_{0}^{\sin^2 \lambda} \left( 1 - {t} \right)^{\frac{p + \alpha_2 - 1}{2}} {t}^{\frac{n}{2} - 1} \mathrm{d}{t} \right]^{\frac{1}{p}}.
				\end{align*}
				By substituting ${t} = x \sin^2 \lambda$, we get
				\begin{align}
					\label{LPMassSphereZero2}
					\| f_{\lambda} \|_{{L^p_{\alpha_1, \alpha_2} \left( \mathbb{S}^n \right)}} &= C \left[ \int\limits_{0}^{1} \left( 1 - x \sin^2 \lambda \right)^{\frac{p + \alpha_2 - 1}{2}} \sin^{n - 2} \lambda \ x^{\frac{n}{2} - 1} \sin^2 \lambda \ \mathrm{d}x \right]^{\frac{1}{p}} = C \sin^{\frac{n}{p}} \lambda \left[ {}_2F_1 \left( \frac{1 - p - \alpha_2}{2}, \frac{n}{2}; 1 + \frac{n}{2}; \sin^2 \lambda \right) \right]^{\frac{1}{p}}.
				\end{align}
				In the last equality, we have used the integral form of the hypergeometric function, given in Equation \eqref{IntegralForm2F1}.
		
				On the other hand, using Equation \eqref{NormRkSn}, we have,
				\begin{align*}
					\| R_kf_{{\lambda}} \|_{L^r_{\beta_1, \beta_2} \left( \Xi_k \left( \mathbb{S}^n \right) \right)} &= C \left[ \int\limits_{0}^{1} \left| \int\limits_{0}^{s} \frac{\chi_{\left( \cos \lambda, 1 \right)} \left( t \right) t \left( 1 - t^2 \right)^{- \frac{\alpha_1}{2p}} \left( s^2 - t^2 \right)^{\frac{k}{2} - 1}}{s^{k - 1}} \right|^r s^{\beta_2 + k} \left( 1 - s^2 \right)^{\frac{\beta_1 + n - k}{2} - 1} \mathrm{d}s \right]^{\frac{1}{r}} \\
					&= C \left[ \int\limits_{\cos \lambda}^{1} \left| \int\limits_{\cos \lambda}^{s} t \left( 1 - t^2 \right)^{- \frac{\alpha_1}{2p}} \left( s^2 - t^2 \right)^{\frac{k}{2} - 1} \mathrm{d}t \right|^r s^{\beta_2 + k - r \left( k - 1 \right)} \left( 1 - s^2 \right)^{\frac{\beta_1 + n - k}{2} - 1} \mathrm{d}s \right]^{\frac{1}{r}}.
				\end{align*}
				We now substitute $t^2 = s^2 - x \left( s^2 - \cos^2 \lambda \right)$ in the inner integral. After simplification, we arrive at
				\begin{align*}
					&\| R_kf_{{\lambda}} \|_{L^r_{\beta_1, \beta_2} \left( \Xi_k \left( \mathbb{S}^n \right) \right)} \\
					&= C \left[ \int\limits_{\cos \lambda}^{1} \left| \int\limits_{0}^{1} \left( 1 + \left( \frac{s^2 - \cos^2 \lambda}{1 - s^2} \right) x \right)^{- \frac{\alpha_1}{2p}} x^{\frac{k}{2} - 1} {\mathrm{d}x} \right|^r s^{\beta_2 + k - r \left( k - 1 \right)} \left( 1 - s^2 \right)^{\frac{\beta_1 + n - k}{2} - \frac{\alpha_1 r}{2p} - 1} \left( s^2 - \cos^2 \lambda \right)^{\frac{rk}{2}} \mathrm{d}s \right]^{\frac{1}{r}} \\
					&= C \left[ \int\limits_{\cos \lambda}^{1} \left| {}_2F_1 \left( \frac{\alpha_1}{2p}, \frac{k}{2}, 1 + \frac{k}{2}; - \left( \frac{s^2 - \cos^2 \lambda}{1 - s^2} \right) \right) \right|^r s^{\beta_2 + k - r \left( k - 1 \right)} \left( 1 - s^2 \right)^{\frac{\beta_1 + n - k}{2} - \frac{\alpha_1 r}{2p} - 1} \left( s^2 - \cos^2 \lambda \right)^{\frac{rk}{2}} \mathrm{d}s \right]^{\frac{1}{r}}.
				\end{align*}
				In the last equality, we have used the integral form of the hypergeometric function given in Equation \eqref{IntegralForm2F1}. Now, using the transformation of Equation \eqref{Transformation2F1} {\eqref{2F1Transformation1}}, we get
				\begin{align*}
					&\| R_kf_{{\lambda}} \|_{L^r_{\beta_1, \beta_2} \left( \Xi_k \left( \mathbb{S}^n \right) \right)} \\
					&= C \sin^{- \frac{\alpha_1}{p}} \lambda \left[ \int\limits_{\cos \lambda}^{1} \left| {}_2F_1 \left( \frac{\alpha_1}{2p}, 1; 1 + \frac{k}{2}; \frac{s^2 - \cos^2 \lambda}{\sin^2 \lambda} \right) \right|^r s^{\beta_2 + k - r \left( k - 1 \right)} \left( 1 - s^2 \right)^{\frac{\beta_1 + n - k}{2} - 1} \left( s^2 - \cos^2 \lambda \right)^{\frac{rk}{2}} \mathrm{d}s \right]^{\frac{1}{r}}.
				\end{align*}
				{Next}, we put $s^2 = 1 - x \sin^2 \lambda$ {and simplify} to get,
				\begin{align}
					\label{Beta1SnCondition2}
					&\| R_kf_{{\lambda}} \|_{L^r_{\beta_1, \beta_2} \left( \Xi_k \left( \mathbb{S}^n \right) \right)} \\
					&= C \sin^{{\frac{\beta_1 + n - k}{r} - \frac{\alpha_1}{p} + k}} \lambda \left[ \int\limits_{0}^{1} \left| {}_2F_1 \left( \frac{\alpha_1}{2p}, 1; 1 + \frac{k}{2}; 1 - x \right) \right|^r \left( 1 - x \sin^2 \lambda \right)^{\frac{\beta_2 - \left( r - 1 \right) \left( k - 1 \right)}{2}} x^{\frac{\beta_1 + n - k}{2} - 1} \left( 1 - x \right)^{\frac{rk}{2}} \mathrm{d}x \right]^{\frac{1}{r}}. \nonumber
				\end{align}
				{Again, from Equation \eqref{Beta1SnCondition2}, it is clear that $\beta_1 > k - n$ is necessary, for otherwise, $R_kf_{\lambda} \notin L^r_{\beta_1, \beta_2} \left( \Xi_k \left( \mathbb{S}^n \right) \right)$.} {Further, we}  notice that {since $\alpha_1 \leq -n < 0$, by Remark \ref{HypergemetricFunctionBoundedBelow}}, the hypergoemetric function in Equation \eqref{Beta1SnCondition2} is bounded {below by a positive constant}. Hence, we have using Equation \eqref{IntegralForm2F1},
				\begin{equation}
					\label{RKFNormMixedZero2}
					\| R_kf \|_{L^r_{\beta_1, \beta_2} \left( \Xi_k \left( \mathbb{S}^n \right) \right)} \geq C \sin^{\frac{\beta_1 + n - k}{r} + k - \frac{\alpha_1}{p}} \lambda \left[ {}_2F_1 \left( \frac{\left( r - 1 \right) \left( k - 1 \right) - \beta_2}{2}, \frac{\beta_1 + n - k}{2}; 1 + \frac{rk + \beta_1 + n - k}{2}; \sin^2 \lambda \right) \right]^{\frac{1}{r}}.
				\end{equation}
				Using Equations \eqref{LPMassSphereZero2} and \eqref{RKFNormMixedZero2} in Inequality \eqref{RequiredLPLRSn}, we get {for ``small" $\lambda > 0$},
				\begin{align*}
					&\sin^{\frac{\beta_1 + n - k}{r} + k - \frac{\alpha_1}{p}} \lambda \left[ {}_2F_1 \left( \frac{\left( r - 1 \right) \left( k - 1 \right) - \beta_2}{2}, \frac{\beta_1 + n - k}{2}; 1 + \frac{rk + \beta_1 + n - k}{2}; \sin^2 \lambda \right) \right]^{\frac{1}{r}} \\
					&\leq C \sin^{\frac{n}{p}} \lambda \left[ {}_2F_1 \left( \frac{1 - p - \alpha_2}{2}, \frac{n}{2}; 1 + \frac{n}{2}; \sin^2 \lambda \right) \right]^{\frac{1}{p}}.
				\end{align*}
				As $\lambda \rightarrow 0$, the above inequality forces ${\frac{\beta_1 + n - k}{r} + k \geq \frac{\alpha_1 + n}{p}}$.
			
				\textbf{{Case} II:} To get other necessary conditions {mentioned in Equation \eqref{LPLRSNNecessaryConditions}}, we now look at functions concentrated on the equator of $\mathbb{S}^n$. Here, we consider $\lambda > 0$ to be ``near" $\frac{\pi}{2}$. {In what follows, we consider the set $E_{\lambda} := \mathbb{S}^n_+ \setminus B \left( 0, \lambda \right)$.} The analysis is again split into two cases.
				
				\textbf{Case {II(A)}:} First, let us assume $p + \alpha_2 + 1 > 0$. In this case, we consider the function {$f_{\lambda}: \mathbb{S}^n_+ \rightarrow \mathbb{C}$, defined as} $f_{\lambda} \left( x \right) = \chi_{{E_{\lambda}}} \left( x \right) \cos d \left( 0, x \right)$. We have using Equation \eqref{LpNormRadialSn2},
				\begin{align*}
					\| f_{\lambda} \|_{L^p_{\alpha_1, \alpha_2} \left( \mathbb{S}^n \right)} &= C \left[ \int\limits_{0}^{\cos^2 \lambda} {t}^{\frac{\alpha_2 + p - 1}{2}} \left( 1 - {t} \right)^{\frac{\alpha_1 + n}{2} - 1} \mathrm{d}{t} \right]^{\frac{1}{p}}.
				\end{align*}
				By substituting ${t} = x \cos^2 \lambda$, we get
				\begin{align}
					\| f_{\lambda} \|_{L^p_{\alpha_1, \alpha_2} \left( \mathbb{S}^n \right)} &= C \left[ \int\limits_{0}^{1} \cos^{\alpha_2 + p - 1} \lambda \ x^{\frac{\alpha_2 + p - 1}{2}} \left( 1 - x \cos^2 \lambda \right)^{\frac{\alpha_1 + n}{2} - 1} \cos^2 \lambda \ \mathrm{d}x \right]^{\frac{1}{p}} \nonumber \\
					\label{LPMassSphereEquator1}
					&= C \cos^{1 + \frac{\alpha_2 + 1}{p}} \lambda \left[ {}_2F_1 \left( 1 - \frac{\alpha_1 + n}{2}, \frac{\alpha_2 + p + 1}{2}; \frac{\alpha_2 + p + 3}{2}; \cos^2 \lambda \right) \right]^{\frac{1}{p}}.
				\end{align}
				In the last step, we have used the integral form of the hypergeometric function given in Equation \eqref{IntegralForm2F1}.
				
				For the weighted norm of the $k$-plane transform, we have from Equation \eqref{NormRkSn},
				\begin{align*}
					\| R_kf_{\lambda} \|_{L^r_{\beta_1, \beta_2} \left( \Xi_k \left( \mathbb{S}^n \right) \right)} &= C \left[ \int\limits_{0}^{1} \left| \int\limits_{0}^{\min \left\lbrace s, \cos \lambda \right\rbrace} \frac{t \left( s^2 - t^2 \right)^{\frac{k}{2} - 1}}{s^{k - 1}} \mathrm{d}t \right|^r s^{\beta_2 + k} \left( 1 - s^2 \right)^{\frac{\beta_1 + n - k}{2} - 1} \mathrm{d}s \right]^{\frac{1}{r}} \\
					&= C \left[ \int\limits_{0}^{1} \left( s^k - \left( s^2 - \left( \min \left\lbrace s, \cos \lambda \right\rbrace \right)^2 \right)^{\frac{k}{2}} \right)^r s^{\beta_2 + k - r \left( k - 1 \right)} \left( 1 - s^2 \right)^{\frac{\beta_1 + n - k}{2} - 1} \mathrm{d}s \right]^{\frac{1}{r}} \\
					&= C \left[ \int\limits_{0}^{\cos \lambda} s^{rk + \beta_2 + k - r \left( k - 1 \right)} \left( 1 - s^2 \right)^{\frac{\beta_1 + n - k}{2} - 1} \mathrm{d}s \right. \\
					&\left. + \int\limits_{\cos \lambda}^{1} \left( s^k - \left( s^2 - \cos^2 \lambda \right)^{\frac{k}{2}} \right)^r s^{\beta_2 + k - r \left( k - 1 \right)} \left( 1 - s^2 \right)^{\frac{\beta_1 + n - k}{2} - 1} \mathrm{d}s \right]^{\frac{1}{r}}.
				\end{align*}
				Let
				\begin{equation}
					\label{I1Case1}
					I_1 = \int\limits_{0}^{\cos \lambda} s^{rk + \beta_2 + k - r \left( k - 1 \right)} \left( 1 - s^2 \right)^{\frac{\beta_1 + n - k}{2} - 1} \mathrm{d}s,
				\end{equation}
				and
				\begin{equation}
					\label{I2Case1}
					I_2 = \int\limits_{\cos \lambda}^{1} \left( s^k - \left( s^2 - \cos^2 \lambda \right)^{\frac{k}{2}} \right)^r s^{\beta_2 + k - r \left( k - 1 \right)} \left( 1 - s^2 \right)^{\frac{\beta_1 + n - k}{2} - 1} \mathrm{d}s.
				\end{equation}
				We now substitute $s^2 = y \cos^2 \lambda$ in Equation \eqref{I1Case1}, and obtain
				\begin{align*}
					I_1 &= \cos^{\beta_2 + k + r + 1} \lambda \int\limits_{0}^{1} y^{\frac{\beta_2 + k + r - 1}{2}} \left( 1 - y \cos^2 \lambda \right)^{\frac{\beta_1 + n - k}{2} - 1} \mathrm{d}y \\
					&= \cos^{\beta_2 + k + r + 1} \lambda \ {}_2F_1 \left( 1 - \frac{\beta_1 + n - k}{2}, \frac{\beta_2 + k + r + 1}{2}; \frac{\beta_2 + k + r + 3}{2}; \cos^2 \lambda \right),
				\end{align*}
				where, the last equality follows from Equation \eqref{IntegralForm2F1}.
				
				Now, by substituting $s^2 = \left( 1 + y \right) \cos^2 \lambda$ in Equation \eqref{I2Case1}, we get
				\begin{align}
					I_2 &= \cos^{\beta_2 + r + k + 1} \lambda \int\limits_{0}^{\tan^2 \lambda} \left( \left( 1 + y \right)^{\frac{k}{2}} - y^{\frac{k}{2}} \right)^r \left( 1 + y \right)^{\frac{\beta_2 - \left( r - 1 \right) \left( k - 1 \right)}{2}} \left( \sin^2 \lambda - y \cos^2 \lambda \right)^{\frac{\beta_1 + n - k}{2} - 1} \mathrm{d}y \nonumber \\
					\label{I2LastEqualityCase1}
					&= \cos^{\beta_1 + \beta_2 + r + n - 1} \lambda \int\limits_{0}^{\tan^2 \lambda} \left( \left( 1 + y \right)^{\frac{k}{2}} - y^{\frac{k}{2}} \right)^r \left( 1 + y \right)^{\frac{\beta_2 - \left( r - 1 \right) \left( k - 1 \right)}{2}} \left( \tan^2 \lambda - y \right)^{\frac{\beta_1 + n - k}{2} - 1} \mathrm{d}y.
				\end{align}
				To get an estimate on $I_2$, we consider two cases.
				
				\textbf{Case {II(A.1)}:} First, let us assume $k \geq 2$. Then, using the mean value theorem, we get
				$$\left( 1 + y \right)^{\frac{k}{2}} - y^{\frac{k}{2}} = \frac{k}{2} \left( \eta \left( y \right) \right)^{\frac{k}{2} - 1},$$
				for some $\eta \left( y \right) \in \left( y, y + 1 \right)$. Since $k \geq 2$, we have,
				$$\left( 1 + y \right)^{\frac{k}{2}} - y^{\frac{k}{2}} \geq \frac{k}{2} \ y^{\frac{k}{2} - 1}.$$
				Thus, Equation \eqref{I2LastEqualityCase1} becomes
				\begin{align*}
					I_2 &\geq C \cos^{\beta_1 + \beta_2 + r + n - 1} \lambda \int\limits_{0}^{\tan^2 \lambda} y^{r \left( \frac{k}{2} - 1 \right)} \left( 1 + y \right)^{\frac{\beta_2 - \left( r - 1 \right) \left( k - 1 \right)}{2}} \left( \tan^2 \lambda - y \right)^{\frac{\beta_1 + n - k}{2} - 1} \mathrm{d}y.
				\end{align*}
				Substituting $y = u \tan^2 \lambda$, we get from Equations \eqref{IntegralForm2F1} and \eqref{Transformation2F1} {\eqref{2F1Transformation1}},
				\begin{align*}
					I_2 &\geq C \cos^{\beta_1 + \beta_2 + r + n - 1} \lambda \tan^{2r \left( \frac{k}{2} - 1 \right) + \beta_1 + n - k} \lambda \int\limits_{0}^{1} \frac{u^{r \left( \frac{k}{2} - 1 \right)} \left( 1 - u \right)^{\frac{\beta_1 + n - k}{2} - 1}}{\left( 1 + u \tan^2 \lambda \right)^{\frac{\left( r - 1 \right) \left( k - 1 \right) - \beta_2}{2}}} \mathrm{d}u \\
					&= \cos^{\beta_2 + 3r + k - rk - 1} \lambda \sin^{2r \left( \frac{k}{2} - 1 \right) + \beta_1 + n - k} \lambda \ {}_2F_1 \left( \frac{\left( r - 1 \right) \left( k - 1 \right) - \beta_2}{2}, 1 + r \left( \frac{k}{2} - 1 \right); 1 - r + \frac{\beta_1 + n - k + rk}{2}; - \tan^2 \lambda \right) \\
					&= \cos^{2r} \lambda \sin^{2r \left( \frac{k}{2} - 1 \right) + \beta_1 + n - k} \lambda \ {}_2F_1 \left( \frac{\left( r - 1 \right) \left( k - 1 \right) - \beta_2}{2}, \frac{\beta_1 + n - k}{2}; 1 + r \left( \frac{k}{2} - 1 \right) + \frac{\beta_1 + n - k}{2}; \sin^2 \lambda \right).
				\end{align*}
				That is for $k \geq 2$, we have,
				\begin{align}
					\label{RKFNormMixedEquator1}
					&\| R_kf_{\lambda} \|_{{L^r_{\beta_1, \beta_2} \left( \Xi_k \right)}} \\
					&\geq C \left[ \cos^{\beta_2 + k + r + 1} \lambda \ {}_2F_1 \left( 1 - \frac{\beta_1 + n - k}{2}, \frac{\beta_2 + k + r + 1}{2}; \frac{\beta_2 + k + r + 3}{2}; \cos^2 \lambda \right) \right. \nonumber \\
					&\left. + \cos^{2r} \lambda \sin^{rk - 2r + \beta_1 + n - k} \lambda \ {}_2F_1 \left( \frac{\left( r - 1 \right) \left( k - 1 \right) - \beta_2}{2}, \frac{\beta_1 + n - k}{2}; 1 - r + \frac{\beta_1 + n - k + rk}{2}; \sin^2 \lambda \right) \right]^{\frac{1}{r}}. \nonumber
				\end{align}
				Using Equations \eqref{RKFNormMixedEquator1} and \eqref{LPMassSphereEquator1} in Inequality \eqref{RequiredLPLRSn}, we get
				\begin{equation}
					\label{NecessaryInequalitySnCase1}
					\begin{aligned}
						&\left[ \cos^{\beta_2 + k + r + 1} \lambda \ {}_2F_1 \left( 1 - \frac{\beta_1 + n - k}{2}, \frac{\beta_2 + k + r + 1}{2}; \frac{\beta_2 + k + r + 3}{2}; \cos^2 \lambda \right) \right. \\
						&\left. + \cos^{2r} \lambda \sin^{rk - 2r + \beta_1 + n - k} \lambda \ {}_2F_1 \left( \frac{\left( r - 1 \right) \left( k - 1 \right) - \beta_2}{2}, \frac{\beta_1 + n - k}{2}; 1 - r + \frac{\beta_1 + n - k + rk}{2}; \sin^2 \lambda \right) \right]^{\frac{1}{r}} \\
						&\leq C \cos^{1 + \frac{\alpha_2 + 1}{p}} \lambda \left[ {}_2F_1 \left( 1 - \frac{\alpha_1 + n}{2}, \frac{\alpha_2 + p + 1}{2}; \frac{\alpha_2 + p + 3}{2}; \cos^2 \lambda \right) \right]^{\frac{1}{p}}.
					\end{aligned}
				\end{equation}
				First, we notice that as $\lambda \rightarrow \frac{\pi}{2}$, $\cos \lambda \rightarrow 0$. Therefore, the hypergeometric function{s} in the first term on the left hand side and that on the right hand side of Inequality \eqref{NecessaryInequalitySnCase1} are bounded (above and below). The analysis now rests on the behaviour of the hypergeometric function that involves $\sin^2 \lambda$ in its argument. We wish to employ Theorem \ref{Behaviour2F1} to get the necessary conditions.
				
				First, if $\beta_2 > r - k - 1$, then from Equation \eqref{Behaviour2F1Case1}, we have as $\lambda \rightarrow \frac{\pi}{2}$,
				$$\left( \cos^{\beta_2 + k + 1 + r} \lambda + \cos^{2r} \lambda \right)^{\frac{1}{r}} \leq C \cos^{1 + \frac{\alpha_2 + 1}{p}} \lambda.$$
				This is same as
				$$\left( 1 + \cos^{\beta_2 - r + k + 1} \lambda \right)^{\frac{1}{r}} \leq C \cos^{\frac{\alpha_2 + 1}{p} - 1} \lambda.$$
				This, however, gives us $\frac{\alpha_2 + 1}{p} \leq 1$, i.e., $p \geq 1 + \alpha_2$. We notice that this condition is redundant{, since it is required for the existence of the $k$-plane transform (see Theorem \ref{ExistenceSN})}.
				
				In the case when $\beta_2 = r - k - 1$, we have from Equation \eqref{Behaviour2F1Case2}, {as $\lambda \rightarrow \frac{\pi}{2}$}
				$$\left( \cos^{2r} \lambda - \cos^{2r} \lambda \log \left( \cos^2 \lambda \right) \right)^{\frac{1}{r}} \leq C \cos^{1 + \frac{\alpha_2 + 1}{p}} \lambda,$$
				or, what is the same,
				$$\left( 1 - \log \left( \cos^2 \lambda \right) \right)^{\frac{1}{r}} \leq C \cos^{\frac{\alpha_2 + 1}{p} - 1} \lambda.$$
				This gives us the restriction $p > 1 + \alpha_2$. {This is again redundant, except when $\alpha_2 = 0$.}
				
				In the case when $\beta_2 <r - k - 1$, we have from Equation \eqref{Behaviour2F1Case3}, {as $\lambda \rightarrow \frac{\pi}{2}$}
				$$\cos^{\frac{\beta_2 + k + 1}{r} + 1} \lambda \leq C \cos^{1 + \frac{\alpha_2 + 1}{p}} \lambda,$$
				which forces $\frac{\beta_2 + k + 1}{r} \geq \frac{\alpha_2 + 1}{p}$.
				
				\textbf{Case {II(A.2)}:} We now turn to the case when $k = 1$. We again estimate the integral $I_2$ of Equation \eqref{I2Case1}. Here, by using the mean value theorem, we get
				$$\left( 1 + y \right)^{\frac{1}{2}} - y^{\frac{1}{2}} = \frac{1}{2} \left( \eta \left( y \right) \right)^{- \frac{1}{2}} \geq \frac{1}{2} \left( 1 + y \right)^{- \frac{1}{2}}.$$
				Thus, Equation \eqref{I2LastEqualityCase1} becomes
				\begin{align*}
					I_2 &\geq \cos^{\beta_2 + r + n - 1} \lambda \int\limits_{0}^{\tan^2 \lambda} \left( 1 + y \right)^{\frac{\beta_2 - r}{2}} \left( \tan^2 \lambda - y \right)^{\frac{\beta_1 + n - {1}}{2} - 1} \mathrm{d}y.
				\end{align*}
				By substituting $y = z \tan^2 \lambda$, and using Equations \eqref{IntegralForm2F1} and \eqref{Transformation2F1} {\eqref{2F1Transformation1}}, we get
				\begin{align*}
					I_2 &\geq \cos^{\beta_1 + \beta_2 + r + n - 1} \lambda \tan^{\beta_1 + n - 1} \lambda \int\limits_{0}^{1} \frac{\left( 1 - z \right)^{\frac{\beta_1 + n - 1}{2} - 1}}{\left( 1 + z \tan^2 \lambda \right)^{\frac{r - \beta_2}{2}}} \mathrm{d}z \\
					&= \cos^{\beta_2 + r} \lambda \sin^{\beta_1 + n - 1} \lambda \ {}_2F_1 \left( \frac{r - \beta_2}{2}, 1; 1 + \frac{\beta_1 + n - 1}{2}; - \tan^2 \lambda \right) \\
					&= \cos^{2r} \lambda \sin^{\beta_1 + n - 1} \lambda \ {}_2F_1 \left( \frac{r - \beta_2}{2}, \frac{\beta_1 + n - 1}{2}; 1 + \frac{\beta_1 + n - 1}{2}; \sin^2 \lambda \right).
				\end{align*}
				That is, we have,
				\begin{align}
					\label{RKFNormMixedEquator2}
					\| R_1f_{\lambda} \|_{L^r_{\beta_1, \beta_2} \left( \Xi_k \left( \mathbb{S}^n \right) \right)} &\geq C \left[ \cos^{\beta_2 + r + 2} \lambda \ {}_2F_1 \left( 1 - \frac{\beta_1 + n - 1}{2}, \frac{\beta_2 + r + 2}{2}; \frac{\beta_2 + r + 4}{2}; \cos^2 \lambda \right) \right. \\
					&\left. + \cos^{2r} \lambda \sin^{\beta_1 + n - 1} \lambda \ {}_2F_1 \left( \frac{r - \beta_2}{2}, \frac{\beta_1 + n - 1}{2}; 1 + \frac{\beta_1 + n - 1}{2}; \sin^2 \lambda \right) \right]^{\frac{1}{r}}. \nonumber
				\end{align}
				Using Equations \eqref{LPMassSphereEquator1} and \eqref{RKFNormMixedEquator2} in Inequality \eqref{RequiredLPLRSn}, we get
				\begin{equation}
					\label{NecessaryInequalityCase1B}
					\begin{aligned}
						&\left[ \cos^{\beta_2 + r + 2} \lambda {}_2F_1 \left( 1 - \frac{\beta_1 + n - 1}{2}, \frac{\beta_2 + r + 2}{2}; \frac{\beta_2 + r + 4}{2}; \cos^2 \lambda \right) \right. \\
						&\left. + \cos^{2r} \lambda \sin^{\beta_1 + n - 1} \lambda {}_2F_1 \left( \frac{r - \beta_2}{2}, \frac{\beta_1 + n - 1}{2}; 1 + \frac{\beta_1 + n - 1}{2}; \sin^2 \lambda \right) \right]^{\frac{1}{r}} \\
						&\leq C \cos^{1 + \frac{\alpha_2 + 1}{p}} \lambda \left[ {}_2F_1 \left( 1 - \frac{\alpha_1 + n}{2}, \frac{\alpha_2 + p + 1}{2}; \frac{\alpha_2 + p + 3}{2}; \cos^2 \lambda \right) \right]^{\frac{1}{p}}.
					\end{aligned}
				\end{equation}
				As with the Case II(A.1), the analysis rests on the behaviour of the hypergeometric function that involves the sine term. We argue as before.
				
				First, if $\beta_2 > r - 2$, we get from Equation \eqref{Behaviour2F1Case1}, {as $\lambda \rightarrow \frac{\pi}{2}$},
				$$\left( \cos^{\beta_2 + r + 2} \lambda + \cos^{2r} \lambda \right)^{\frac{1}{r}} \leq C \cos^{\frac{\alpha_2 + 1}{p} + 1} \lambda,$$
				which is equivalent to
				$$\left( 1 + \cos^{\beta_2 - r + 2} \lambda \right)^{\frac{1}{r}} \leq C \cos^{\frac{\alpha_2 + 1}{p} - 1} \lambda.$$
				This forces $p \geq 1 + \alpha_2$. This condition is again redundant.
				
				If $\beta_2 = r - 2$, we have from Equation \eqref{Behaviour2F1Case2}, {as $\lambda \rightarrow \frac{\pi}{2}$},
				$$\left( \cos^{2r} \lambda - \cos^{2r} \lambda \log \left( \cos^2 \lambda \right) \right)^{\frac{1}{r}} \leq C \cos^{1 + \frac{\alpha_2 + 1}{p}} \lambda.$$
				Equivalently, we have
				$$\left( 1 - \log \left( \cos^2 \lambda \right) \right)^{\frac{1}{r}} \leq C \cos^{\frac{\alpha_2 + 1}{p} - 1} \lambda.$$
				This gives us $p > 1 + \alpha_2$.
				
				Lastly, if $\beta_2 < r - 2$, we have from Equation \eqref{Behaviour2F1Case3}, {as $\lambda \rightarrow \frac{\pi}{2}$},
				$$\cos^{\frac{\beta_2 + 2}{r} + 1} \lambda \leq C \cos^{\frac{\alpha_2 + 1}{p} + 1} \lambda,$$
				which gives $\frac{\beta_2 + 2}{r} + 1 \geq \frac{\alpha_2 + 1}{p}$.
				
				This completes the analysis of the case when $p > -1 - \alpha_2$.
				
				\textbf{Case {II(B)}:} Let us now consider $p \leq -1 - \alpha_2$. Here, we look at the function {$f_{\lambda}: \mathbb{S}^n_+ \rightarrow \mathbb{C}$, defined as} $f_{\lambda} \left( x \right) = \chi_{{E_{\lambda}}} \left( x \right) \cos^{1 - \frac{\alpha_2}{p}} d \left( 0, x \right)$. We have using Equation \eqref{LpNormRadialSn2},
				\begin{align*}
					\| f_{\lambda} \|_{L^p_{\alpha_1, \alpha_2} \left( \mathbb{S}^n \right)} &= C \left[ \int\limits_{0}^{\cos^2 \lambda} {t}^{\frac{p - 1}{2}} \left( 1 - {t} \right)^{\frac{\alpha_1 + n}{2} - 1} \mathrm{d}{t} \right]^{\frac{1}{p}}.
				\end{align*}
				Substituting ${t} = x \cos^2 \lambda$, we have,
				\begin{align}
				\label{LpMassSphereEquator2}
					\| f_{\lambda} \|_{L^p_{\alpha_1, \alpha_2} \left( \mathbb{S}^n \right)} &= C \cos^{1 + \frac{1}{p}} \lambda \left[ \int\limits_{0}^{1} x^{\frac{p - 1}{\lambda}} \left( 1 - x \cos^2 \lambda \right)^{\frac{\alpha_1 + n}{2} - 1} \mathrm{d}x \right]^{\frac{1}{p}} = C \cos^{1 + \frac{1}{p}} \lambda \left[ {}_2F_1 \left( 1 - \frac{\alpha_1 + n}{2}, \frac{p  + 1}{2}; \frac{p + 3}{2}; \cos^2 \lambda \right) \right]^{\frac{1}{p}}.
				\end{align}
				In the last step, we have used Equation \eqref{IntegralForm2F1}.
				
				For the weighted norm of the $k$-plane transform, we have from Equation \eqref{NormRkSn},
				\begin{align*}
					\| R_kf_{\lambda} \|_{L^r_{\beta_1, \beta_2} \left( \Xi_k \left( \mathbb{S}^n \right) \right)} &= C \left[ \int\limits_{0}^{1} \left| \int\limits_{0}^{\min \left\lbrace s, \cos \lambda \right\rbrace} \frac{t^{1 - \frac{\alpha_2}{p}} \left( s^2 - t^2 \right)^{\frac{k}{2} - 1}}{s^{k - 1}} \mathrm{d}t \right|^r s^{\beta_2 + k} \left( 1 - s^2 \right)^{\frac{\beta_1 + n - k}{2} - 1} \mathrm{d}t \right]^{\frac{1}{r}} \\
					&= C \left[ \int\limits_{0}^{\cos \lambda} \left| \int\limits_{0}^{s} t^{1 - \frac{\alpha_2}{p}} \left( s^2 - t^2 \right)^{\frac{k}{2} - 1} \mathrm{d}t \right|^r s^{\beta_2 + k - r \left( k - 1 \right)} \left( 1 - s^2 \right)^{\frac{\beta_1 + n - k}{2} - 1} \mathrm{d}s \right. \\
					&\left. + \int\limits_{\cos \lambda}^{1} \left| \int\limits_{0}^{\cos \lambda} t^{1 - \frac{\alpha_2}{p}} \left( s^2 - t^2 \right)^{\frac{k}{2} - 1} \mathrm{d}t \right|^r s^{\beta_2 + k - r \left( k - 1 \right)} \left( 1 - s^2 \right)^{\frac{\beta_1 + n - k}{2} - 1} \mathrm{d}s \right]^{\frac{1}{r}} \\
					&= C \left[ I_1 + I_2 \right]^{\frac{1}{r}},
				\end{align*}
				where,
				\begin{equation}
					\label{I1Case2}
					I_1 = \int\limits_{0}^{\cos \lambda} \left| \int\limits_{0}^{s} t^{1 - \frac{\alpha_2}{p}} \left( s^2 - t^2 \right)^{\frac{k}{2} - 1} \mathrm{d}t \right|^r s^{\beta_2 + k - r \left( k - 1 \right)} \left( 1 - s^2 \right)^{\frac{\beta_1 + n - k}{2} - 1} \mathrm{d}s,
				\end{equation}
				and
				\begin{equation}
					\label{I2Case2}
					I_2 = \int\limits_{\cos \lambda}^{1} \left| \int\limits_{0}^{\cos \lambda} t^{1 - \frac{\alpha_2}{p}} \left( s^2 - t^2 \right)^{\frac{k}{2} - 1} \mathrm{d}t \right|^r s^{\beta_2 + k - r \left( k - 1 \right)} \left( 1 - s^2 \right)^{\frac{\beta_1 + n - k}{2} - 1} \mathrm{d}s.
				\end{equation}
				We begin by substituting $t^2 = xs^2$ in {the inner integral of} Equation \eqref{I1Case2}. Then, we get,
				\begin{align*}
					I_1 &= \int\limits_{0}^{\cos \lambda} \left| \int\limits_{0}^{1} x^{- \frac{\alpha_2}{2p}} \left( 1 - x \right)^{\frac{k}{2} - 1} \mathrm{d}x \right|^r s^{{\beta_2 + k + r - \frac{r \alpha_1}{p}}} \left( 1 - s^2 \right)^{\frac{\beta_1 + n - k}{2} - 1} \mathrm{d}s = C \int\limits_{0}^{\cos \lambda} s^{\beta_2 + k + r - \frac{\alpha_2 r}{p}} \left( 1 - s^2 \right)^{\frac{\beta_1 + n - k}{2} - 1} \mathrm{d}s.
				\end{align*}
				The constant $C$ in the last equation is $\left| \Beta \left( 1 - \frac{\alpha_2}{2p}, \frac{k}{2} \right) \right|^r$, where $\Beta \left( \cdot, \cdot \right)$ is the Beta function. In further calculations, we do not write this constant explicitly.
				
				Now, we substitute $s^2 = y \cos^2 \lambda$. Then, using Equation \eqref{IntegralForm2F1}, we have,
				\begin{align*}
					I_1 &= C \int\limits_{0}^{1} \cos^{\beta_2 + k + r - 1 - \frac{\alpha_2 r}{p}} \lambda \ y^{\frac{\beta_2 + k - 1}{2} + \frac{r}{2} \left( 1 - \frac{\alpha_2}{p} \right)} \left( 1 - y \cos^2 \lambda \right)^{\frac{\beta_1 + n - k}{2} - 1} \cos^2 \lambda \mathrm{d}y \\
					&= C \cos^{\beta_2 + k + 1 + r \left( 1 - \frac{\alpha_2}{p} \right)} \lambda \ {}_2F_1 \left( 1 - \frac{\beta_1 + n - k}{2}, \frac{\beta_2 + k + 1}{2} + \frac{r}{2} \left( 1 - \frac{\alpha_2}{p} \right); \frac{\beta_2 + k + 3}{2} + \frac{r}{2} \left( 1 - \frac{\alpha_2}{p} \right); \cos^2 \lambda \right).
				\end{align*}
				Next, we substitute $t^2 = x \cos^2 \lambda$ in {in the inner integral of} Equation \eqref{I2Case2}, and use Equation \eqref{IntegralForm2F1} to obtain,
				\begin{align}
					I_2 &= \int\limits_{\cos \lambda}^{1} \left| \int\limits_{0}^{1} \cos^{- \frac{\alpha_2}{p}} \lambda \ x^{- \frac{\alpha_2}{2p}} \left( s^2 - x \cos^2 \lambda \right)^{\frac{k}{2} - 1} \cos^2 \lambda \ \mathrm{d}x \right|^r s^{\beta_2 + k - r \left( k - 1 \right)} \left( 1 - s^2 \right)^{\frac{\beta_1 + n - k}{2} - 1} \mathrm{d}s \nonumber \\
					\label{I2IntermediateCaseIIBSnNecessary}
					&= \cos^{r \left( 2 - \frac{\alpha_2}{p} \right) \lambda} \int\limits_{\cos \lambda}^{1} \left| {}_2F_1 \left( 1 - \frac{k}{2}, 1 - \frac{\alpha_2}{2p}; 2 - \frac{\alpha_2}{2p}; \frac{\cos^2 \lambda}{s^2} \right) \right|^r s^{\beta_2 + k - r} \left( 1 - s^2 \right)^{\frac{\beta_1 + n - k}{2} - 1} \mathrm{d}s.
				\end{align}
				{We notice by Equation \eqref{Behaviour2F1Case1} of Theorem \ref{Behaviour2F1}, that}
				$${\lim\limits_{s \rightarrow \cos \lambda} {}_2F_1 \left( 1 - \frac{k}{2}, 1 - \frac{\alpha_2}{2p}; 2 - \frac{\alpha_2}{2p}; \frac{\cos^2 \lambda}{s^2} \right) = \frac{\Gamma \left( 2 - \frac{\alpha_2}{2p} \right) \Gamma \left( \frac{k}{2} \right)}{\Gamma \left( 1 + \frac{k}{2} - \frac{\alpha_2}{2p} \right)} > 0.}$$
				{Also, we know (see \cite{Olver}) that}
				$${\frac{d}{dz} {}_2F_1 \left( 1 - \frac{k}{2}, 1 - \frac{\alpha_2}{2p}; 2 - \frac{\alpha_2}{2p}; z \right) = \frac{\left( 1 - \frac{k}{2} \right) \left( 1 - \frac{\alpha_2}{2p} \right)}{\left( 2 - \frac{\alpha_2}{2p} \right)} {}_2F_1 \left( 2 - \frac{k}{2}, 2 - \frac{\alpha_2}{2p}; 3 - \frac{\alpha_2}{2p}; z \right),}$$
				{for $0 < z < 1$. By the integral formula for hypergeometric function given in Equation \eqref{IntegralForm2F1}, it is easy to see that for $0 < z < 1$, we have ${}_2F_1 \left( 1 - \frac{k}{2}, 1 - \frac{\alpha_2}{2p}; 2 - \frac{\alpha_2}{2p}; z \right) \geq 0$. Hence, ${}_2F_1 \left( 1 - \frac{k}{2}, 1 - \frac{\alpha_2}{2p}; 2 - \frac{\alpha_2}{2p}; z \right)$ is non-increasing in $0 < z < 1$. Consequently, for any $s \in \left( \cos \lambda, 1 \right)$, we have}
				$${{}_2F_1 \left( 1 - \frac{k}{2}, 1 - \frac{\alpha_2}{2p}; 2 - \frac{\alpha_2}{2p}; \frac{\cos^2 \lambda}{s^2} \right) \geq \frac{\Gamma \left( 2 - \frac{\alpha_2}{2p} \right) \Gamma \left( \frac{k}{2} \right)}{\Gamma \left( 1 + \frac{k}{2} - \frac{\alpha_2}{2p} \right)} > 0.}$$
				{That is, the hypergeometric function in Equation \eqref{I2IntermediateCaseIIBSnNecessary} is bounded away from $0$}, so that we have,
				\begin{align*}
					I_2 &\geq C \cos^{r \left( 2 - \frac{\alpha_2}{p} \right)} \lambda \int\limits_{\cos \lambda}^{1} s^{\beta_2 + k - r} \left( 1 - s^2 \right)^{\frac{\beta_1 + n - k}{2} - 1} \mathrm{d}s.
				\end{align*}
				By using $s^2 = 1 - z \sin^2 \lambda$ and again using Equation \eqref{IntegralForm2F1}, we obtain,
				\begin{align*}
					I_2 &\geq C \cos^{r \left( 2 - \frac{\alpha_2}{p} \right)} \lambda \sin^{\beta_1 + n - k} \lambda \int\limits_{0}^{1} \left( 1 - z \sin^2 \lambda \right)^{\frac{\beta_2 + k - r - 1}{2}} z^{\frac{\beta_1 + n - k}{2} - 1} \mathrm{d}z \\
					&= C \cos^{r \left( 2 - \frac{\alpha_2}{p} \right)} \lambda \sin^{\beta_1 + n - k} \lambda \ {}_2F_1 \left( \frac{r + 1 - k - \beta_2}{2}, \frac{\beta_1 + n - k}{2}; 1 + \frac{\beta_1 + n - k}{2}; \sin^2 \lambda \right).
				\end{align*}
				That is, we have,
				\begin{equation}
					\label{RKFNormMixedEequator3}
					\begin{aligned}
						&\| R_kf_{\lambda} \|_{L^r_{\beta_1, \beta_2} \left( \Xi_k \left( \mathbb{S}^n \right) \right)} \geq \\
						&C \left[ \cos^{\beta_2 + k + 1 + r \left( 1 - \frac{\alpha_2}{p} \right)} \lambda \ {}_2F_1 \left( 1 - \frac{\beta_1 + n - k}{2}, \frac{\beta_2 + k + 1}{2} + \frac{r}{2} \left( 1 - \frac{\alpha_2}{p} \right); \frac{\beta_2 + k + 3}{2} + \frac{r}{2} \left( 1 - \frac{\alpha_2}{p} \right); \cos^2 \lambda \right)\right. \\
						&\left. + \cos^{r \left( 2 - \frac{\alpha_2}{p} \right)} \lambda \sin^{\beta_1 + n - k} \lambda \ {}_2F_1 \left( \frac{r + 1 - k - \beta_2}{2}, \frac{\beta_1 + n - k}{2}; 1 + \frac{\beta_1 + n - k}{2}; \sin^2 \lambda \right) \right]^{\frac{1}{r}}.
					\end{aligned}
				\end{equation}
				Using Equation \eqref{LpMassSphereEquator2} and {Inequality} \eqref{RKFNormMixedEequator3} in Inequality \eqref{RequiredLPLRSn}, we get {for $\lambda$ ``near" $\frac{\pi}{2}$},
				\begin{align*}
					&\left[ \cos^{\beta_2 + k + 1 + r \left( 1 - \frac{\alpha_2}{p} \right)} \lambda \ {}_2F_1 \left( 1 - \frac{\beta_1 + n - k}{2}, \frac{\beta_2 + k + 1}{2} + \frac{r}{2} \left( 1 - \frac{\alpha_2}{p} \right); \frac{\beta_2 + k + 3}{2} + \frac{r}{2} \left( 1 - \frac{\alpha_2}{p} \right); \cos^2 \lambda \right)\right. \nonumber \\
					&\left. + \cos^{r \left( 2 - \frac{\alpha_2}{p} \right)} \lambda \sin^{\beta_1 + n - k} \lambda \ {}_2F_1 \left( \frac{r + 1 - k - \beta_2}{2}, \frac{\beta_1 + n - k}{2}; 1 + \frac{\beta_1 + n - k}{2}; \sin^2 \lambda \right) \right]^{\frac{1}{r}} \\
					&\leq C \cos^{1 + \frac{1}{p}} \lambda \left[ {}_2F_1 \left( 1 - \frac{\alpha_1 + n}{2}, \frac{p  + 1}{2}; \frac{p + 3}{2}; \cos^2 \lambda \right) \right]^{\frac{1}{p}}.
				\end{align*}
				Again, depending on the value of $\beta_2$, we have the following cases from Theorem \ref{Behaviour2F1}. 
				
				First, if $\beta_2 > r - k - 1$, we have from Equation \eqref{Behaviour2F1Case1}, {as $\lambda \rightarrow \frac{\pi}{2}$},
				$$\left( \cos^{\beta_2 + k + 1 + r \left( 1 - \frac{\alpha_2}{p} \right)} \lambda + \cos^{r \left( 2 - \frac{\alpha_2}{p} \right)} \lambda \right)^{\frac{1}{r}} \leq C \cos^{1 + \frac{1}{p}} \lambda.$$
				{The above inequality} is equivalent to
				$$\left( 1 + \cos^{\beta_2 + k + 1 - r} \lambda \right)^{\frac{1}{r}} \leq C \cos^{\frac{\alpha_2 + 1}{p} - 1} \lambda.$$
				This forces $p \geq 1 + \alpha_2$, which is redundant. 
				
				When $\beta_2 = r - k - 1$, we have from Equation \eqref{Behaviour2F1Case2}, {as $\lambda \rightarrow \frac{\pi}{2}$},
				$$\left( \cos^{r \left( 2 - \frac{{\alpha_2}}{p} \right)} \lambda - \cos^{r \left( 2 - \frac{{\alpha_2}}{p} \right)} \lambda \ \log \left( \cos^2 \lambda \right) \right)^{\frac{1}{r}} \leq C \cos^{1 + \frac{1}{p}} \lambda.$$
				Equivalently,
				$$\left( 1 - \log \cos^2 \lambda \right)^{\frac{1}{r}} \leq C \cos^{\frac{\alpha_2 + 1}{p} - 1} \lambda,$$
				which gives $p > 1 + \alpha_2$. This condition {is redundant except when $\alpha_2 = 0$.}
				
				Lastly, when $\beta_2 < r - k - 1$, we have from Equation \eqref{Behaviour2F1Case3}, {as $\lambda \rightarrow \frac{\pi}{2}$},
				$$\cos^{\frac{\beta_2 + k + 1}{r} + 1 - \frac{\alpha_2}{p}} \lambda \leq C \cos^{1 + \frac{1}{p}} \lambda.$$
				This inequality gives $\frac{\beta_2 + k + 1}{r} \geq \frac{1 + \alpha_2}{p}$.
			\end{proof}
			We now proceed to check the sufficiency of the conditions mentioned in Theorem \ref{NecessaryMixed}.
			\begin{theorem}
				\label{SufficientMixed}
				Let $\alpha_1, \alpha_2, {\beta_1}, \beta_2 \in \mathbb{R}$, {and $p \geq \max \left\lbrace 1, 1 + \alpha_2 \right\rbrace$ with $p > 1 + \alpha_2$ when $\alpha_2 > 0$}. {Let} $r \geq p$.
				\begin{enumerate}
					\item[\mylabel{LPLRSnMixedWeightA}{(A)}] For $k \geq 2$, Inequality \eqref{RequiredLPLRSn} holds for all radial functions in $L^p_{\alpha_1, \alpha_2} \left( \mathbb{S}^n \right)$ {if and only if $\beta_1$ and $\beta_2$ satisfy the conditions mentioned in Equation \eqref{LPLRSNNecessaryConditions}}.
					\item[\mylabel{LPLRSnMixedWeightB}{(B)}] For $k = 1$, {we assume in addition that $\frac{1}{p} - \frac{1}{r} \leq \frac{1}{2}$ when $p > 1$ and $r < 2$ when $p = 1$. Then,} Inequality \eqref{RequiredLPLRSn} holds for all radial functions in $L^p_{\alpha_1, \alpha_2} \left( \mathbb{S}^n \right)$ {if and only if $\beta_1$ and $\beta_2$ satisfy the conditions mentioned in Equation \eqref{LPLRSNNecessaryConditions}}.
				\end{enumerate}
			\end{theorem}
			\begin{proof}~
				\begin{enumerate}
					\item[\mylabel{LPLRSnMixedWeightAProof}{(A)}] Let us first prove the result for $k \geq 2$. The necessary conditions of Theorem \ref{NecessaryMixed} require $\frac{\beta_1 + n - k}{r} \geq \frac{\alpha_1 + n}{p} - k$. Since $\beta_1$ is a power of sine, which is a bounded function, it is enough to prove the result for $\frac{\beta_1 + n - k}{r} = \frac{\alpha_1 + n}{p} - k$. From Equation \eqref{NormRkSn}, we have
					\begin{align*}
						\| R_kf \|_{L^r_{\beta_1, \beta_2} \left( \Xi_k \left( \mathbb{S}^n \right) \right)} &= C \left[ \int\limits_{0}^{1} \left| \int\limits_{0}^{s} \frac{\tilde{f} \left( t \right) \left( s^2 - t^2 \right)^{\frac{k}{2} - 1}}{s^{k - 1}} \mathrm{d}t \right|^r  s^{\beta_2 + k} \left( 1 - s^2 \right)^{\frac{\beta_1 + n - k}{2} - 1} \mathrm{d}s \right]^{\frac{1}{r}}.
					\end{align*}
					Substituting $s = \left( 1 + u \right)^{-\frac{1}{2}}$ and $t = \left( 1 + v \right)^{- \frac{1}{2}}$ and simplifying, we get
					\begin{align}
						\| R_kf \|_{L^r_{\beta_1, \beta_2} \left( \Xi_k \left( \mathbb{S}^n \right) \right)} &= C \left[ \int\limits_{0}^{\infty} \left| \int\limits_{u}^{\infty} \tilde{f} \left( \frac{1}{\sqrt{1 + v}} \right) v^{\frac{k}{2}} \left( 1 - \frac{u}{v} \right)^{\frac{k}{2} - 1} \frac{\left( 1 + u \right)^{\frac{1}{2} - \frac{\beta_1 + \beta_2 + n + 1}{2r}}}{\left( 1 + v \right)^{\frac{k + 1}{2}}} u^{\frac{\beta_1 + n - k}{2r}} \ \frac{\mathrm{d}v}{v} \right|^r \frac{\mathrm{d}u}{u} \right]^{\frac{1}{r}} \nonumber \\
						\label{RkEstimate}
						&= C \left[ \int\limits_{0}^{\infty} \left| \int\limits_{u}^{\infty} \tilde{f} \left( \frac{1}{\sqrt{1 + v}} \right) \frac{v^{\frac{\alpha_1 + n}{2p}}}{\left( 1 + v \right)^{\frac{\alpha_1 + \alpha_2 + n + 1}{2p}}} \frac{\left( 1 + u \right)^{\frac{1}{2} - \frac{\beta_1 + n - k}{2r} - \frac{\beta_2 + k + 1}{2r}}}{\left( 1 + v \right)^{\frac{k + 1}{2} - \frac{\alpha_1 + \alpha_2 + n + 1}{2p}}} \times \right. \right. \\
						&\left. \left.  \left( \frac{u}{v} \right)^{\frac{\beta_1 + n - k}{2r}} \left( 1 - \frac{u}{v} \right)^{\frac{k}{2} - 1} \frac{\mathrm{d}v}{v} \right|^r \frac{\mathrm{d}u}{u} \right]^{\frac{1}{r}}. \nonumber
					\end{align}
					Using $\frac{\beta_1 + n - k}{r} = \frac{\alpha_1 + n}{p} - k$ and $\frac{\beta_2 + k + 1}{2r} \geq \frac{1 + \alpha_2}{2p}$, and $1 + u > 1$, we have
					$$\left( 1 + u \right)^{\frac{1}{2} - \frac{\beta_1 + n - k}{2r} - \frac{\beta_2 + k + 1}{2r}} \leq \left( 1 + u \right)^{\frac{1}{2} + \frac{k}{2} - \frac{\alpha_1 + n}{2p} - \frac{\alpha_2 + 1}{2p}} = \left( 1 + u \right)^{\frac{k + 1}{2} - \frac{\alpha_1 + \alpha_2 + n + 1}{2p}}.$$
					This gives,
					\begin{equation}
						\label{RkEstimate1}
						\| R_kf \|_{L^r_{\beta_1, \beta_2} \left( \Xi_k \left( \mathbb{S}^n \right) \right)} \leq C \left[ \int\limits_{0}^{\infty} \left| \int\limits_{u}^{\infty} \tilde{f} \left( \frac{1}{\sqrt{1 + v}} \right) \frac{v^{\frac{\alpha_1 + n}{2p}}}{\left( 1 + v \right)^{\frac{\alpha_1 + \alpha_2 + n + 1}{2p}}} \left( \frac{1 + u}{1 + v} \right)^{\frac{k + 1}{2} - \frac{\alpha_1 + \alpha_2 + n + 1}{2p}} \left( \frac{u}{v} \right)^{\frac{\beta_1 + n - k}{2r}} \left( 1 - \frac{u}{v} \right)^{\frac{k}{2} - 1} \frac{\mathrm{d}v}{v} \right|^r \frac{\mathrm{d}u}{u} \right]^{\frac{1}{r}}.
					\end{equation}
					We divide the analysis into three cases.
					
					\textbf{Case I:} First, let us assume that $p \geq \frac{\alpha_1 + \alpha_2 + n + 1}{k + 1}$. In this case, we have,
					$$\left( \frac{1 + u}{1 + v} \right)^{\frac{k + 1}{2} - \frac{\alpha_1 + \alpha_2 + n + 1}{2p}} \leq 1,$$
					so that from {Inequality} \eqref{RkEstimate1}, we get
					\begin{align*}
						\| R_kf \|_{L^r_{\beta_1, \beta_2} \left( \Xi_k \left( \mathbb{S}^n \right) \right)} &\leq C \left[ \int\limits_{0}^{\infty} \left| \int\limits_{u}^{\infty} \tilde{f} \left( \frac{1}{\sqrt{1 + v}} \right) \frac{v^{\frac{\alpha_1 + n}{2p}}}{\left( 1 + v \right)^{\frac{\alpha_1 + \alpha_2 + n + 1}{2p}}} \left( \frac{u}{v} \right)^{\frac{\beta_1 + n - k}{2r}} \left( 1 - \frac{u}{v} \right)^{\frac{k}{2} - 1} \frac{\mathrm{d}v}{v} \right|^r \frac{\mathrm{d}u}{u} \right]^{\frac{1}{r}} \\
						&= \| F * G \|_{L^r \left( \left( 0, \infty \right), \frac{\mathrm{d}t}{t} \right)},
					\end{align*}
					where,
					\begin{equation}
						\label{ConvolutionFSn1}
						F \left( t \right) = \tilde{f} \left( \frac{1}{\sqrt{1 + {t}}} \right) \frac{{t}^{\frac{\alpha_1 + n}{2p}}}{\left( 1 + {t} \right)^{\frac{\alpha_1 + \alpha_2 + n + 1}{2p}}},
					\end{equation}
					and
					\begin{equation}
						\label{ConvolutionGSn1}
						G \left( t \right) = t^{\frac{\beta_1 + n - k}{2r}} \left( 1 - t \right)^{\frac{k}{2} - 1} \chi_{\left( 0, 1 \right)} \left( t \right).
					\end{equation}
					Then, by the Young's convolution inequality, we have,
					$$\| R_kf \|_{L^r_{\beta_1, \beta_2} \left( \Xi_k \left( \mathbb{S}^n \right) \right)} \leq C \| F \|_{L^p \left( \left( 0, \infty \right), \frac{\mathrm{d}t}{t} \right)} \| G \|_{L^q \left( \left( 0, \infty \right), \frac{\mathrm{d}t}{t} \right)},$$
					where, $1 + \frac{1}{r} = \frac{1}{p} + \frac{1}{q}$. From Equation \eqref{LpNormRadialSn4} it is clear that $\| F \|_{L^p \left( \left( 0, \infty \right), \frac{\mathrm{d}t}{t} \right)} = C \| f \|_{L^p_{\alpha_1, \alpha_2} \left( \mathbb{S}^n \right)}$. Also, since $\beta_1 > k - n$ and $k \geq 2$, we have $G \in L^q \left( \left( 0, \infty \right), \frac{\mathrm{d}t}{t} \right)$, for any $q \geq 1$. Therefore, we get
					$$\| R_kf \|_{L^r_{\beta_1, \beta_2} \left( \Xi_k \left( \mathbb{S}^n \right) \right)} \leq C \| f \|_{L^p_{\alpha_1, \alpha_2} \left( \mathbb{S}^n \right)}.$$
					
					\textbf{Case II:} Now, we assume $1 + \alpha_2 < p < \frac{\alpha_1 + \alpha_2 + n + 1}{k + 1}$. Then, we have,
					$$\left( \frac{1 + u}{1 + v} \right)^{\frac{k + 1}{2} - \frac{\alpha_1 + \alpha_2 + n + 1}{2p}} \leq \left( \frac{u}{v} \right)^{\frac{k + 1}{2} - \frac{\alpha_1 + \alpha_2 + n + 1}{2p}},$$
					so that, from Equation \eqref{RkEstimate1}, we obtain
					\begin{align*}
						\| R_kf \|_{L^r_{\beta_1, \beta_2} \left( \Xi_k \left( \mathbb{S}^n \right) \right)} &\leq C \left[ \int\limits_{0}^{\infty} \left| \int\limits_{u}^{\infty} \tilde{f} \left( \frac{1}{\sqrt{1 + v}} \right) \frac{v^{\frac{\alpha_1 + n}{2p}}}{\left( 1 + v \right)^{\frac{\alpha_1 + \alpha_2 + n + 1}{2p}}} \left( \frac{u}{v} \right)^{\frac{1}{2} - \frac{\alpha_2 + 1}{2p}} \left( 1 - \frac{u}{v} \right)^{\frac{k}{2} - 1} \frac{\mathrm{d}v}{v} \right|^r \frac{\mathrm{d}u}{u} \right]^{\frac{1}{r}} \\
						&= \| F * H \|_{L^r \left( \left( 0, \infty \right), \frac{\mathrm{d}t}{t} \right)},
					\end{align*}
					where $F$ is given by Equation \eqref{ConvolutionFSn1}, and
					\begin{equation}
						\label{ConvolutionHSn1}
						H \left( t \right) = t^{\frac{1}{2} - \frac{\alpha_2 + 1}{2p}} \left( 1 - t \right)^{\frac{k}{2} - 1} \chi_{\left( 0, 1 \right)} \left( t \right).
					\end{equation}
					Then, by the Young's convolution inequality, we have,
					$$\| R_kf \|_{L^r_{\beta_1, \beta_2} \left( \Xi_k \left( \mathbb{S}^n \right) \right)} \leq C \| F \|_{L^p \left( \left( 0, \infty \right), \frac{\mathrm{d}t}{t} \right)} \| H \|_{L^q \left( \left( 0, \infty \right), \frac{\mathrm{d}t}{t} \right)},$$
					where, $1 + \frac{1}{r} = \frac{1}{p} + \frac{1}{q}$. Since $k \geq 2$ and $p > 1 + \alpha_2$, it is clear that $H \in L^q \left( \left( 0, \infty \right), \frac{\mathrm{d}t}{t} \right)$ for any $q \geq 1$. Using Equation \eqref{LpNormRadialSn4}, we now get
					$$\| R_kf \|_{L^r_{\beta_1, \beta_2} \left( \Xi_k \left( \mathbb{S}^n \right) \right)} \leq C \| f \|_{L^p_{\alpha_1, \alpha_2} \left( \mathbb{S}^n \right)}.$$
					
					\textbf{Case III:} It only remains to check the result for $\alpha_2 {=} 0$ and $p = 1 < \frac{\alpha_1 {+ \alpha_2} + n + 1}{k + 1}$. First, we observe that {as mentioned in Equation \eqref{LPLRSNNecessaryConditions}} of Theorem \ref{NecessaryMixed}, we must have $\beta_2 > r - k - 1$. From Equation \eqref{NormRkSn}, we have,
					\begin{align*}
						\| R_kf \|_{L^r_{\beta_1, \beta_2} \left( \Xi_k \left( \mathbb{S}^n \right) \right)} &= C \left[ \int\limits_{0}^{1} \left| \int\limits_{0}^{s} \frac{\tilde{f} \left( t \right) \left( s^2 - t^2 \right)^{\frac{k}{2} - 1}}{s^{k - 1}} s^{\frac{\beta_2 + k}{r}} \left( 1 - s^2 \right)^{\frac{\beta_1 + n - k}{2r} - \frac{1}{r}} \mathrm{d}t \right|^r {\mathrm{d}s} \right]^{\frac{1}{r}}.
					\end{align*}
					Now, by using Minkowski inequality, we have,
					\begin{align*}
						\| R_kf \|_{L^r_{\beta_1, \beta_2} \left( \Xi_k \left( \mathbb{S}^n \right) \right)} \leq C \int\limits_{0}^{1} \left| \tilde{f} \left( t \right) \right| \left[ \int\limits_{t}^{1} \left( s^2 - t^2 \right)^{r \left( \frac{k}{2} - 1 \right)} s^{\beta_2 + k - r \left( k - 1 \right)} \left( 1 - s^2 \right)^{\frac{\beta_1 + n - k}{2} - 1} \mathrm{d}s \right]^{\frac{1}{r}} \mathrm{d}t.
					\end{align*}
					Upon substituting $s^2 = v \left( 1 - t^2 \right) + t^2$ {in the inner integral} and simplifying, we get
					\begin{align*}
						\| R_kf \|_{L^r_{\beta_1, \beta_2} \left( \Xi_k \left( \mathbb{S}^n \right) \right)} \leq C \int\limits_{0}^{1} \left| \tilde{f} \left( t \right) \right| \left[ \frac{\left( 1 - t^2 \right)^{r \left( \frac{k}{2} + \frac{\beta_1 + n - k}{2r} - 1 \right)}}{t^{\left( r - 1 \right) \left( k - 1 \right) - \beta_2}} \int\limits_{0}^{1} \frac{v^{r \left( \frac{k}{2} - 1 \right)} \left( 1 - v \right)^{\frac{\beta_1 + n - k}{2} - 1}}{\left( 1 + \left( \frac{1 - t^2}{t^2} \right) v \right)^{\frac{\left( r - 1 \right) \left( k - 1 \right)}{2}}} \mathrm{d}v \right]^{\frac{1}{r}} \mathrm{d}t.
					\end{align*}
					From Equation \eqref{IntegralForm2F1}, we get
					\begin{align*}
						\| R_kf \|_{L^r_{\beta_1, \beta_2} \left( \Xi_k \left( \mathbb{S}^n \right) \right)} &\leq C \int\limits_{0}^{1} \left| \tilde{f} \left( t \right) \right| \left( 1 - t^2 \right)^{\frac{k}{2} + \frac{\beta_1 + n - k}{2{r}} - 1} t^{\frac{\beta_2 - \left( r - 1 \right) \left( k - 1 \right)}{r}} {\times} \\
						&\left[ {}_2F_1 \left( \frac{\left( r - 1 \right) \left( k - 1 \right) - \beta_2}{2}, 1 + r \left( \frac{k}{2} - 1 \right); 1 + r \left( \frac{k}{2} - 1 \right) + \frac{\beta_1 + n - k}{2}; \frac{t^2 - 1}{t^2} \right) \right]^{\frac{1}{r}} \mathrm{d}t.
					\end{align*}
					{Now, using Equation \eqref{Transformation2F1} \eqref{2F1Transformation1}}, we get
					\begin{align*}
						&\| R_kf \|_{L^r_{\beta_1, \beta_2} \left( \Xi_k \left( \mathbb{S}^n \right) \right)} \\
						&\leq C \int\limits_{0}^{1} \left| \tilde{f} \left( t \right) \right| \left( 1 - t^2 \right)^{\frac{k}{2} + \frac{\beta_1 + n - k}{2{r}} - 1} \left[ {}_2F_1 \left( \frac{\left( r - 1 \right) \left( k - 1 \right) - \beta_2}{2}, \frac{\beta_1 + n - k}{2}; 1 + r \left( \frac{k}{2} - 1 \right) + \frac{\beta_1 + n - k}{2}; 1 - t^2 \right) \right]^{\frac{1}{r}} \mathrm{d}t.
					\end{align*}
					It is clear from {Equation \eqref{Behaviour2F1Case1} of} Theorem \ref{Behaviour2F1} that the hypergeometric function is bounded as $t$ varies over $\left( 0, 1 \right)$. Moreover, from the necessary conditions of Theorem \ref{NecessaryMixed}, we have $\frac{k}{2} + \frac{\beta_1 + n - k}{2r} \geq \frac{\alpha_1 + n}{2}$. Combining this fact with $1 - t^2 \leq 1$, we get
					$$\left( 1 - t^2 \right)^{\frac{k}{2} + \frac{\beta_1 + n - k}{2r} - 1} \leq \left( 1 - t^2 \right)^{\frac{\alpha_1 + n}{2} - 1}.$$
					Hence,
					\begin{align*}
						{\| R_kf \|_{L^r_{\beta_1, \beta_2} \left( \Xi_k \left( \mathbb{S}^n \right) \right)} \leq C \int\limits_{0}^{1} \left| \tilde{f} \left( t \right) \right| \left( 1 - t^2 \right)^{\frac{\alpha_1 + n}{2} - 1} \mathrm{d}t = C \| f \|_{L^1_{\alpha_1, 0} \left( \mathbb{S}^n \right)}}.
					\end{align*}
					
					This completes the proof for $k \geq 2$.
					\item[\mylabel{LPLRSnMixedWeightBProof}{(B)}] Now, we prove the result for $k = 1$. First, we observe that for $\frac{1}{p} - \frac{1}{r} < \frac{1}{2}$, the proof is verbatim to that of \ref{LPLRSnMixedWeightAProof}. The condition $\frac{1}{p} - \frac{1}{r} < \frac{1}{2}$ is required for the functions $G$ and $H$ of Equations \eqref{ConvolutionGSn1} and \eqref{ConvolutionHSn1} to be in $L^q \left( \left( 0, \infty \right), \frac{\mathrm{d}t}{t} \right)$, where $1 + \frac{1}{r} = \frac{1}{p} + \frac{1}{q}$. The case $p = 1$ is also covered in Case III of \ref{LPLRSnMixedWeightAProof}.
					
					The proof will be complete once we check the case $\frac{1}{p} - \frac{1}{r} = \frac{1}{2}$. Indeed, we want to go in the setting of Theorem \ref{LPLQBoundednessFIFiniteOmega}. From Equation \eqref{NormRkSn}, we get
					\begin{align*}
						\| R_1f \|_{L^r_{\beta_1, \beta_2} \left( \Xi_k \left( \mathbb{S}^n \right) \right)} &= C \left[ \int\limits_{0}^{1} \left| \int\limits_{0}^{s} \tilde{f} \left( t \right) \left( s^2 - t^2 \right)^{- \frac{1}{2}} \mathrm{d}t \right|^r s^{\beta_2 + 1} \left( 1 - s^2 \right)^{\frac{\beta_1 + n - 1}{2} - 1} \mathrm{d}s \right]^{\frac{1}{r}}.
					\end{align*}
					By substituting $s^2 = u$ and $t^2 = v$, we arrive at
					\begin{align}
						\label{RKFISNEquation}
						\| R_kf \|_{L^r_{\beta_1, \beta_2} \left( \Xi_k \left( \mathbb{S}^n \right) \right)} &= C \left[ \int\limits_{0}^{1} \left| \int\limits_{0}^{u} \frac{\tilde{f} \left( \sqrt{v} \right) v^{- \frac{1}{2}}}{\left( u - v \right)^{\frac{1}{2}}} u^{\frac{\beta_2}{2r}} \left( 1 - u \right)^{\frac{\beta_1 + n - 1}{2r} - \frac{1}{r}} \mathrm{d}v \right|^r \mathrm{d}u \right]^{\frac{1}{r}} = C \| \rho_+ I_{0+}^{\frac{1}{2}} \varphi \|_{L^r \left( \Omega \right)},
					\end{align}
					where, in the notation of Theorem \ref{LPLQBoundednessFIFiniteOmega}, we have $\Omega = \left[ 0, 1 \right]$, and
					$$\varphi \left( x \right) = \tilde{f} \left( \sqrt{x} \right) x^{- \frac{1}{2}}, \text{ and } \rho_+ \left( x \right) = \rho \left( x \right) = x^{\frac{\beta_2}{2r}} \left( 1 - x \right)^{\frac{\beta_1 + n - 1}{2r} - \frac{1}{r}}.$$
					By applying Theorem \ref{LPLQBoundednessFIFiniteOmega}, we get
					\begin{align*}
						\| R_kf \|_{L^r_{\beta_1, \beta_2} \left( \Xi_k \left( \mathbb{S}^n \right) \right)} &\leq C \| \rho \varphi \|_{p} = C \left[ \int\limits_{0}^{1} \left| \tilde{f} \left( \sqrt{x} \right) \right|^p x^{- \frac{p}{2}} x^{\frac{p \beta_2}{2r}} \left( 1 - x \right)^{p \left( \frac{\beta_1 + n - k}{2r} - \frac{1}{r} \right)} \mathrm{d}x \right]^{\frac{1}{p}}.
					\end{align*}
					Let us substitute $x = \cos^2 t$. Then, we get
					\begin{align*}
						\| R_kf \|_{r, \beta_1, \beta_2} &\leq C \left[ \int\limits_{0}^{\frac{\pi}{2}} \left| \tilde{f} \left( \cos t \right) \right|^p \cos^{p \left( \frac{\beta_2}{r} - 1 \right)} t \sin^{p \left( \frac{\beta_1 + n - 1}{r} - \frac{2}{r} \right)} t \cos t \sin t \ \mathrm{d}t \right]^{\frac{1}{p}} \\
						&= C \left[ \int\limits_{0}^{\frac{\pi}{2}} \left| \tilde{f} \left( \cos t \right) \right|^p \cos^{\alpha_2} t \sin^{\alpha_1 + n - 1} t \cos^{p \left( \frac{\beta_2}{r} - 1 \right) - \alpha_2 + 1} t \sin^{p \left( \frac{\beta_1 + n - 1}{r} - \frac{2}{r} \right) + 2 - \alpha_1 - n} t \ \mathrm{d}t \right]^{\frac{1}{p}}.
					\end{align*}
					First, we notice from {Equation \eqref{LPLRSNNecessaryConditions}} of Theorem \ref{NecessaryMixed} that
					$$p \left( \frac{\beta_1 + n - 1}{r} - \frac{2}{r} \right) + 2 - \alpha_1 - n = p \left( \frac{\beta_1 + n - 1}{r} + 1 - \frac{\alpha_1 + n}{p} \right) \geq 0.$$
					{Similarly}, we deduce that
					$$p \left( \frac{\beta_2}{r} - 1 \right) + 1 - \alpha_2 = p \left( \frac{\beta_2}{r} - 1 + \frac{1 - \alpha_2}{p} \right) = p \left( \frac{\beta_2}{r} - \frac{\alpha_2 + 1}{p} + \frac{2}{p} - 1 \right) = p \left( \frac{\beta_2 + 2}{r} - \frac{\alpha_2 + 1}{p} \right) \geq 0.$$
					Hence,
					$$\| R_kf \|_{L^r_{\beta_1, \beta_2} \left( \Xi_k \left( \mathbb{S}^n \right) \right)} \leq C \left[ \int\limits_{0}^{\frac{\pi}{2}} \left| \tilde{f} \left( \cos t \right) \right|^p \cos^{\alpha_2} t \sin^{\alpha_1 + n - 1} t \ \mathrm{d}t \right]^{\frac{1}{p}} = C \| f \|_{L^p_{\alpha_1, \alpha_2} \left( \mathbb{S}^n \right)}.$$
				\end{enumerate}
				This completes the proof.
			\end{proof}
			\begin{remark}
				\normalfont
				Theorems \ref{NecessaryMixed} and \ref{SufficientMixed} give the answer to Question \ref{LPLRSnQuestion}. {We observe that for} $k \geq 2$, the question is completely answered. {On the other hand, for $k = 1$, we have additionally assumed that} $\frac{1}{p} - \frac{1}{r} \leq \frac{1}{2}$. For the case $\frac{\alpha_1 + n}{p} - \frac{\beta_1 + n - k}{r} = k$ and $\frac{\beta_2 + k + 1}{r} = \frac{\alpha_2 + 1}{p}$, using arguments similar to those mentioned in \cite{KumarRayWE}, we get the necessity of $\frac{1}{p} - \frac{1}{r} \leq \frac{1}{2}$ when $p > 1$. However, when $\frac{\alpha_1 + n}{p} - \frac{\beta_1 + n - k}{r} < k$ or $\frac{\beta_2 + k + 1}{r} > \frac{\alpha_2 + 1}{p}$, we do have the necessity of this condition.
				
				{On the other hand, it is readily seen from Equation \eqref{RKFISNEquation} and Example \ref{L1L2NotPossibleIHalfPlus} that we cannot expect a (weighted) $L^1$-$L^2$ inequality for the X-Ray transform on $\mathbb{S}^n$.}
			\end{remark}
			Next, we answer Question \ref{LPLInfinitySnQuestion} about weighted $L^p$-$L^{\infty}$ estimates for the $k$-plane transform of radial functions. As before, we first wish to get the necessary conditions required for Inequality \eqref{RequiredLPLInfinitySn} to hold. 
			\begin{theorem}
				\label{NecessaryLPLInfinity}
				Let $\alpha_1, \alpha_2 \in \mathbb{R}$ and $p \geq {\max \left\lbrace 1, 1 + \alpha_2 \right\rbrace}$ {with $p > 1 + \alpha_2$ when $\alpha_2 = 0$}. Then, Inequality \eqref{RequiredLPLInfinitySn} holds only if
				\begin{equation}
					\label{LPLInfinitySNNecessaryConditions}
					\gamma_1 \geq \max \left\lbrace 0, \frac{\alpha_1 + n}{p} - k \right\rbrace \text{ and } \gamma_2 \geq \frac{1 + \alpha_2}{p}.
				\end{equation}
			\end{theorem}
			\begin{proof}
				We get the necessary conditions {by considering} two {cases}, similar to the proof of Theorem \ref{NecessaryMixed}.
				
				\textbf{{Case} I:} First, we consider functions centered at the origin of the sphere. {Here,} we consider two {situations, depending on the value of $\alpha_1$}. 
				
				\textbf{Case I{(A)}:} {For} $\alpha_1 > -n${, and}  $\lambda > 0$ (small), {let us consider} ${f_{\lambda}: \mathbb{S}^n_+ \rightarrow \mathbb{C}}${, defined as} $f_{\lambda} \left( x \right) = \chi_{B \left( 0, \lambda \right)} \left( x \right) \cos d \left( 0, x \right)$. {We know from Equation \eqref{LPMassSphereOrigin1} that $f_{\lambda} \in L^p_{\alpha_1, \alpha_2} \left( \mathbb{S}^n \right)$.} {Now}, from Equation \eqref{KPlaneBallSnEquation}, we get
				\begin{equation}
					\label{RkfNearZero}
					\sin^{\gamma_1} d \left( 0, \xi \right) \cos^{\gamma_2} d \left( 0, \xi \right) R_kf_{\lambda} \left( \xi \right) = C \sin^{\gamma_1} d \left( 0, \xi \right) \cos^{\gamma_2 - k + 1} d \left( 0, \xi \right) \left( \sin^2 \lambda - \sin^2 d \left( 0, \xi \right) \right)^{\frac{k}{2}} \chi_{\left( 0, \lambda \right)} \left( d \left( 0, \xi \right) \right).
				\end{equation}
				It is clear from Equation \eqref{RkfNearZero} that if $\gamma_1 < 0$, then $\| \sin^{\gamma_1} d \left( 0, \cdot \right) \cos^{\gamma_2} d \left( 0, \cdot \right) R_kf \|_{L^{\infty} \left( \Xi_k \left( \mathbb{S}^n \right) \right)} = + \infty$, and hence Inequality \eqref{RequiredLPLInfinitySn} cannot hold. In other words, a necessary condition for the required inequality is $\gamma_1 \geq 0$. Also, by choosing $d \left( 0, \xi \right)$ such that $\sin^2 d \left( 0, \xi \right) = \frac{1}{2} \sin^2 \lambda$, we have $\cos^2 d \left( 0, \xi \right) = \frac{1 + \cos^2 \lambda}{2}$, and hence
				\begin{equation}
					\label{EndPointZero1}
					\| \sin^{\gamma_1} d \left( 0, \cdot \right) \cos^{\gamma_2} d \left( 0, \cdot \right) R_kf \|_{L^{\infty} \left( \Xi_k \right)} \geq C \sin^{\gamma_1 + k} \lambda \left( 1 + \cos^2 \lambda \right)^{\frac{\gamma_2 - k + 1}{2}}.
				\end{equation}
				Using Equations \eqref{EndPointZero1} and \eqref{LPMassSphereOrigin1} in Inequality \eqref{RequiredLPLInfinitySn}, we get
				$$\sin^{\gamma_1 + k} \lambda \left( 1 + \cos^2 \lambda \right)^{\frac{\gamma_2 - k + 1}{2}} \leq C \sin^{\frac{\alpha_1 + n}{p}} \lambda \left[ {}_2F_1 \left( \frac{1 - p - \alpha_2}{2}, \frac{\alpha_1 + n}{2}; 1 + \frac{\alpha_1 + n}{2}; \sin^2 \lambda \right) \right]^{\frac{1}{p}}.$$
				As $\lambda \rightarrow 0$, the above inequality forces $\gamma_1 \geq \frac{\alpha_1 + n}{p} - k$.
				
				\textbf{Case {I(B)}:} Now, we consider the case when $\alpha_1 \leq -n$. Here, we consider {$f_{\lambda}: \mathbb{S}^n_+ \rightarrow \mathbb{C}$, defined as} 
				$$f_{\lambda} \left( x \right) = \chi_{B \left( 0, \lambda \right)} \left( x \right) \cos d \left( 0, x \right) \sin^{- \frac{\alpha_1}{p}} d \left( 0, x \right).$$ 
				We know from Equation \eqref{LPMassSphereZero2} that $f_{\lambda} \in L^p_{\alpha_1, \alpha_2} \left( \mathbb{S}^n \right)$. Now, using Equation {\eqref{KPlaneOriginSn2Equation}}, we have,
				\begin{align}
					&\sin^{\gamma_1} d \left( 0, \xi \right) \cos^{\gamma_2} d \left( 0, \xi \right) R_kf_{\lambda} \left( \xi \right) \nonumber \\
					\label{CaseIBLPLInfinitySnNecessaryEquationRkf}
					&= C \sin^{\gamma_1} d \left( 0, \xi \right) \cos^{\gamma_2 - k + 1} d \left( 0, \xi \right) \sin^{k - \frac{\alpha_1}{p}} \lambda \left( 1 - \frac{\sin^2 d \left( 0, \xi \right)}{\sin^2 \lambda} \right)^{\frac{k}{2}} {}_2F_1 \left( \frac{\alpha_1}{2p}, 1; 1 + \frac{k}{2}; 1 - \frac{\sin^2 d \left( 0, \xi \right)}{\sin^2 \lambda} \right) {\chi_{\left( 0, \lambda \right)} \left( d \left( 0, \xi \right) \right)}.
				\end{align}
				{Since $\alpha_1 \leq -n < 0$, from Remark \ref{HypergemetricFunctionBoundedBelow}, we have that the hypergeometric function in Equation \eqref{CaseIBLPLInfinitySnNecessaryEquationRkf} is bounded below by a positive constant.} It is, {therefore}, clear at once that $\gamma_1 \geq 0$ is {necessary for Inequality \eqref{RequiredLPLInfinitySn} to hold,} for otherwise, $\| \sin^{\gamma_1} d \left( 0, \cdot \right) \cos^{\gamma_2} d \left( 0, \cdot \right) R_kf \|_{L^{\infty} \left( \Xi_k \left( \mathbb{S}^n \right) \right)} = + \infty$. Now, by choosing $d \left( 0, \xi \right)$ such that $\sin^2 d \left( 0, \xi \right) = \frac{1}{2} \sin^2 \lambda$, we get
				\begin{equation}
					\label{EndPointZero2}
					\| \sin^{\gamma_1} d \left( 0, \cdot \right) \cos^{\gamma_2} d \left( 0, \cdot \right) R_kf \|_{L^{\infty} \left( \Xi_k \right)} \geq C \sin^{\gamma_1 + k - \frac{\alpha_1}{p}} \lambda \left( 1 + \cos^2 \lambda \right)^{\frac{\gamma_2 - k + 1}{2}}.
				\end{equation}
				Using Inequality \eqref{EndPointZero2} and Equation \eqref{LPMassSphereZero2} in Inequality \eqref{RequiredLPLInfinitySn}, we get
				\begin{align*}
					\sin^{\gamma_1 + k - \frac{\alpha_1}{p}} \lambda \left( 1 + \cos^2 \lambda \right)^{\frac{\gamma_2 - k + 1}{2}} \leq C \sin^{\frac{n}{p}} \lambda \left[ {}_2F_1 \left( \frac{1 - p - \alpha_2}{2}, \frac{n}{2}; 1 + \frac{n}{2}; \sin^2 \lambda \right) \right]^{\frac{1}{p}}.
				\end{align*}
				As $\lambda \rightarrow 0$, the above inequality forces $\gamma_1 \geq \frac{\alpha_1 + n}{p} - k$.
				
				\textbf{{Case} II:} We now consider functions centered at the equator of $\mathbb{S}^n$. {In what follows, we consider $E_{\lambda}: = \mathbb{S}^n_+ \setminus B \left( 0, \lambda \right)$, for $0 < \lambda < \frac{\pi}{2}$.} Again, we consider two cases.
				
				\textbf{Case {II(A)}:} First, we take up the case $\alpha_2 > -1 - p$. {Here, we consider the function $f_{\lambda}: \mathbb{S}^n_+ \rightarrow \mathbb{C}$, defined as $f_{\lambda} \left( x \right) = \chi_{E_{\lambda}} \left( x \right) \cos d \left( 0, x \right)$.} {We know from Equation \eqref{LPMassSphereEquator1} that $f_{\lambda} \in L^p_{\alpha_1, \alpha_2} \left( \mathbb{S}^n \right)$.} {Using} Equation \eqref{KPlaneEquatorSnEquation}, we have,
				\begin{equation}
					\label{RkfNearEquator}
					\begin{aligned}
						&\sin^{\gamma_1} d \left( 0, \xi \right) \cos^{\gamma_2} d \left( 0, \xi \right) R_kf_{\lambda} \left( \xi \right) \\
						&= C \begin{cases}
									\sin^{\gamma_1} d \left( 0, \xi \right) \cos^{\gamma_2 + 1} d \left( 0, \xi \right), & d \left( 0, \xi \right) > \lambda. \\
									\sin^{\gamma_1} d \left( 0, \xi \right) \left[ \cos^{\gamma_2 + 1} d \left( 0, \xi \right) - \cos^{\gamma_2 - k + 1} d \left( 0, \xi \right) \left( \sin^2 \lambda - \sin^2 d \left( 0, \xi \right) \right)^{\frac{k}{2}} \right], & d \left( 0, \xi \right) \leq \lambda.
								\end{cases}
					\end{aligned}
				\end{equation}
				From Equation \eqref{RkfNearEquator}, we deduce that if $\gamma_2 < -1$, then $\| \sin^{\gamma_1} d \left( 0, \cdot \right) \cos^{\gamma_2} d \left( 0, \cdot \right) R_kf \|_{L^{\infty} \left( \Xi_k \left( \mathbb{S}^n \right) \right)} = + \infty$. Hence, a necessary condition condition for Inequality \eqref{RequiredLPLInfinitySn} to hold is $\gamma_2 \geq -1$. However, we see later that this restriction is redundant.
				
				Further, by choosing $d \left( 0, \xi \right) = \lambda$, we have
				\begin{equation}
					\label{RkfNearEquator1Estimate}
					\| \sin^{\gamma_1} d \left( 0, \cdot \right) \cos^{\gamma_2} d \left( 0, \cdot \right) R_kf \|_{L^{\infty} \left( \Xi_k \left( \mathbb{S}^n \right) \right)} \geq C \sin^{\gamma_1} \lambda \cos^{\gamma_2 + 1} \lambda.
				\end{equation}
				Using Equations \eqref{LPMassSphereEquator1} and \eqref{RkfNearEquator1Estimate} in Inequality \eqref{RequiredLPLInfinitySn}, we get
				$$C \sin^{\gamma_1} \lambda \cos^{\gamma_2 + 1} \lambda \leq C \cos^{1 + \frac{\alpha_2 + 1}{p}} \lambda \left[ {}_2F_1 \left( 1 - \frac{\alpha_1 + n}{2}, \frac{\alpha_2 + p + 1}{2}; \frac{\alpha_2 + p + 3}{2}; \cos^2 \lambda \right) \right]^{\frac{1}{p}}.$$
				Hence, as $\lambda \rightarrow \frac{\pi}{2}$, we require $\gamma_2 \geq \frac{1 + \alpha_2}{p}$. We notice that $\frac{1 + \alpha_2}{p} \geq -1$, making the condition $\gamma_2 \geq -1$ redundant.
				
				\textbf{Case II{(B)}:} Next, we see $\alpha_2 \leq  -1 - p$. Here, we {consider} {$f_{\lambda}: \mathbb{S}^n_+ \rightarrow \mathbb{C}$, given by} $f_{\lambda} \left( x \right) = \chi_{{E_{\lambda}}} \left( x \right) \cos^{1 - \frac{\alpha_2}{p}} d \left( 0, x \right)$. We know from Equation \eqref{LpMassSphereEquator2} that $f_{\lambda} \in L^p_{\alpha_1, \alpha_2} \left( \mathbb{S}^n \right)$. Now, using Equation \eqref{KPlaneEquatorSn2Equation}, we get
				\begin{align*}
					&\sin^{\gamma_1} d \left( 0, \xi \right) \cos^{\gamma_2} d \left( 0, \xi \right) R_kf \left( \xi \right) \\
					&= C \begin{cases}
												\sin^{\gamma_1} d \left( 0, \xi \right) \cos^{\gamma_2 + 1 - \frac{\alpha_2}{p}} d \left( 0, \xi \right), & d \left( 0, \xi \right) > \lambda. \\
												\sin^{\gamma_1} d \left( 0, \xi \right) \cos^{\gamma_2 - 1} d \left( 0, \xi \right) \cos^{2 - \frac{\alpha_2}{2p}} \lambda \ {}_2F_1 \left( 1 - \frac{k}{2}, 1 - \frac{\alpha_2}{2p}; 2 - \frac{\alpha_2}{2p}; \frac{\cos^2 \lambda}{\cos^2 d \left( 0, \xi \right)} \right), & d \left( 0, \xi \right) \leq \lambda.
											\end{cases}
				\end{align*}
				By choosing $d \left( 0, \xi \right) = \lambda$, we get
				\begin{equation}
					\label{RkFNearEquator2Estimate}
					\| \sin^{\gamma_1} d \left( 0, \cdot \right) \cos^{\gamma_2} d \left( 0, \cdot \right) R_kf \|_{L^{\infty} \left( \Xi_k \left( \mathbb{S}^n \right) \right)} \geq C \sin^{\gamma_1} \lambda \cos^{\gamma_2 + 1 - \frac{\alpha_2}{p}} \lambda.
				\end{equation}
				{Further, by} using Inequality \eqref{RkFNearEquator2Estimate} and Equation \eqref{LpMassSphereEquator2} in Inequality \eqref{RequiredLPLInfinitySn}, we get
				$$\sin^{\gamma_1} \lambda \cos^{\gamma_2 + 1 - \frac{\alpha_2}{p}} \lambda \leq C \cos^{1 + \frac{1}{p}} \lambda \left[ {}_2F_1 \left( 1 - \frac{\alpha_1 + n}{2}, \frac{p  + 1}{2}; \frac{p + 3}{2}; \cos^2 \lambda \right) \right]^{\frac{1}{p}}.$$
				As $\lambda \rightarrow \frac{\pi}{2}$, we get $\gamma_2 \geq \frac{1 + \alpha_2}{p}$.
				
				This completes the proof!
			\end{proof}
			We now check the sufficiency of the conditions stated in {Equation \eqref{LPLInfinitySNNecessaryConditions} of} Theorem \ref{NecessaryLPLInfinity}.
			\begin{theorem}
				\label{LPLInfinitySufficient}
				Let $\alpha_1, \alpha_2{, \gamma_1, \gamma_2} \in \mathbb{R}$ and $p { \geq \max \left\lbrace 1, 1 + \alpha_2 \right\rbrace}${, with $p > 1 + \alpha_2$ whenever $\alpha_2 > 0$}. Further, assume that {if $p = \frac{\alpha_1 + n}{k}$, then $\gamma_1 > 0$.}
				\begin{enumerate}
					\item[\mylabel{LPLInfinitySnWeightedA}{(A)}] For $k \geq 2$, Inequality \eqref{RequiredLPLInfinitySn} is valid for all radial functions in $L^p_{\alpha_1, \alpha_2} \left( \mathbb{S}^n \right)$ {if and only if $\gamma_1$ and $\gamma_2$ satisfy Equation \eqref{LPLInfinitySNNecessaryConditions}}.
					\item[\mylabel{LPLInfinitySnWeightedB}{(B)}] For $k = 1$, Inequality \eqref{RequiredLPLInfinitySn} is valid for all radial functions in $L^p_{\alpha_1, \alpha_2} \left( \mathbb{S}^n \right)$ {if and only if $\gamma_1$ and $\gamma_2$ satisfy Equation \eqref{LPLInfinitySNNecessaryConditions} and} $p > 2$.
				\end{enumerate}
			\end{theorem}
			\begin{proof}~
				\begin{enumerate}
					\item[\mylabel{LPLInfinitySnWeightedAProof}{(A)}] We first prove the result for $k \geq 2$. This is done in two cases.
					
					\textbf{Case I:} First, we assume $p > 1$. Here, using Equation \eqref{KPlaneTransformRadialSnEquation2}, we have
					\begin{align*}
						\left| \sin^{\gamma_1} d \left( 0, \xi \right) \cos^{\gamma_2} d \left( 0, \xi \right) R_kf \left( \xi \right) \right| &\leq C \sin^{\gamma_1} d \left( 0, \xi \right) \cos^{\gamma_2 - 1} d \left( 0, \xi \right) \int\limits_{d \left( 0, \xi \right)}^{\frac{\pi}{2}} \left| \tilde{f} \left( \cos t \right) \right| \left( 1 - \frac{\tan^2 d \left( 0, \xi \right)}{\tan^2 t} \right)^{\frac{k}{2} - 1} \sin^{k - 1} t \ \mathrm{d}t \\
						&= C \sin^{\gamma_1} d \left( 0, \xi \right) \cos^{\gamma_2 - 1} d \left( 0, \xi \right) \int\limits_{0}^{\frac{\pi}{2}} \left| \tilde{f} \left( \cos t \right) \right| \sin^{\frac{\alpha_1 + n - 1}{p}} t \cos^{\frac{\alpha_2}{p}} t \left( 1 - \frac{\tan^2 d \left( 0, \xi \right)}{\tan^2 t} \right)^{\frac{k}{2} - 1} \times \\
						&\sin^{k - 1 - \frac{\alpha_1 + n - 1}{p}} t \cos^{- \frac{\alpha_2}{p}} t \ \chi_{\left( d \left( 0, \xi \right), \frac{\pi}{2} \right)} \left( t \right) \mathrm{d}t \\
						&\leq C \sin^{\gamma_1} d \left( 0, \xi \right) \cos^{\gamma_2 - 1} d \left( 0, \xi \right) \| f \|_{L^p_{\alpha_1, \alpha_2} \left( \mathbb{S}^n \right)} I^{\frac{1}{p'}},
					\end{align*}
					by using H\"{o}lder's inequality. Here,
					\begin{align*}
						I &= \int\limits_{d \left( 0, \xi \right)}^{\frac{\pi}{2}} \sin^{p' \left( k - 1 - \frac{\alpha_1 + n - 1}{p} \right)} t \cos^{-p' \frac{\alpha_2}{p}} t \left( 1 - \frac{\tan^2 d \left( 0, \xi \right)}{\tan^2 t} \right)^{p' \left( \frac{k}{2} - 1 \right)} \mathrm{d}t. \\
						&= \cos^{-2p' \left( \frac{k}{2} - 1 \right)} d \left( 0, \xi \right) \int\limits_{d \left( 0, \xi \right)}^{\frac{\pi}{2}} \sin^{p' \left( 1 - \frac{\alpha_1 + n - 1}{p} \right)} t \cos^{-p' \frac{\alpha_2}{p}} t \left( \cos^2 d \left( 0, \xi \right) - \cos^2 t \right)^{p' \left( \frac{k}{2} - 1 \right)} \mathrm{d}t.
					\end{align*}
					By substituting $\cos^2 t = x \cos^2 d \left( 0, \xi \right)$, we get
					\begin{align}
						\label{LPLInfinitySneightedCase1I}
						I &= C \cos^{1 - p' \frac{\alpha_2}{p}} d \left( 0, \xi \right) \int\limits_{0}^{1} \left( 1 - x \cos^2 d \left( 0, \xi \right) \right)^{\frac{1}{2} \left( p' \left( 1 - \frac{\alpha_1 + n - 1}{p} \right) - 1 \right)} x^{- \frac{1}{2} \left( 1 + \frac{p' \alpha_2}{p} \right)} \left( 1 - x \right)^{p' \left( \frac{k}{2} - 1 \right)} \mathrm{d}t \\
						&= C \cos^{1 - \frac{p' \alpha_2}{p}} d \left( 0, \xi \right) \times \nonumber \\
						&{}_2F_1 \left( \frac{1}{2} \left( 1 - p' \left( 1 - \frac{\alpha_1 + n - 1}{p} \right) \right), 1 - \frac{1}{2} \left( \frac{p' \alpha_2}{p} + 1 \right); 2 - \frac{1}{2} \left( \frac{p' \alpha_2}{p} + 1 \right) + p' \left( \frac{k}{2} - 1 \right); \cos^2 d \left( 0, \xi \right) \right). \nonumber
					\end{align}
					In the last equality, we have used the integral form of hypergeometric function given in Equation \eqref{IntegralForm2F1}. Hence,
					\begin{equation}
						\label{IntermediateInequalityLPLInfinitySnWeightedCase1}
						\begin{aligned}
							&\left| \sin^{\gamma_1} d \left( 0, \xi \right) \cos^{\gamma_2} d \left( 0, \xi \right) R_kf \left( \xi \right) \right| \leq C \| f \|_{p, \alpha_1, \alpha_2} \sin^{\gamma_1} d \left( 0, \xi \right) \cos^{\gamma_2 - \frac{1 + \alpha_2}{p}} d \left( 0, \xi \right) \times \\
							&\left[ {}_2F_1 \left( \frac{1}{2} \left( 1 - p' \left( 1 - \frac{\alpha_1 + n - 1}{p} \right) \right), 1 - \frac{1}{2} \left( \frac{p' \alpha_2}{p} + 1 \right); 2 - \frac{1}{2} \left( \frac{p' \alpha_2}{p} + 1 \right) + p' \left( \frac{k}{2} - 1 \right); \cos^2 d \left( 0, \xi \right) \right) \right]^{\frac{1}{{p'}}}.
						\end{aligned}
					\end{equation}
					First, we notice that for $d \left( 0, \xi \right) \in \left( \frac{\pi}{4}, \frac{\pi}{2} \right)$, {owing to the fact that $\gamma_1 \geq 0$ and $\gamma_2 \geq \frac{1 + \alpha_2}{p}$,} all the terms in the right hand side of the above inequality are bounded. That is for any $\xi \in \Xi_k \left( \mathbb{S}^n \right)$ such that $d \left( 0, \xi \right) \in \left( \frac{\pi}{4}, \frac{\pi}{2} \right)$, we have
					$$\left| \sin^{\gamma_1} d \left( 0, \xi \right) \cos^{\gamma_2} d \left( 0, \xi \right) R_kf \left( \xi \right) \right| \leq C \| f \|_{{L^p_{\alpha_1, \alpha_2} \left( \mathbb{S}^n \right)}}.$$
					We only need to check the {validity of the desired} inequality for $d \left( 0, \xi \right) \in \left( 0, \frac{\pi}{4} \right)$. Depending on the value of $p$, we have various cases of the behaviour of the hypergeometric function as given in Theorem \ref{Behaviour2F1}.
					
					First, let us consider $p > \frac{\alpha_1 + n}{k}$. In this case, using Equation \eqref{Behaviour2F1Case1}, we see that the hypergeometric term {in Inequality \eqref{IntermediateInequalityLPLInfinitySnWeightedCase1}} is bounded. Moreover, $\gamma_1 \geq 0$ and $\gamma_2 \geq \frac{1 + \alpha_2}{p}$ imply that $\sin^{\gamma_1} d \left( 0, \xi \right) \leq 1$ and $\cos^{\gamma_2 - \frac{1 + \alpha_2}{p}} d \left( 0, \xi \right) \leq 1$. Hence, we have,
					$$\left| \sin^{\gamma_1} d \left( 0, \xi \right) \cos^{\gamma_2} d \left( 0, \xi \right) R_kf \left( \xi \right) \right| \leq C \| f \|_{{L^p_{\alpha_1, \alpha_2} \left( \mathbb{S}^n \right)}}.$$
					
					Let us now assume $p = \frac{\alpha_1 + n}{k}$. From the assumptions of the theorem, we must have $\gamma_1 > 0$. Therefore, from Equation \eqref{Behaviour2F1Case2} {and Inequality \eqref{IntermediateInequalityLPLInfinitySnWeightedCase1}}, we have
					\begin{align}
						\label{EndPointMotivation1}
						\left| \sin^{\gamma_1} d \left( 0, \xi \right) \cos^{\gamma_2} d \left( 0, \xi \right) R_kf \left( \xi \right) \right| &\leq C \| f \|_{{L^p_{\alpha_1, \alpha_2} \left( \mathbb{S}^n \right)}} \sin^{\gamma_1} d \left( 0, \xi \right) \cos^{\gamma_2 - \frac{1 + \alpha_2}{p}} d \left( 0, \xi \right) \left( - \ln \sin^2 d \left( 0, \xi \right) \right)^{\frac{1}{p'}} \leq C \| f \|_{{L^p_{\alpha_1, \alpha_2} \left( \mathbb{S}^n \right)}},
					\end{align}
					
					In the case when $p < \frac{\alpha_1 + n}{k}$, first we observe that $\gamma_1 \geq \frac{\alpha_1 + n}{p} - k > 0$. Also, from Equation \eqref{Behaviour2F1Case3}, we have,
					\begin{align*}
						\left| \sin^{\gamma_1} d \left( 0, \xi \right) \cos^{\gamma_2} d \left( 0, \xi \right) R_kf \left( \xi \right) \right| &\leq C \| f \|_{{L^p_{\alpha_1, \alpha_2} \left( \mathbb{S}^n \right)}} \sin^{\gamma_1 - \frac{\alpha_1 + n}{p} + k} d \left( 0, \xi \right) \cos^{\gamma_2 - \frac{1 + \alpha_2}{p}} d \left( 0, \xi \right) \leq C \| f \|_{{L^p_{\alpha_1, \alpha_2} \left( \mathbb{S}^n \right)}}.
					\end{align*}
					
					\textbf{Case II:} Now, let us assume $p = 1$. Notice that from the existence conditions, this can only happen when $\alpha_2 \leq 0$. Again, using Equation \eqref{KPlaneTransformRadialSnEquation2}, we get
					\begin{align*}
						\left| \sin^{\gamma_1} d \left( 0, \xi \right) \cos^{\gamma_2} d \left( 0, \xi \right) R_kf \left( \xi \right) \right| &\leq C \sin^{\gamma_1} d \left( 0, \xi \right) \cos^{\gamma_2 - 1} d \left( 0, \xi \right) \int\limits_{d \left( 0, \xi \right)}^{\frac{\pi}{2}} \left| \tilde{f} \left( \cos t \right) \right| \sin^{k - 1} t \left( 1 - \frac{\tan^2 d \left( 0, \xi \right)}{\tan^2 t} \right)^{\frac{k}{2} - 1} \\
						&\leq C \sin^{\gamma_1} d \left( 0, \xi \right) \cos^{\gamma_2 - 1} d \left( 0, \xi \right) \int\limits_{d \left( 0, \xi \right)}^{\frac{\pi}{2}} \left| \tilde{f} \left( \cos t \right) \right| \sin^{\alpha_1 + n - 1} t \cos^{\alpha_2} t \sin^{k - \alpha_1 - n} t \cos^{- \alpha_2} t \ \mathrm{d}t \\
						&\leq C \sin^{\gamma_1} d \left( 0, \xi \right) \cos^{\gamma_2 - 1 - \alpha_2} d \left( 0, \xi \right) \int\limits_{d \left( 0, \xi \right)}^{\frac{\pi}{2}} \left| \tilde{f} \left( \cos t \right) \right| \sin^{\alpha_1 + n - 1} t \cos^{\alpha_2} t \sin^{k - \alpha_1 - n} t \ \mathrm{d}t \\
						&\leq C \sin^{\gamma_1} d \left( 0, \xi \right) \int\limits_{d \left( 0, \xi \right)}^{\frac{\pi}{2}} \left| \tilde{f} \left( \cos t \right) \right| \sin^{\alpha_1 + n - 1} t \cos^{\alpha_2} t \sin^{k - \alpha_1 - n} t \ \mathrm{d}t,
					\end{align*}
					since $\gamma_2 \geq 1 + \alpha_2$. We now have the following cases. 
					
					\textbf{Case II(A):} First, let us consider $\alpha_1 > k - n$. Then, $\sin^{k - \alpha_1 - n} t \leq \sin^{k - \alpha_1 - n} d \left( 0, \xi \right)$. Consequently,
					\begin{align*}
						\left| \sin^{\gamma_1} d \left( 0, \xi \right) \cos^{\gamma_2} d \left( 0, \xi \right) R_kf \left( \xi \right) \right| &\leq C \sin^{\gamma_1 + k - \alpha_1 - n} d \left( 0, \xi \right) \int\limits_{d \left( 0, \xi \right)}^{\frac{\pi}{2}} \left| \tilde{f} \left( \cos t \right) \right| \sin^{\alpha_1 + n - 1} t \cos^{\alpha_2} t \ \mathrm{d}t \leq C \| f \|_{L^1_{\alpha_1, \alpha_2} \left( \mathbb{S}^n \right)},
					\end{align*}
					since $\gamma_1 \geq \alpha_1 + n - k$.
					
					\textbf{Case II(B):} Now, let us consider $\alpha_1 \leq k - n$,. Here, we have $\sin^{k - \alpha_1 - n} t \leq 1$, so that
					\begin{align*}
						\left| \sin^{\gamma_1} d \left( 0, \xi \right) \cos^{\gamma_2} d \left( 0, \xi \right) R_kf \left( \xi \right) \right| &\leq C \sin^{\gamma_1} d \left( 0, \xi \right) \int\limits_{d \left( 0, \xi \right)}^{\frac{\pi}{2}} \left| \tilde{f} \left( \cos t \right) \right| \sin^{\alpha_1 + n - 1} t \cos^{\alpha_2} t \ \mathrm{d}t \leq C \| f \|_{L^1_{\alpha_1, \alpha_2} \left( \mathbb{S}^n \right)},
					\end{align*}
					since $\gamma_1 \geq 0$.
					
					This completes the proof for $k \geq 2$.
					\item[\mylabel{LPLInfinitySnWeightedBProof}{(B)}] For $k = 1$ and $p > 2$, the proof is exactly the same as that of Case I of {\ref{LPLInfinitySnWeightedAProof}, where the condition that $p > 2$ is used in guaranteeing the convergence of the integral $I$ mentioned in Equation \eqref{LPLInfinitySneightedCase1I}. The necessity of the condition $p > 2$ when $k = 1$ is given in Example \ref{XRayPAtLeast2} and Remark \ref{L2LinfinityImpossibleSn}.}
				\end{enumerate}
			\end{proof}
			We now see through an example that for $k = 1$, we cannot expect Inequality \eqref{RequiredLPLInfinitySn} when $p < 2$.
			\begin{example}
				\label{XRayPAtLeast2}
				\normalfont
				Consider the function $f: \mathbb{S}^n \rightarrow \mathbb{C}$, defined as $f \left( x \right) = \chi_{A} \left( x \right)$, where $A := \left\lbrace x \in \mathbb{S}^n | \theta_1 < d \left( 0, x \right) < \theta_2 \right\rbrace$ is an annulus {in} $\mathbb{S}^n$. Let $a = \cos \theta_2$ and $b = \cos \theta_1$. Then, we can write $f \left( x \right) = \chi_{\left( a, b \right)} \left( \cos d \left( 0, x \right) \right)$. {We remark that to obtain a $L^p$-$L^{\infty}$ norm inequality as mentioned in Inequality \eqref{RequiredLPLInfinitySn}, we must necessarily have the following Lorentz-norm inequality.
				\begin{equation}
					\label{LPLInfinityRequiredSnLorentz}
					\| \sin^{\gamma_1} d \left( 0, \cdot \right) \cos^{\gamma_2} d \left( 0, \cdot \right) R_1f \|_{{L^{\infty} \left( \Xi_k \left( \mathbb{S}^n \right) \right)}} \leq C \| f \|_{{L^{p, 1}_{\alpha_1, \alpha_2} \left( \mathbb{S}^n \right)}},
				\end{equation}
				where, the constant $C$ is independent of $f$.} 
				{Therefore, let us first} consider the weighted Lorentz norm of the function $f$. We have,
				\begin{equation}
					\label{WeightedLorentzNormAnnulusSphere}
					\| f \|_{L^{p, 1}_{\alpha_1, \alpha_2} \left( \mathbb{S}^n \right)} = C \left[ \int\limits_{a}^{b} \left( 1 - x^2 \right)^{\frac{\alpha_1 + n}{2} - 1} x^{\alpha_2} \mathrm{d}x \right]^{\frac{1}{p}}.
				\end{equation}
				On the other hand, from Equation \eqref{KPlaneTransformRadialSnEquation}, the $X$-Ray transform is given by
				\begin{align*}
					R_1f \left( \xi \right) &= C \int\limits_{d \left( 0, \xi \right)}^{\frac{\pi}{2}} {\chi_{\left( a, b \right)}} \left( \cos t \right) \left( \cos^2 d \left( 0, \xi \right) - \cos^2 t \right)^{- \frac{1}{2}} \sin t \ \mathrm{d}t.
				\end{align*}
				By substituting $\cos t = x \cos d \left( 0, \xi \right)$, we get
				\begin{align}
					R_1f \left( \xi \right) &= C \int\limits_{0}^{1} {\chi_{\left( a, b \right)}} \left( x \cos d \left( 0, \xi \right) \right) \left( 1 - x^2 \right)^{- \frac{1}{2}} \mathrm{d}x. \nonumber \\
					\label{XRayTransformAnnulus}
					&= \begin{cases}
							0, & \cos d \left( 0, \xi \right) < a. \\
							\int\limits_{\frac{a}{\cos d \left( 0, \xi \right)}}^{1} \left( 1 - x^2 \right)^{- \frac{1}{2}} \mathrm{d}x, & a \leq \cos d \left( 0, \xi \right) < b. \\
							\int\limits_{\frac{a}{\cos d \left( 0, \xi \right)}}^{\frac{b}{\cos d \left( 0, \xi \right)}} \left( 1 - x^2 \right)^{- \frac{1}{2}} \mathrm{d}x, & b \leq \cos d \left( 0, \xi \right) \leq 1.
						\end{cases}
				\end{align}
				By choosing $\cos d \left( 0, \xi \right) = b$ in Equation \eqref{XRayTransformAnnulus}, we have,
				\begin{equation}
					\label{EstimateAnnulusXRay}
					\| \sin^{\gamma_1} d \left( 0, \cdot \right) \cos^{\gamma_2} d \left( 0, \cdot \right) R_1f \|_{{L^{\infty} \left( \Xi_k \left( \mathbb{S}^n \right) \right)}} \geq C \left( 1 - b^2 \right)^{\frac{\gamma_1}{2}} b^{\gamma_2} \int\limits_{\frac{a}{b}}^{1} \left( 1 - x^2 \right)^{- \frac{1}{2}} \mathrm{d}x.
				\end{equation}
				Using Equations \eqref{WeightedLorentzNormAnnulusSphere} and \eqref{EstimateAnnulusXRay} in Inequality \eqref{LPLInfinityRequiredSnLorentz}, we get
				\begin{equation}
					\label{LPLInfinityLorentzNecessarySn}
					\left( 1 - b^2 \right)^{\frac{p \gamma_1}{2}} b^{p \gamma_2} \left( \int\limits_{\frac{a}{b}}^{1} \left( 1 - x^2 \right)^{- \frac{1}{2}} \mathrm{d}x \right)^p \leq C \int\limits_{a}^{b} \left( 1 - x^2 \right)^{\frac{\alpha_1 + n}{2} - 1} x^{\alpha_2} \ \mathrm{d}x.
				\end{equation}
				For a fixed $b \in \left( 0, 1 \right)$, we consider the function
				$$H \left( a \right) = \frac{\left( \int\limits_{\frac{a}{b}}^{1} \left( 1 - x^2 \right)^{- \frac{1}{2}} \mathrm{d}x \right)^p}{\int\limits_{a}^{b} \left( 1 - x^2 \right)^{\frac{\alpha_1 + n}{2} - 1} x^{\alpha_2} \ \mathrm{d}x},$$
				defined for $0 < a < b$. For Inequality \eqref{LPLInfinityLorentzNecessarySn} to hold, we want $H$ to be bounded above. Particularly, $\lim\limits_{a \rightarrow b} H \left( a \right)$ should be bounded. We notice that this limit is in an indeterminate form, and therefore, we apply the L'H\^{o}spital rule.
				\begin{align*}
					\lim\limits_{a \rightarrow b} H \left( a \right) &= \lim\limits_{a \rightarrow b} \frac{p \left( \int\limits_{\frac{a}{b}}^{1} \left( 1 - x^2 \right)^{- \frac{1}{2}} \mathrm{d}x \right)^{p - 1} \left( - \left( 1 - \frac{a^2}{b^2} \right)^{- \frac{1}{2}} \right) \left( \frac{1}{b} \right)}{- \left( 1 - a^2 \right)^{\frac{\alpha_1 + n}{2} - 1} a^{\alpha_2}} \\
					&= \lim\limits_{a \rightarrow b} \frac{p \left( \int\limits_{\frac{a}{b}}^{1} \left( 1 - x^2 \right)^{- \frac{1}{2}} \mathrm{d}x \right)^{p - 1}}{\sqrt{b^2 - a^2} \left( 1 - a^2 \right)^{\frac{\alpha_1 + n}{2} - 1} a^{\alpha_2}} \\
					&= \lim\limits_{a \rightarrow b} \frac{p \left( p - 1 \right) \left( \int\limits_{\frac{a}{b}}^{1} \left( 1 - x^2 \right)^{- \frac{1}{2}} \mathrm{d}x \right)^{p - 2}}{\left( 1 - a^2 \right)^{\frac{\alpha_1 + n}{2} - 1} a^{\alpha_2 + 1} + \left( \alpha_1 + n - 2 \right) \left( b^2 - a^2 \right) \left( 1 - a^2 \right)^{\frac{\alpha_1 + n}{2} - 2} a^{\alpha_2 + 1} - \alpha_2 \left( b^2 - a^2 \right) \left( 1 - a^2 \right)^{\frac{\alpha_1 + n}{2} - 1} a^{\alpha_2 - 1}}.
				\end{align*}
				Clearly, this limit is infinite for $p < 2$. That is, we cannot expect the $L^p$-$L^{\infty}$ boundedness of the $X$-Ray transform for $p < 2$.
		\end{example}
		\begin{remark}
			\label{L2LinfinityImpossibleSn}
			\normalfont
			{We now see that, as in the case of the Euclidean and the hyperbolic space, we cannot expect a weighted $L^2$-$L^{\infty}$ norm inequality for the $X$-ray transform on $\mathbb{S}^n$. Again, the idea is motivated from Hardy and Littlewood's work in \cite{HardyLittlewoodFI} and the theory of fractional integrals given in Subsection \ref{FIBoundednessSection}. As we have done before, we first convert Inequality \eqref{RequiredLPLInfinitySn} to an equivalent inequality concerning certain fractional integrals.}
			
			{To do so, let us begin with the left-hand-side of Inequality \eqref{RequiredLPLInfinitySn}. By substituting $\cos^2 t = u$ in Equation \eqref{KPlaneTransformRadialSnEquation}, we have,
			\begin{align*}
				R_1f \left( \xi \right) &= C \int\limits_{0}^{\cos^2 d \left( 0, \xi \right)} \tilde{f} \left( \sqrt{u} \right) \left( \cos^2 d \left( 0, \xi \right) - u \right)^{- \frac{1}{2}} u^{- \frac{1}{2}} \mathrm{d}u = C I^{\frac{1}{2}}_+ \varphi \left( \cos^2 d \left( 0, \xi \right) \right),
			\end{align*}
			where,
			\begin{equation}
				\label{PhiEquationL2LInfinitySnImpossible}
				\varphi \left( t \right) = \frac{\tilde{f} \left( \sqrt{t} \right)}{\sqrt{t}}.
			\end{equation}
			Further,by writing $\cos^2 d \left( 0, \xi \right) = x$, the left-hand-side of Inequality \eqref{RequiredLPLInfinitySn} becomes
			\begin{equation}
				\label{RequiredLHSL2LInfinitySn}
				\| \sin^{\gamma_1} d \left( 0, \cdot \right) \cos^{\gamma_2} d \left( 0, \cdot \right) R_1f \|_{L^{\infty} \left( \Xi_k \left( \mathbb{S}^n \right) \right)} = C \| \left( 1 - x \right)^{\frac{\gamma_1}{2}} x^{\frac{\gamma_2}{2}} I^{\frac{1}{2}}_+ \varphi \|_{L^{\infty} \left( 0, 1 \right)}.
			\end{equation}
			On the other hand, we also have,
			$$\| f \|_{L^2_{\alpha_1, \alpha_2} \left( \mathbb{S}^n \right)} = C \left[ \int\limits_{0}^{\frac{\pi}{2}} \left| \tilde{f} \left( \cos t \right) \right|^2 \sin^{\alpha_1 + n - 1} t \cos^{\alpha_2} t \ \mathrm{d}t \right]^{\frac{1}{2}}.$$
			By substituting $\cos^2 t = u$, we get
			\begin{align}
				\| f \|_{L^2_{\alpha_1, \alpha_2} \left( \mathbb{S}^n \right)} &= C \left[ \int\limits_{0}^{1} \left| \tilde{f} \left( \sqrt{u} \right) \right|^2 \left( 1 - u \right)^{\frac{\alpha_1 + n}{2} - 1} u^{\frac{\alpha_2 - 1}{2}} \mathrm{d}u \right]^{\frac{1}{2}} \nonumber \\
				\label{L2NormfSnFI}
				&= C \left[ \int\limits_{0}^{1} \left| \frac{\tilde{f} \left( \sqrt{u} \right)}{\sqrt{u}} \right|^2 \left( 1 - u \right)^{\frac{\alpha_2 + n}{2} - 1} u^{\frac{\alpha_2 + 1}{2}} \mathrm{d}u \right]^{\frac{1}{2}} = C \| \varphi \|_{L^2_{\frac{\alpha_1 + n}{2} - 1, \frac{\alpha_2 + 1}{2}} \left( 0, 1 \right)},
			\end{align}
			where $\varphi$ is given in Equation \eqref{PhiEquationL2LInfinitySnImpossible}. Now, using Equations \eqref{RequiredLHSL2LInfinitySn} and Equation \eqref{L2NormfSnFI}, it is readily seen that Inequality \eqref{RequiredLPLInfinitySn} for $p = 2$ is equivalent to
			\begin{equation}
				\label{RequiredL2LInfinitySnFIEquivalent}
				\| \left( 1 - x \right)^{\frac{\gamma_1}{2}} x^{\frac{\gamma_2}{2}} I^{\frac{1}{2}}_+ \varphi \|_{L^{\infty} \left( 0, 1 \right)} \leq C \| \varphi \|_{L^2_{\frac{\alpha_1 + n}{2} - 1, \frac{\alpha_2 + 1}{2}} \left( 0, 1 \right)},
			\end{equation}
			for all $\varphi \in L^2_{\frac{\alpha_1 + n}{2} - 1, \frac{\alpha_2 + 1}{2}} \left( 0, 1 \right)$. However, we have seen in Proposition \ref{L2LInfinityIHalfPlusNotPossible} that this is not possible.}
		\end{remark}
		\begin{remark}
			\label{RemarkEndPointsSn}
			\normalfont
			Theorem \ref{LPLInfinitySufficient} leaves out a few ``end-point" cases. We mention them here for reference.
			\begin{enumerate}
				\item The end-point for existence of the $k$-plane transform, $p = 1 + \alpha_2$.
				\item The end-point for Inequality \eqref{RequiredLPLInfinitySn} to hold, namely $p = \frac{\alpha_1 + n}{k}$ and $\gamma_1 = 0$.
			\end{enumerate}
			These cases are dealt with in Section \ref{EndPointSection}, where we shall see a Lorentz norm estimate.
		\end{remark}
	\section{{Certain} end-point estimates for radial functions}
		\label{EndPointSection}
		In this section, we look for the ``end-point" estimates for the $k$-plane transform on the hyperbolic space and the sphere. We have two type of end-points at hand. Either they occur as an end-point of existence (given by Theorems \ref{ExistenceMixed} and \ref{ExistenceSN}), or we have encountered them in proving the weighted $L^p$-$L^{\infty}$ estimates (Theorems \ref{SufficientLPLInfinityHn} and \ref{LPLInfinitySufficient}). {It is to be noticed that at the end-points of existence, we cannot expect a norm inequality on the Lebesgue space for the $k$-plane transform, since the $k$-plane transform is not well-defined here. Nonetheless, in this section we get similar estimates for weighted Lorentz spaces $L^{p, 1}$.}
		
		Before we begin formalizing our ideas, we would like to mention a few end-point results available in literature {that serve as the motivation to this section}.
		\begin{theorem}[Kumar and Ray (\cite{KumarRayWE}, \cite{KumarRay})]
			\label{KumarRayEndPointResult}
			We have the following:
			\begin{enumerate}
				\item Let $\alpha > k - n$, $p = \frac{\alpha + n}{k}$ and $f \in L^{p, 1}_{{\alpha}} \left( \mathbb{R}^n \right)$ {be radial}. Then, for $k \geq 2$ there is a constant (independent of $f$) $C > 0$ such that,
				\begin{equation}
					\label{WeightedEuclideanEndPoint1}
					\| R_kf \|_{L^{\infty} \left( \Xi_k \left( \mathbb{R}^n \right) \right)} \leq C \| f \|_{L^{p, 1}_{\alpha} \left( \mathbb{R}^n \right)}.
				\end{equation}
				Also, for $k = 1$, Inequality \eqref{WeightedEuclideanEndPoint1} holds if $p = \alpha + n \geq 2$.
				\item For $k \geq 2$ and $p = \frac{n - 1}{k - 1}$, there is some constant $C > 0$ such that for all {radial functions} $f \in L^{p, 1} \left( \mathbb{H}^n \right)$, we have
				\begin{equation}
					\label{HyperbolicUnweightedEndPoint}
					\| \cosh d \left( 0, \cdot \right) R_kf \left( \cdot \right) \|_{\infty} \leq C \| f \|_{p, 1}.
				\end{equation}
				\item For $k \geq 1$ and $p = \frac{n}{k}$, there is some constant $C > 0$ such that for all {radial functions} $f \in L^{p, 1} \left( \mathbb{S}^n \right)$, we have
				\begin{equation}
					\label{SphereUnWeightedEndPoint}
					\| \cos d \left( 0, \cdot \right) R_kf \left( \cdot \right) \|_{\infty} \leq C \| f \|_{p, 1}.
				\end{equation}
			\end{enumerate}
		\end{theorem}
		\begin{remark}
			\normalfont
			We remark that a unified proof of Theorem \ref{KumarRayEndPointResult} is recently given in \cite{UnifiedEndPointEstimatesSelf}.
		\end{remark}
		We notice that the end-point estimates of Inequalities \eqref{WeightedEuclideanEndPoint1} and \eqref{HyperbolicUnweightedEndPoint} are the end-points of existence of the $k$-plane transform. That is, we know that the $k$-plane transform on $\mathbb{R}^n$ is only defined on $L^p_{\alpha} \left( \mathbb{R}^n \right)$ for $p < \frac{\alpha + n}{k}$, and that on $\mathbb{H}^n$ is defined on $L^p \left( \mathbb{H}^n \right)$ provided $p < \frac{n - 1}{k - 1}$ {(see for instance Theorem \ref{ExistenceMixed})}. On the other hand, there is no restriction {on $p$} (at least in the unweighted case) for the existence of the $k$-plane transform on $\mathbb{S}^n$. The end-point obtained in Inequality \eqref{SphereUnWeightedEndPoint} is a result of asking for {an} $L^{p, 1}$-$L^{\infty}$ estimate (see \cite{KumarRay}).
		
		We dedicate this section to {improve} the results of \cite{KumarRay} on hyperbolic space and the sphere. {We} divide this section into two {subsections}: first, we deal with the hyperbolic case, and then the sphere. In proving {the main results of this section}, the following lemma of \cite{KumarRayWE} will be essential. We state it here for the reader's reference.
		\begin{lemma}[\cite{KumarRayWE}]
			\label{KumarRayMainLemma}
			Let $n \in \mathbb{N}$ and $\gamma \geq 1$. For real numbers $x_1 \geq x_2 \geq \cdots \geq x_n \geq 0$, we have
			\begin{equation}
				\label{AnnulusEstimateLemma}
				\left( \sum\limits_{i = 1}^{n} \left( -1 \right)^{i - 1} x_i \right)^{\gamma} \leq \sum\limits_{i = 1}^{n} \left( -1 \right)^{i - 1} x_i^{\gamma}.
			\end{equation}
		\end{lemma}
		The following consequence of Lemma \ref{KumarRayMainLemma} is immediate, and can be found in \cite{KumarRay}. We again state it as a lemma since it {is} of use in the sequel.
		\begin{lemma}
			\label{ConsequenceMainLemma}
			For characteristic functions of measurable sets (of finite measure) in $\left( 0, \infty \right)$, we have for $p \geq 1$,
			\begin{equation}
				\label{ConsequenceAnnulusInequality}
				\int\limits_{0}^{\infty} f \left( t \right) \mathrm{d}t \leq p^{\frac{1}{p}} \left( \int\limits_{0}^{\infty} f \left( t \right) t^{p - 1} \mathrm{d}t \right)^{\frac{1}{p}}.
			\end{equation}
			We also have for characteristic functions $f$ of measurable sets of finite measure in $\mathbb{R}$, and $\gamma > 0$,
			\begin{equation}
				\label{ExponentialInequalityLemma}
				\int\limits_{\mathbb{R}} f \left( t \right) e^{\gamma t} \mathrm{d}t \leq \gamma^{-\frac{1}{p'}} p^{\frac{1}{p}} \left( \int\limits_{\mathbb{R}} f \left( t \right) e^{p \gamma t} \mathrm{d}t \right)^{\frac{1}{p}}.
			\end{equation}
			Moreover, for $\gamma > 0$ and characteristic functions $f$ of measurable sets of finite measure in $\left( 0, \infty \right)$, we have
			\begin{equation}
				\label{PolynomialInequalityLemma}
				\int\limits_{0}^{\infty} f \left( t \right) t^{\gamma - 1} \mathrm{d}t \leq p^{\frac{1}{p}} \gamma^{-\frac{1}{p'}} \left( \int\limits_{0}^{\infty} f \left( t \right) t^{p \gamma - 1} \mathrm{d}t \right)^{\frac{1}{p}}.
			\end{equation}
		\end{lemma}
		\begin{proof}
			To get the estimate in Equation \eqref{ConsequenceAnnulusInequality}, it is enough to look at a disjoint union of intervals, since any measurable set (with finite measure) can be approximated by such sets, and any measurable set can be approximated by measurable sets of finite measure. Therefore, we begin with a characteristic function $\chi_{\bigcup\limits_{i = 1}^{n} \left[ a_i, b_i \right]}$, where $b_1 \geq a_1 \geq b_2 \geq a_2 {\geq} \cdots \geq b_n \geq a_n > 0$. Then, we have using Inequality \eqref{AnnulusEstimateLemma},
			\begin{align*}
				\left( \int\limits_{0}^{\infty} \chi_{\bigcup\limits_{i = 1}^{n} \left[ a_i, b_i \right]} \left( t \right) \mathrm{d}t \right)^p &= \left( \sum\limits_{i = 1}^{n} \left( b_i - a_i \right) \right)^p \leq \sum\limits_{i = 1}^{n} \left( b_i^p - a_i^p \right) = p \sum\limits_{i = 1}^{n} \int\limits_{a_i}^{b_i} t^{p - 1} \mathrm{d}t = p \int\limits_{0}^{\infty} \chi_{\bigcup\limits_{i = 1}^{n} \left[ a_i, b_i \right]} \left( t \right) t^{p - 1} \mathrm{d}t.
			\end{align*}
			This proves {Inequality} \eqref{ConsequenceAnnulusInequality}. To get Inequality \eqref{ExponentialInequalityLemma}, we replace $t$ with $e^{\gamma t}$ in {Inequality} \eqref{ConsequenceAnnulusInequality}. Finally, by replacing $e^t$ with $t$ in {Inequality} \eqref{ExponentialInequalityLemma}, we get {Inequality} \eqref{PolynomialInequalityLemma}.
		\end{proof}
		We are now ready to begin our analysis for the {end-point estimates for the $k$-plane transform on the} hyperbolic space and the sphere.
		\subsection{The Hyperbolic Space \texorpdfstring{$\mathbb{H}^n$}{}}
			In this section, we address the following questions about the end-point behaviour of the $k$-plane transform on $\mathbb{H}^n$. {We first recall that in the proof of Theorem \ref{SufficientLPLInfinityHn}, we had to assume that if $p = \frac{\alpha_1 + n}{k}$, then $\gamma_1 > 0$. This assumption was needed in the proof to handle a logarithmic term. We start by asking whether we can have a Lorentz norm estimate when $p = \frac{\alpha_1 + n}{k}$ and $\gamma_1 = 0$. Particularly, we pose the following question.}
			\begin{question}
				\label{EndPointHnQuestion1}
				{Given $\alpha_1, \alpha_2 \in \mathbb{R}$ with $\alpha_1 + \alpha_2 \geq k - n$,} what are the admissible values of $\gamma \in \mathbb{R}$ such that the following inequality holds for all radial functions in $L^{p, 1}_{\alpha_1, \alpha_2} \left( \mathbb{H}^n \right)$, when $p = \frac{\alpha_1 + n}{k}$?
				\begin{equation}
					\label{RequiredEndPointHn2}
					\| \cosh^{\gamma} d \left( 0, \cdot \right) R_kf \|_{L^{\infty} \left( \Xi_k \left( \mathbb{H}^n \right) \right)} \leq C \| f \|_{L^{p, 1}_{\alpha_1, \alpha_2} \left( \mathbb{H}^n \right)}.
				\end{equation}
			\end{question}
			Also, we recall that the point $p = \frac{\alpha_1 + \alpha_2 + n - 1}{k - 1}$ is the end-point of existence of the $k$-plane transform. Therefore, the following question, analogous to the one asked by Kumar and Ray in \cite{KumarRay}, arises.
			\begin{question}
				\label{EndPointHnQuestion2}
				{Given $\alpha_1, \alpha_2 \in \mathbb{R}$ with $\alpha_1 + \alpha_2 > k - n$ and $p = \frac{\alpha_1 + \alpha_2 + n - 1}{k - 1}$,} what are the admissible values of $\gamma_1, \gamma_2 \in \mathbb{R}$ such that the following inequality holds for all radial functions in $L^{p, 1}_{\alpha_1, \alpha_2} \left( \mathbb{H}^n \right)$, where $p = \frac{\alpha_1 + \alpha_2 + n - 1}{k - 1}$?
				\begin{equation}
					\label{RequiredEndPointHn1}
					\| \sinh^{\gamma_1} d \left( 0, \cdot \right) \cosh^{\gamma_2} d \left( 0, \cdot \right) R_kf \|_{L^{\infty} \left( \Xi_k \left( \mathbb{H}^n \right) \right)} \leq C \| f \|_{L^{p, 1}_{\alpha_1, \alpha_2} \left( \mathbb{H}^n \right)}.
				\end{equation}
			\end{question}
			{By answering Question \ref{EndPointHnQuestion2}, we improve the result stated in Inequality \eqref{HyperbolicUnweightedEndPoint} of Theorem \ref{KumarRayEndPointResult}.} To prove the end-point result stated in {Inequality \eqref{HyperbolicUnweightedEndPoint}}, the authors {of \cite{KumarRay}} use Lemma \ref{ConsequenceMainLemma}. We wish to prove a mixed weighted end-point estimate. For this, we require a hyperbolic version of Lemma \ref{ConsequenceMainLemma}, with weight {involving powers} of $\sinh$ and $\cosh$.
		\begin{lemma}
			\label{HyperbolicMainLemmaCosh}
			For characteristic functions $f$ of measurable sets of finite measure in $\left( 0, \infty \right)$, $\gamma > 0$, $\alpha \in \mathbb{R}$, and $p \geq \max \left\lbrace 1, 1 - \alpha \right\rbrace$, we have,
			\begin{equation}
				\label{HyperbolicMainLemmaCoshEquation}
				\int\limits_{0}^{\infty} f \left( t \right) \sinh^{\gamma} t \ \mathrm{d}t \leq C_{p, \gamma, n} \left( \int\limits_{0}^{\infty} f \left( t \right) \cosh^{\alpha} t \sinh^{p \gamma - \alpha} t \mathrm{d}t \right)^{\frac{1}{p}},
			\end{equation}
			where,
			$$C_{p, \gamma, n} = p^{\frac{1}{p}} \left[ \left( \gamma + 1 \right)^{- \frac{1}{p'}} \min \left\lbrace 1, \sinh^{\gamma - \frac{\alpha}{p}} 1 \right\rbrace \max \left\lbrace 1, \cosh^{- \frac{\alpha}{p}} 1 \right\rbrace + \frac{\gamma^{- \frac{1}{p'}}}{2^{\gamma}} \max \left\lbrace 2^{\gamma - \frac{\alpha}{p}}, \left( \frac{e}{\sinh 1} \right)^{\frac{\alpha}{p}} \right\rbrace \max \left\lbrace 2^{\frac{\alpha}{p}}, 1 \right\rbrace \right].$$
		\end{lemma}
		\begin{proof}
			First, we observe that for $t \in \left[ 0, 1 \right]$, we have
			\begin{equation}
				\label{SinhNear0}
				t \leq \sinh t \leq \left( \sinh 1 \right) t,
			\end{equation}
			and for $t \geq 1$, we have
			\begin{equation}
				\label{SinhAwayFrom0}
				\left( \dfrac{\sinh 1}{e} \right) e^t \leq \sinh t \leq {\frac{e^t}{2}}.
			\end{equation}
			We also have for any $t > 0$,
			\begin{equation}
				\label{CoshEstimate}
				\frac{e^t}{2} \leq \cosh t \leq e^t.
			\end{equation}
			Thus, {by using Equations \eqref{SinhNear0} and \eqref{SinhAwayFrom0}, we have for characteristic function of measurable sets in $\left( 0, \infty \right)$,}
			\begin{align*}
				\int\limits_{0}^{\infty} f \left( t \right) \sinh^{\gamma} t \ \mathrm{d}t &= \int\limits_{0}^{1} f \left( t \right) \sinh^{\gamma} t \ \mathrm{d}t + \int\limits_{1}^{\infty} f \left( t \right) \sinh^{\gamma} t \ \mathrm{d}t \\
				&\leq \left( \sinh 1 \right)^{\gamma} \int\limits_{\mathbb{R}} f \left( t \right) \chi_{\left( 0, 1 \right)} \left( t \right) t^{\gamma} \ \mathrm{d}t + \dfrac{1}{2^{\gamma}} \int\limits_{\mathbb{R}} f \left( t \right) \chi_{\left[ 1, \infty \right)} \left( t \right) e^{{\gamma}t} \ \mathrm{d}t.
			\end{align*}
			Now, using Equations \eqref{ExponentialInequalityLemma} and \eqref{PolynomialInequalityLemma}, we get
			\begin{align*}
				\int\limits_{0}^{\infty} f \left( t \right) \sinh^{\gamma} t \ \mathrm{d}t  &\leq p^{\frac{1}{p}} \left( \gamma + 1 \right)^{-\frac{1}{p'}} \left( \sinh 1 \right)^{\gamma} \left( \int\limits_{\mathbb{R}} f \left( t \right) \chi_{\left( 0, 1 \right)} \left( t \right) t^{p \gamma + p - 1} \ \mathrm{d}t \right)^{\frac{1}{p}} + \dfrac{p^{\frac{1}{p}} \gamma^{-\frac{1}{p'}}}{2^{\gamma}} \left( \int\limits_{\mathbb{R}} f \left( t \right) \chi_{\left[ 1, \infty \right)} \left( t \right) e^{p \gamma t} \mathrm{d}t \right)^{\frac{1}{p}} \\
				&= p^{\frac{1}{p}} \left( \gamma + 1 \right)^{-\frac{1}{p'}} \left( \sinh 1 \right)^{\gamma} \left( \int\limits_{0}^{1} f \left( t \right) t^{p \gamma - \alpha} t^{\alpha + p - 1} \mathrm{d}t \right)^{\frac{1}{p}} + \frac{p^{\frac{1}{p}} \gamma^{-\frac{1}{p'}}}{2^{\gamma}} \left( \int\limits_{1}^{\infty} f \left( t \right) e^{\left( p \gamma - \alpha \right) t} e^{\alpha t} \mathrm{d}t \right)^{\frac{1}{p}}.
			\end{align*}
			Since $p \geq 1 - \alpha$, for $t \in \left( 0, 1 \right)$, we have $t^{\alpha + p - 1} \leq 1$. Hence, from Equations \eqref{SinhNear0}, \eqref{SinhAwayFrom0}, and \eqref{CoshEstimate}, {we get}
			\begin{align*}
				\int\limits_{0}^{\infty} f \left( t \right) \sinh^{\gamma} t \ \mathrm{d}t &\leq p^{\frac{1}{p}} \left( \gamma + 1 \right)^{-\frac{1}{p'}} \left( \sinh 1 \right)^{\gamma} \left( \int\limits_{0}^{1} f \left( t \right) t^{p \gamma - \alpha} \mathrm{d}t \right)^{\frac{1}{p}} + \frac{p^{\frac{1}{p}} \gamma^{-\frac{1}{p'}}}{2^{\gamma}} \left( \int\limits_{1}^{\infty} f \left( t \right) e^{\left( p \gamma - \alpha \right) t} e^{\alpha t} \mathrm{d}t \right)^{\frac{1}{p}} \\
				&\leq p^{\frac{1}{p}} \left[ \left( \gamma + 1 \right)^{- \frac{1}{p'}} \min \left\lbrace 1, \sinh^{\gamma - \frac{\alpha}{p}} 1 \right\rbrace \max \left\lbrace 1, \cosh^{- \frac{\alpha}{p}} 1 \right\rbrace + \frac{\gamma^{- \frac{1}{p'}}}{2^{\gamma}} \max \left\lbrace 2^{\gamma - \frac{\alpha}{p}}, \left( \frac{e}{\sinh 1} \right)^{\frac{\alpha}{p}} \right\rbrace \max \left\lbrace 2^{\frac{\alpha}{p}}, 1 \right\rbrace \right] \times \\
				&\left( \int\limits_{0}^{\infty} f \left( t \right) \sinh^{p \gamma - \alpha} t \cosh^{\alpha} t \ \mathrm{d}t \right)^{\frac{1}{p}}.
			\end{align*}
			This completes the proof!
		\end{proof}
		\begin{remark}
			\normalfont
			{We remark here that the constant $C_{p, \gamma, n}$ mentioned in Lemma \ref{HyperbolicMainLemmaCosh} might not be the best possible one, and in our main result, we do not use it explicitly.}
		\end{remark}
		We are now in a position to answer Question \ref{EndPointHnQuestion1}. Particularly, we have the following result.
		\begin{theorem}
			{Let $\alpha_1, \alpha_2 \in \mathbb{R}$ be such that $\alpha_1 + \alpha_2 > k - n$, and $1 \leq p = \frac{\alpha_1 + n}{k} < \frac{\alpha_1 + \alpha_2 + n - 1}{k - 1}$ or $\alpha_1 + \alpha_2 = k - n$ and $\frac{\alpha_1 + n}{k} = p = 1$.
			\begin{enumerate}
				\item[\mylabel{EndPointLPLInfinityHnWeightedA}{(A)}] For $k \geq 2$, Inequality \eqref{RequiredEndPointHn2} holds for all radial functions $f \in L^{p, 1}_{\alpha_1, \alpha_2} \left( \mathbb{H}^n \right)$ if and only if $\gamma \leq \frac{\alpha_1 + \alpha_2 + n - 1}{p}$.
				\item[\mylabel{EndPointLPLInfinityHnWeightedB}{(B)}] For $k = 1$, Inequality \eqref{RequiredEndPointHn2} holds for all radial functions in $L^{p, 1}_{\alpha_1, \alpha_2} \left( \mathbb{H}^n \right)$ if and only if $\gamma \leq \frac{\alpha_1 + \alpha_2 + n - 1}{p}$ and $p \geq 2$.
			\end{enumerate}}			
		\end{theorem}
		\begin{proof}
			We first remark that due to Theorem \ref{LorentzSpaceBoundedness}, it is sufficient to get Inequality \eqref{RequiredEndPointHn2} for characteristic functions of radial sets of $\mathbb{H}^n$.
			\begin{enumerate}
				\item[\mylabel{EndPointLPLInfinityHnWeightedAProof}{(A)}] We begin with proving the result for $k \geq 2$.
				
				Using Equation \eqref{KPlaneTransformRadialHnEquation3}, we have
				\begin{align}
					\left| \cosh^{\gamma} d \left( 0, \xi \right) R_kf \left( \xi \right) \right| &\leq C \cosh^{\gamma - k + 1} d \left( 0, \xi \right) \int\limits_{d \left( 0, \xi \right)}^{\infty} \left| \tilde{f} \left( \cosh t \right) \right| \left( 1 - \frac{\sinh^2 d \left( 0, \xi \right)}{\sinh^2 t} \right)^{\frac{k}{2} - 1} \sinh^{k - 1} t \ \mathrm{d}t \nonumber \\
					\label{EndPointLPLInfinityHnInitialEstimate}
					&\leq C \cosh^{\gamma - k + 1} d \left( 0, \xi \right) \int\limits_{d \left( 0, \xi \right)}^{\infty} \left| \tilde{f} \left( \cosh t \right) \right| \sinh^{k - 1} t \ \mathrm{d}t.
				\end{align}
				We now consider two cases.
			
				\textbf{Case I:} First, let us consider $d \left( 0, \xi \right) \geq 1$. {We notice that since $\frac{\alpha_1 + n}{k} = p \geq 1$ and $\alpha_1 + \alpha_2 \geq k - n$, we have $p \geq 1 - \alpha_2$.} {Hence, by using Lemma \ref{HyperbolicMainLemmaCosh}}, we get
				\begin{align*}
					\left| \cosh^{\gamma} d \left( 0, \xi \right) R_kf \left( \xi \right) \right| &\leq C \cosh^{\gamma - k + 1} d \left( 0, \xi \right) \left( \int\limits_{d \left( 0, \xi \right)}^{\infty} \left| \tilde{f} \left( \cosh t \right) \right| \cosh^{\alpha_2} t \sinh^{\alpha_1 + n - p - \alpha_2} t \ \mathrm{d}t \right)^{\frac{1}{p}} \\
					&= C \cosh^{\gamma - k + 1} d \left( 0, \xi \right) \left( \int\limits_{d \left( 0, \xi \right)}^{\infty} \left| \tilde{f} \left( \cosh t \right) \right| \sinh^{\alpha_1 + n - 1} t \cosh^{\alpha_2} t \sinh^{1 - p - \alpha_2} t \ \mathrm{d}t \right)^{\frac{1}{p}}.
				\end{align*}
				{If $\alpha_1 + \alpha_2 > k - n$, then} $p = \frac{\alpha_1 + n}{k} < \frac{\alpha_1 + \alpha_2 + n - 1}{k - 1}$, {and hence} $1 - p - \alpha_2 < 0$. {On the other hand, if $\alpha_1 + \alpha_2 = k - n$, then $p = 1$ and we must have $\alpha_2 = 0$, so that $1 - p - \alpha_2 = 0$.} Also, from the necessary condition of Theorem \ref{NecessaryLPLInfinityHn}, we must have $\gamma \leq \frac{\alpha_1 + \alpha_2 + n - 1}{p} = k \left( 1 {+} \frac{\alpha_2 - 1}{\alpha_1 + n} \right)$. Hence, $\gamma - k + 1 \leq 1 + \frac{\alpha_2 - 1}{p}$. Since $\cosh u \geq 1$ for any $u \in \mathbb{R}$, we have $\cosh^{\gamma - k + 1} d \left( 0, \xi \right) \leq \cosh^{1 + \frac{\alpha_2 - 1}{p}} d \left( 0, \xi \right)$. Consequently, we get
				\begin{align*}
					\left| \cosh^{\gamma} d \left( 0, \xi \right) R_kf \left( \xi \right) \right| &\leq C  \cosh^{1 + \frac{\alpha_2 - 1}{p}} d \left( 0, \xi \right) \left( \int\limits_{d \left( 0, \xi \right)}^{\infty} \left| \tilde{f} \left( \cosh t \right) \right| \sinh^{\alpha_1 + n - 1} t \cosh^{\alpha_2} t \sinh^{1 - p - \alpha_2} t \ \mathrm{d}t \right)^{\frac{1}{p}} \\
					&\leq C \left( \int\limits_{d \left( 0, \xi \right)}^{\infty} \left| \tilde{f} \left( \cosh t \right) \right| \sinh^{\alpha_1 + n - 1} t \cosh^{\alpha_2} t \tanh^{1 - p - \alpha_2} t \ \mathrm{d}t \right)^{\frac{1}{p}},
				\end{align*}
				since $\cosh d \left( 0, \xi \right) \leq \cosh t$ for $t \geq d \left( 0, \xi \right)$ and $p + \alpha_2 - 1 > 0$. Since $t \geq d \left( 0, \xi \right) \geq 1$, we have $\tanh^{1 - p - \alpha_2} t \leq \tanh^{1 - p - \alpha_2} 1$. Therefore, in this case, we have,
				\begin{equation}
					\label{EndPointEstimateLPLInfinityHn1}
					\left| \cosh^{\gamma} d \left( 0, \xi \right) R_kf \left( \xi \right) \right| \leq C \left( \int\limits_{d \left( 0, \xi \right)}^{\infty} \left| \tilde{f} \left( \cosh t \right) \right| \sinh^{\alpha_1 + n - 1} t \cosh^{\alpha_2} t \ \mathrm{d}t \right)^{\frac{1}{p}} = C \| f \|_{L^{p, 1}_{\alpha_1, \alpha_2} \left( \mathbb{H}^n \right)}.
				\end{equation}
			
				\textbf{Case II:} Let us now consider $d \left( 0, \xi \right) < 1$. Here, we have from Equation \eqref{EndPointLPLInfinityHnInitialEstimate},
				\begin{align}
					\left| \cosh^{\gamma} d \left( 0, \xi \right) R_kf \left( \xi \right) \right| &\leq C \cosh^{\gamma - k + 1} d \left( 0, \xi \right) \int\limits_{d \left( 0, \xi \right)}^{\infty} \left| \tilde{f} \left( \cosh t \right) \right| \sinh^{k - 1} t \ \mathrm{d}t \nonumber \\
					\label{EndPointLPLInfinityHnInitialEstimate2}
					&= C \cosh^{\gamma - k + 1} d \left( 0, \xi \right) \left[ \int\limits_{d \left( 0, \xi \right)}^{1} \left| \tilde{f} \left( \cosh t \right) \right| \sinh^{k - 1} t \ \mathrm{d}t + \int\limits_{1}^{\infty} \left| \tilde{f} \left( \cosh t \right) \right| \sinh^{k - 1} t \ \mathrm{d}t \right].
				\end{align}
				Let us consider,
				\begin{equation}
					\label{I1EndPoint1Hn}
					I_1 = \int\limits_{d \left( 0, \xi \right)}^{1} \left| \tilde{f} \left( \cosh t \right) \right| \sinh^{k - 1} t \ \mathrm{d}t,
				\end{equation}
				and
				\begin{equation}
					\label{I2EndPoint1Hn}
					I_2 = \int\limits_{1}^{\infty} \left| \tilde{f} \left( \cosh t \right) \right| \sinh^{k - 1} t \ \mathrm{d}t.
				\end{equation}
				Using Equation \eqref{SinhNear0} in Equation \eqref{I1EndPoint1Hn}, we get
				\begin{align*}
					I_1 &\leq C \int\limits_{d \left( 0, \xi \right)}^{1} \left| \tilde{f} \left( \cosh t \right) \right| t^{k - 1} \mathrm{d}t.
				\end{align*}
				Further, from Equation \eqref{PolynomialInequalityLemma} with $\gamma = k$ and $p = \frac{\alpha_1 + n}{k}$, we get
				\begin{align*}
					I_1 &\leq C \left( \int\limits_{d \left( 0, \xi \right)}^{1} \left| \tilde{f} \left( \cosh t \right) \right| t^{\alpha_1 + n - 1} \mathrm{d}t \right)^{\frac{1}{p}}.
				\end{align*}
				Again using Equation \eqref{SinhNear0} and the fact that for small $u \in \mathbb{R}$, $\cosh u$ behaves like a constant, we get
				\begin{equation}
					\label{I1HnFinalEstimate}
					I_1 \leq C \left( \int\limits_{d \left( 0, \xi \right)}^{1} \left| \tilde{f} \left( \cosh t \right) \right| \sinh^{\alpha_1 + n - 1} t \cosh^{\alpha_2} t \mathrm{d}t \right)^{\frac{1}{p}} = C \| f \|_{L^{p, 1}_{\alpha_1, \alpha_2} \left( \mathbb{H}^n \right)}.
				\end{equation}
				For the integral $I_2$, we use the analysis done in Case I (by replacing $d \left( 0, \xi \right)$ with $1$) to obtain
				\begin{equation}
					\label{I2HnFinalEstimate}
					{I_2 \leq C \| f \|_{L^{p, 1}_{\alpha_1, \alpha_2} \left( \mathbb{H}^n \right)}},
				\end{equation}
				since $\sinh^{1 - p - \alpha_2} t \leq \sinh^{1 - p - \alpha_2} 1$. Using Equations \eqref{I1HnFinalEstimate} and \eqref{I2HnFinalEstimate} in Equation \eqref{EndPointLPLInfinityHnInitialEstimate2}, we get
				\begin{equation}
					\label{EndPointEstimateLPLInfinityHn2}
					\left| \cosh^{\gamma} d \left( 0, \xi \right) R_kf \left( \xi \right) \right| \leq C \cosh^{\gamma - k + 1} d \left( 0, \xi \right) \| f \|_{L^{p, 1}_{\alpha_1, \alpha_2} \left( \mathbb{H}^n \right)} \leq C \| f \|_{L^{p, 1}_{\alpha_1, \alpha_2} \left( \mathbb{H}^n \right)}.
				\end{equation}
				In the last inequality, we have used the fact that since $d \left( 0, \xi \right) < 1$, any power of $\cosh d \left( 0, \xi \right)$ behaves like a constant.
				
				This completes the proof of \ref{EndPointLPLInfinityHnWeightedAProof}.
				\item[\mylabel{EndPointLPLInfinityHnWeightedBProof}{(B)}] {Now, we prove the result for $k = 1$. Here, we have assumed, additionally, that $\alpha_1 + n = p \geq 2$. We have, from Equation \eqref{KPlaneTransformRadialHnEquation}, for any characteristic function $f$ of measurable subset of $\mathbb{H}^n$ with $f \in L^{p, 1}_{\alpha_1, \alpha_2} \left( \mathbb{H}^n \right)$,
				$$\cosh^{\gamma} d \left( 0, \xi \right) R_1f \left( \xi \right) = C \cosh^{\gamma} d \left( 0, \xi \right) \int\limits_{d \left( 0, \xi \right)}^{\infty} \tilde{f} \left( \cosh t \right) \left( \cosh^2 t - \cosh^2 d \left( 0, \xi \right) \right)^{- \frac{1}{2}} \sinh t \ \mathrm{d}t.$$
				Let us substitute $\cosh t = \cosh x \cosh d \left( 0, \xi \right)$. After a simplification, we get
				\begin{align*}
					\cosh^{\gamma} d \left( 0, \xi \right) R_1f \left( \xi \right) &= C \cosh^{\gamma} d \left( 0, \xi \right) \int\limits_{0}^{\infty} \tilde{f} \left( \cosh x \cosh d \left( 0, \xi \right) \right) \mathrm{d}x \\
					&= C \cosh^{\gamma} d \left( 0, \xi \right) \left[ \int\limits_{0}^{1} \tilde{f} \left( \cosh x \cosh d \left( 0, \xi \right) \right) \mathrm{d}x + \int\limits_{1}^{\infty} \tilde{f} \left( \cosh x \cosh d \left( 0, \xi \right) \right) \mathrm{d}x \right].
				\end{align*}
				We know that $\alpha_1 + \alpha_2 + n - 1 \geq 0$, and for $x \geq 1$, $\sinh x \geq \sinh 1$. Hence, $1 \leq \frac{\sinh^{\frac{\alpha_1 + \alpha_2 + n - 1}{p}} x}{\sinh^{\frac{\alpha_1 + \alpha_2 + n - 1}{p}} 1}$. Therefore, we have,
				\begin{align*}
					\cosh^{\gamma} d \left( 0, \xi \right) R_1f \left( \xi \right) &\leq C \cosh^{\gamma} d \left( 0, \xi \right) \left[ \int\limits_{0}^{1} \tilde{f} \left( \cosh x \cosh d \left( 0, \xi \right) \right) \mathrm{d}x + \int\limits_{1}^{\infty} \tilde{f} \left( \cosh x \cosh d \left( 0, \xi \right) \right) \sinh^{\frac{\alpha_1 + \alpha_2 + n - 1}{p}} x \ \mathrm{d}x \right].
				\end{align*}
				We use Inequality \eqref{ConsequenceAnnulusInequality} in the first integral, and Inequality \eqref{HyperbolicMainLemmaCoshEquation} in the second integral. Then, we get,
				\begin{align*}
					\cosh^{\gamma} d \left( 0, \xi \right) R_1f \left( \xi \right) &\leq C \cosh^{\gamma} d \left( 0, \xi \right) \left[ \left( \int\limits_{0}^{1} \tilde{f} \left( \cosh x \cosh  d \left( 0, \xi \right) \right) x^{p - 1} \mathrm{d}x \right)^{\frac{1}{p}} \right. \\
					&\left. + \left( \int\limits_{1}^{\infty} \tilde{f} \left( \cosh x \cosh d \left( 0, \xi \right) \right) \sinh^{\alpha_1 + n - 1} x \cosh^{\alpha_2} x \ \mathrm{d}x \right)^{\frac{1}{p}} \right] \\
					&\leq C \cosh^{\gamma} \left[ \left( \int\limits_{0}^{1} \tilde{f} \left( \cosh x \cosh d \left( 0, \xi \right) \right) \sinh^{p - 1} x \ \mathrm{d}x \right)^{\frac{1}{p}} \right. \\
					&\left. + \left( \int\limits_{1}^{\infty} \tilde{f} \left( \cosh x \cosh d \left( 0, \xi \right) \right) \sinh^{\alpha_1 + n - 1} x \cosh^{\alpha_2} x \ \mathrm{d}x \right)^{\frac{1}{p}} \right],
				\end{align*}
				where, we have used the fact that $x \leq \sinh x$ and $p \geq 1$. Further, we observe that for $0 \leq x \leq 1$, we have $\cosh^{\alpha_2} x \geq \min \left\lbrace 1, \cosh^{\alpha_2} 1 \right\rbrace$. Also, we have that $p = \alpha_1 + n$. Therefore, we get,
				\begin{align*}
					\cosh^{\gamma} d \left( 0, \xi \right) R_1f \left( \xi \right) &\leq C \cosh^{\gamma} d \left( 0, \xi \right) \left[ \left( \int\limits_{0}^{1} \tilde{f} \left( \cosh x \cosh d \left( 0, \xi \right) \right) \sinh^{\alpha_1 + n - 1} x \cosh^{\alpha_2} x \ \mathrm{d}x \right)^{\frac{1}{p}} \right. \\
					&\left. + \left( \int\limits_{1}^{\infty} \tilde{f} \left( \cosh x \cosh  d \left( 0, \xi \right) \right) \sinh^{\alpha_1 + n - 1} x \ \cosh^{\alpha_2} x \ \mathrm{d}x \right)^{\frac{1}{p}} \right] \\
					&\leq C \cosh^{\gamma} d \left( 0, \xi \right) \left( \int\limits_{0}^{\infty} \tilde{f} \left( \cosh x \cosh d \left( 0, \xi \right) \right) \sinh^{\alpha_1 + n - 1} x \ \cosh^{\alpha_2} x \ \mathrm{d}x \right)^{\frac{1}{p}}.
				\end{align*}
				Now, we substitute $\cosh x \cosh d \left( 0, \xi \right) = \cosh t$ and obtain
				$$\cosh^{\gamma} d \left( 0, \xi \right) R_1f \left( \xi \right) \leq C \cosh^{\gamma} d \left( 0, \xi \right) \left( \int\limits_{0}^{\infty} \tilde{f} \left( \cosh t \right) \left( \frac{\cosh^2 t}{\cosh^2 d \left( 0, \xi \right)} - 1 \right)^{\frac{\alpha_1 + n}{2} - 1} \left( \frac{\cosh t}{\cosh d \left( 0, \xi \right)} \right)^{\alpha_2} \frac{\sinh t}{\cosh d \left( 0, \xi \right)} \mathrm{d}t \right)^{\frac{1}{p}}.$$
				Since from our assumptions, we have $p = \alpha_1 + n \geq 2$, we get $\left( \frac{\cosh^2 t}{\cosh^2 d \left( 0, \xi \right)} - 1 \right)^{\frac{\alpha_1 + n}{2} - 1} \leq \frac{\sinh^{\alpha_1 + n - 2} t}{\cosh^{\alpha_1 + n - 1} d \left( 0, \xi \right)}$. Therefore,
				$$\cosh^{\gamma} d \left( 0, \xi \right) R_1f \left( \xi \right) \leq C \cosh^{\gamma - \frac{\alpha_1 + \alpha_2 + n - 1}{p}} d \left( 0, \xi \right) \left( \int\limits_{0}^{\infty} \tilde{f} \left( \cosh t \right) \sinh^{\alpha_1 + n - 1} t \cosh^{\alpha_2} t \ \mathrm{d}t \right)^{\frac{1}{p}} \leq C \| f \|_{L^{p, 1}_{\alpha_1, \alpha_2} \left( \mathbb{H}^n \right)}.$$
				In the final inequality, we have used the fact that $\gamma \leq \frac{\alpha_1 + \alpha_2 + n - 1}{p}$ and hence $\cosh^{\gamma - \frac{\alpha_1 + \alpha_2 + n - 1}{p}} d \left( 0, \xi \right) \leq 1$, for any $\xi \in \Xi_k \left( \mathbb{H}^n \right)$. The necessity of $p \geq 2$ when $k = 1$ is shown in Example \ref{CounterExampleLpLinfinityHn}.}
			\end{enumerate}
			{The necessity of the condition $\gamma \leq \frac{\alpha_1 + \alpha_2 + n - 1}{p}$ is shown in Remark \ref{NecessaryDoubleAnnulusHn}.}
		\end{proof}
		\begin{remark}
			\label{NecessaryDoubleAnnulusHn}
			\normalfont
			{The necessity of the condition $\gamma \leq \frac{\alpha_1 + \alpha_2 + n - 1}{p}$ was shown in Theorem \ref{NecessaryLPLInfinityHn} for $L^p$-$L^{\infty}$ norm inequality. Let us now see that the same is necessary for the Lorentz norm estimate of Inequality \eqref{RequiredEndPointHn2}. Here, we fix $\alpha_1, \alpha_2 \in \mathbb{R}$ such that $\alpha_1 + \alpha_2 > k - n$ and $1 \leq p = \frac{\alpha_1 + n}{k} < \frac{\alpha_1 + \alpha_2 + n - 1}{k - 1}$ or $\alpha_1 + \alpha_2 = k - n$ and $p = 1$. To do so, let us consider the function $f_{a, b}: \mathbb{H}^n \rightarrow \mathbb{C}$, given by $f_{a, b} \left( x \right) = \chi_{a, b} \left( \cosh d \left( 0, x \right) \right)$, for $b > a > 1$. Here, we consider the case $b = 2a$, and denote $f_a = f_{a, 2a}$. It is clear from Equation \eqref{LPMassWeightedAnnulusHn}, by substituting $x = a \sqrt{u}$,
			\begin{align*}
				\| f_a \|_{L^{p, 1}_{\alpha_1, \alpha_2} \left( \mathbb{H}^n \right)} &= C a^{\frac{\alpha_2 + 1}{p}} \left[ \int\limits_{1}^{4} \left( a^2u - 1 \right)^{\frac{\alpha_1 + n}{2} - 1} u^{\frac{\alpha_2 - 1}{2}} \mathrm{d}u \right]^{\frac{1}{p}}.
			\end{align*}
			For ``large" $a > 1$, say $a \geq 2$, we have that $a^2u - 1$ behaves like $a^2$, whenever $1 \leq u \leq 4$. Hence,
			\begin{equation}
				\label{LorentzNormDoubleAnnulusHnEstimate}
				\| f_a \|_{L^{p, 1}_{\alpha_1, \alpha_2} \left( \mathbb{H}^n \right)} \leq C a^{\frac{\alpha_2 + 1}{p} + \frac{2}{p} \left( \frac{\alpha_1 + n}{2} - 1 \right)} = C a^{\frac{\alpha_2 - 1}{p} + k}.
			\end{equation}
			As for the $k$-plane transform, we use $\cosh t = x \cosh d \left( 0, \xi \right)$ in Equation \eqref{KPlaneTransformRadialHnEquation} and simply to obtain,
			\begin{align}
				R_kf_a \left( \xi \right) &= C \int\limits_{1}^{\infty} \chi_{\left( a, 2a \right)} \left( x \cosh d \left( 0, \xi \right) \right) \left( x^2 - 1 \right)^{\frac{k}{2} - 1} \mathrm{d}x \nonumber \\
				\label{KPlaneDoubleAnnulusHnEquation}
				&= C \begin{cases}
							0, & \cosh d \left( 0, \xi \right) \geq 2a. \\
							\int\limits_{1}^{\frac{2a}{\cosh d \left( 0, \xi \right)}} \left( x^2 - 1 \right)^{\frac{k}{2} - 1} \mathrm{d}x, & a \leq \cosh d \left( 0, \xi \right) \leq 2a. \\
							\int\limits_{\frac{a}{\cosh d \left( 0, \xi \right)}}^{\frac{2a}{\cosh d \left( 0, \xi \right)}} \left( x^2 - 1 \right)^{\frac{k}{2} - 1} \mathrm{d}x, & \cosh d \left( 0, \xi \right) < a.
						\end{cases}
			\end{align}
			By choosing $\xi \in \Xi_k \left( \mathbb{H}^n \right)$ such that $\cosh d \left( 0, \xi \right) = a$, it is clear that
			\begin{equation}
				\label{RkfEstimateLorentzDoubleAnnulusHn}
				\| \cosh^{\gamma} d \left( 0, \cdot \right) R_kf \|_{L^{\infty} \left( \Xi_k \left( \mathbb{H}^n \right) \right)} \geq C a^{\gamma} \int\limits_{1}^{2} \left( x^2 - 1 \right)^{\frac{k}{2} - 1} \mathrm{d}x = C a^{\gamma}.
			\end{equation}
			Using Inequalities \eqref{RkfEstimateLorentzDoubleAnnulusHn} and \eqref{LorentzNormDoubleAnnulusHnEstimate} in Inequality \eqref{RequiredEndPointHn2}, we get for $a \geq 2$,
			$$a^{\gamma} \leq C a^{\frac{\alpha_2 - 1}{p} + k}.$$
			As $a \rightarrow \infty$, this forces $\gamma \leq \frac{\alpha_2 - 1}{p} + k = \frac{\alpha_1 + \alpha_2 + n - 1}{p}$, since $p = \frac{\alpha_1 + n}{k}$.}
		\end{remark}
		Theorem \ref{RequiredEndPointHn2} answers Question \ref{EndPointHnQuestion1}. We first get the necessary conditions on $\gamma_1$ and $\gamma_2$ for Inequality \eqref{RequiredEndPointHn1}.
		\begin{theorem}
			\label{EndPointMixedWeightNecessary}
			Let $k \geq 2$ and $\alpha_1 + \alpha_2 > k - n$. For $p = \frac{\alpha_1 + \alpha_2 + n - 1}{k - 1}$, Inequality \eqref{RequiredEndPointHn1} holds only if 
			\begin{equation}
				\label{EndPointHn1NecessaryConditions}
				{\gamma_1 \geq \max \left\lbrace 0, \frac{\alpha_1 + n}{k} \right\rbrace} \text{ and } {\gamma_1 + \gamma_2 \leq k - 1}.
			\end{equation}
		\end{theorem}
		\begin{proof}
			{We prove this result in two cases.}
			
			{\textbf{Case I:}} First, we assume that $\alpha_1 > - n$. Let us consider the function {$f_{\lambda}: \mathbb{H}^n \rightarrow \mathbb{C}$, defined as} $f_{\lambda} \left( x \right) = \chi_{B \left( 0, \lambda \right)} \left( x \right)${, for $\lambda > 0$}. We have,
			\begin{align}
				\label{LorentzNormDoubleAnnulusHnInit}
				\| f_{\lambda} \|_{L^{p, 1}_{\alpha_1, \alpha_2} \left( \mathbb{H}^n \right)} &= C \left[ \int\limits_{0}^{\lambda} \cosh^{\alpha_2} t \sinh^{\alpha_1 + n - 1} t \mathrm{d}t \right]^{\frac{1}{p}}.
			\end{align}
			By substituting $\sinh^2 t = u$, we get from Equation \eqref{IntegralForm2F1},
			\begin{align*}
				\| f_{\lambda} \|_{L^{p, 1}_{\alpha_1, \alpha_2} \left( \mathbb{H}^n \right)} &= C \left[ \int\limits_{0}^{\sinh^{{2}} \lambda} \left( 1 + u \right)^{\frac{\alpha_2 - 1}{2}} u^{\frac{\alpha_1 + n}{2} - 1} \mathrm{d}u \right]^{\frac{1}{p}} = C \sinh^{\frac{\alpha_1 + n}{p}} \lambda \left[ {}_2F_1 \left( \frac{1 - \alpha_2}{2}, \frac{\alpha_1 + n}{2}; 1 + \frac{\alpha_1 + n}{2}; - \sinh^2 \lambda \right) \right]^{\frac{1}{p}}.
			\end{align*}
			By using the transformation given in Equation \eqref{Transformation2F1} \eqref{2F1Transformation1}, we get
			\begin{equation}
				\label{FBallNormMixed}
				\| f_{\lambda} \|_{L^{p, 1}_{\alpha_1, \alpha_2} \left( \mathbb{H}^n \right)} = C \sinh^{\frac{\alpha_1 + n}{p}} \lambda \cosh^{\frac{\alpha_2 - 1}{p}} \lambda \left[ {}_2F_1 \left( \frac{1 - \alpha_2}{2}, 1; 1 + \frac{\alpha_1 + n}{2}; \tanh^2 \lambda \right) \right]^{\frac{1}{p}}.
			\end{equation}
			It is clear from Theorem \ref{Behaviour2F1} that the hypergeometric factor is bounded. Hence, we have,
			\begin{equation}
				\label{FBallNormMixedEstimate}
				\| f_{\lambda} \|_{L^{p, 1}_{\alpha_1, \alpha_2} \left( \mathbb{H}^n \right)} \leq C \sinh^{\frac{\alpha_1 + n}{p}} \lambda \cosh^{\frac{\alpha_2 - 1}{p}} \lambda.
			\end{equation}
			The $k$-plane transform of $f_{\lambda}$ is given in Equation \eqref{KPlaneballHnEquation}. {By choosing $d \left( 0, \xi \right)$ such that $\sinh^2 d \left( 0, \xi \right) = \frac{1}{2} \sinh^2 \lambda$, we have $\cosh^2 d \left( 0, \xi \right) = \frac{1}{2} \left( 1 + \cosh^2 \lambda \right)$ and $\tanh^2 d \left( 0, \xi \right) = \frac{\sinh^2 \lambda}{1 + \cosh^2 \lambda}$.} Hence from Equation \eqref{KPlaneBallHn}, we get
			\begin{equation}
				\label{SupNormEstimateRkFCosh}
				\| {\sinh^{\gamma_1} d \left( 0, \xi \right)} \cosh^{\gamma_{{2}}} d \left( 0, \cdot \right) R_kf_{\lambda} \left( \cdot \right) \|_{\infty} \geq C {\sinh^{\gamma_1} \lambda} \left( 1 + \cosh^2 \lambda \right)^{\frac{\gamma_{{2}}}{2}} \tanh^k \lambda \ {}_2F_1 \left( \frac{k + 1}{2}, \frac{k}{2}; 1 + \frac{k}{2}; \frac{\sinh^2 \lambda}{1 + \sinh^2 \lambda + \cosh^2 \lambda} \right).
			\end{equation}
			Using Inequalities \eqref{FBallNormMixedEstimate} and \eqref{SupNormEstimateRkFCosh} in Inequality \eqref{RequiredEndPointHn1}, we have
			$${\sinh^{\gamma_1} \lambda} \left( 1 + \cosh^2 \lambda \right)^{\frac{\gamma_{{2}}}{2}} \tanh^k \lambda \ {}_2F_1 \left( \frac{k + 1}{2}, \frac{k}{2}; 1 + \frac{k}{2}; \frac{\sinh^2 \lambda}{1 + \sinh^2 \lambda + \cosh^2 \lambda} \right) \leq C \sinh^{\frac{\alpha_1 + n}{p}} \lambda \cosh^{\frac{\alpha_2 - 1}{p}} \lambda.$$
			As $\lambda \rightarrow 0$, we get
			$$\lambda^{{\gamma_1 +} k} \leq C \lambda^{\frac{\alpha_1 + n}{p}}.$$
			The above inequality forces ${\gamma_1 \geq \frac{\alpha_1 + n}{p} - k}$. {The necessity of $\gamma_1 \geq 0$ is also clear, for otherwise we have 
			$$\| \sinh^{\gamma_1} d \left( 0, \cdot \right) \cosh^{\gamma_2} d \left( 0, \cdot \right) R_kf_{\lambda} \|_{L^{\infty} \left( \Xi_k \left( \mathbb{H}^n \right) \right)} = + \infty.$$}
			On the other hand, for $\lambda \rightarrow \infty$, we have
			$$e^{{\left( \gamma_1 + \gamma_2 \right)} \lambda} \leq C e^{\frac{\alpha_1 + \alpha_2 + n - 1}{p}\lambda} = C e^{\left( k - 1 \right) \lambda}.$$
			Hence, we must have ${\gamma_1 + \gamma_2} \leq k - 1$.
			
			{\textbf{Case II:} Now, we consider $\alpha_1 \leq -n$. To see the necessity of the condition $\gamma_1 + \gamma_2 \leq k - 1$, we consider the function $f_a$ described in Remark \ref{NecessaryDoubleAnnulusHn}. From Equation \eqref{LorentzNormDoubleAnnulusHnEstimate} {and} the fact {that} $p = \frac{\alpha_1 + \alpha_2 + n - 1}{k - 1}$, {we have}
			\begin{equation}
				\label{LorentzNormDoubleAnnulusHnEndPoint2}
				\| f_a \|_{L^{p, 1}_{\alpha_1, \alpha_2} \left( \mathbb{H}^n \right)} \leq C a^{\frac{\alpha_1 + \alpha_2 + n - 1}{p}} = C a^{k - 1},
			\end{equation}
			for large $a > 1$. Now, by choosing $\xi \in \Xi_k \left( \mathbb{H}^n \right)$ such that $d \left( 0, \xi \right) = a$, we get
			\begin{equation}
				\label{LInfinityNormWeightedKPlaneDoubleAnnulusHn}
				\| \sinh^{\gamma_1} d \left( 0, \cdot \right) \cosh^{\gamma_2} d \left( 0, \cdot \right) R_kf_{a} \|_{L^{\infty} \left( \Xi_k \left( \mathbb{H}^n \right) \right)} \geq C \left( a^2 - 1 \right)^{\frac{\gamma_1}{2}} a^{\gamma_2} \geq C a^{\gamma_1 + \gamma_2},
			\end{equation}
			for large $a > 1$. Using Inequalities \eqref{LorentzNormDoubleAnnulusHnEndPoint2} and \eqref{LInfinityNormWeightedKPlaneDoubleAnnulusHn} in Inequality \eqref{RequiredEndPointHn1}, we get
			$$a^{\gamma_1 + \gamma_2} \leq C a^{k - 1},$$
			for large $a > 1$. As $a \rightarrow \infty$, this forces $\gamma \leq k - 1$.}
			
			{It is also clear from Equation \eqref{LorentzNormDoubleAnnulusHnInit} that
			\begin{equation}
				\label{LorentzNormDoubleAnnulusHnNearZero}
				\| f_a \|_{L^{p, 1}_{\alpha_1, \alpha_2} \left( \mathbb{H}^n \right)} \leq C a^{\frac{\alpha_2 + 1}{p}} \left( 4a^2 - 1 \right)^{\frac{1}{p} \left( \frac{\alpha_1 + n}{2} - 1 \right)}.
			\end{equation}
			Using Inequalities \eqref{LorentzNormDoubleAnnulusHnNearZero} and \eqref{LInfinityNormWeightedKPlaneDoubleAnnulusHn} in Inequality \eqref{RequiredEndPointHn1}, we get for any $a > 1$,
			$$\left( a^2 - 1 \right)^{\frac{\gamma_1}{2}} a^{\gamma_2} \leq C a^{\frac{\alpha_2 + 1}{p}} \left( 4a^2 - 1 \right)^{\frac{1}{p} \left( \frac{\alpha_1 + n}{2} - 1 \right)}.$$
			As $a \rightarrow 0$, the above inequality gives us $\gamma_1 \geq 0 = \max \left\lbrace 0, \frac{\alpha_1 + n}{p} - k \right\rbrace$.}
		\end{proof}
		We now see that the conditions mentioned in Theorem \ref{EndPointMixedWeightNecessary} are also sufficient.
		\begin{theorem}
			\label{EndPointMixedWeightFull}
			{Let $k \geq 2$ and $\alpha_1 + \alpha_2 > k - n$. For $p = \frac{\alpha_1 + \alpha_2 + n - 1}{k - 1} > 1$}, Inequality \eqref{RequiredEndPointHn1} holds if and only if {$\gamma_1$ and $\gamma_2$ satisfy the conditions of Equation \eqref{EndPointHn1NecessaryConditions}}.
		\end{theorem}
		\begin{proof}
			We have seen the necessity of the conditions mentioned in Equation \eqref{EndPointHn1NecessaryConditions} in Theorem \ref{EndPointMixedWeightNecessary}. To see the sufficiency of these conditions, {we consider the following}. {Due to Theorem \ref{LorentzSpaceBoundedness}, we only consider characteristic functions $f$ of radially symmetric sets for which $f \in L^{p, 1}_{\alpha_1, \alpha_2} \left( \mathbb{H}^n \right)$.}
			
			{\textbf{Case I}: First, let us assume that $p \geq 1 - \alpha_2$. In this case, we may use Lemma \ref{HyperbolicMainLemmaCosh}. We start with Equation \eqref{KPlaneTransformRadialHnEquation3} and use the fact that $k \geq 2$ to obtain
			$$\| \sinh^{\gamma_1} d \left( 0, \cdot \right) \cosh^{\gamma_2} d \left( 0, \cdot \right) R_kf \|_{L^{\infty} \left( \Xi_k \left( \mathbb{H}^n \right) \right)} \leq C \sinh^{\gamma_1} d \left( 0, \xi \right) \cosh^{\gamma_2 - k + 1} d \left( 0, \xi \right) \int\limits_{d \left( 0, \xi \right)}^{\infty} \left| \tilde{f} \left( \cosh t \right) \right| \sinh^{k - 1} t \ \mathrm{d}t.$$
			Since $\gamma_1 \geq 0$ and $\gamma_1 + \gamma_2 \leq k - 1$, we have $\sinh^{\gamma_1} d \left( 0, \xi \right) \cosh^{\gamma_2} d \left( 0, \xi \right) \leq C e^{\left( \gamma_1 + \gamma_2 - k + 1 \right) d \left( 0, \xi \right)} \leq C$. Hence, by the use of Equation \eqref{HyperbolicMainLemmaCoshEquation} of Lemma \ref{HyperbolicMainLemmaCosh}, we get
			$$\| \sinh^{\gamma_1} d \left( 0, \cdot \right) \cosh^{\gamma_2} d \left( 0, \cdot \right) R_kf \|_{L^{\infty} \left( \Xi_k \left( \mathbb{H}^n \right) \right)} \leq C \left( \int\limits_{0}^{\infty} \tilde{f} \left( \cosh t \right) \sinh^{\alpha_1 + n - 1} t \cosh^{\alpha_2} t \ \mathrm{d}t \right)^{\frac{1}{p}} = C \| f \|_{L^{p, 1}_{\alpha_1, \alpha_2} \left( \mathbb{H}^n \right)}.$$}
			
			{\textbf{Case II}: We now consider $1 \leq p < 1 - \alpha_2$. We notice that in this case, we cannot use Lemma \ref{HyperbolicMainLemmaCosh} to produce $\cosh^{\alpha_2} t$ and therefore the above method fails. To prove the desired result, we now divide our analysis into two parts.}
			
			{\textbf{Case II(A)}: First, let us consider $d \left( 0, \xi \right) \geq 1$. In this case, using Equation \eqref{KPlaneTransformRadialHnEquation3} and the fact that $k \geq 2$, we get
			$$\| \sinh^{\gamma_1} d \left( 0, \cdot \right) \cosh^{\gamma_2} d \left( 0, \cdot \right) R_kf \|_{L^{\infty} \left( \Xi_k \left( \mathbb{H}^n \right) \right)} \leq C \sinh^{\gamma_1} d \left( 0, \xi \right) \cosh^{\gamma_2 - k + 1} d \left( 0, \xi \right) \int\limits_{d \left( 0, \xi \right)}^{\infty} \left| \tilde{f} \left( \cosh t \right) \right| \sinh^{k - 1} t \ \mathrm{d}t.$$
			Now, we use Equation \eqref{HyperbolicMainLemmaCoshEquation} with $\alpha = 0$ of Lemma \ref{HyperbolicMainLemmaCosh} to obtain
			\begin{align*}
				&\| \sinh^{\gamma_1} d \left( 0, \cdot \right) \cosh^{\gamma_2} d \left( 0, \cdot \right) R_kf \|_{L^{\infty} \left( \Xi_k \left( \mathbb{H}^n \right) \right)} \\
				&\leq C \sinh^{\gamma_1} d \left( 0, \xi \right) \cosh^{\gamma_2 - k + 1} d \left( 0, \xi \right) \left( \int\limits_{d \left( 0, \xi \right)}^{\infty} \left| \tilde{f} \left( \cosh t \right) \right| \sinh^{\alpha_1 + \alpha_2 + n - 1} t \ \mathrm{d}t \right)^{\frac{1}{p}} \\
				&= C \sinh^{\gamma_1} d \left( 0, \xi \right) \cosh^{\gamma_2} d \left( 0, \xi \right) \left( \int\limits_{d \left( 0, \xi \right)}^{\infty} \tilde{f} \left( \cosh t \right) \sinh^{\alpha_1 + n - 1} t \cosh^{\alpha_2} t \tanh^{\alpha_2} t \ \mathrm{d}t \right)^{\frac{1}{p}}.
			\end{align*}
			Now, we observe that since $1 \leq p < 1 - \alpha_2$, we must have $\alpha_2 < 0$. Therefore, for $t \geq d \left( 0, \xi \right) \geq 1$, we have $\tanh^{\alpha_2} t \leq \tanh^{\alpha_2} d \left( 0, \xi \right) \leq \tanh^{\alpha_2} 1$. Also, since $\gamma_1 \geq 0$ and $\gamma_1 + \gamma_2 \leq k - 1$, we have as in Case I, $\sinh^{\gamma_1} d \left( 0, \xi \right) \cosh^{\gamma_2 - k + 1} d \left( 0, \xi \right) \leq C$. That is, we get
			$$\| \sinh^{\gamma_1} d \left( 0, \cdot \right) \cosh^{\gamma_2} d \left( 0, \cdot \right) R_kf \|_{L^{\infty} \left( \Xi_k \left( \mathbb{H}^n \right) \right)} \leq C \left( \int\limits_{0}^{\infty} \tilde{f} \left( \cosh t \right) \sinh^{\alpha_1 + n - 1} t \cosh^{\alpha_2} t \ \mathrm{d}t \right)^{\frac{1}{p}} = C \| f \|_{L^{p, 1}_{\alpha_1, \alpha_2} \left( \mathbb{H}^n \right)}.$$}
			
			{\textbf{Case II(B)}: Let us now consider the case when $d \left( 0, \xi \right) < 1$. Here, we write
			\begin{align*}
				\| \sinh^{\gamma_1} d \left( 0, \cdot \right) \cosh^{\gamma_2} d \left( 0, \cdot \right) R_kf \|_{L^{\infty} \left( \Xi_k \left( \mathbb{H}^n \right) \right)} \leq C \sinh^{\gamma_1} d \left( 0, \xi \right) \cosh^{\gamma_2 - k + 1} d \left( 0, \xi \right) &\left[ \int\limits_{d \left( 0, \xi \right)}^{1} \left| \tilde{f} \left( \cosh t \right) \right| \sinh^{k - 1} t \ \mathrm{d}t \right. \\
				&\left. + \int\limits_{1}^{\infty} \left| \tilde{f} \left( \cosh t \right) \right| \sinh^{k - 1} t \ \mathrm{d}t \right].
			\end{align*}
			Let
			$$I_1 = \int\limits_{d \left( 0, \xi \right)}^{1} \left| \tilde{f} \left( \cosh t \right) \right| \sinh^{k - 1} t \ \mathrm{d}t,$$
			and
			$$I_2 = \int\limits_{1}^{\infty} \left| \tilde{f} \left( \cosh t \right) \right| \sinh^{k - 1} t \ \mathrm{d}t.$$
			Then, we know from Case II(A) that
			$$\sinh^{\gamma_1} 1 \cosh^{\gamma_2} 1 I_2 \leq C \| f \|_{L^{p, 1}_{\alpha_1, \alpha_2} \left( \mathbb{H}^n \right)}.$$
			We only need to estimate $I_2$. To do so, we start by observing that when $t \leq 1$, we have $\sinh^{k - 1} t \leq C t^{k - 1}$. Therefore,
			$$I_1 \leq C \int\limits_{d \left( 0, \xi \right)}^{1} \tilde{f} \left( \cosh t \right) t^{k - 1} \mathrm{d}t.$$
			Now, using Equation \eqref{PolynomialInequalityLemma} of Lemma \ref{AnnulusEstimateLemma}, we get
			$$I_1 \leq C \left( \int\limits_{d \left( 0, \xi \right)}^{\infty} \tilde{f} \left( \cosh t \right) t^{pk - 1} \mathrm{d}t \right)^{\frac{1}{p}}.$$
			Now, we see that $pk - 1 = \alpha_1 + \alpha_2 + n - 1 + p - 1$. Also, for $t \in \left( d \left( 0, \xi \right), 1 \right)$, $\sinh t$ behaves like $t$. Therefore, we have,
			$$I_1 \leq C \left( \int\limits_{d \left( 0, \xi \right)}^{1} \tilde{f} \left( \cosh t \right) \sinh^{\alpha_1 + n - 1} t \sinh^{\alpha_2 + p - 1} t \ \mathrm{d}t \right)^{\frac{1}{p}}.$$
			From our assumption, we have $\alpha_2 + p - 1 < 0$. Therefore, $\sinh^{\alpha_2 + p - 1} t \leq \sinh^{\alpha_2 + p - 1} d \left( 0, \xi \right)$, for $t \geq d \left( 0, \xi \right)$. Now, using the fact that for small values of $t$, $\cosh t$ behaves as a constant, we get
			$$I_1 \leq C \sinh^{1 + \frac{\alpha_2 - 1}{p}} d \left( 0, \xi \right) \| f \|_{L^{p, 1}_{\alpha_1, \alpha_2} \left( \mathbb{H}^n \right)}.$$
			Combining the analysis done on the integrals $I_1$ and $I_2$, we conclude that
			\begin{align*}
				\| \sinh^{\gamma_1} d \left( 0, \cdot \right) \cosh^{\gamma_2} d \left( 0, \cdot \right) R_kf \|_{L^{\infty} \left( \Xi_k \left( \mathbb{H}^n \right) \right)} &\leq C \sinh^{\gamma_1 + 1 + \frac{\alpha_2 - 1}{p}} d \left( 0, \xi \right) \cosh^{\gamma_2 - k + 1} d \left( 0, \xi \right) \| f \|_{L^{p, 1}_{\alpha_1, \alpha_2} \left( \mathbb{H}^n \right)} \\
				&+ \sinh^{\gamma_1} d \left( 0, \xi \right) \cosh^{\gamma_2 - k + 1} d \left( 0, \xi \right) \| f \|_{L^{p, 1}_{\alpha_1, \alpha_2} \left( \mathbb{H}^n \right)}.
			\end{align*}
			In both the terms in the above inequality, since $d \left( 0, \xi \right) < 1$, the hyperbolic cosine behaves like a constant. Also, owing to the fact that $\gamma_1 \geq 0$, the hyperbolic sine in the second term can be dominated by a constant. Now, from the necessary conditions, we know that we require $\gamma_1 \geq \frac{\alpha_1 + n}{p} - k = -1 - \frac{\alpha_2 - 1}{p}$. Therefore, the hyperbolic sine in the first term can also be dominated by a constant. Consequently, we have,
			$$\| \sinh^{\gamma_1} d \left( 0, \cdot \right) \cosh^{\gamma_2} d \left( 0, \cdot \right) R_kf \|_{L^{\infty} \left( \Xi_k \left( \mathbb{H}^n \right) \right)} \leq C \| f \|_{L^{p, 1}_{\alpha_1, \alpha_2} \left( \mathbb{H}^n \right)}.$$}
			This completes the proof!
		\end{proof}
		\subsection{The Sphere \texorpdfstring{$\mathbb{S}^n$}{}}
			In this section, we consider the end-point estimates for the ``end-points" {mentioned} in {Remark \ref{RemarkEndPointsSn}}. Particularly, we ask the following questions.
			\begin{question}
				\label{EndPointSnQuestion1}
				{Given $\alpha_1, \alpha_2 \in \mathbb{R}$ and $p \geq \max \left\lbrace 1, 1 + \alpha_2 \right\rbrace$ with $p > 1 + \alpha_2$ when $\alpha_2 > 0$,} what are the admissible values of $\gamma \in \mathbb{R}$ such that the following inequality holds for all radial functions in $L^{p, 1}_{\alpha_1, \alpha_2} \left( \mathbb{S}^n \right)$, when $p = \frac{\alpha_1 + n}{k}$?
				\begin{equation}
					\label{RequiredEndPointSn1}
					\| \cos^{\gamma} d \left( 0, \cdot \right) R_kf \|_{L^{\infty} \left( \Xi_k \left( \mathbb{S}^n \right) \right)} \leq C \| f \|_{L^{p, 1}_{\alpha_1, \alpha_2} \left( \mathbb{S}^n \right)}.
				\end{equation}
			\end{question}
			We have {also} seen that due to the existence conditions of Theorem \ref{ExistenceSN}, $p = 1 + \alpha_2$ is an end-point for $\alpha_2 > 0$. Therefore, the following question is natural.
			\begin{question}
				\label{EndPointSnQuestion2}
				{Given $\alpha_1, \alpha_2 \in \mathbb{R}$ with $\alpha_2 > 0$ and $p = 1 + \alpha_2$,} what are the admissible values of $\gamma_1, \gamma_2 \in \mathbb{R}$ such that the following inequality holds {for} all radial functions in $L^{p, 1}_{\alpha_1, \alpha_2} \left( \mathbb{S}^n \right)$?
				\begin{equation}
					\label{RequiredEndPointSn2}
					\| \sin^{\gamma_1} d \left( 0, \cdot \right) \cos^{\gamma_2} d \left( 0, \cdot \right) R_kf \|_{L^{\infty} \left( \Xi_k \left( \mathbb{S}^n \right) \right)} \leq C \| f \|_{L^{p, 1}_{\alpha_1, \alpha_2} \left( \mathbb{S}^n \right)}.
				\end{equation}
			\end{question}
			First, we notice that Question \ref{EndPointSnQuestion1} arises only when $\alpha_1 \geq k - n$. Also, the case when $\alpha_1 = k - n$, and hence $p = 1$, is already settled in Theorem \ref{LPLInfinitySufficient}. Therefore, to answer Question \ref{EndPointSnQuestion1}, we only consider $\alpha_1 > k - n$, i.e., $p > 1$. To answer our question, we require an analogue of Lemma \ref{ConsequenceMainLemma} suitable for the sphere $\mathbb{S}^n$.
			\begin{lemma}
				\label{MainLemmaSphere}
				Let $\eta_1, \eta_2 > 0$ and $p \geq 1$. Then, for characteristic functions $f$ of measurable sets in $\left( 0, \frac{\pi}{2} \right)$, we have
				\begin{equation}
					\label{MainLemmaEquationSphere}
					\int\limits_{0}^{\frac{\pi}{2}} f \left( t \right) \sin^{\eta_1 - 1} t \cos^{\eta_2 - 1} t \ \mathrm{d}t \leq C_{p, \eta_1, \eta_2} \left( \int\limits_{0}^{\frac{\pi}{2}} f \left( t \right) \sin^{p \eta_1 - 1} t \cos^{p \eta_2 - 1} t \ \mathrm{d}t \right)^{\frac{1}{p}},
				\end{equation}
				where,
				\begin{equation}
					\label{ConstantLemmaSphere}
					\begin{aligned}
						C_{p, \eta_1, \eta_2} = p^{\frac{1}{p}} &\left[ {\eta_1^{-\frac{1}{p'}} {\max \left\lbrace \frac{2^{\frac{3}{2} \left( \eta_1 - 1 \right)}}{\pi^{\eta_1 - 1}}, 1 \right\rbrace \max \left\lbrace 2^{- \frac{\eta_2 - 1}{2}}, 1 \right\rbrace} \max \left\lbrace \frac{\pi^{p \eta_1 - 1}}{2^{\frac{3}{2} \left( p \eta_1 - 1 \right)}}, 1 \right\rbrace \max \left\lbrace 2^{\frac{p \eta_2 - 1}{2}}, 1 \right\rbrace} \right. \\
						&\left.{ + \eta_2^{-\frac{1}{p'}} {\max \left\lbrace 2^{- \frac{\eta_1 - 1}{2}}, 1 \right\rbrace \max \left\lbrace \frac{2^{\frac{3}{2} \left( \eta_2 - 1 \right)}}{\pi^{\eta_2 - 1}}, 1 \right\rbrace} \max \left\lbrace 2^{\frac{p \eta_1 - 1}{2}}, 1 \right\rbrace \max \left\lbrace \frac{\pi^{p \eta_2 - 1}}{2^{\frac{3}{2} \left( p \eta_2 - 1 \right)}}, 1 \right\rbrace} \right].
					\end{aligned}
				\end{equation}
			\end{lemma}
			\begin{proof}
				We first observe some simple properties of sine and cosine functions. We notice that in the interval $\left( 0, \frac{\pi}{4} \right)$, we have
				\begin{equation}
					\label{SineCosNearZero}
					\frac{2 \sqrt{2}}{\pi} t \leq \sin t \leq t, \text{ and } \frac{1}{\sqrt{2}} \leq \cos t \leq 1.
				\end{equation}
				{On the other hand, when} ${t \in} \left( \frac{\pi}{4}, \frac{\pi}{2} \right)$, we have
				\begin{equation}
					\label{SineCosNearEquator}
					\frac{2 \sqrt{2}}{\pi} \left( \frac{\pi}{2} - t \right) \leq \cos t \leq \frac{\pi}{2} - t, \text{ and } \frac{1}{\sqrt{2}} \leq \sin t \leq 1.
				\end{equation}
				Therefore, from Equations \eqref{SineCosNearZero} and \eqref{SineCosNearEquator}, {we get}
				\begin{align*}
					\int\limits_{0}^{\frac{\pi}{2}} f \left( t \right) \sin^{\eta_1 - 1} t \cos^{\eta_2 - 1} t \ \mathrm{d}t &= \int\limits_{0}^{\frac{\pi}{4}} f \left( t \right) \sin^{\eta_1 - 1} t \cos^{\eta_2 - 1} t \ \mathrm{d}t + \int\limits_{\frac{\pi}{4}}^{\frac{\pi}{2}} f \left( t \right) \sin^{\eta_1 - 1} t \cos^{\eta_2 - 1} t \ \mathrm{d}t \\
					&\leq {\max \left\lbrace \frac{2^{\frac{3}{2} \left( \eta_1 - 1 \right)}}{\pi^{\eta_1 - 1}}, 1 \right\rbrace \max \left\lbrace 2^{- \frac{\eta_2 - 1}{2}}, 1 \right\rbrace} \int\limits_{0}^{\frac{\pi}{4}} f \left( t \right) t^{\eta_1 - 1} \mathrm{d}t \\
					&+ {\max \left\lbrace 2^{- \frac{\eta_1 - 1}{2}}, 1 \right\rbrace \max \left\lbrace \frac{2^{\frac{3}{2} \left( \eta_2 - 1 \right)}}{\pi^{\eta_2 - 1}}, 1 \right\rbrace} \int\limits_{\frac{\pi}{4}}^{\frac{\pi}{2}} f \left( t \right) \left( \frac{\pi}{2} - t \right)^{\eta_2 - 1} \mathrm{d}t.
				\end{align*}
				Changing the variables $\frac{\pi}{2} - t = s$ in the second integral, we have,
				\begin{align*}
					\int\limits_{0}^{\frac{\pi}{2}} f \left( t \right) \sin^{\eta_1 - 1} t \cos^{\eta_2 - 1} t \ \mathrm{d}t &\leq {\max \left\lbrace \frac{2^{\frac{3}{2} \left( \eta_1 - 1 \right)}}{\pi^{\eta_1 - 1}}, 1 \right\rbrace \max \left\lbrace 2^{- \frac{\eta_2 - 1}{2}}, 1 \right\rbrace} \int\limits_{0}^{\infty} f \left( t \right) \chi_{\left( 0, \frac{\pi}{4} \right)} \left( t \right) t^{\eta_1 - 1} \mathrm{d}t \\
					&+ {\max \left\lbrace 2^{- \frac{\eta_1 - 1}{2}}, 1 \right\rbrace \max \left\lbrace \frac{2^{\frac{3}{2} \left( \eta_2 - 1 \right)}}{\pi^{\eta_2 - 1}}, 1 \right\rbrace} \int\limits_{0}^{\infty} f \left( \frac{\pi}{2} - s \right) \chi_{\left( 0, \frac{\pi}{4} \right)} \left( s \right) s^{\eta_2 - 1} \mathrm{d}s.
				\end{align*}
				Now, we use Inequality \eqref{PolynomialInequalityLemma} on both integrals. Then, we have,
				\begin{align*}
					\int\limits_{0}^{\frac{\pi}{2}} f \left( t \right) \sin^{\eta_1 - 1} t \cos^{\eta_2 - 1} t \ \mathrm{d}t \leq p^{\frac{1}{p}} &\left[ \eta_1^{-\frac{1}{p'}} {\max \left\lbrace \frac{2^{\frac{3}{2} \left( \eta_1 - 1 \right)}}{\pi^{\eta_1 - 1}}, 1 \right\rbrace \max \left\lbrace 2^{- \frac{\eta_2 - 1}{2}}, 1 \right\rbrace} \left( \int\limits_{0}^{\infty} f \left( t \right) \chi_{\left( 0, \frac{\pi}{4} \right)} \left( s \right) t^{p\eta_1 - 1} \mathrm{d}t \right)^{\frac{1}{p}} \right. \\
					&\left. + \eta_2^{-\frac{1}{p'}} {\max \left\lbrace 2^{- \frac{\eta_1 - 1}{2}}, 1 \right\rbrace \max \left\lbrace \frac{2^{\frac{3}{2} \left( \eta_2 - 1 \right)}}{\pi^{\eta_2 - 1}}, 1 \right\rbrace} \left( \int\limits_{0}^{\infty} f \left( \frac{\pi}{2} - s \right) \chi_{\left( 0, \frac{\pi}{4} \right)} \left( s \right) s^{p \eta_2 - 1} \mathrm{d}s \right)^{\frac{1}{p}} \right].
				\end{align*}
				Again changing the variables $\frac{\pi}{2} - s = t$ in the second integral in the right hand side of the above inequality, we get
				\begin{align*}
					\int\limits_{0}^{\frac{\pi}{2}} f \left( t \right) \sin^{\eta_1 - 1} t \cos^{\eta_2 - 1} t \ \mathrm{d}t \leq p^{\frac{1}{p}} &\left[ \eta_1^{-\frac{1}{p'}} {\max \left\lbrace \frac{2^{\frac{3}{2} \left( \eta_1 - 1 \right)}}{\pi^{\eta_1 - 1}}, 1 \right\rbrace \max \left\lbrace 2^{- \frac{\eta_2 - 1}{2}}, 1 \right\rbrace} \left( \int\limits_{0}^{\frac{\pi}{4}} f \left( t \right) t^{p\eta_1 - 1} \mathrm{d}t \right)^{\frac{1}{p}} \right. \\
					&\left. + \eta_2^{-\frac{1}{p'}} {\max \left\lbrace 2^{- \frac{\eta_1 - 1}{2}}, 1 \right\rbrace \max \left\lbrace \frac{2^{\frac{3}{2} \left( \eta_2 - 1 \right)}}{\pi^{\eta_2 - 1}}, 1 \right\rbrace} \left( \int\limits_{\frac{\pi}{4}}^{\frac{\pi}{2}} f \left( t \right) \left( \frac{\pi}{2} - t \right)^{p \eta_2 - 1} \mathrm{d}{t} \right)^{\frac{1}{p}} \right].
				\end{align*}
				{Using} Equations \eqref{SineCosNearZero} and \eqref{SineCosNearEquator}, we get
				\begin{align*}
					&\int\limits_{0}^{\frac{\pi}{2}} f \left( t \right) \sin^{\eta_1 - 1} t \cos^{\eta_2 - 1} t \ \mathrm{d}t \leq C_{p, \eta_1, \eta_2} \left( \int\limits_{0}^{\frac{\pi}{2}} f \left( t \right) \sin^{p \eta_1 - 1} t \cos^{p \eta_2 - 1} t \ {\mathrm{d}t} \right)^{\frac{1}{p}}.
				\end{align*}
				Here, $C_{p, \eta_1, \eta_2}$ is as given in Equation \eqref{ConstantLemmaSphere}. This completes the proof!
			\end{proof}
			\begin{remark}
				\normalfont
				With a simple change of variables (replacing $\cos t$ by $t$) in Equation \eqref{MainLemmaEquationSphere}, Lemma \ref{MainLemmaSphere} can be written in the following form for characteristic functions $f$ of measurable sets in $\left( 0, 1 \right)$.
				\begin{equation}
					\label{MainLemmaSphereEquation2}
					\int\limits_{0}^{1} f \left( t \right) \left( 1 - t^2 \right)^{\frac{\eta_1}{2} - 1} t^{\eta_2 - 1} \mathrm{d}t \leq C \left( \int\limits_{0}^{1} f \left( t \right) \left( 1 - t^2 \right)^{\frac{p \eta_1}{2} - 1} t^{p \eta_2 - 1} \mathrm{d}t \right)^{\frac{1}{p}},
				\end{equation}
				for $\eta_1, \eta_2 > 0$ and $p \geq 1$. This form of Lemma \ref{MainLemmaSphere} {is} of interest to us in the analysis of $X$-Ray transform.
			\end{remark}
			We are now in a position to answer Question \ref{EndPointSnQuestion1}.
			\begin{theorem}
				\label{EndPointResult1}
				Let {$\alpha_1, \alpha_2 \in \mathbb{R}$ with} $\alpha_1 \geq k - n$, $p = \frac{\alpha_1 + n}{k} {\geq 1 + \alpha_2}$, {with a strict inequality when $\alpha_2 > 0$}.
				\begin{enumerate}
					\item[\mylabel{EndPointLPLInfinitySnWeightedA}{(A)}] For $k \geq 2$, the {Inequality \eqref{RequiredEndPointSn1}} holds for all radial functions $f \in L^{p, 1}_{\alpha_1, \alpha_2} \left( \mathbb{S}^n \right)$ {if and only if $\gamma \geq \frac{1 + \alpha_2}{p}$}.
					\item[\mylabel{EndPointLPLInfinitySnWeightedB}{(B)}] For $k = 1$, Inequality \eqref{RequiredEndPointSn1} holds for all radial functions $f \in L^{p, 1}_{\alpha_1, \alpha_2} \left( \mathbb{S}^n \right)$ {if and only if $\gamma \geq \frac{1 + \alpha_2}{p}$ and} $p \geq 2$.
				\end{enumerate}
			\end{theorem}
			\begin{proof}~
				\begin{enumerate}
					\item[\mylabel{EndPointLPLInfinitySnWeightedAProof}{(A)}] We first prove the result for $k \geq 2$. For characteristic functions of measurable radial sets in $\mathbb{S}^n$, from Equation \eqref{KPlaneTransformRadialSnEquation2}, {we have}
					\begin{align*}
						\left| \cos^{\gamma} d \left( 0, \xi \right) R_kf \left( \xi \right) \right| &\leq C \cos^{\gamma - 1} d \left( 0, \xi \right) \int\limits_{d \left( 0, \xi \right)}^{\frac{\pi}{2}} \left| \tilde{f} \left( \cos t \right) \right| \sin^{k - 1} t \left( 1 - \frac{\tanh^2 d \left( 0, \xi \right)}{\tan^2 t} \right)^{\frac{k}{2} - 1} \mathrm{d}t \\
						&\leq C \cos^{\gamma - 1} d \left( 0, \xi \right) \int\limits_{d \left( 0, \xi \right)}^{\frac{\pi}{2}} \left| \tilde{f} \left( \cos t \right) \right| \sin^{k - 1} t \ \mathrm{d}t.
					\end{align*}
					Using {Inequality \eqref{MainLemmaEquationSphere}} Lemma \ref{MainLemmaSphere} with $\eta_1 = k$ and $\eta_2 = 1$, we get
					\begin{align*}
						\left| \cos^{\gamma} d \left( 0, \xi \right) R_kf \left( \xi \right) \right| &\leq C \cos^{\gamma - 1} d \left( 0, \xi \right) \left( \int\limits_{d \left( 0, \xi \right)}^{\frac{\pi}{2}} \left| \tilde{f} \left( \cos t \right) \right| \sin^{\alpha_1 + n - 1} t \cos^{p - 1} t \ \mathrm{d}t \right)^{\frac{1}{p}} \\
						&= C \cos^{\gamma - 1} d \left( 0, \xi \right) \left( \int\limits_{d \left( 0, \xi \right)}^{\frac{\pi}{2}} \left| \tilde{f} \left( \cos t \right) \right| \sin^{\alpha_1 + n - 1} t \cos^{\alpha_2} t \cos^{p - 1 - \alpha_2} t \ \mathrm{d}t \right)^{\frac{1}{p}}.
					\end{align*}
					We know from the existence conditions of Theorem \ref{ExistenceSN}, that $p \geq 1 + \alpha_2$. Thus, we get
					\begin{align*}
						\left| \cos^{\gamma} d \left( 0, \xi \right) R_kf \left( \xi \right) \right| \leq C \cos^{\gamma - \frac{1 + \alpha_2}{p}} d \left( 0, \xi \right) \| f \|_{L^{p, 1}_{\alpha_1, \alpha_2} \left( \mathbb{S}^n \right)} {\leq C \| f \|_{L^{p, 1}_{\alpha_1, \alpha_2} \left( \mathbb{S}^n \right)}},
					\end{align*}
					{since} $\gamma \geq \frac{1 + \alpha_2}{p}$.
					\item[\mylabel{EndPointLPLInfinitySnWeightedBProof}{(B)}] For the case $k = 1$, {the proof is} a bit more involved since $\left( 1 - \frac{\tan^2 d \left( 0, \xi \right)}{\tan^2 t} \right)^{- \frac{1}{2}}$ is not bounded. Nonetheless, we can use Lemma \ref{MainLemmaSphere} in the form of Inequality \eqref{MainLemmaSphereEquation2} to obtain the end-point result. We have, for characteristic functions $f$ of measurable radially symmetric sets in $\mathbb{S}^n$, by Equation \eqref{KPlaneTransformRadialSnEquation},
					\begin{align*}
						R_1f \left( \xi \right) &= C \int\limits_{d \left( 0, \xi \right)}^{\frac{\pi}{2}} \tilde{f} \left( t \right) \left( \cos^2 d \left( 0, \xi \right) - \cos^2 t \right)^{- \frac{1}{2}} \sin t \ \mathrm{d}t.
					\end{align*}
					By substituting $\cos t = x \cos d \left( 0, \xi \right)$, we get
					\begin{align}
						\label{XRayNewFormula}
						R_1 f \left( \xi \right) &= C \int\limits_{0}^{1} \tilde{f} \left( x \cos d \left( 0, \xi \right) \right) \left( \cos^2 d \left( 0, \xi \right) - x^{{2}} \cos^2 d \left( 0, \xi \right) \right)^{- \frac{1}{2}} \cos d \left( 0, \xi \right) \mathrm{d}x = C \int\limits_{0}^{1} \tilde{f} \left( x \cos d \left( 0, \xi \right) \right) \left( 1 - x^2 \right)^{- \frac{1}{2}} \mathrm{d}x.
					\end{align}
					Therefore, we have, by using Inequality \eqref{MainLemmaSphereEquation2} with $\eta_1 = 1 = \eta_2$, {and the fact that $p = \alpha_1 + n$},
					\begin{align*}
						\left| \cos^{\gamma} d \left( 0, \xi \right) R_1f \left( \xi \right) \right| &= C \cos^{\gamma} d \left( 0, \xi \right) \left( \int\limits_{0}^{1} \tilde{f} \left( x \cos d \left( 0, \xi \right) \right) \left( 1 - x^2 \right)^{\frac{\alpha_1 + n}{2} - 1} x^{p - 1} \ \mathrm{d}x \right)^{\frac{1}{p}} \\
						&\leq C \cos^{\gamma} d \left( 0, \xi \right) \left( \int\limits_{0}^{1} \tilde{f} \left( x \cos d \left( 0, \xi \right) \right) \left( 1 - x^2 \right)^{\frac{\alpha_1 + n}{2} - 1} x^{\alpha_2} \ \mathrm{d}x \right)^{\frac{1}{p}},
					\end{align*}
					since $p \geq 1 + \alpha_2$ and $x \leq 1$. Now, we substitute $x \cos d \left( 0, \xi \right) = u$. Then, we have,
					\begin{align*}
						\left| \cos^{\gamma} d \left( 0, \xi \right) R_1f \left( \xi \right) \right| &\leq C \cos^{\gamma} d \left( 0, \xi \right) \left( \int\limits_{0}^{\cos d \left( 0, \xi \right)} \tilde{f} \left( u \right) \left( 1 - \frac{u^2}{\cos^2 d \left( 0, \xi \right)} \right)^{\frac{\alpha_1 + n}{2} - 1} \left( \frac{u}{\cos d \left( 0, \xi \right)} \right)^{\alpha_2} \frac{\mathrm{d}u}{\cos d \left( 0, \xi \right)} \right)^{\frac{1}{p}}.
					\end{align*}
					From our assumptions, we have $p = \alpha_1 + n \geq 2$. Hence,
					$$\left( 1 - \frac{u^2}{\cos^2 d \left( 0, \xi \right)} \right)^{\frac{\alpha_1 + n}{2} - 1} \leq \left( 1 - u^2 \right)^{\frac{\alpha_1 + n}{2} - 1}.$$
					Then,
					\begin{align*}
						\left| \cos^{\gamma} d \left( 0, \xi \right) R_1f \left( \xi \right) \right| &\leq C \cos^{\gamma - \frac{1 + \alpha_2}{p}} d \left( 0, \xi \right) \left( \int\limits_{0}^{\cos d \left( 0, \xi \right)} \tilde{f} \left( u \right) \left( 1 - u^2 \right)^{\frac{\alpha_1 + n}{2} - 1} u^{\alpha_2} \mathrm{d}u \right)^{\frac{1}{p}}.
					\end{align*}
					Again, from the necessary conditions, we know that $\gamma \geq \frac{1 + \alpha_2}{p}$. Therefore, $\cos^{\gamma - \frac{1 + \alpha_2}{p}} d \left( 0, \xi \right) \leq 1$. Also, by substituting $u = \cos t$, we get
					\begin{align*}
						\left| \cos^{\gamma} d \left( 0, \xi \right) R_1f \left( \xi \right) \right| &\leq C \left( \int\limits_{0}^{1} \tilde{f} \left( \cos t \right) \sin^{\alpha_1 + n - 1} t \cos^{\alpha_2} t \ \mathrm{d}t \right)^{\frac{1}{p}} = C \| f \|_{L^{p, 1}_{\alpha_1, \alpha_2} \left( \mathbb{S}^n \right)}.
					\end{align*}
					{The necessity of the condition $p \geq 2$ when $k = 1$ is shown in Example \ref{XRayPAtLeast2}.}
				\end{enumerate}
				{The necessity of the condition $\gamma \geq \frac{1 + \alpha_2}{p}$ is shown in Remark \ref{RemarkGammaNecessryEndPoint1Sn}.}
			\end{proof}
			\begin{remark}
				\label{RemarkGammaNecessryEndPoint1Sn}
				\normalfont
				{We now show that the condition on $\gamma$ mentioned in Theorem \ref{EndPointResult1} is necessary. To this end, we consider the function $f_{a, b}: \mathbb{S}^n_+ \rightarrow \mathbb{C}$ defined as $f_{a, b} \left( x \right) = \chi_{\left( a, b \right)} \left( \cos d \left( 0, x \right) \right)$, where $0 < a < b < 1$. We are interested particularly when $b = 2a$, and we denote $f_a = f_{a, 2a}$. This forces $a < \frac{1}{2}$. In fact, in what follows, we consider $a < \frac{1}{4}$ for certain computational convenience. Next, by substituting $x = at$ in Equation \eqref{WeightedLorentzNormAnnulusSphere}, and using the fact that $b = 2a$, we get
				\begin{align*}
					\| f_a \|_{L^{p, 1}_{\alpha_1, \alpha_2} \left( \mathbb{S}^n \right)} &= C a^{\frac{\alpha_2 + 1}{p}} \left[ \int\limits_{1}^{2} \left( 1 - a^2t^2 \right)^{\frac{\alpha_1 + n}{2} - 1} t^{\alpha_2} \mathrm{d}t \right]^{\frac{1}{2}}.
				\end{align*}
				Now, we observe that since $a < \frac{1}{4}$ and $1 \leq t \leq 2$, we have $\frac{3}{4} \leq 1 - a^2t^2 \leq 1$. Therefore, we have,
				\begin{equation}
					\label{LorentzNormEstimateDoubleAnnulusSn}
					\| f \|_{L^p_{\alpha_1, \alpha_2} \left( \mathbb{S}^n \right)} \leq C a^{\frac{1 + \alpha_2}{p}}.
				\end{equation}
				On the other hand, to estimate the $k$-plane transform of $f_a$, we use Equation \eqref{KPlaneTransformRadialSnEquation}. Particularly, by substituting $\cos t = x \cos d \left( 0, \xi \right)$ in Equation \eqref{KPlaneTransformRadialSnEquation} and simplifying, we get
				\begin{align*}
					R_kf_a \left( \xi \right) &= C \int\limits_{0}^{1} \chi_{\left( a, 2a \right)} \left( x \cos d \left( 0, \xi \right) \right) \left( 1 - x^2 \right)^{\frac{k}{2} - 1} \mathrm{d}x \\
					&= C \begin{cases}
								0, & \cos d \left( 0, \xi \right) < a. \\
								\int\limits_{\frac{a}{\cos d \left( 0, \xi \right)}}^{1} \left( 1 - x^2 \right)^{\frac{k}{2} - 1} \mathrm{d}x, & a \leq \cos d \left( 0, \xi \right) < 2a. \\
								\int\limits_{\frac{a}{\cos d \left( 0, \xi \right)}}^{\frac{2a}{\cos d \left( 0, \xi \right)}} \left( 1 - x^2 \right)^{\frac{k}{2} - 1} \mathrm{d}x, & \cos d \left( 0, \xi \right) \geq 2a.
							\end{cases}
				\end{align*}
				By choosing $\xi \in \Xi_k \left( \mathbb{S}^n \right)$ such that $\cos d \left( 0, \xi \right) = 2a$, it is clear at once that
				\begin{equation}
					\label{RkfEstimateDountableAnnulusSn}
					\| \cos^{\gamma} d \left( 0, \cdot \right) R_kf \|_{L^{\infty} \left( \Xi_k \left( \mathbb{S}^n \right) \right)} \geq C a^{\gamma} \int\limits_{\frac{1}{2}}^{1} \left( 1 - x^2 \right)^{\frac{k}{2} - 1} \mathrm{d}x = C a^{\gamma}.
				\end{equation}
				We now use Inequalities \eqref{LorentzNormEstimateDoubleAnnulusSn} and \eqref{RkfEstimateDountableAnnulusSn} in Inequality \eqref{RequiredEndPointSn1} to obtain
				$$a^{\gamma} \leq C a^{\frac{1 + \alpha_2}{p}}.$$
				As $a \rightarrow 0$, the above inequality forces $\gamma \geq \frac{1 + \alpha_2}{p}$.}
			\end{remark}
			Next, we move on to answer Question \ref{EndPointSnQuestion2}. {As a first step, we get the necessary conditions for Inequality \eqref{RequiredEndPointSn2} to hold with $p = 1 + \alpha_2$ and $\alpha_2 > 0$. We do so by testing the inequality against the characteristic functions of certain ``nice" radially symmetric sets.} 
			\begin{theorem}
				\label{EndPointSn2NecessityTheorem}
				{Let $\alpha_2 > 0$, $\alpha_1 \in \mathbb{R}$ and $p = 1 + \alpha_2$. Then, Inequality \eqref{RequiredEndPointSn2} holds only if 
				\begin{equation}
					\label{EndPointSn2NecessaryConditions}
					\gamma_1 \geq \max \left\lbrace 0, \frac{\alpha_1 + n}{p} - k \right\rbrace \text{ and } \gamma_2 \geq 1.
				\end{equation}}
			\end{theorem}
			\begin{proof}
				{We prove this result in two cases.}
				
				{\textbf{Case I}: First, we consider the functions centered at the origin of the sphere. Here, we further divide our analysis in two parts, depending on the value of $\alpha_1$.}
				
				{\textbf{Case I(A)}: First, let us consider $\alpha_1 > - n$. Then, for $\lambda > 0$, we consider $f_{\lambda}: \mathbb{S}^n_+ \rightarrow \mathbb{C}$, defined as $f_{\lambda} \left( x \right) = \chi_{B \left( 0, \lambda \right)} \left( x \right)$, where $B \left( 0, \lambda \right)$ is the ball of radius $\lambda$ in $\mathbb{S}^n_+$, centered at the origin. Then, we have,
				$$\| f_{\lambda} \|_{L^{p, 1}_{\alpha_1, \alpha_2} \left( \mathbb{S}^n \right)} = C \left[ \int\limits_{0}^{\lambda} \sin^{\alpha_1 + n - 1} t \cos^{\alpha_2} t \ \mathrm{d}t \right]^{\frac{1}{p}}.$$
				By substituting $\sin^2 t = x \sin^2 \lambda$ and simplifying, we get
				$$\| f_{\lambda} \|_{L^{p, 1}_{\alpha_1, \alpha_2} \left( \mathbb{S}^n \right)} = C \sin^{\frac{\alpha_1 + n}{p}} \lambda \left[ \int\limits_{0}^{1} x^{\frac{\alpha_1 + n}{2} - 1} \left( 1 - x \sin^2 \lambda \right)^{\frac{\alpha_2 - 1}{2}} \mathrm{d}x \right]^{\frac{1}{p}}.$$
				Now, by using Equation \eqref{IntegralForm2F1}, we have,
				\begin{equation}
					\label{LorentzNormBallSn}
					\| f_{\lambda} \|_{L^{p, 1}_{\alpha_1, \alpha_2} \left( \mathbb{S}^n \right)} = C \sin^{\frac{\alpha_1 + n}{p}} \lambda \left[ {}_2F_1 \left( \frac{1 - \alpha_2}{2}, \frac{\alpha_1 + n}{2}; 1 + \frac{\alpha_1 + n}{2}; \sin^2 \lambda \right) \right]^{\frac{1}{p}}.
				\end{equation}
				The $k$-plane transform of $f_{\lambda}$ is evaluated in Equation \eqref{KPlaneBallSnWithoutMassEquation} of Theorem \ref{KPlaneBallSnWithoutMass}. Therefore, we have,
				\begin{align*}
					&\sin^{\gamma_1} d \left( 0, \xi \right) \cos^{\gamma_2} d \left( 0, \xi \right) R_kf \left( \xi \right) \\
					&= C \cos^{\gamma_2 - k - 1} d \left( 0, \xi \right) \sin^{\gamma_1} d \left( 0, \xi \right) \sin^k \lambda \left( 1 - \frac{\sin^2 d \left( 0, \xi \right)}{\sin^2 \lambda} \right)^{\frac{k}{2}} {}_2F_1 \left( \frac{1}{2}, \frac{k}{2}; 1 + \frac{k}{2}; 1 - \frac{\cos^2 \lambda}{\cos^2 d \left( 0, \xi \right)} \right) \chi_{\left( 0, \lambda \right)} \left( d \left( 0, \xi \right) \right).
				\end{align*}
				First, we observe that if $\gamma_1 < 0$, then the above quantity is unbounded as $d \left( 0, \xi \right) \rightarrow 0$. Therefore, a necessary requirement for Inequality \eqref{RequiredEndPointSn2} is that $\gamma_1 \geq 0$. Further, by choosing $\xi \in \Xi_k \left( \mathbb{S}^n \right)$ such that $\sin^2 d \left( 0, \xi \right) = \frac{1}{2} \sin^2 \lambda$, it is clear that
				\begin{equation}
					\label{RkfLambdaLInfinityEstimateEndPointSn2}
					\| \sin^{\gamma_1} d \left( 0, \cdot \right) \cos^{\gamma_2} d \left( 0, \cdot \right) R_kf_{\lambda} \|_{L^{\infty} \left( \Xi_k \left( \mathbb{S}^n \right) \right)} \geq C \left( 1 + \cos^2 \lambda \right)^{\frac{\gamma_2 - k - 1}{2}} \sin^{\gamma_1 + k} \lambda \ {}_2F_1 \left( \frac{1}{2}, \frac{k}{2}; 1 + \frac{k}{2}; \frac{\sin^2 \lambda}{1 + \cos^2 \lambda} \right).
				\end{equation}
				Using Inequality \eqref{RkfLambdaLInfinityEstimateEndPointSn2} and Equation \eqref{LorentzNormBallSn} in Inequality \eqref{RequiredEndPointSn2}, we get
				$$\left( 1 + \cos^2 \lambda \right)^{\frac{\gamma_2 - k - 1}{2}} \sin^{\gamma_1 + k} \lambda \ {}_2F_1 \left( \frac{1}{2}, \frac{k}{2}; 1 + \frac{k}{2}; \frac{\sin^2 \lambda}{1 + \cos^2 \lambda} \right) \leq C \sin^{\frac{\alpha_1 + n}{p}} \lambda \left[ {}_2F_1 \left( \frac{1 - \alpha_2}{2}, \frac{\alpha_1 + n}{2}; 1 + \frac{\alpha_1 + n}{2}; \sin^2 \lambda \right) \right]^{\frac{1}{p}}.$$
				As $\lambda \rightarrow 0$, the hypergeometric terms in the above inequality approach $1$. Similarly, the cosine terms approach $1$. Therefore, as $\lambda \rightarrow 0$, the above inequality forces $\gamma_1 \geq \frac{\alpha_1 + n}{p} - k$.}
				
				{\textbf{Case I(B)}: Now, let us consider $\alpha_1 \leq -n$. In this case, we cannot use the characteristic function of a ball centered at origin, since the integral that leads to Equation \eqref{LorentzNormBallSn} is not convergent. In this case, therefore, let us consider the set $E_{\lambda} := B \left( 0, 2 \lambda \right) \setminus B \left( 0, \lambda \right)$, which is an annulus on $\mathbb{S}^n_+$ of inner radius $\lambda$ and outer radius $2 \lambda$. For the sake of computational convenience, we consider $\lambda < \frac{\pi}{8}$, and define the function $f_{\lambda}: \mathbb{S}^n_+ \rightarrow \mathbb{C}$ as $f_{\lambda} \left( x \right) = \chi_{E_{\lambda}} \left( x \right)$. We begin by estimating the Lorentz norm of $f_{\lambda}$. We have,
				$$\| f_{\lambda} \|_{L^{p, 1}_{\alpha_1, \alpha_2} \left( \mathbb{S}^n \right)} = C \left[ \int\limits_{\lambda}^{2 \lambda} \sin^{\alpha_1 + n - 1} t \cos^{\alpha_2} t \ \mathrm{d}t \right]^{\frac{1}{p}}.$$
				By substituting $\sin^2 t = u$ and simplifying, we obtain
				$$\| f_{\lambda} \|_{L^{p, 1}_{\alpha_1, \alpha_2} \left( \mathbb{S}^n \right)} = C \left[ \int\limits_{\sin^2 \lambda}^{\sin^2 2 \lambda} u^{\frac{\alpha_1 + n}{2} - 1} \left( 1 - u \right)^{\frac{\alpha_2 - 1}{2}} \mathrm{d}u \right]^{\frac{1}{p}}.$$
				Since $\alpha_1 \leq -n$, we have $u^{\frac{\alpha_1 + n}{2} - 1} \leq \sin^{\alpha_1 + n - 2} \lambda$, whenever $u \geq \sin^2 \lambda$. Also, since $\lambda < \frac{\pi}{8}$, we have $\frac{1}{2} \leq 1 - u \leq 1$, for $\sin^2 \lambda \leq u \leq \sin^2 2 \lambda$. Consequently, we have,
				\begin{equation}
					\label{LorentzNormEstimateDoubleAnnulusSnNearZero}
					\| f_{\lambda} \|_{L^{p, 1}_{\alpha_1, \alpha_2} \left( \mathbb{S}^n \right)} \leq C \sin^{\frac{\alpha_1 + n - 2}{p}} \lambda \left( \sin^2 2 \lambda - \sin^2 \lambda \right)^{\frac{1}{p}} = C \sin^{\frac{\alpha_1 + n}{p}} \lambda \left( 4 \cos^2 \lambda - 1 \right)^{\frac{1}{p}} \leq C \sin^{\frac{\alpha_1 + n}{p}} \lambda.
				\end{equation}
				Next, we estimate the $k$-plane transform of $f_{\lambda}$. By substituting $\cos t = x \cos d \left( 0, \xi \right)$ in Equation \eqref{KPlaneTransformRadialSnEquation}, we get
				\begin{align*}
					R_kf_{\lambda} \left( \xi \right) = C \int\limits_{0}^{1} \chi_{\left( \cos 2 \lambda, \cos \lambda \right)} \left( x \cos d \left( 0, \xi \right) \right) \left( 1 - x^2 \right)^{\frac{k}{2} - 1} \mathrm{d}x = C \begin{cases}
														0, & d \left( 0, \xi \right) \geq 2 \lambda. \\
														\int\limits_{\frac{\cos 2 \lambda}{\cos d \left( 0, \xi \right)}}^{1} \left( 1 - x^2 \right)^{\frac{k}{2} - 1} \mathrm{d}x, & \lambda < d \left( 0, \xi \right) < 2 \lambda. \\
														\int\limits_{\frac{\cos 2 \lambda}{\cos d \left( 0, \xi \right)}}^{\frac{\cos \lambda}{\cos d \left( 0, \xi \right)}} \left( 1 - x^2 \right)^{\frac{k}{2} - 1} \mathrm{d}x, & d \left( 0, \xi \right) \leq \lambda.
													\end{cases}
				\end{align*}
				It is again clear that for the left-hand-side of Inequality \eqref{RequiredEndPointSn2} to be finite, we must have $\gamma_1 \geq 0$. Now, for $\xi \in \Xi_k \left( \mathbb{S}^n \right)$ such that $d \left( 0, \xi \right) = \lambda$, we have,
				$$R_kf_{\lambda} \left( \xi \right) = \int\limits_{\frac{\cos 2 \lambda}{\cos \lambda}}^{1} \left( 1 - x^2 \right)^{\frac{k}{2} - 1} \mathrm{d}x \geq \int\limits_{\frac{\cos 2 \lambda}{\cos \lambda}}^{1} \left( 1 - x \right)^{\frac{k}{2} - 1} \mathrm{d}x = C \left( 1 - \frac{\cos 2 \lambda}{\cos \lambda} \right)^{\frac{k}{2}} = C \frac{\left( \cos \lambda - \cos^2 \lambda + \sin^2 \lambda \right)^{\frac{k}{2}}}{\cos^{\frac{k}{2}} \lambda} \geq C \frac{\sin^k \lambda}{\cos^{\frac{k}{2}} \lambda}.$$
				That is, we have,
				\begin{equation}
					\label{RkfEstimateDoubleAnnulusSnNearZero}
					\| \sin^{\gamma_1} d \left( 0, \cdot \right) \cos^{\gamma_2} d \left( 0, \cdot \right) R_kf_{\lambda} \|_{L^{\infty} \left( \Xi_k \left( \mathbb{S}^n \right) \right)} \geq C \sin^{\gamma_1 + k} \lambda \cos^{\gamma_2 - \frac{k}{2}} \lambda.
				\end{equation}
				Using Inequalities \eqref{LorentzNormEstimateDoubleAnnulusSnNearZero} and \eqref{RkfEstimateDoubleAnnulusSnNearZero} in Inequality \eqref{RequiredEndPointSn2}, we get, for $\lambda < \frac{\pi}{8}$,
				$$\sin^{\gamma_1 + k} \lambda \cos^{\gamma_2 - \frac{k}{2}} \lambda \leq C \sin^{\frac{\alpha_1 + n}{p}} \lambda.$$
				As $\lambda \rightarrow 0$, this forces $\gamma_1 \geq \frac{\alpha_1 + n}{p} - k$.}
				
				{Thus, by considering functions that are concentrated near the origin of the sphere, we get a necessary condition on $\gamma_1$. To get the necessary condition on $\gamma_2$, we now consider functions concentrated near the equator of $\mathbb{S}^n$.}
				
				{\textbf{Case II}: For $0 < \lambda < \frac{\pi}{2}$, let us consider the set $E_{\lambda} = \mathbb{S}^n_+ \setminus B \left( 0, \lambda \right)$ and the function $f_{\lambda}: \mathbb{S}^n_+ \rightarrow \mathbb{C}$ defined as $f_{\lambda} \left( x \right) = \chi_{E_{\lambda}} \left( x \right)$. Here, $B \left( 0, \lambda \right)$ is the ball centered at origin of $\mathbb{S}^n$ of radius $\lambda$. Then, we have,
				$$\| f_{\lambda} \|_{L^{p, 1}_{\alpha_1, \alpha_2} \left( \mathbb{S}^n \right)} = C \left[ \int\limits_{\lambda}^{\frac{\pi}{2}} \sin^{\alpha_1 + n - 1} t \cos^{\alpha_2} t \ \mathrm{d}t \right]^{\frac{1}{p}}.$$
				By substituting $\cos t = \cos \lambda \sqrt{x}$ and using the fact that $p = 1 + \alpha_2$, we get
				$$\| f_{\lambda} \|_{L^{p, 1}_{\alpha_1, \alpha_2} \left( \mathbb{S}^n \right)} = C \cos \lambda \left[ \int\limits_{0}^{1} x^{\frac{\alpha_2 - 1}{2}} \left( 1 - x \cos^2 \lambda \right)^{\frac{\alpha_1 + n}{2} - 1} \mathrm{d}x \right]^{\frac{1}{p}}.$$
				Finally, using the integral form of the hypergeometric function given in Equation \eqref{IntegralForm2F1}, we conclude that
				\begin{equation}
					\label{LorentzNormEquatorStripSn}
					\| f_{\lambda} \|_{L^{p, 1}_{\alpha_1, \alpha_2} \left( \mathbb{S}^n \right)} = C \cos \lambda \left[ {}_2F_1 \left( 1 - \frac{\alpha_1 + n}{2}, \frac{\alpha_2 + 1}{2}; \frac{\alpha_2 + 3}{2}; \cos^2 \lambda \right) \right]^{\frac{1}{p}}.
				\end{equation}
				The $k$-plane transform of $f_{\lambda}$ was evaluated in Equation \eqref{KPlaneEquatorSnWithoutMassEquation} of Theorem \ref{KPlaneEquatorSnWithoutMass}. By choosing $\xi \in \Xi_k \left( \mathbb{S}^n \right)$ such that $d \left( 0, \xi \right) = \lambda$, it is readily seen that
				\begin{equation}
					\label{RkfEstimateEquatorStripSn}
					\| \sin^{\gamma_1} d \left( 0, \cdot \right) \cos^{\gamma_2} d \left( 0, \cdot \right) R_kf_{\lambda} \|_{L^{\infty} \left( \Xi_k \left( \mathbb{S}^n \right) \right)} \geq C \sin^{\gamma_1} \lambda \cos^{\gamma_2} \lambda.
				\end{equation}
				Further, by using Inequalities \eqref{LorentzNormEquatorStripSn} and \eqref{RkfEstimateEquatorStripSn} in Inequality \eqref{RequiredEndPointSn2}, we get
				$$\sin^{\gamma_1} \lambda \cos^{\gamma_2} \lambda \leq C \cos \lambda \left[ {}_2F_1 \left( 1 - \frac{\alpha_1 + n}{2}, \frac{\alpha_2 + 1}{2}; \frac{\alpha_2 + 3}{2}; \cos^2 \lambda \right) \right]^{\frac{1}{p}}.$$
				As $\lambda \rightarrow \frac{\pi}{2}$, the above inequality forces $\gamma_2 \geq 1$.}
			\end{proof}
			{We now prove the sufficiency of the conditions stated in Theorem \ref{EndPointSn2NecessityTheorem}.} Again, we employ Lemma \ref{MainLemmaSphere} in the proof.
			\begin{theorem}
				\label{EndPointSn2}
				Let {$\alpha_1, \alpha_2 \in \mathbb{R}$ with $\alpha_2 > 0$  and $p = 1 + \alpha_2$.}
				\begin{enumerate}
					\item[\mylabel{EndPoint2SnWeightedA}{(A)}] For $k \geq 2$, the Inequality \eqref{RequiredEndPointSn2} holds for all radial functions $f \in L^{p, 1}_{\alpha_1, \alpha_2} \left( \mathbb{S}^n \right)$ {if and only if $\gamma_1$ and $\gamma_2$ satisfy the conditions of Equation \eqref{EndPointSn2NecessaryConditions}}.
					\item[\mylabel{EndPoint2SnWeightedB}{(B)}] For $k = 1$, Inequality \eqref{RequiredEndPointSn2} holds {for all radial functions $f \in L^{p, 1}_{\alpha_1, \alpha_2} \left( \mathbb{S}^n \right)$ if and only if $\gamma_1$ and $\gamma_2$ satisfy the conditions of Equation \eqref{EndPointSn2NecessaryConditions} and} $p \geq 2$.
				\end{enumerate}
			\end{theorem}
			\begin{proof}~
				\begin{enumerate}
					\item[\mylabel{EndPoint2SnWeightedAProof}{(A)}] Let us begin with $k \geq 2$. As {mentioned earlier}, we wish to use Lemma \ref{MainLemmaSphere}. We have from Equation \eqref{KPlaneTransformRadialSnEquation2},
					\begin{align*}
						\left| \sin^{\gamma_1} d \left( 0, \xi \right) \cos^{\gamma_2} d \left( 0, \xi \right) R_kf \left( \xi \right) \right| &\leq C \sin^{\gamma_1} d \left( 0, \xi \right) \cos^{\gamma_2 - 1} d \left( 0, \xi \right) \int\limits_{0}^{\frac{\pi}{2}} \left| \tilde{f} \left( \cos t \right) \right| \left( 1 - \frac{\tan^2 d \left( 0, \xi \right)}{\tan^2 t} \right)^{\frac{k}{2} - 1} \sin^{k - 1} t \ \mathrm{d}t \\
						&\leq C \sin^{\gamma_1} d \left( 0, \xi \right) \cos^{\gamma_2 - 1} d \left( 0, \xi \right) \int\limits_{d \left( 0, \xi \right)}^{\frac{\pi}{2}} \left| \tilde{f} \left( \cos t \right) \right| \sin^{k - 1} t \ \mathrm{d}t.
					\end{align*}
					Using {Inequality \eqref{MainLemmaEquationSphere} of} Lemma \ref{MainLemmaSphere} with $\eta_1 = k$ and $\eta_2 = 1$, we get
					\begin{align}
						\left| \sin^{\gamma_1} d \left( 0, \xi \right) \cos^{\gamma_2} d \left( 0, \xi \right) R_kf \left( \xi \right) \right| &\leq C \sin^{\gamma_1} d \left( 0, \xi \right) \cos^{\gamma_2 - 1} d \left( 0, \xi \right) \left( \int\limits_{d \left( 0, \xi \right)}^{\frac{\pi}{2}} \left| \tilde{f} \left( \cos t \right) \right| \sin^{pk - 1} t \cos^{p - 1} t \ \mathrm{d}t \right)^{\frac{1}{p}} \nonumber \\
						\label{KPlaneEstimateSnEndPoint}
						&= C \sin^{\gamma_1} d \left( 0, \xi \right) \cos^{\gamma_2 - 1} d \left( 0, \xi \right) \left( \int\limits_{d \left( 0, \xi \right)}^{\frac{\pi}{2}} \left| \tilde{f} \left( \cos t \right) \right| \sin^{\alpha_1 + n - 1} t \cos^{\alpha_2} t \sin^{pk - \alpha_1 - n} t \ \mathrm{d}t \right)^{\frac{1}{p}}.
					\end{align}
					We now divide our analysis in two cases.
					
					\textbf{Case I:} Let us assume that $p \geq \frac{\alpha_1 + n}{k}$. Then, $\sin^{pk - \alpha_1 - n} t \leq 1$, for all $t \in \left( 0, \frac{\pi}{2} \right)$. From Equation \eqref{KPlaneEstimateSnEndPoint} and using the fact that $\gamma_1 \geq 0$ and $\gamma_2 \geq 1$, we get
					\begin{align*}
						\left| \sin^{\gamma_1} d \left( 0, \xi \right) \cos^{\gamma_2} d \left( 0, \xi \right) R_kf \left( \xi \right) \right| &\leq C \| f \|_{L^{p, 1}_{\alpha_1, \alpha_2} \left( \mathbb{S}^n \right)}.
					\end{align*}
					
					\textbf{Case II:} Now, we consider $p < \frac{\alpha_1 + n}{k}$, we have that $\sin^{pk - \alpha_1 - n} t \leq \sin^{pk - \alpha_1 - n} d \left( 0, \xi \right)$. Keeping in mind that $\gamma_1 \geq \frac{\alpha_1 + n}{p} - k$ and $\gamma_2 \geq 1$, we get,
					\begin{align*}
						\left| \sin^{\gamma_1} d \left( 0, \xi \right) \cos^{\gamma_2} d \left( 0, \xi \right) R_kf \left( \xi \right) \right| &\leq C \sin^{\gamma_1 + k - \frac{\alpha_1 + n}{p}} d \left( 0, \xi \right) \cos^{\gamma_2 - 1} d \left( 0, \xi \right) \left( \int\limits_{0}^{\frac{\pi}{2}} \left| \tilde{f} \left( \cos t \right) \right| \sin^{\alpha + n - 1} t \cos^{\alpha_2} t \ \mathrm{d}t \right)^{\frac{1}{p}} \\
						&\leq C \| f \|_{L^{p, 1}_{\alpha_1, \alpha_2} \left( \mathbb{S}^n \right)}.
					\end{align*}
					This completes the proof of \ref{EndPoint2SnWeightedAProof}.
					\item[\mylabel{EndPoint2SnWeightedBProof}{(B)}] We now look at the case $k = 1$. We have, from Equation \eqref{XRayNewFormula}, 
					\begin{align*}
						\left| \sin^{\gamma_1} d \left( 0, \xi \right) \cos^{\gamma_2} d \left( 0, \xi \right) R_1f \left( \xi \right) \right| &= C \sin^{\gamma_1} d \left( 0, \xi \right) \cos^{\gamma_2} d \left( 0, \xi \right) \int\limits_{0}^{1} \tilde{f} \left( x \cos d \left( 0, \xi \right) \right) \left( 1 - x^2 \right)^{- \frac{1}{2}} \mathrm{d}x.
					\end{align*}
					Using Inequality \eqref{MainLemmaSphereEquation2} with $\eta_1 = 1 = \eta_2$, {and the fact that $p = 1 + \alpha_2$,} we get
					\begin{align*}
						\left| \sin^{\gamma_1} d \left( 0, \xi \right) \cos^{\gamma_2} d \left( 0, \xi \right) R_1f \left( \xi \right) \right| &\leq C \sin^{\gamma_1} d \left( 0, \xi \right) \cos^{\gamma_2} d \left( 0, \xi \right) \left( \int\limits_{0}^{1} \tilde{f} \left( x \cos d \left( 0, \xi \right) \right) \left( 1 - x^2 \right)^{\frac{p}{2} - 1} x^{{\alpha_2}} \mathrm{d}x \right)^{\frac{1}{p}}.
					\end{align*}
					By substituting $x \cos d \left( 0, \xi \right) = u$, we get
					\begin{align*}
						&\left| \sin^{\gamma_1} d \left( 0, \xi \right) \cos^{\gamma_2} d \left( 0, \xi \right) R_1f \left( \xi \right) \right| \\
						&\leq C \sin^{\gamma_1} d \left( 0, \xi \right) \cos^{\gamma_2} d \left( 0, \xi \right) \left( \int\limits_{0}^{\cos d \left( 0, \xi \right)} \tilde{f} \left( u \right) \left( 1 - \frac{u^2}{\cos^2 d \left( 0, \xi \right)} \right)^{\frac{p}{2} - 1} \left( \frac{u}{\cos d \left( 0, \xi \right)} \right)^{\alpha_2} \frac{\mathrm{d}u}{\cos d \left( 0, \xi \right)} \right)^{\frac{1}{p}}.
					\end{align*}
					Since $p \geq 2$, we have
					$$\left( 1 - \frac{u^2}{\cos^2 d \left( 0, \xi \right)} \right)^{\frac{p}{2} - 1} \leq \left( 1 - u^2 \right)^{\frac{p}{2} - 1}.$$
					Hence, we get
					\begin{align}
						\left| \sin^{\gamma_1} d \left( 0, \xi \right) \cos^{\gamma_2} d \left( 0, \xi \right) R_1f \left( \xi \right) \right| &\leq C \sin^{\gamma_1} d \left( 0, \xi \right) \cos^{\gamma_2 - 1} d \left( 0, \xi \right) \left( \int\limits_{0}^{\cos d \left( 0, \xi \right)} \tilde{f} \left( u \right) \left( 1 - u^2 \right)^{\frac{p}{2} - 1} u^{\alpha_2} \ \mathrm{d}u \right)^{\frac{1}{p}} \nonumber \\
						\label{XRayEndPointSnEstimate}
						&= C \sin^{\gamma_1} d \left( 0, \xi \right) \cos^{\gamma_2 - 1} d \left( 0, \xi \right) \left( \int\limits_{0}^{\cos d \left( 0, \xi \right)} \tilde{f} \left( u \right) \left( 1 - u^2 \right)^{\frac{\alpha_1 + n}{2} - 1} u^{\alpha_2} \left( 1 - u^2 \right)^{\frac{p - \alpha_1 - n}{2}} \mathrm{d}u \right)^{\frac{1}{p}}.
					\end{align}
					As in the proof of \ref{EndPoint2SnWeightedAProof}, we divide our analysis in two steps.
					
					\textbf{Case I:} First, let us assume $p \geq \alpha_1 + n$. Then, $\left( 1 - u^2 \right)^{\frac{p - \alpha_1 - n}{2}} \leq 1$, for all $u \in \left( 0, \cos d \left( 0, \xi \right) \right)$. Since we have $\gamma_1 \geq 0$ and $\gamma_2 \geq 1$, from Equation \eqref{XRayEndPointSnEstimate},
					\begin{align*}						
						\left| \sin^{\gamma_1} d \left( 0, \xi \right) \cos^{\gamma_2} d \left( 0, \xi \right) R_1f \left( \xi \right) \right| &\leq C \| f \|_{L^{p, 1}_{\alpha_1, \alpha_2} \left( \mathbb{S}^n \right)}.
					\end{align*}
					
					\textbf{Case II:} Now, we move to the case when $\alpha_1 + n > p$. Here, we have $\left( 1 - u^2 \right)^{\frac{p - \alpha_1 - n}{2}} \leq \sin^{p - \alpha_1 - n} d \left( 0, \xi \right)$, for $u \in \left( 0, \cos d \left( 0, \xi \right) \right)$. Hence, from Equation \eqref{XRayEndPointSnEstimate} and the fact that $\gamma_1 \geq \frac{\alpha_1 + n}{p} - 1$ and $\gamma_2 \geq 1$, we get
					\begin{align*}
						\left| \sin^{\gamma_1} d \left( 0, \xi \right) \cos^{\gamma_2} d \left( 0, \xi \right) R_1f \left( \xi \right) \right| &\leq C \sin^{\gamma_1 + 1 - \frac{\alpha_1 + n}{p}} d \left( 0, \xi \right) \cos^{\gamma_2 - 1} d \left( 0, \xi \right) \left( \int\limits_{0}^{\cos d \left( 0, \xi \right)} \tilde{f} \left( u \right) \left( 1 - u^2 \right)^{\frac{\alpha_1 + n}{2} - 1} u^{\alpha_2} \mathrm{d}u \right)^{\frac{1}{p}} \\
						&\leq C \| f \|_{L^{p, 1}_{\alpha_1, \alpha_2} \left( \mathbb{S}^n \right)}.
					\end{align*}
					{The necessity of the condition $p \geq 2$ when $k = 1$ is shown in Example \ref{XRayPAtLeast2}.}
				\end{enumerate}
				This completes the proof!
			\end{proof}

	\section{Formula for  \texorpdfstring{$k$}{k}-plane transform for general functions}
		\label{FormulaSection}
		In this section, we derive a formula for the totally-geodesic $k$-plane transform for functions defined on $\mathbb{R}^n$, $\mathbb{H}^n$, and $\mathbb{S}^n$. To this end, we use Kurusa's technique from \cite{KurusaHyperbolicSpace} and \cite{KurusaSphere}. The crux of this derivation is that the totally-geodesic submanifolds in spaces of constant curvature have rotational symmetry. Particularly, given a $k$-dimensional totally-geodesic submanifold $\xi$ in $X$ (where, $X = \mathbb{R}^n, \mathbb{H}^n$, or $\mathbb{S}^n$), and a geodesic $\gamma$ in $\xi$, we can cover all of $\xi$ by rotating the geodesic $\gamma$ inside $\xi$. {We employ this observation to get the desired formulae.}
		\subsection{The Euclidean Space \texorpdfstring{$\mathbb{R}^n$}{}}
			Here we get a formula for the $k$-plane transform in $\mathbb{R}^n$. This formula is not new; in fact, it a restatement of Rubin's formula in \cite{RubinKPlane}. However, we derive this formula using Kurusa's technique. This helps us to understand the workings, and get the formulae for the $k$-plane transform in the hyperbolic space and the sphere.
			
			{To begin with, let us consider the $X$-ray transform in $\mathbb{R}^2$.} This is the simplest case and {forms} the base of the general formula for the $k$-plane transform. {A line in $\mathbb{R}^2$ is parametrized by a direction, and its distance from the origin.} We denote by $\xi \left( \bar{\omega}, h \right)$, the line perpendicular to $\bar{\omega} \in \mathbb{S}^1$ at a distance $h > 0$ from the origin.
			
\begin{figure}[ht!]
				\centering
				\includegraphics[scale=1.0]{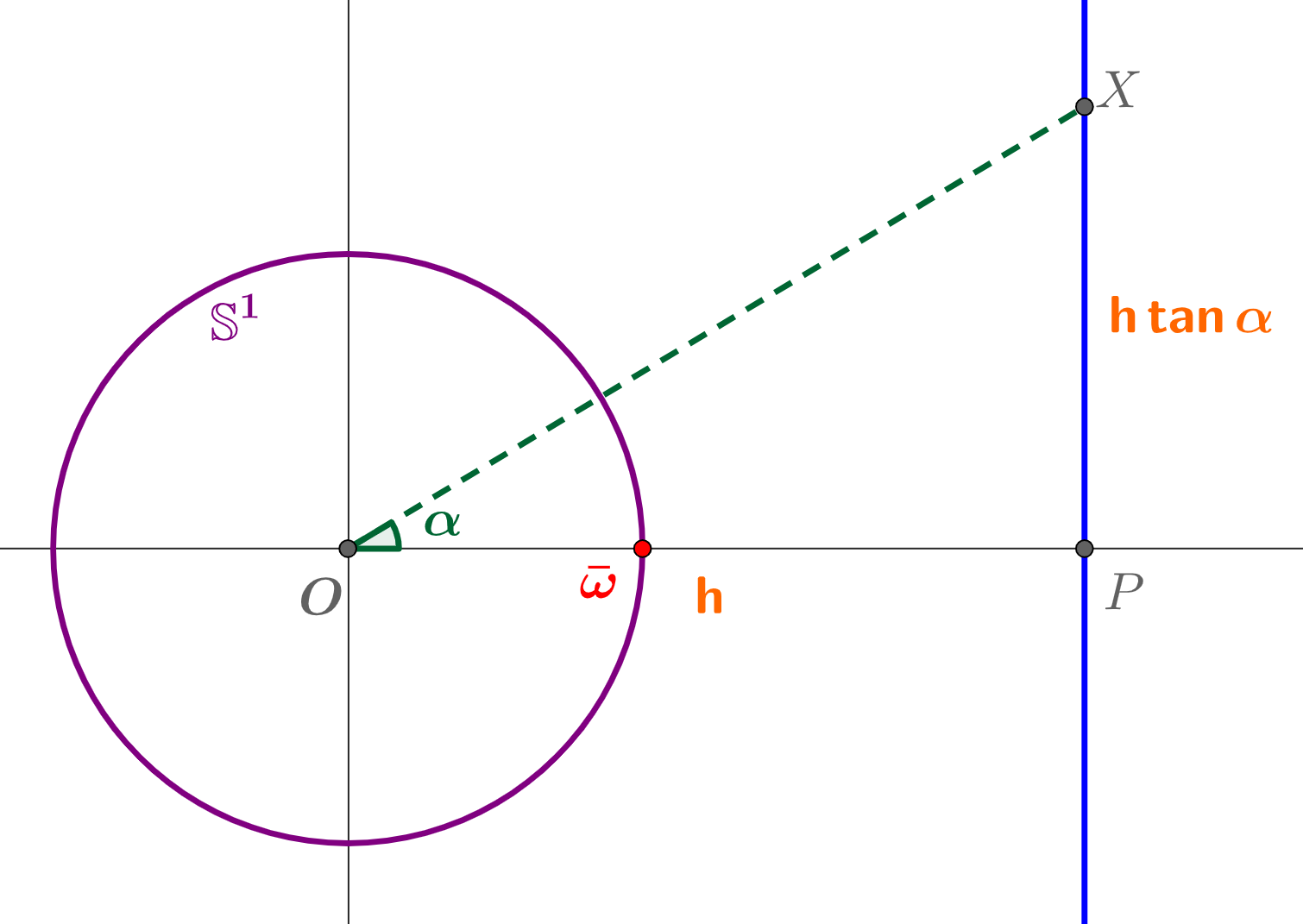}
				\caption{Parameterizing points on a line using angles at origin.}
				\label{LineAngles}
			\end{figure}

			Let $P$ be the point on the line closest to $0$, and $X$ be any point on the line. The triangle $OPX$ is a right-angled triangle, and the point $X$ can be determined in two ways: one, by the {signed} distance from $P$, and another, by the {directed} angle made with the vector $\bar{\omega}$. This situation is shown in Figure \ref{LineAngles}. The notation made in the figure is followed in the sequel. \\
			
			Using elementary trigonometry, it is clear that the {signed} distance $PX$ is given by $h \tan \alpha$, where $\alpha$ is the {directed} angle made by $X$ with $\bar{\omega}$. Suppose the angle of $\bar{\omega}$ with the $X$-axis is $\bar{\alpha}$. Then, the point $X$ can be given (in the polar coordinates of $\mathbb{R}^2$) by $\left( \alpha + \bar{\alpha}, h \sec \alpha \right)$. On the other hand, in the polar coordinates of the line $\xi \left( \bar{\omega}, h \right)$ with $P$ as its center, the point $X$ is given by $d \left( P, X \right) = h \tan \alpha$, where $\alpha \in \left( - \frac{\pi}{2}, \frac{\pi}{2} \right)$.
		
		Hence, we have
		\begin{align*}
			Rf \left( \bar{\omega}, h \right) &= \int\limits_{\xi \left( \bar{\omega}, h \right)} f \left( x \right) \mathrm{d}x = \int\limits_{- \frac{\pi}{2}}^{\frac{\pi}{2}} f \left( \alpha + \bar{\alpha}, {h} \sec \alpha \right) \frac{\mathrm{d}}{\mathrm{d}\alpha} d \left( P, x \right) \mathrm{d}\alpha = \int\limits_{- \frac{\pi}{2}}^{\frac{\pi}{2}} f \left( \alpha + \bar{\alpha}, {h} \sec \alpha \right) {h} \sec^2 \alpha \ \mathrm{d}\alpha.
		\end{align*}
		This pretty much gives the formula for the X-ray transform in $\mathbb{R}^2$. However, in this form, the formula is not quite useful since the notion of angles in higher dimensions is not easy to handle. However, we notice that the angle $\alpha + \bar{\alpha}$  {defines} a point $\omega$ on the circle $\mathbb{S}^1$. Moreover, since $\alpha \in \left( - \frac{\pi}{2}, \frac{\pi}{2} \right)$ {is the angle made by the point with $\bar{\omega}$}, we must have that $\langle \omega, \bar{\omega} \rangle > 0$. Converting the quantities involving the angles $\alpha$ and $\bar{\alpha}$ into inner products, we get the following formula.
		\begin{theorem}
			\label{XRayR2Formula}
			Let $f: \mathbb{R}^2 \rightarrow \mathbb{C}$ be a function that is integrable on almost every line in $\mathbb{R}^2$. When a line $\xi$ is parametrized by $\left( \bar{\omega}, h \right) \in \mathbb{S}^1 \times \left( 0, \infty \right)$, we have
			\begin{equation}
				\label{XRayR2}
				Rf \left( \bar{\omega}, h \right) = \int\limits_{\mathbb{S}^1_{\langle \omega, \bar{\omega} \rangle > 0}} f \left( \omega, \frac{{h}}{\langle \omega, \bar{\omega} \rangle} \right) \frac{{h}}{\langle \omega, \bar{\omega} \rangle^2} \ \mathrm{d}\omega,
			\end{equation}
			where, $\mathbb{S}^1_{\langle \omega, \bar{\omega} \rangle > 0} = \left\lbrace \omega \in \mathbb{S}^1 | \langle \omega, \bar{\omega} \rangle > 0 \right\rbrace$.
		\end{theorem}
		We are now ready to get the formula for the $k$-plane transform. The idea is to use rotational symmetry of a $k$-plane. That is, if we fix a line inside a $k$-plane $\eta$, and perform all ``$k$-dimensional rotations" inside $\eta$, we can cover whole of $\eta$. In this process, we consider the point in $\eta$ that is closest to the origin as our center for rotations.
		
		A $k$-plane $\eta$ in $\mathbb{R}^n$ can be parametrized by a $k$-dimensional subspace $\xi$ and a point $u \in \xi^{\perp}$, which is the ``foot of the perpendicular" from the origin {on} the plane $\eta$. The point $u \in \xi^{\perp}$ is the point in $\eta$ closest to $0$, and in our considerations is treated as the ``center" of the plane $\eta$.
		
		Any point on the plane $\eta$ can be {written as} $x + u$, for a unique $x \in \xi$. Further, using the polar decomposition in $\xi$, we have $x + u = t \sigma + u$, for $\sigma \in \mathbb{S}^{n - 1} \cap \xi$ and $t \in \left( 0, \infty \right)$. Therefore, we have
		$$R_kf \left( \xi, u \right) = \int\limits_{\xi} f \left( x + u \right) \mathrm{d}x = \int\limits_{\mathbb{S}^{n - 1} \cap \xi} \int\limits_{0}^{\infty} f \left( t \sigma + u \right) t^{k - 1} \mathrm{d}t \mathrm{d}\sigma.$$
		
\begin{figure}[ht!]
			\centering
			\includegraphics[scale=0.5]{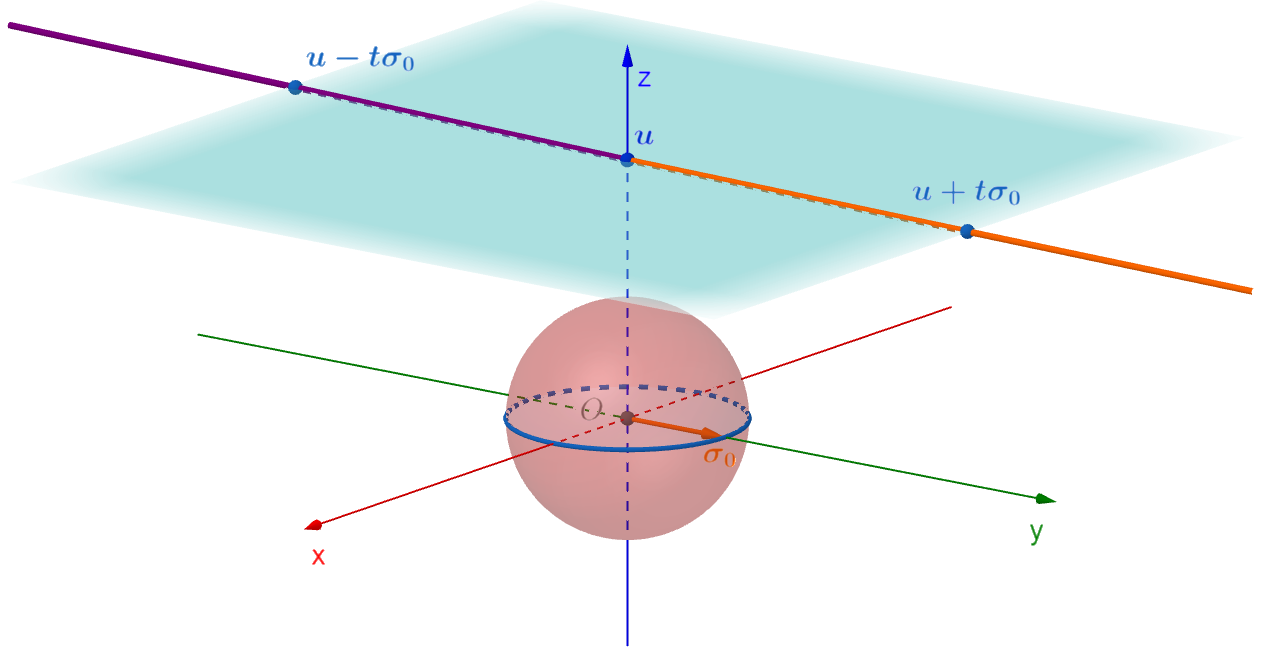}
			\caption{{A line in a plane obtained from a direction of $\mathbb{S}^{n - 1}$.}}
			\label{Plane2}
		\end{figure}

		Let us fix $\sigma_0 \in \mathbb{S}^{n - 1} \cap \xi$. Then, the set $u + \left\lbrace t \sigma_0 | t > 0 \right\rbrace$ is a ray in the plane $\eta$ starting at $u$ and in the ``positive" direction of $\sigma_0$. However, there is a rotation in $SO_{\xi} \left( k \right)$ (the rotation group of $\xi$), {say $R_{\sigma_0}$}, that takes $\sigma_0$ to $- \sigma_0$, and as $t \in \left( 0, \infty \right)$ varies, we get the whole line in $\eta$ through $u$ in the direction $\sigma_0$. Now, the entire plane is covered (upto a double repetition) by $k$-dimensional rotations in $\eta$ of the line $u + \left\lbrace t \sigma_0 | t \in \mathbb{R} \right\rbrace$. An illustration of this concept is shown in Figure \ref{Plane2}.
		
		Since integration over $\mathbb{S}^{n - 1} \cap \xi$ can be written as an integration over $SO_{\xi} \left( k \right)$, upto an appropriate normalization, we have
		\begin{align*}
			R_kf \left( \xi, u \right) &= \int\limits_{\mathbb{S}^{n - 1} \cap \xi} \int\limits_{0}^{\infty} f \left( t \sigma + u \right) t^{k - 1} \ \mathrm{d}t \mathrm{d}\sigma \\
			&= |\mathbb{S}^{n - 1} \cap \xi| \int\limits_{SO_{\xi} \left( k \right)} \int\limits_{0}^{\infty} f \left( t g \sigma_0 + u \right) t^{k - 1} \ \mathrm{d}t \mathrm{d}g \\
			&= |\mathbb{S}^{n - 1} \cap \xi| \int\limits_{\faktor{SO_{\xi} \left( k \right)}{\left\lbrace {I, R_{\sigma_0}} \right\rbrace}} \int\limits_{-\infty}^{\infty} f_{g} \left( t \sigma_0 + u \right) |t|^{k - 1} \ \mathrm{d}t \mathrm{d}g.
		\end{align*}
		{Here, $f_g \left( x \right) = f \left( gx \right)$, for any $g \in SO(n)$.} However, the inner integral represents the $X$-ray transform {of} the function $f_g \left( \cdot \right) \left| d \left( u, \cdot \right) \right|^{k - 1}$ in the $\left\lbrace u, \sigma_0 \right\rbrace$-plane along the line through $u$ in the direction $\sigma_0$. {We consider $u \neq 0$ so that $\left\lbrace u, \sigma_0 \right\rbrace$ indeed spans a plane.} Let us denote by $\mathbb{S}^1_{\left\lbrace u, \sigma_0 \right\rbrace}$, the unit circle in the plane spanned by $\left\lbrace u, \sigma_0 \right\rbrace$. Then, using Theorem \ref{XRayR2Formula}, we get
		\begin{align*}
			R_kf \left( \xi, u \right) &= |\mathbb{S}^{n - 1} \cap \xi| \int\limits_{\faktor{SO_{\xi} \left( k \right)}{\left\lbrace \pm I \right\rbrace}} \int\limits_{\mathbb{S}^{1}_{\left\lbrace u, \sigma_0 \right\rbrace, \langle \omega, \frac{u}{\| u \|} \rangle > 0}} f_g \left( \omega, \frac{\| u \|}{\langle \omega, \frac{u}{\| u \|} \rangle} \right) \left( \frac{\| u \| \sqrt{1 - \langle \omega, \frac{u}{\| u \|} \rangle^2}}{\langle \omega, \frac{u}{\| u \|} \rangle} \right)^{k - 1} \frac{\| u \|}{\langle \omega, \frac{u}{\| u \|} \rangle^2} \ \mathrm{d}\omega \mathrm{d}g.
		\end{align*}
		Now, we notice that as $g$ varies over $SO_{\xi} \left( k \right)$, $g\omega$ varies over $\mathbb{S}^{n - 1} \cap \xi$. On the other hand, since $u \in \xi^{\perp}$, $\langle \omega, \frac{u}{\| u \|} \rangle$ is invariant under the action of $SO_{\xi} \left( k \right)$, and gives the cosine of the angle between $\omega$ and $\frac{u}{\| u \|}$. Consequently, we have using the polar decomposition of the $k$-dimensional sphere,
		\begin{align*}
			R_kf \left( \xi, u \right) &= \int\limits_{\mathbb{S}^{n - 1} \cap \xi} \int\limits_{0}^{\frac{\pi}{2}} f \left( \omega, \frac{\| u \|}{\cos \theta} \right) \frac{\| u \|^{k}}{\cos^{k + 1} \theta} \sin^{k - 1} \theta \ \mathrm{d}\theta \mathrm{d} \omega = \int\limits_{\mathbb{S}^{n - 1}_{\langle \omega, \frac{u}{\| u \|} \rangle > 0} \cap \left( \xi \oplus u \right)} f \left( \omega, \frac{\| u \|}{\langle \omega, \frac{u}{\| u \|} \rangle} \right) \frac{\| u \|^{k}}{\langle \omega, \frac{u}{\| u \|} \rangle^{k + 1}} \ \mathrm{d}\omega,
		\end{align*}
		where, by $\xi \oplus u$, we mean the direct sum of $\xi$ and the subspace spanned by $u$. That is, we have the following formula for the $k$-plane transform in $\mathbb{R}^n$.
		\begin{theorem}
			\label{KPlaneRNFormula}
			Let $f: \mathbb{R}^n \rightarrow \mathbb{C}$ be integrable over almost every $k$-dimensional plane. {By parametrizing a $k$-plane $\eta$ by $\left( \xi, u \right)$}, where $\xi$ is a $k$-dimensional subspace and $u \in \xi^{\perp} {\setminus \left\lbrace 0 \right\rbrace}$, we have
			\begin{equation}
				\label{KPlaneRN}
				R_kf \left( \xi, u \right) = \int\limits_{\mathbb{S}^{n - 1}_{\langle \omega, \frac{u}{\| u \|} \rangle > 0} \cap \left( \xi \oplus u \right)} f \left( \omega, \frac{\| u \|}{\langle \omega, \frac{u}{\| u \|} \rangle} \right) \frac{\| u \|^k}{\langle \omega, \frac{u}{\| u \|} \rangle^{k + 1}} \mathrm{d}\omega.
			\end{equation}
		\end{theorem}
		We remark here that Equation \eqref{KPlaneRN} {is a different version of} Rubin's {formula} in \cite{RubinKPlane}. {In \cite{RubinKPlane}, the author uses the formula for the $k$-plane transform to get weighted $L^p$-$L^p$ estimates, and obtain the operator norm of the operator $R_k$.} Motivated from this work, we {derive certain} weighted $L^p$-$L^p$ estimates for the $k$-plane transform on the hyperbolic space and the sphere. For the same, we derive the respective formulae in the next two subsections.
		\subsection{The Hyperbolic Space \texorpdfstring{$\mathbb{H}^n$}{}}
\begin{figure}[ht!]
			\centering
			\includegraphics[scale=0.9]{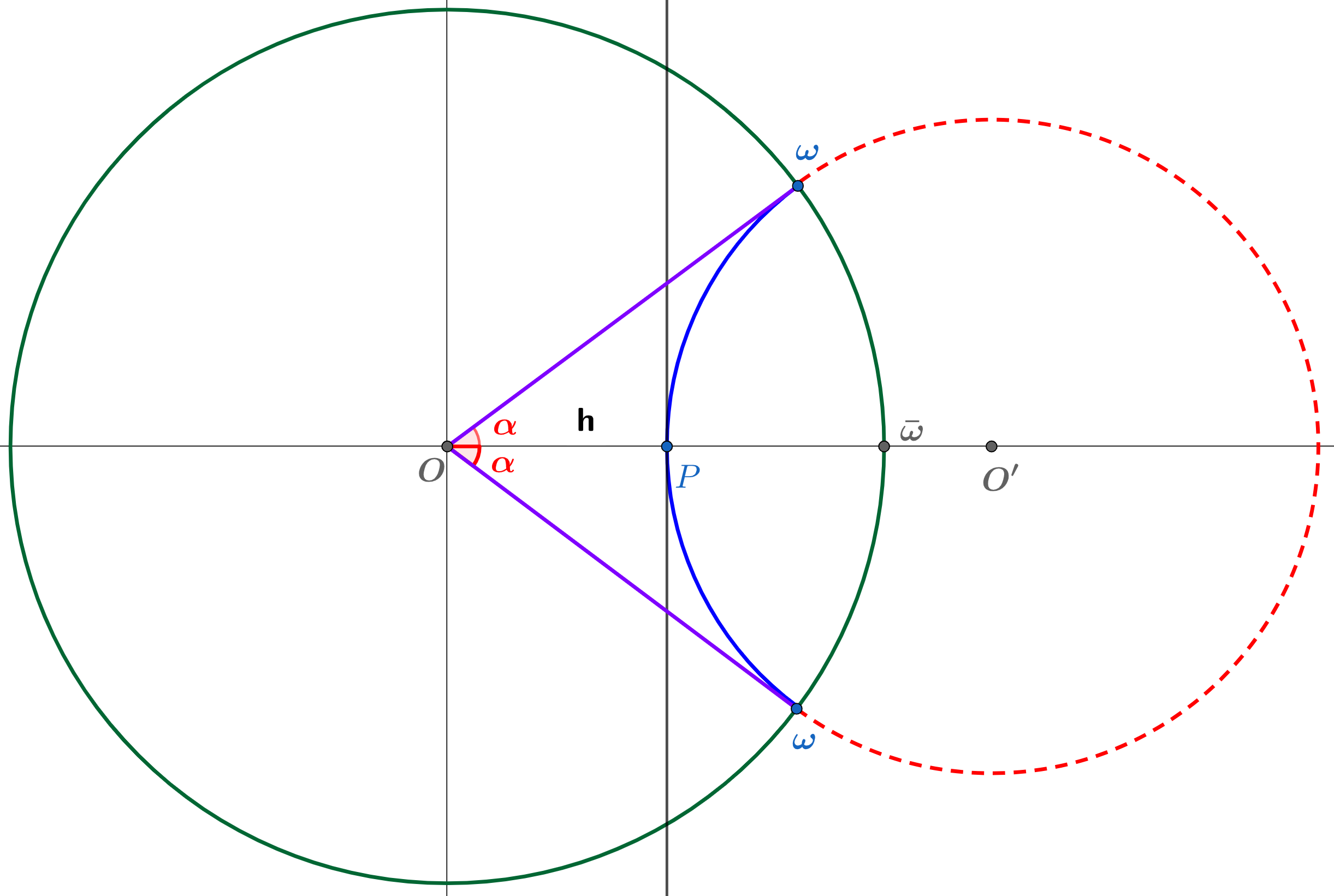}
			\caption{For a given geodesic at a distance $h$, the angles made at the origin cannot exceed $\alpha = \cos^{-1} \left( \tanh h \right)$.}
			\label{TGS}
		\end{figure}

			The formula for the Radon transform for the real Hyperbolic space was found by Kurusa (\cite{KurusaHyperbolicSpace}). We use Kurusa's techniques to get a formula for the $k$-plane transform on the real hyperbolic space $\mathbb{H}^n$. We recall that we have fixed the Poincar\'{e} ball model as the description of the real hyperbolic space $\mathbb{H}^n$. Before we begin with the formal calculations to obtain the formula for the $k$-plane transform on $\mathbb{H}^n$, we state the formula for X-ray transform in $\mathbb{H}^2$, obtained by Kurusa in \cite{KurusaHyperbolicSpace}.
		\begin{theorem}[\cite{KurusaHyperbolicSpace}]
			\label{XRayH2Formula}
			Let $f: \mathbb{H}^2 \rightarrow \mathbb{C}$ be a function that is integrable on almost every geodesic in $\mathbb{H}^2$, where $\mathbb{H}^2$ is the real hyperbolic space of dimension $2$ described by the Poincar\'{e} ball model. We parametrize a geodesic $\xi$ by $\left( \bar{\omega}, h \right) \in \mathbb{S}^{1} \times \left( 0, \infty \right)$, where $\bar{\omega}$ {is} the direction of the foot of the perpendicular to $\xi$, and $h$ is the hyperbolic distance of $\xi$ from $0$, the hyperbolic origin. Then, the X-ray transform of $f$ is given by
			\begin{equation}
				\label{XRayH2}
				Rf \left( \bar{\omega}, h \right) = \int\limits_{{\mathbb{S}^{1}_{\langle \omega, \bar{\omega} \rangle > \tanh h}}} f \left( \omega, \frac{1}{2} \bar{\mu} \left( \frac{\tanh h}{\langle \omega, \bar{\omega} \rangle} \right) \right) \frac{\left( \coth^2 h \ \langle \omega, \bar{\omega} \rangle^2 - 1 \right)^{-1}}{\sinh h} \ \mathrm{d}\omega,
			\end{equation}
			where $\bar{\mu} \left( x \right) = \ln \left( \frac{1 + x}{1 - x} \right) = 2 \tanh^{-1} x$, for $0 < x < 1$, and $\mathbb{S}^{1}_{\langle \omega, \bar{\omega} \rangle > \tanh h} = \left\lbrace \omega \in \mathbb{S}^1 | \langle \omega, \bar{\omega} \rangle > \tanh h \right\rbrace$.
		\end{theorem}
		It is to be noticed that Equation \eqref{XRayH2} is much different from the Euclidean case. Particularly, depending on the distance {$h$ of the geodesic $\left( \bar{\omega}, h \right)$} from the origin, we have a restricted spherical cap on which we integrate. In contrast, for a Euclidean geodesic (line), we get the domain of integration as a semicircle, irrespective of the distance of the line from the origin. Figure \ref{TGS} illustrates that the restriction is natural due to the geometry of the space and its geodesics. \\
		
		{We quickly see that for a given geodesic at a hyperbolic distance $h > 0$, the angle $\alpha$ mentioned in Figure \ref{TGS} is $\cos^{-1} \left( \tanh h \right)$. To do so, we first observe that $O\omega$ is perpendicular to $O'\omega$. Also, the Euclidean length of $OP$ is $\tanh \left( \frac{h}{2} \right)$. Therefore, if the length of the segment joining $O'$ and $\omega$ is $r$, then by the Pythagoras' theorem for the Euclidean space, we have $\left( r + \tanh \left( \frac{h}{2} \right) \right)^2 = 1 + r^2$. Upon simplifying, we get $r = \text{cosech } h$. Therefore, the Euclidean length of $OO'$ is $\tanh \left( \frac{h}{2} \right) + \text{cosech } h = \coth h$. It is now easy to see that $\cos \alpha = \tanh h$.}
		
		We use Theorem \ref{XRayH2Formula} and the parametriation described in Theorem \ref{TGSParametrization} to get the formula for the $k$-plane transform.
				
		{We use the strategy as in the Euclidean case, to decompose the integrals into polar form, and realise the inner integral as the $X$-ray transform in two dimensions. By using polar decomposition on the $k$-dimensional totally-geodesic submanifold $\left( \xi, u \right)$ in $\mathbb{H}^n$, we have}
		\begin{align*}
			R_kf \left( \xi, u \right) &= \int\limits_{\mathbb{S}^{n - 1} \cap \xi} \int\limits_{0}^{\infty} f \left( \exp_u \left( t\sigma \right) \right) \sinh^{k - 1} t \ \mathrm{d}t \mathrm{d}\sigma \\
			&= \left| \mathbb{S}^{n - 1} \cap \xi \right| \int\limits_{SO_{\xi} \left( k \right)} \int\limits_{0}^{\infty} f \left( \exp_u \left( gt\sigma_0 \right) \right) \sinh^{k - 1} t \ \mathrm{d}t \mathrm{d}g \\
			&= \left| \mathbb{S}^{n - 1} \cap \xi \right| \int\limits_{\faktor{SO_{\xi} \left( k \right)}{\left\lbrace {I, R_{\sigma_0}} \right\rbrace}} \int\limits_{- \infty}^{\infty} f_g \left( \exp_{u} \left( t \sigma_0 \right) \right) \sinh^{k - 1} |t| \ \mathrm{d}t \mathrm{d}g,
		\end{align*}
		{where, $R_{\sigma_0} \in SO(n)$ is a rotation that takes $\sigma_0$ to $- \sigma_0$. Now, the inner integral can be realised as the $X$-ray transform of the function $f_g \left( \cdot \right) \sinh^{k - 1} d \left( 0, \cdot \right)$ along the geodesic through $u$ in the direction of $\sigma_0$. Hence, by using Equation \eqref{XRayH2Formula}, we get}
		\begin{align*}
			R_kf \left( \xi, u \right) = \left| \mathbb{S}^{n - 1} \cap \xi \right| \int\limits_{\faktor{SO_{\xi} \left( k \right)}{\left\lbrace {I, R_{\sigma_0}} \right\rbrace}} \int\limits_{\mathbb{S}^1_{\left\lbrace u, \sigma_0 \right\rbrace \cap \left\langle \omega, \frac{u}{\| u \|} \right\rangle > \tanh h}} &f_g \left( \omega, \frac{1}{2} \bar{\mu} \left( \frac{\tanh h}{\langle \omega, \bar{\omega} \rangle} \right) \right) \frac{\left( \coth^2 h \langle \omega, \frac{u}{\| u \|} \rangle^2 - 1 \right)^{-1}}{\sinh h} \times \\
			&\sinh^{k - 1} \left( \frac{1}{2} \bar{\mu} \left( \frac{\tanh h}{\langle \omega, \frac{u}{\| u \|} \rangle} \right) \right) \mathrm{d} \omega \mathrm{d}g
		\end{align*}
		{We observe that as $g$ varies over $SO_{\xi} \left( k \right)$, the inner integral varies over a spherical cap centered at $\frac{u}{\| u \|}$ and of spherical radius $\tanh h$, in the subspace containing both $\xi$ and $u$. That is, we have,}
		\begin{align*}
			R_kf \left( \xi, u \right) &= \int\limits_{\mathbb{S}^{n - 1}_{\left\langle \omega, \frac{u}{\| u \|} \right\rangle > {\tanh h}} \cap \left( \xi \oplus \frac{u}{\| u \|} \right)} f \left( \omega, \frac{1}{2} \bar{\mu} \left( \frac{\tanh h}{\langle \omega, \frac{u}{\| u \|} \rangle} \right) \right) \frac{\left( \coth^2 h \langle \omega, \frac{u}{\| u \|} \rangle - 1 \right)^{- \frac{\left( k + 1 \right)}{2}}}{\sinh h} \ \mathrm{d}\omega.
		\end{align*}
		Thus, we have the formula for the $k$-plane transform on $\mathbb{H}^n$.
		\begin{theorem}
			\label{KPlaneHNFormula}
			Let $f: \mathbb{H}^n \rightarrow \mathbb{C}$ be integrable on every $k$-dimensional totally-geodesic submanifold. Let us parametrize a $k$-dimensional totally-geodesic submanifold $\eta$ by $\left( \xi, u \right)$, where $\xi \in G_{n, k}$ and $u \in \xi^{\perp} {\setminus \left\lbrace 0 \right\rbrace}$, as given in Theorem \ref{TGSParametrization}. Then,
			\begin{equation}
				\label{KPlaneHN}
				R_kf \left( \xi, u \right) = \int\limits_{\mathbb{S}^{n - 1}_{\left\langle \omega, \frac{u}{\| u \|} \right\rangle > {\tanh h}} \cap \left( \xi \oplus \frac{u}{\| u \|} \right)} f \left( \omega, \frac{1}{2} \bar{\mu} \left( \frac{\tanh h}{\langle \omega, \frac{u}{\| u \|} \rangle} \right) \right) \frac{\left( \coth^2 h \langle \omega, \frac{u}{\| u \|} \rangle^2 - 1 \right)^{- \frac{\left( k + 1 \right)}{2}}}{\sinh h} \ \mathrm{d}\omega.
			\end{equation}
			We also have
			\begin{equation}
				\label{KPlaneHNGroup}
				R_kf \left( g \cdot \mathbb{R}^k, h g \cdot e_{k + 1} \right) = \int\limits_{\mathbb{S}^{k}_{\langle \omega, e_{k + 1} \rangle > \tanh h}} f_g \left( \omega, \frac{1}{2} \bar{\mu} \left( \frac{\tanh h}{\langle \omega, e_{k + 1} \rangle} \right) \right) \frac{\left( \coth^2 h \langle \omega, e_{k + 1} \rangle^2 - 1 \right)^{- \frac{\left( k + 1 \right)}{2}}}{\sinh h} \ \mathrm{d}\omega,
			\end{equation}
			where $g \in SO \left( n \right)$ is such that $g \cdot \mathbb{R}^k = \xi$ and $g \cdot e_{k + 1} = \frac{u}{\| u \|}$, in Equation \eqref{KPlaneHN}{, and $h > 0$ is the distance of the totally-geodesic submanifold $\left( \xi, u \right)$ from $0 \in \mathbb{H}^n$}.
		\end{theorem}
		\subsection{The Sphere \texorpdfstring{$\mathbb{S}^n$}{}}
\begin{figure}[ht!]
				\centering
				\begin{subfigure}[t]{0.45\textwidth}
					\centering
					\includegraphics[scale=0.5]{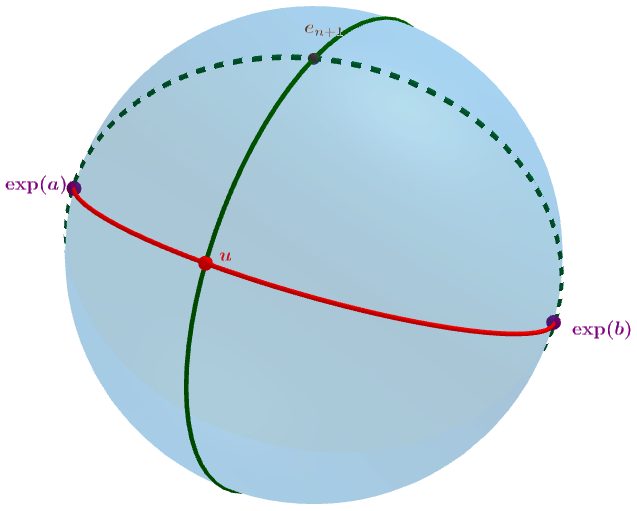}
					\caption{Half of geodesic ({red}) on $\mathbb{S}^n$.}
					\label{HalfGeodesicSn}
				\end{subfigure}
				~~~~~~~~
				\begin{subfigure}[t]{0.45\textwidth}
					\centering
					\includegraphics[scale=1.2]{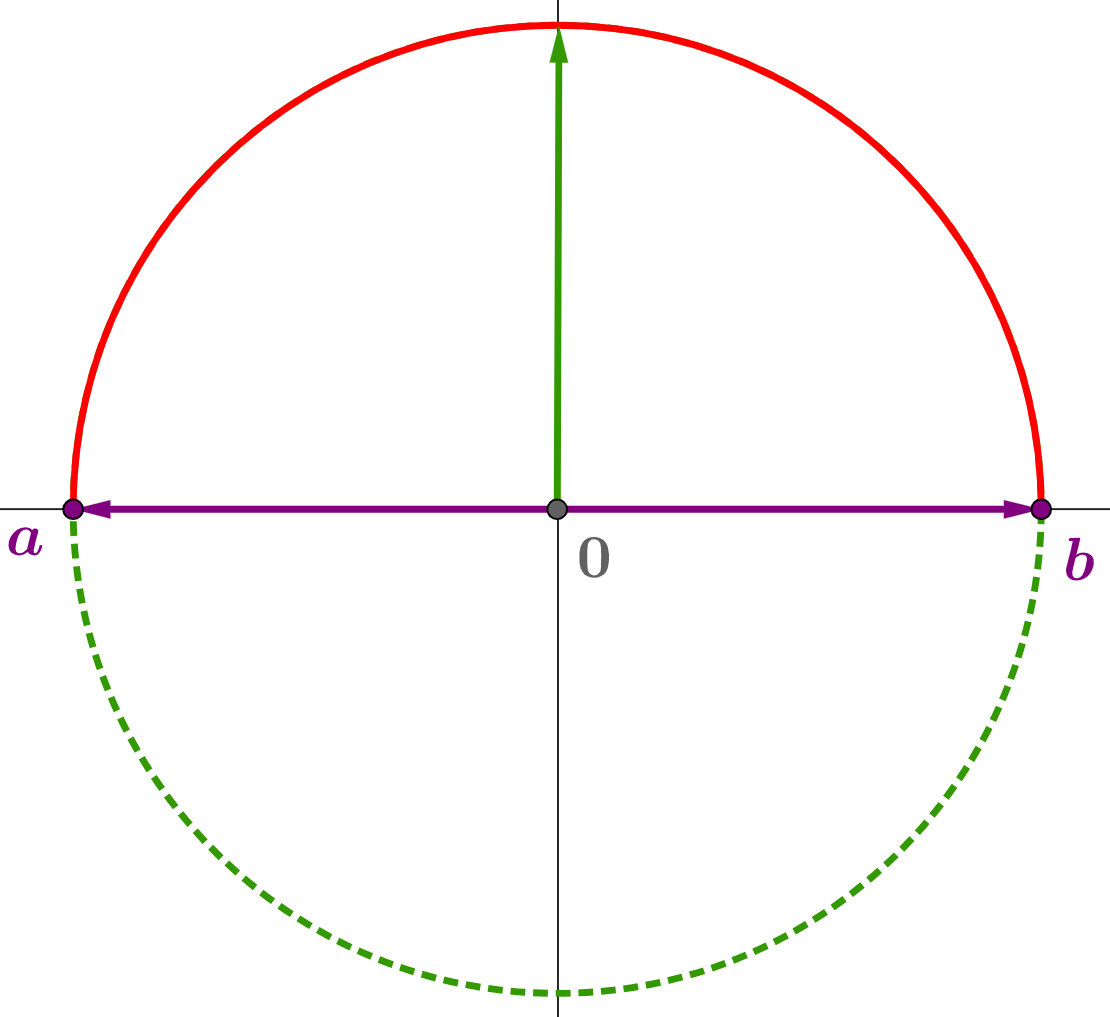}
					\caption{A semi-circle in the tangent space of $\mathbb{S}^n$ at $e_{n + 1}$.}
					\label{TangentSpaceSemiCircle}
				\end{subfigure}
				\caption{As a point moves along half of a geodesic on the sphere, shown in \ref{HalfGeodesicSn}, we get a semicircle in the tangent space of $\mathbb{S}^n$ at $e_{n + 1}$, shown in \ref{TangentSpaceSemiCircle}.}
				\label{TGSSnParameter}
			\end{figure}

			Similar to the case of the hyperbolic space, $\mathbb{H}^n$, Kurusa derived the formula for the Radon transform on the sphere, $\mathbb{S}^n$, in \cite{KurusaSphere}. An important observation in the case of the sphere is that the study of Radon (and $k$-plane) transforms is sufficient for even functions, since odd functions are in their null space. In what follows, we consider only even functions on $\mathbb{S}^n$. Kurusa's idea to get the formula for the Radon transform came from the fact that the codimension-$1$ totally-geodesic submanifolds of $\mathbb{S}^n$ are precisely lower dimensional spheres. Therefore, it was sufficient to get a formula for the $X$-ray transform on $\mathbb{S}^2$, and then by rotation symmetry, we get the formula for the Radon transform. We use this idea to get the formula for the $k$-plane transform. We observe that since we only consider even functions, to study the $X$-Ray transform along a geodesic, we only need to consider half of the given geodesic. As seen in Figure \ref{TGSSnParameter}, it is evident that as a point moves along half of a geodesic, its pre-image under the exponential map covers a semicircle on the tangent space of $e_{n + 1}$, the origin of $\mathbb{S}^n$.
			
			We now state Kurusa's formula for the $X$-Ray transform in $\mathbb{S}^2$, given in \cite{KurusaSphere}.
		\begin{theorem}[\cite{KurusaSphere}]
			\label{XRayS2}
			Let $f: \mathbb{S}^2 \rightarrow \mathbb{C}$ be a function integrable on almost every geodesic in $\mathbb{S}^2$. Upon parametrizing the geodesics in $\mathbb{S}^2$ by $\left( \bar{\omega}, h \right)$, {where $\bar{\omega}  \in \mathbb{S}^1$ is the direction from $e_{n + 1}$ to the point on $\left( \bar{\omega}, h \right)$ closest to $e_{n + 1}$, and $h > 0$ is the spherical distance of the geodesic from $e_{n + 1}$,} the $X$-ray transform of $f$ is given by
			\begin{equation}
				\label{XRayS2Formula}
				Rf \left( \bar{\omega}, h \right) = \int\limits_{\mathbb{S}^1_{\langle \omega, \bar{\omega} \rangle > 0}} f \left( \omega, \tan^{-1} \left( \frac{\tan h}{\langle \omega, \bar{\omega} \rangle} \right) \right) \frac{\left( \cot^2 h \ \langle \omega, \bar{\omega} \rangle^2 + 1 \right)^{-1}}{\sin h} \ \mathrm{d}\omega,
			\end{equation}
			where $\mathbb{S}^1_{\langle \omega, \bar{\omega} \rangle > 0} = \left\lbrace \omega \in \mathbb{S}^1 | \langle \omega, \bar{\omega} \rangle > 0 \right\rbrace$.
		\end{theorem}
		
		With the parametrization of {$k$-dimensional totally-geodesic submanifolds of $\mathbb{S}^n$ given in} Theorem \ref{TGSParametrization}, we now get the formula for the $k$-plane transform on $\mathbb{S}^n$. {Analogous to the Euclidean and Hyperbolic spaces,} {using the polar decomposition on the totally-geodesic submanifold $\left( \xi, u \right)$ in $\mathbb{S}^n$, we get}
		\begin{align*}
			R_kf \left( \xi, u \right) &= \int\limits_{\mathbb{S}^{n - 1} \cap \xi} \int\limits_{0}^{\frac{\pi}{2}} f \left( \exp_u \left( t\sigma \right) \right) \sinh^{k - 1} t \ \mathrm{d}t \mathrm{d}\sigma \\
			&= \left| \mathbb{S}^{n - 1} \cap \xi \right| \int\limits_{SO_{\xi} \left( k \right)} \int\limits_{0}^{\frac{\pi}{2}} f \left( \exp_u \left( gt\sigma_0 \right) \right) \sin^{k - 1} t \ \mathrm{d}t \mathrm{d}g \\
			&= \left| \mathbb{S}^{n - 1} \cap \xi \right| \int\limits_{\faktor{SO_{\xi} \left( k \right)}{\left\lbrace {I, R_{\sigma_0}} \right\rbrace}} \int\limits_{- \frac{\pi}{2}}^{\frac{\pi}{2}} f_g \left( \exp_{u} \left( t \sigma_0 \right) \right) \sin^{k - 1} |t| \ \mathrm{d}t \mathrm{d}g,
			\end{align*}
			{where, $R_{\sigma_0} \in SO_{\xi}(k)$ is such that $R_{\sigma_0} \left( \sigma_0 \right) = - \sigma_0$. We now notice that the inner integral is the $X$-ray transform of the function $f_g \left( \cdot \right) \sin^{k - 1} d \left( 0, \cdot \right)$ along the geodesic through $u$ in the direction of $\sigma_0$. Therefore, by using Equation \eqref{XRayS2Formula}, we get}
			\begin{align*}
				R_kf \left( \xi, u \right) = \left| \mathbb{S}^{n - 1} \cap \xi \right| \int\limits_{\faktor{SO_{\xi} \left( k \right)}{\left\lbrace I, R_{\sigma_0} \right\rbrace}} \int\limits_{\mathbb{S}^1_{\left\lbrace u, \sigma_0 \right\rbrace \cap \langle \omega, \frac{u}{\| u \|} \rangle > 0}} &f_g \left( \omega, \tan^{-1} \left( \frac{\tan h}{\langle \omega, \bar{\omega} \rangle} \right) \right) \frac{\left( \cot^2 h \langle \omega, \frac{u}{\| u \|} \rangle^2 + 1 \right)^{-1}}{\sin h} \times \\
				&\sin^{k - 1} \left( \tan^{-1} \left( \frac{\tan h}{\langle \omega, \frac{u}{\| u \|} \rangle} \right) \right) \mathrm{d}\omega \mathrm{d}g.
			\end{align*}
			{It is now to be noticed that as $g$ varies over all of $SO_{\xi} \left( k \right)$, the inner integral gives a spherical cap of the sphere of the subspace containing $\xi$ and $u$, with the center $\frac{u}{\| u \|}$. That is, we have,}
			\begin{align*}
				R_kf \left( \xi, u \right) &= \int\limits_{\mathbb{S}^{n - 1}_{\langle \omega, \frac{u}{\| u \|} \rangle > 0} \cap \left( \xi \oplus \frac{u}{\| u \|} \right)} f \left( \omega, \tan^{-1} \left( \frac{\tan h}{\langle \omega, \frac{u}{\| u \|} \rangle} \right) \right) \frac{\left( \cot^2 h \ \langle \omega, \frac{u}{\| u \|} \rangle {+} 1 \right)^{- \frac{\left( k + 1 \right)}{2}}}{\sin h} \ \mathrm{d}\omega.
		\end{align*}
		Thus, we have the following formula for the $k$-plane transform on $\mathbb{S}^n$.
		\begin{theorem}
			\label{KPlaneSNFormula}
			Let $f: \mathbb{S}^n \rightarrow \mathbb{C}$ be integrable on every $k$-dimensional totally-geodesic submanifold of $\mathbb{S}^n$. Upon parametrizing the set $\Xi_k$ of all $k$-dimensional totally-geodesic submanifolds of $\mathbb{S}^n$ by $\left( \xi, u \right)$, where $\xi \in G_{n, k}$ and $u \in \left( \xi^{\perp} \cap B \left( 0, \frac{\pi}{2} \right) \right) {\setminus \left\lbrace 0 \right\rbrace}$, as given in Theorem \ref{TGSParametrization}, we have
			\begin{equation}
				\label{KPlaneSN}
				R_kf \left( \xi, u \right) = \int\limits_{\mathbb{S}^{n - 1}_{\langle \omega, \frac{u}{\| u \|} \rangle > 0} \cap \left( \xi \oplus \frac{u}{\| u \|} \right)} f \left( \omega, \tan^{-1} \left( \frac{\tan h}{\langle \omega, \frac{u}{\| u \|} \rangle} \right) \right) \frac{\left( \cot^2 h \ \langle \omega, \frac{u}{\| u \|} \rangle {+} 1 \right)^{- \frac{\left( k + 1 \right)}{2}}}{\sin h} \ \mathrm{d}\omega.
			\end{equation}
			We also have
			\begin{equation}
				\label{KPlaneSNGroup}
				R_kf \left( g \cdot \mathbb{R}^k, h g \cdot e_{k + 1} \right) = \int\limits_{\mathbb{S}^k_{\langle \omega, e_{k + 1} \rangle > 0}} f_g \left( \omega, \tan^{-1} \left( \frac{\tan h}{\langle \omega, e_n \rangle} \right) \right) \frac{\left( \cot^2 h \ \langle \omega, e_{k + 1} \rangle^2 {+} 1 \right)^{- \frac{\left( k + 1 \right)}{2}}}{\sin h} \ \mathrm{d}\omega,
			\end{equation}
			where, $g \in SO \left( n \right)$ is such that $g \cdot \mathbb{R}^k = \xi$ and $g \cdot e_{k + 1} = \frac{u}{\| u \|}$, in Equation \eqref{KPlaneSN}{, and $0< h < \frac{\pi}{2}$ is the distance of the totally-geodesic submanifold $\left( \xi, u \right)$ from $e_{n + 1} \in \mathbb{S}^n$}.
		\end{theorem}
		
	\section{ \texorpdfstring{$L^p$}{Lp}-\texorpdfstring{$L^p$}{Lp} estimates for general functions}
		\label{LPLPEstimatesSection}
		In this section, we get weighted $L^p$-$L^p$ estimates for the $k$-plane transform in $\mathbb{H}^n$ and $\mathbb{S}^n$. Similar questions have been answered by Rubin in \cite{RubinKPlane} for the Euclidean case. Particularly, in the Euclidean case, we have the following result.
		\begin{theorem}[\cite{RubinKPlane}]
			\label{RubinEuclideanLPLPResult}
			For $\alpha > k - n$, $1 \leq p < \frac{\alpha + n}{k}$, and $\beta > k - n$, the $k$-plane transform $R_k: L^p_{\alpha} \left( \mathbb{R}^n \right) \rightarrow L^p_{\beta} \left( G \left( n, k \right) \right)$ is a bounded operator if and only if $\beta = \alpha - k \left( p - 1 \right)$.
		\end{theorem}
		{We proceed to get analogous results in the hyperbolic space and the sphere with radial power weights. We start with the hyperbolic space $\mathbb{H}^n$.}
		\subsection{The Hyperbolic Space \texorpdfstring{$\mathbb{H}^n$}{}}
			In this section, we look at the boundedness of $k$-plane transform from the weighted Lebesgue space $L^p_{\alpha_1, \alpha_2} \left( \mathbb{H}^n \right)$	to $L^p_{\beta_1, \beta_2} \left( \Xi_k \left( \mathbb{H}^n \right) \right)$. The necessary conditions {given in Equation \eqref{NecessaryMixedWeightInequality}} of Theorem \ref{MixedWeightNecessary} for $p = r$ become
			\begin{equation}
				\label{MixedWeightLPLPNecessary}
				{\beta_1 > k - n \text{ and }} \beta_2 - \alpha_2 \leq \alpha_1 - \beta_1 \leq k \left( p - 1 \right).
			\end{equation}
			Thus, we ask the following question.
			\begin{question}
				\label{LPLPHyperbolicQuestion}
				Is the $k$-plane transform $R_k: L^p_{\alpha_1, \alpha_2} \left( \mathbb{H}^n \right) \rightarrow L^p_{\beta_1, \beta_2} \left( \Xi_k \left( \mathbb{H}^n \right) \right)$ a bounded operator, when ${p,} \alpha_1, \alpha_2, \beta_1, \beta_2 {\in \mathbb{R}}$ {are such that $\alpha_1 + \alpha_2 > k - n$, $1 \leq p < \frac{\alpha_1 + \alpha_2 + n - 1}{k - 1}$ and $\beta_1$ and $\beta_2$} satisfy the {conditions mentioned in Equation} \eqref{MixedWeightLPLPNecessary}?
			\end{question}
			In answering Question \ref{LPLPHyperbolicQuestion}, we obtain the boundedness for all functions (not just radial, which have been dealt with in Section \ref{RadialFunctionSection}) for a restricted range of $p$. Particularly, we have the following result.
			\begin{theorem}
				\label{LPLPBoundednessMixedWeightHn}
				Let $\alpha_1, \alpha_2 \in \mathbb{R}$ be such that $\alpha_1 + \alpha_2 > k - n$. Let $\beta_1, \beta_2 \in \mathbb{R}$ {satisfy the conditions of Equation \eqref{MixedWeightLPLPNecessary}}. {Then, for any} $1 \leq p \leq \frac{\alpha_1 + \alpha_2 + n + 1}{k + 1}$, {the $k$-plane transform is bounded from $L^p_{\alpha_1, \alpha_2} \left( \mathbb{H}^n \right)$ to $L^p_{\beta_1, \beta_2} \left( \Xi_k \left( \mathbb{H}^n \right) \right)$ if $\beta_1 + \beta_2 < p \left( k + 1 \right) - \left( n + 1 \right)$.}
			\end{theorem}
			\begin{proof}
				{To obtain an estimate on the $L^p_{\beta_1, \beta_2} \left( \Xi_k \left( \mathbb{H}^n \right) \right)$-norm of $R_kf$, we employ Equation \eqref{KPlaneHNGroup}. We have,
				\begin{align*}
					&\| R_kf \|_{L^p_{\beta_1, \beta_2} \left( \Xi_k \left( \mathbb{H}^n \right) \right)} \\
					&= \left[ \int\limits_{SO(n)} \int\limits_{0}^{\infty} \left| \int\limits_{\mathbb{S}^k_{\langle \omega, e_{k + 1} \rangle > \tanh h}} f_g \left( \omega, \tanh^{-1} \left( \frac{\tanh h}{\langle \omega, e_{k + 1} \rangle} \right) \right) \frac{\left( \coth^2 h \ \langle \omega, e_{k + 1} \rangle^2 - 1 \right)^{- \left( \frac{k + 1}{2} \right)}}{\sinh h} \right|^p  \times \right. \\
					&\left. \sinh^{\beta_1 + n - k - 1} h \cosh^{\beta_2 + k} h \ \mathrm{d}h \mathrm{d}g \right]^{\frac{1}{p}}.
				\end{align*}
				By using the Minkowski integral inequality, we get
				\begin{align*}
					&\| R_kf \|_{L^p_{\beta_1, \beta_2} \left( \Xi_k \left( \mathbb{H}^n \right) \right)} \\
					&\leq \int\limits_{\mathbb{S}^k_{\lambda \omega, e_{k + 1} \rangle > 0}} \left[ \int\limits_{SO(n)} \int\limits_{0}^{\tanh^{-1} \langle \omega, e_{k + 1} \rangle} \left| f_g \left( \omega, \tanh^{-1} \left( \frac{\tanh h}{\langle \omega, e_{k + 1} \rangle} \right) \right) \right|^p \left( \coth^2 h \ \langle \omega, e_{k + 1} \rangle^2 - 1 \right)^{- p\left( \frac{k + 1}{2} \right)} \times \right. \\
					&\left.  \sinh^{\beta_1 + n - k - 1 - p} h \cosh^{\beta_2 + k} h \ \mathrm{d}h \mathrm{d}g \right]^{\frac{1}{p}} \mathrm{d} \omega.
				\end{align*}
				Now, we substitute $\tanh h = \langle \omega, e_{k + 1} \rangle \tanh t$. After a simplification, we obtain,
				\begin{align*}
					&\| R_kf \|_{L^p_{\beta_1, \beta_2} \left( \Xi_k \left( \mathbb{H}^n \right) \right)} \\
					&\leq \int\limits_{\mathbb{S}^k_{\langle \omega, e_{k + 1} \rangle > 0}} \langle \omega, e_{k + 1} \rangle^{\frac{\beta_1 + n - k}{p} - 1} \times \\
					&\left[ \int\limits_{SO(n)} \int\limits_{0}^{\infty} \left| f_g \left( \omega, t \right) \right|^p \sinh^{\beta_1 + k \left( p - 1 \right) + n - 1} t \cosh^{p + k - \beta_1 - n - 1} t \left( 1 - \langle \omega, e_{k + 1} \rangle^2 \tanh^2 t \right)^{\frac{p - \beta_1 - \beta_2 - n - 1}{2}} \mathrm{d}t \mathrm{d}g \right]^{\frac{1}{p}} \mathrm{d}\omega \\
					&= \int\limits_{\mathbb{S}^k_{\langle \omega, e_{k + 1} \rangle > 0}} \langle \omega, e_{k + 1} \rangle^{\frac{\beta_1 + n - k}{p} - 1} \times \\
					&\left[ \int\limits_{SO(n)} \int\limits_{0}^{\infty} \left| f_g \left( \omega, t \right) \right|^p \sinh^{\alpha_1 + n - 1} t \cosh^{\alpha_2} t \tanh^{\beta_1 - \alpha_1 + k \left( p - 1 \right)} t \cosh^{p \left( k + 1 \right) - \left( \alpha_1 + \alpha_2 + n + 1 \right)} t \times \right. \\
					&\left. \left( 1 - \langle \omega, e_{k + 1} \rangle^2 \tanh^2 t \right)^{\frac{p - \beta_1 - \beta_2 - n - 1}{2}} \mathrm{d}t \mathrm{d}g \right]^{\frac{1}{p}} \mathrm{d}\omega.
				\end{align*}
				From Equation \eqref{MixedWeightLPLPNecessary}, we see that $\beta_1 - \alpha_1 + k \left( p - 1 \right) \geq 0$. Also, from the hypothesis of this theorem, we have $p \left( k + 1 \right) - \left( \alpha_1 + \alpha_2 + n + 1 \right) \leq 0$. Hence, we get,
				\begin{align*}
					\| R_kf \|_{L^p_{\beta_1, \beta_2} \left( \Xi_k \left( \mathbb{H}^n \right) \right)} &\leq C \int\limits_{\mathbb{S}^k_{\langle \omega, e_{k + 1} \rangle > 0}} \langle \omega, e_{k + 1} \rangle^{\frac{\beta_1 + n - k}{p} - 1} \times \\
					&\left[ \int\limits_{SO(n)} \int\limits_{0}^{\infty} \left| f_g \left( \omega, t \right) \right|^p \sinh^{\alpha_1 + n - 1} t \cosh^{\alpha_2} t \left( 1 - \langle \omega, e_{k + 1} \rangle^2 \tanh^2 t \right)^{\frac{p - \beta_1 - \beta_2 - n - 1}{2}} \mathrm{d}t \mathrm{d}g \right]^{\frac{1}{p}} \mathrm{d}\omega.
				\end{align*}
				We now have two cases.}
				
				{Case I: Let us assume that $\beta_1 + \beta_2 > p - n - 1$. That is, $p - \beta_1 - \beta_2 - n - 1 < 0$. Since $\tanh^2 t \leq 1$, we have $\left( 1 - \langle \omega, e_{k + 1} \rangle^2 \tanh^2 t \right)^{\frac{p - \beta_1 - \beta_2 - n - 1}{2}} \leq \left( 1 - \langle \omega, e_{k + 1} \rangle^2 \right)^{\frac{p - \beta_1 - \beta_2 - n - 1}{2}}$. Hence, we get,
				\begin{align*}
					\| R_kf \|_{L^p_{\beta_1, \beta_2} \left( \Xi_k \left( \mathbb{H}^n \right) \right)} &\leq C \int\limits_{\mathbb{S}^k_{\langle \omega, e_{k + 1} \rangle > 0}} \langle \omega, e_{k + 1} \rangle^{\frac{\beta_1 + n - k}{p} - 1} \left( 1 - \langle \omega, e_{k + 1} \rangle^2 \right)^{\frac{1}{2} - \frac{\beta_1 + \beta_2 + n + 1}{2p}} \times \\
					&\left[ \int\limits_{SO(n)} \int\limits_{0}^{\infty} \left| f_g \left( \omega, t \right) \right|^p \sinh^{\alpha_1 + n - 1} t \cosh^{\alpha_2} t \ \mathrm{d}t \mathrm{d}g \right]^{\frac{1}{p}} \mathrm{d}\omega \\
					&= C \| f \|_{L^p_{\alpha_1, \alpha_2} \left( \mathbb{H}^n \right)} \int\limits_{\mathbb{S}^k_{\langle \omega, e_{k + 1} \rangle > 0}} \langle \omega, e_{k + 1} \rangle^{\frac{\beta_1 + n - k}{p} - 1} \left( 1 - \langle \omega, e_{k + 1} \rangle^2 \right)^{\frac{1}{2} - \frac{\beta_1 + \beta_2 + n + 1}{2p}} \mathrm{d}\omega.
				\end{align*}
				Using the polar decomposition on $\mathbb{S}^k$, we have that,
				\begin{align*}
					\int\limits_{\mathbb{S}^k_{\langle \omega, e_{k + 1} \rangle > 0}} \langle \omega, e_{k + 1} \rangle^{\frac{\beta_1 + n - k}{p} - 1} \left( 1 - \langle \omega, e_{k + 1} \rangle^2 \right)^{\frac{1}{2} - \frac{\beta_1 + \beta_2 + n + 1}{2p}} \mathrm{d}\omega &= C \int\limits_{0}^{\frac{\pi}{2}} \cos^{\frac{\beta_1 + n - k}{p} - 1} \theta \sin^{k + 1 - \frac{\beta_1 + \beta_2 + n + 1}{p} - 1} \theta \ \mathrm{d}\theta < +\infty,
				\end{align*}
				since $\beta_1 > k - n$ and $\beta_1 + \beta_2 < p \left( k + 1 \right) - \left( n + 1 \right)$. This proves the boundedness of the $k$-plane transform in this case.}
				
				{Case II: We now consider $\beta_1 + \beta_2 \leq p - n - 1$. That is, $p - \beta_1 - \beta_2 - n - 1 \geq 0$. Here, we have,
				\begin{align*}
					\| R_kf \|_{L^p_{\beta_1, \beta_2} \left( \Xi_k \left( \mathbb{H}^n \right) \right)} &\leq C \int\limits_{\mathbb{S}^k_{\langle \omega, e_{k + 1} \rangle > 0}} \langle \omega, e_{k + 1} \rangle^{\frac{\beta_1 + n - k}{p} - 1} \left[ \int\limits_{SO(n)} \int\limits_{0}^{\infty} \left| f_g \left( \omega, t \right) \right|^p \sinh^{\alpha_1 + n - 1} t \cosh^{\alpha_2} t \ \mathrm{d}t \mathrm{d}g \right]^{\frac{1}{p}} \mathrm{d}\omega \\
					&= C \| f \|_{L^p_{\alpha_1, \alpha_2} \left( \mathbb{H}^n \right)} \int\limits_{\mathbb{S}^k_{\langle \omega, e_{k + 1} \rangle > 0}} \langle \omega, e_{k + 1} \rangle^{\frac{\beta_1 + n - k}{p} - 1} \mathrm{d}\omega \leq C \| f \|_{L^p_{\alpha_1, \alpha_2} \left( \mathbb{H}^n \right)},
				\end{align*}
				since $\beta_1 > k - n$. This completes the proof!}
			\end{proof}
		\subsection{The Sphere \texorpdfstring{$\mathbb{S}^n$}{}}
			In this section, we see the continuity of the $k$-plane transform $R_k$ from the weighted {Lebesgue} space $L^p_{\alpha_1, \alpha_2} \left( \mathbb{S}^n \right)$ to $L^p_{\beta_1, \beta_2} \left( \Xi_k \left( \mathbb{S}^n \right) \right)$. The necessary conditions {given in Equation \eqref{LPLRSNNecessaryConditions}} of Theorem \ref{NecessaryMixed} for $p = r$ become,
			\begin{equation}
				\label{NecessaryLPLPSnEquation}
				{\beta_1 > k - n,} \ \beta_1 \geq \alpha_1 - k \left( p - 1 \right) \text{ and } 	\beta_2 \geq \alpha_2 - k {\text{ with strict inequality if } p = 1 + \alpha_2}.
			\end{equation}
			Now, we ask the following question.
			\begin{question}
				\label{LPLPSnQuestion}
				{Let $\alpha_1, \alpha_2 \in \mathbb{R}$ and $p \geq \max \left\lbrace 1, 1 + \alpha_2 \right\rbrace$ with $p > 1 + \alpha_2$ when $\alpha_2 > 0$. Is the $k$-plane transform a bounded operator from $L^p_{\alpha_1, \alpha_2} \left( \mathbb{S}^n \right)$ to $L^p_{\beta_1, \beta_2} \left( \Xi_k \left( \mathbb{S}^n \right) \right)$, whenever $\beta_1, \beta_2 \in \mathbb{R}$ satisfy the conditions of Equation \eqref{NecessaryLPLPSnEquation}?}
			\end{question}
			We have the following boundedness result for the $k$-plane transform that completely answers Question \ref{LPLPSnQuestion}.
			\begin{theorem}
				\label{LPLPMixedWeightSphere}
				Let $\alpha_1, \alpha_2 \in \mathbb{R}$.
				\begin{enumerate}
					\item[\mylabel{LPLPSnWeightedA}{(A)}] For $p > 1 + \alpha_2$, the $k$-plane transform $R_k$ is a bounded linear operator from $L^p_{\alpha_1, \alpha_2} \left( \mathbb{S}^n \right)$ to $L^p_{\beta_1, \beta_2} \left( \Xi_k \left( \mathbb{S}^n \right) \right)$ {if and only if $\beta_1, \beta_2 \in \mathbb{R}$ satisfy the conditions of Equation \eqref{NecessaryLPLPSnEquation}.}.
					\item[\mylabel{LPLPSnWeightedB}{(B)}] For $\alpha_2 = 0$ and $p = 1$, the $k$-plane transform $R_k$ is a bounded linear operator from $L^1_{\alpha_1, 0} \left( \mathbb{S}^n \right)$ to $L^{{1}}_{\beta_1, \beta_2} \left( \Xi_k \left( \mathbb{S}^n \right) \right)$ {if and only if $\beta_1, \beta_2 \in \mathbb{R}$ satisfy the conditions of Equation \eqref{NecessaryLPLPSnEquation}.}.
				\end{enumerate} 
			\end{theorem}
			\begin{proof}~
				\begin{enumerate}
					\item[\mylabel{LPLPSnWeightedAProof}{(A)}] First, we look at the case when $p > 1 + \alpha_2$. From Equations \eqref{KPlaneSNGroup} and \eqref{MeasureXiKSN}, we get
					\begin{align*}
						\| R_kf \|_{L^p_{\beta_1, \beta_2} \left( \Xi_k \left( \mathbb{S}^n \right) \right)} = C \Bigg[ \int\limits_{SO(n)} \int\limits_{0}^{\frac{\pi}{2}} \Bigg| \int\limits_{\mathbb{S}^k_{\langle \omega, e_{k + 1} \rangle > 0}} &f_g \left( \omega, \tan^{-1} \left( \frac{\tan h}{\langle \omega, e_{k + 1} \rangle} \right) \right) \frac{\left( \cot^2 h \langle \omega, e_{k + 1} \rangle^2 + 1 \right)^{- \frac{k + 1}{2}}}{\sin h} \mathrm{d}\omega \Bigg|^p \times \\
						&\sin^{\beta_1 + n - k - 1} h \cos^{\beta_2 + k} h \ \mathrm{d}h \mathrm{d}g \Bigg]^{\frac{1}{p}}.
					\end{align*}
					By using Minkowski inequality, we have,
					\begin{align*}
						\| R_kf \|_{L^p_{\beta_1, \beta_2} \left( \Xi_k \left( \mathbb{S}^n \right) \right)} \leq C \int\limits_{\mathbb{S}^k_{{\langle} \omega, e_{k + 1} \rangle > 0}} \Bigg[ \int\limits_{SO(n)} \int\limits_{0}^{\frac{\pi}{2}} &\left| {f} \left( \omega, \tan^{-1} \left( \frac{\tan h}{\langle \omega, e_{k + 1} \rangle} \right) \right) \right|^p \left( \cot^2 h \langle \omega, e_{k + 1} \rangle^2 + 1 \right)^{- p \left( \frac{k + 1}{2} \right)} \times \\
						&\sin^{\beta_1 + n - k - 1 - p} h \cos^{\beta_2 + k} h \ \mathrm{d}h \mathrm{d}g \Bigg]^{\frac{1}{p}} \mathrm{d}\omega.
					\end{align*}
					We now substitute $\tan h = \langle \omega, e_{k + 1} \rangle \tan t$ to obtain,
					\begin{align*}
						\| R_kf \|_{L^p_{\beta_1, \beta_2} \left( \Xi_k \left( \mathbb{S}^n \right) \right)} \leq C \int\limits_{\mathbb{S}^k_{\langle \omega, e_{k + 1} \rangle > 0}} \Bigg[ \int\limits_{SO(n)} \int\limits_{0}^{\frac{\pi}{2}} &\left| f_g \left( \omega, t \right) \right|^p \left( \cot^2 t + 1 \right)^{- p \left( \frac{k + 1}{2} \right)} \times \\
						&\left( \frac{\langle \omega, e_{k + 1} \rangle \tan t}{\sqrt{1 + \langle \omega, e_{k + 1} \rangle^2 \tan^2 t}} \right)^{\beta_1 + n - k - 1 - p} \left( \frac{1}{\sqrt{1 + \langle \omega, e_{k + 1} \rangle^2 \tan^2 t}} \right)^{\beta_2 + k} \times \\
						&\frac{\langle \omega, e_{k + 1} \rangle \sec^2 t}{1 + \langle \omega, e_{k + 1} \rangle^2 \tan^2 t} \ \mathrm{d}t \mathrm{d}g \Bigg]^{\frac{1}{p}} \mathrm{d}\omega. \\
						= C \int\limits_{\mathbb{S}^k_{\langle \omega, e_{k + 1} \rangle > 0}} \langle \omega, e_{k + 1} \rangle^{\frac{\beta_1 + n - k}{p} - 1} \Bigg[ \int\limits_{SO(n)} \int\limits_{0}^{\frac{\pi}{2}} &\left| f_g \left( \omega, t \right) \right|^p \sin^{\beta_1 + k \left( p - 1 \right) + n - 1} t \cos^{p - \beta_1 - n + k - 1} t \times \\
						&\left( 1 + \langle \omega, e_{k + 1} \rangle^2 \tan^2 t \right)^{\frac{p - \beta_1 - \beta_2 - n - 1}{2}} \mathrm{d}t \mathrm{d}g \Bigg]^{\frac{1}{p}} \mathrm{d}\omega \\
						= C \int\limits_{\mathbb{S}^k_{\langle \omega, e_{k + 1} \rangle > 0}} \langle \omega, e_{k + 1} \rangle^{\frac{\beta_1 + n - k}{p} - 1} \Bigg[ \int\limits_{SO(n)} \int\limits_{0}^{\frac{\pi}{2}} &\left| f_g \left( \omega, t \right) \right|^p \sin^{\alpha_1 + n - 1} t \cos^{\alpha_2} t \sin^{\beta + k \left( p - 1 \right) - \alpha_1} t \times \\
						&\cos^{\beta_2 + k - \alpha_2} t \left( \cos^2 t + \langle \omega, e_{k + 1} \rangle^2 \sin^2 t \right)^{\frac{p - \beta_1 - \beta_2 - n - 1}{2}} \mathrm{d}t \mathrm{d}g \Bigg]^{\frac{1}{p}} \mathrm{d}\omega
					\end{align*}
					By the necessary conditions, we know that $\beta_1 + k \left( p - 1 \right) - \alpha_1 \geq 0$, and $\beta_2 + k - \alpha_2 \geq 0$. Hence, we have,
					\begin{align}
						\label{RkEstimateLPLPSn}
						\| R_kf \|_{L^p_{\beta_1, \beta_2} \left( \Xi_k \left( \mathbb{S}^n \right) \right)} \leq C \int\limits_{\mathbb{S}^k_{\langle \omega, e_{k + 1} \rangle > 0}} \langle \omega, e_{k + 1} \rangle^{\frac{\beta_1 + n - k}{p} - 1} \Bigg[ \int\limits_{SO(n)} \int\limits_{0}^{\frac{\pi}{2}} &\left| f_g \left( \omega, t \right) \right|^p \sin^{\alpha_1 + n - 1} t \cos^{\alpha_2} t \times \\
						&\left( \cos^2 t + \langle \omega, e_{k + 1} \rangle^2 \sin^2 t \right)^{\frac{p - \beta_1 - \beta_2 - n - 1}{2}} \mathrm{d}t \mathrm{d}g \Bigg]^{\frac{1}{p}} \mathrm{d}\omega. \nonumber
					\end{align}
					We notice that it is enough to show the boundedness for $\beta_1 + \beta_2 \leq p - n - 1$. This is because for any $\beta_2' > p - n - 1 - \beta_1$ (where, $\beta_1$ is fixed), we have
					$$\| R_kf \|_{L^p_{\beta_1, \beta_2'} \left( \Xi_k \left( \mathbb{S}^n \right) \right)} \leq \| R_kf \|_{{L^p_{\beta_1, p - n - 1 - \beta_1} \left( \Xi_k \left( \mathbb{S}^n \right) \right)}}.$$
					Under the above assumption on $\beta_2$, we have $p - \beta_1 - \beta_2 - n - 1 \geq 0$. Since $\sin t, \cos t$, and $\langle \omega, e_{k + 1} \rangle \leq 1$, we have $\left( \cos^2 t + \langle \omega, e_{k + 1} \rangle^2 \sin^2 t \right)^{\frac{p - \beta_1 - \beta_2 - n - 1}{2}} \leq C$. Hence, from Equation \eqref{RkEstimateLPLPSn}, we get,
					\begin{align*}
						\| R_kf \|_{L^p_{\beta_1, \beta_2} \left( \Xi_k \left( \mathbb{S}^n \right) \right)} &\leq C \int\limits_{\mathbb{S}^k_{\langle \omega, e_{k + 1} \rangle > 0}} \langle \omega, e_{k + 1} \rangle^{\frac{\beta_1 + n - k}{p} - 1} \left[ \int\limits_{SO(n)} \int\limits_{0}^{\frac{\pi}{2}} \left| f_g \left( \omega, t \right) \right|^p \sin^{\alpha + n - 1} t \cos^{\alpha_2} t \ \mathrm{d}t \mathrm{d}g \right]^{\frac{1}{p}} \mathrm{d}\omega \\
						&= C \| f \|_{L^p_{\alpha_1, \alpha_2} \left( \mathbb{S}^n \right)} \int\limits_{\mathbb{S}^k_{\langle \omega, e_{k + 1} \rangle > 0}} \langle \omega, e_{k + 1} \rangle^{\frac{\beta_1 + n - k}{p} - 1} \mathrm{d}\omega.
					\end{align*}
					We notice that by using polar decomposition on $\mathbb{S}^k_{\langle \omega, e_{k + 1} \rangle > 0}$ with center $e_{k + 1}$, we get
					\begin{align*}
						\int\limits_{\mathbb{S}^k_{\langle \omega, e_{k + 1} \rangle > 0}} \langle \omega, e_{k + 1} \rangle^{\frac{\beta_1 + n - k}{p} - 1} \mathrm{d}\omega &= C \int\limits_{0}^{\frac{\pi}{2}} \cos^{\frac{\beta_1 + n - k}{p}} {\theta} \sin^{k - 1} {\theta} \ \mathrm{d}{\theta} < + \infty,
					\end{align*}
					since $\beta_1 > k - n$. {This completes the proof of \ref{LPLPSnWeightedAProof}.}
					\item[\mylabel{LPLPSnWeightedBProof}{(B)}] When $\alpha_2 = 0$ and $p = 1$, we use a duality argument. First, we notice that for any $f \in L^1_{\alpha_1, 0} \left( \mathbb{S}^n \right)$, we have $\left| R_kf \left( \xi \right) \right| \leq R_k\left| f \right| \left( \xi \right)$, for any $\xi \in \Xi_k \left( \mathbb{S}^n \right)$. Now, using the duality of the $k$-plane transform given in Equation \eqref{DualityEquation} and Equation \eqref{DualMixedSphereEquation}, we have,
					\begin{align*}
						\int\limits_{\Xi_k \left( \mathbb{S}^n \right)} \left| R_kf \left( \xi \right) \right| \sin^{\beta_1} d \left( 0, \xi \right) \cos^{\beta_2} d \left( 0, \xi \right) \mathrm{d}\xi &\leq \int\limits_{\Xi_k \left( \mathbb{S}^n \right)} R_k \left| f \right| \left( \xi \right) \sin^{\beta_1} d \left( 0, \xi \right) \cos^{\beta_2} d \left( 0, \xi \right) \mathrm{d}\xi \\
						&= C \int\limits_{\mathbb{S}^n_{{+}}} \left| f \left( x \right) \right| \sin^{\beta_1} d \left( 0, x \right) \cos^{\beta_2} d \left( 0, x \right) {}_2F_1 \left( - \frac{\beta_2}{2}, \frac{k}{2}; \frac{\beta_1 + n}{2}; - \tan^2 d \left( 0, x \right) \right) \mathrm{d}x.
					\end{align*}
					Using Equation \eqref{Transformation2F1}, we get
					\begin{align*}
						&\int\limits_{\Xi_k \left( \mathbb{S}^n \right)} \left| R_kf \left( \xi \right) \right| \sin^{\beta_1} d \left( 0, \xi \right) \cos^{\beta_2} d \left( 0, \xi \right) \mathrm{d}\xi \\
						&\leq C \int\limits_{\mathbb{S}^n_{{+}}} \left| f \left( x \right) \right| \sin^{\alpha_1} d \left( 0, x \right) \sin^{\beta_1 - \alpha_1} d \left( 0, x \right) {}_2F_1 \left( - \frac{\beta_2}{2}, \frac{\beta_1 + n - k}{2}; \frac{\beta_1 + n}{2}; \sin^2 d \left( 0, x \right) \right) \mathrm{d}x.
					\end{align*}
					From Equation \eqref{NecessaryLPLPSnEquation}, we have $\beta_1 \geq \alpha_1$, since $p = 1$, {and} $\beta_2 > {-} k$. Consequently, {$\sin^{\beta_1 - \alpha_2} d \left( 0, x \right) \leq 1$ and} from Equation \eqref{Behaviour2F1Case1} of Theorem \ref{Behaviour2F1}, we conclude that the hypergeometric function is bounded. Hence,
					$$\| R_kf \|_{L^1_{\beta_1, \beta_2} \left( \Xi_k \left( \mathbb{S}^n \right) \right)} = \int\limits_{\Xi_k \left( \mathbb{S}^n \right)} \left| R_kf \left( \xi \right) \right| \sin^{\beta_1} d \left( 0, \xi \right) \cos^{\beta_2} d \left( 0, \xi \right) \mathrm{d}\xi \leq C \int\limits_{\mathbb{S}^n} \left| f \left( x \right) \right| \sin^{\alpha_1} d \left( 0, x \right) \mathrm{d}x = \| f \|_{L^1_{\alpha_1, 0} \left( \mathbb{S}^n \right)}.$$
				\end{enumerate}
			\end{proof}
		
		\begin{remark}
			\normalfont
			It is clear from Theorem \ref{LPLPMixedWeightSphere}, by using $\alpha_1 = \alpha_2 = 0$, that the $k$-plane transform is a bounded operator from $L^p \left( \mathbb{S}^n \right)$ to $L^p \left( \Xi_k \left( \mathbb{S}^n \right) \right)$. This (unweighted) result is well-known (see for instance, \cite{RubinInversionSphere}).
		\end{remark}
	\section{Conclusion}
		\label{ConclusionSection}
		In this article, we have seen weighted $L^p$-improving mapping properties of the $k$-plane transform on constant curvature spaces. The technique used extensively throughout this article is the geodesic correspondence between $\mathbb{R}^n$ and the three spaces of constant curvature. A few questions remain open, and the techniques to answer them elude us. We present {these questions} here.
		\begin{openquestion}
			\label{OQLPDecreasing}
			{Do we have weighted $L^p$-$L^r$ boundedness of the $k$-plane transform on spaces of constant curvature, when $r < p$?}
		\end{openquestion}
		In the Euclidean case, the weighted norm of the $k$-plane transform can be expressed as a norm of certain convolution operator (see \cite{KumarRayWE}). Since convolutions are translation-invariant, we have that the $k$-plane transform is $L^p$-improving (see for instance, \cite{Hormander}). However, we observe that in our work we \textit{\textbf{dominated}} the weighted norm of the $k$-plane transform by certain convolution operator. Therefore, the necessity of the k-plane transform being $L^p$-improving is not immediate. We now try to partially answer the above question.
		
		First, let us observe that for $s > r$ and $\beta_2' > \frac{s}{r} \beta_2 + \left( \beta_1 + n - 1 \right) \left( 1 - \frac{r}{s} \right)$, we have $L^s_{\beta_1, \beta_2'} \left( \Xi_k \left( \mathbb{H}^n \right) \right) \subseteq L^r_{\beta_1, \beta_2} \left( \Xi_k \left( \mathbb{H}^n \right) \right)$. In fact, we have,
		\begin{equation}
			\label{ContainmentInequality}
			\| \varphi \|_{L^r_{\beta_1, \beta_2} \left( \Xi_k \left( \mathbb{H}^n \right) \right)} \leq C \| \varphi \|_{L^s_{\beta_1, \beta_2'} \left( \Xi_k \left( \mathbb{H}^n \right) \right)}.
		\end{equation}
		 Such a containment can be easily seen by an application of H\"{o}lder's inequality with exponents $\frac{s}{r}$ and $\frac{s}{s - r}$. Let us now recall from Theorem \ref{MixedWeightLPLRFullRadial} that we have the following inequality for all radial functions in $L^p_{\alpha_1, \alpha_2} \left( \mathbb{H}^n \right)$.
		\begin{equation}
			\label{LPLPRadialKnown}
			\| R_kf \|_{L^p_{\beta_1, \beta_2} \left( \Xi_k \left( \mathbb{H}^n \right) \right)} \leq C \| f \|_{L^p_{\alpha_1, \alpha_2} \left( \mathbb{H}^n \right)},
		\end{equation}
		where, $\beta_1, \beta_2$ satisfy the conditions of Equation \eqref{NecessaryMixedWeightInequality} with $r = p$. Now, it is easily seen that for $r < p$ and $\beta_2' < \frac{r}{p} \beta_2 - \left( \beta_1 + n - 1 \right) \left( 1 - \frac{r}{p} \right)$, the triple $\left( \beta_1, \beta_2', r \right)$ satisfies Equation \eqref{NecessaryMixedWeightInequality}. Also, from Inequalities \eqref{ContainmentInequality} and \eqref{LPLPRadialKnown}, we get
		\begin{equation}
			\label{LPDecreasingEstimate}
			\| R_kf \|_{L^r_{\beta_1, \beta_2'} \left( \Xi_k \left( \mathbb{H}^n \right) \right)} \leq C \| f \|_{L^p_{\alpha_1, \alpha_2} \left( \mathbb{H}^n \right)},
		\end{equation}
		for all radial functions $f \in L^p_{\alpha_1, \alpha_2} \left( \mathbb{H}^n \right)$. That is, under certain assumptions on the weight $\beta_2$, we can have $L^p$-decreasing estimates for the $k$-plane transform on $\mathbb{H}^n$.
		
		Naturally, we ask if for all admissible values of $\beta_1$ and $\beta_2$ we have $L^p$-decreasing estimates. However, this is answered negatively in the following series of results that are motivated from H\"{o}rmander's work (see \cite{Hormander}). In what follows, we consider a measure space $\left( X, \mu \right)$, where $X$ is a topological space and $\mu$ is a Radon measure on $X$. Let $T: L^p \left( X \right) \rightarrow L^q \left( X \right)$ is a bounded linear operator such that for every $f \in L^p \left( X \right)$, we have $\left| Tf \left( x \right) \right| \leq T \left| f \right| \left( x \right)$, for almost every $x \in X$. Then, it is easy to see that the operator norm $\| T \|$ is given by $\sup \left\lbrace \frac{\| Tf \|_{L^q \left( X \right)}}{\| f \|_{L^p \left( X \right)}} \Bigg| f \in L^p \left( X \right) \text{ with } f \geq 0 \right\rbrace$.
		We now prove the following lemma.
		\begin{lemma}
			\label{TranslationEstimateLPImproving}
			Let $\gamma < 0$ and $L^p_{\gamma} \left( \mathbb{R} \right) := L^p \left( \mathbb{R}, \left( 1 + e^x \right)^{\gamma} \mathrm{d}x \right)$. Then, for any $f \in L^p_{\gamma} \left( \mathbb{R} \right)$ with $f \geq 0$, we have
			\begin{equation}
				\label{TranslationEstimateInequality}
				\liminf\limits_{h \rightarrow - \infty} \| f + \tau_hf \|_{L^p_{\gamma} \left( \mathbb{R} \right)} \geq 2^{\frac{1}{p}} \| f \|_{L^p_{\gamma} \left( \mathbb{R} \right)}.
			\end{equation}
			Here, $\tau_hf \left( x \right) = f \left( x - h \right)$.
		\end{lemma}
		\begin{proof}
			Let us start with a compactly supported $g \in L^p_{\gamma} \left( \mathbb{R} \right)$. Then, there is some $h_0 < 0$ such that for all $h \leq h_0$, we have $\text{supp } g \cap \text{supp } \tau_hg = \emptyset$. Now, for any $h \leq h_0$, we have,
			\begin{align*}
				\| g + \tau_h g \|_{L^p_{\gamma} \left( \mathbb{R} \right)} &= \left[ \int\limits_{\mathbb{R}} \left| g \left( x \right) + \tau_h g \left( x \right) \right|^p \left( 1 + e^x \right)^{\gamma} \mathrm{d}x \right]^{\frac{1}{p}} \\
				&= \left[ \int\limits_{\mathbb{R}} \left| g \left( x \right) \right|^p \left( 1 + e^x \right)^{\gamma} \mathrm{d}x + \int\limits_{\mathbb{R}} \left| g \left( x \right) \right|^p \left( 1 + e^{x + h} \right)^{\gamma} \mathrm{d}x \right]^{\frac{1}{p}} \geq 2^{\frac{1}{p}} \| g \|_{L^p_{\gamma} \left( \mathbb{R} \right)},
			\end{align*}
			since $h < 0$ and $\gamma < 0$. Therefore, for any compactly supported function $g \in L^p_{\gamma} \left( \mathbb{R} \right)$, we have,
			$$\liminf\limits_{h \rightarrow - \infty} \| g + \tau_h g \|_{L^p_{\gamma} \left( \mathbb{R} \right)} \geq 2^{\frac{1}{p}} \| g \|_{L^p_{\gamma} \left( \mathbb{R} \right)}.$$
			Now, let $f \in L^p_{\gamma} \left( \mathbb{R} \right)$ be a non-negative function. By the density of $C_c \left( \mathbb{R} \right)$ in $L^p_{\gamma} \left( \mathbb{R} \right)$, we get a sequence $\left( g_n \right)_{n \in \mathbb{N}}$ of compactly supported functions in $L^p_{\gamma} \left( \mathbb{R} \right)$ such that $\| g_n - f \|_{L^p_{\gamma} \left( \mathbb{R} \right)} \rightarrow 0$ as $n \rightarrow \infty$. In fact, we may choose $0 \leq g_n \leq f$, for every $n \in \mathbb{N}$. Therefore, it is readily seen that for every $n \in \mathbb{N}$, we have,
			$$\liminf\limits_{h \rightarrow -\infty} \| f + \tau_h f \|_{L^p_{\gamma} \left( \mathbb{R} \right)} \geq \liminf\limits_{h \rightarrow - \infty} \| g_n + \tau_h g_n \|_{L^p_{\gamma} \left( \mathbb{R} \right)} \geq 2^{\frac{1}{p}} \| g_n \|_{L^p_{\gamma} \left( \mathbb{R} \right)}.$$
			Since $g_n \rightarrow f$ in $L^p_{\gamma}$-norm, the result is proved!
		\end{proof}
		We now prove a result about translation-invariant operators when the codomain is equipped with a certain weight.
		\begin{theorem}
			\label{LPImprovingWeighted}
			Let $T: L^p \left( \mathbb{R} \right) \rightarrow L^p_{\gamma} \left( \mathbb{R} \right)$ be a translation-invariant bounded linear operator, where $L^p_{\gamma} \left( \mathbb{R} \right)$ is as defined in Lemma \ref{TranslationEstimateLPImproving}. Assume that for any $f \in L^p \left( \mathbb{R} \right)$ and almost every $x \in \mathbb{R}$, we have $\left| Tf \left( x \right) \right| \leq T \left| f \right| \left( x \right)$. Also, assume that whenever $f \geq 0$, we have $Tf \geq 0$. Then, $q \geq p$.
		\end{theorem}
		\begin{proof}
			Let $f \geq 0$ and $f \in L^p \left( \mathbb{R} \right)$. Then, by our assumption, we have $Tf \geq 0$. Since $T$ is translation-invariant, we get for any $h \in \mathbb{R}$, $T \left( f + \tau_h f \right) = Tf + \tau_h Tf \geq 0$. Therefore, from Lemma \ref{TranslationEstimateLPImproving}, we get
			$$2^{\frac{1}{p}} \| Tf \|_{L^q_{\gamma} \left( \mathbb{R} \right)} \leq \liminf\limits_{h \rightarrow - \infty} \| T \left( f + \tau_h f \right) \|_{L^q_{\gamma} \left( \mathbb{R} \right)} \leq \| T \| \liminf\limits_{h \rightarrow - \infty} \| f + \tau_h f \|_{L^p \left( \mathbb{R} \right)} = 2^{\frac{1}{p}} \| T \| \| f \|_{L^p \left( \mathbb{R} \right)},$$
			where, the last equality follows from H\"{o}rmander's result in \cite{Hormander}. It is now readily seen that $\frac{1}{q} \leq \frac{1}{p}$, i.e., $q \geq p$.
		\end{proof}
		We now show that under certain conditions, the $k$-plane transform is an $L^p$-improving operator.
		\begin{corollary}
			\label{KPlaneTransformLPImproving}
			Let $\alpha_1 + \alpha_2 > k - n$ and $1 \leq p < \frac{\alpha_1 + \alpha_2 + n - 1}{k - 1}$. Let $r \geq 1$, and $\beta_1, \beta_2 \in \mathbb{R}$ be such that $\frac{\beta_2 + k - 1}{r} - \frac{\alpha_2 - 1}{p} \leq \frac{\alpha_1 + n}{p} - \frac{\beta_1 + n - k}{r} = k$. If $R_k: L^p_{\alpha_1, \alpha_2} \left( \mathbb{H}^n \right) \rightarrow L^r_{\beta_1, \beta_2} \left( \Xi_k \left( \mathbb{H}^n \right) \right)$ is bounded, then $r \geq p$.
		\end{corollary}
		\begin{proof}
			First, we observe that
			\begin{equation}
				\label{KPlaneConvolutionEquation}
				\begin{aligned}
					&\| R_kf \|_{L^r_{\beta_1, \beta_2} \left( \Xi_k \left( \mathbb{H}^n \right) \right)} \\
					&= \left[ \int\limits_{0}^{\infty} \left| \int\limits_{u}^{\infty} \tilde{f} \left( \sqrt{v + 1} \right) v^{\frac{\alpha_1 + n}{2p}} \left( 1 + v \right)^{\frac{\alpha_2 - 1}{2p}} \left( \frac{u}{v} \right)^{\frac{\beta_1 + n - k}{2r}} \left( 1 - \frac{u}{v} \right)^{\frac{k}{2} - 1} \left( \frac{1 + u}{1 + v} \right)^{\frac{1}{2} + \frac{\alpha_2 - 1}{2p}} \frac{\mathrm{d}v}{v} \right|^r \left( 1 + u \right)^{r \left( \frac{\beta_2 + k - 1}{2r} - \frac{\alpha_2 - 1}{2p} - \frac{k}{2} \right)} \frac{\mathrm{d}u}{u} \right]^{\frac{1}{r}}.
				\end{aligned}
			\end{equation}
			We observe that
			$$\left( \frac{1 + u}{1 + v} \right)^{\frac{1}{2} + \frac{\alpha_2 - 1}{2p}} \geq \begin{cases}
			1, & \text{when } 1 + \frac{\alpha_2 - 1}{p} < 0. \\
			\left( \frac{u}{v} \right)^{\frac{1}{2} + \frac{\alpha_2 - 1}{2p}}, & \text{when } 1 + \frac{\alpha_2 - 1}{p} \geq 0.
		\end{cases}$$
			Therefore, we have,
			\begin{equation}
				\label{KPlaneConvolutionEquationDomination}
				\| R_kf \|_{L^r_{\beta_1, \beta_2} \left( \Xi_k \left( \mathbb{H}^n \right) \right)} \geq \left[ \int\limits_{0}^{\infty} \left| F \left( v \right) G_{p, \alpha_2} \left( \frac{u}{v} \right) \right|^r \left( 1 + u \right)^{r \left( \frac{\beta_2 + k - 1}{2r} - \frac{\alpha_2 - 1}{2p} - \frac{k}{2} \right)} \frac{\mathrm{d}u}{u} \right]^{\frac{1}{r}},
			\end{equation}
			where,
			\begin{equation}
				\label{ConvolutionFEquation}
				F \left( v \right) = \tilde{f} \left( \sqrt{v + 1} \right) v^{\frac{\alpha_1 + n}{2p}} \left( 1 + v \right)^{\frac{\alpha_2 - 1}{2p}},
			\end{equation}
			and,
			\begin{equation}
				\label{ConvolutionGEquation}
				G \left( t \right) = \begin{cases}
											t^{\frac{\beta_1 + n - k}{2r}} \left( 1 - t \right)^{\frac{k}{2} - 1} \chi_{\left( 0, 1 \right)} \left( t \right), & \text{when } 1 + \frac{\alpha_2 - 1}{p} < 0. \\
											t^{\frac{\beta_1 + n - k}{2r} + \frac{1}{2} + \frac{\alpha_2 - 1}{2p}} \left( 1 - t \right)^{\frac{k}{2} - 1} \chi_{\left( 0, 1 \right)} \left( t \right), & \text{when } 1 + \frac{\alpha_2 - 1}{p} \geq 0.
										\end{cases}
			\end{equation}
			Consider the operator $T$ given by $T\tilde{F} \left( x \right) = \left( \tilde{F} * \tilde{G} \right) \left( x \right)$, for $\tilde{F} \left( x \right) = F \left( e^x \right)$ and $\tilde{G} \left( x \right) = G \left( e^x \right)$, for $x \in \mathbb{R}$. From Equation \eqref{KPlaneConvolutionEquationDomination}, it is clear that the boundedness of the $k$-plane transform implies the boundedness of the operator $T: L^p \left( \mathbb{R} \right) \rightarrow L^r_{\gamma} \left( \mathbb{R} \right)$, where $\gamma = r \left( \frac{\beta_2 + k - 1}{2r} - \frac{\alpha_2 - 1}{2p} - \frac{k}{2} \right) \leq 0$. From Theorem \ref{LPImprovingWeighted}, we must have $r \geq p$.
		\end{proof}
		Thus, we have proved, under certain special cases, that the $k$-plane transform is necessarily an $L^p$-improving operator. However, in some other cases, we have also obtained the $L^p$-decreasing nature of the $k$-plane transform on $\mathbb{H}^n$. We remark here that the complete understanding of Open Question \ref{OQLPDecreasing} needs to be explored.
		
		We now end the article with some questions that arise naturally from our study.
		\begin{openquestion}
			\label{OQ1}
			\normalfont
			Can Theorem \ref{LPLPBoundednessMixedWeightHn} be improved to include the full range of $p$? That is, does the weighted $L^p$-$L^p$ boundedness hold for the case $\frac{\alpha + n + 1}{k + 1} < p < \frac{\alpha + n - 1}{k - 1}$, {whenever $\beta_1, \beta_2 \in \mathbb{R}$ satisfy the conditions of Equation \eqref{MixedWeightLPLPNecessary}}?
		\end{openquestion}
		The problem with the current technique is that there are extra factors (particularly $k$-many) of hyperbolic cosine that appear in the weighted $L^p$-norm of $R_kf$. These factors cannot be eliminated easily.
		\begin{openquestion}
			\label{OQ2}
			\normalfont
			What {are} the admissible (necessary and sufficient) values of $\alpha_1, \alpha_2, \beta_1, \beta_2 \in \mathbb{R}$ such that, we have for every $f: \mathbb{S}^n \rightarrow \mathbb{C}$ with $\| \sin^{\alpha_1} d \left( 0, \cdot \right) \cos^{\alpha_2} d \left( 0, \cdot \right) f \|_{\infty} < + \infty$,
			$$\| \sin^{\beta_1} d \left( 0, \cdot \right) \cos^{\beta_2} d \left( 0, \cdot \right) R_kf \|_{\infty} \leq C \| \sin^{\alpha_1} d \left( 0, \cdot \right) \cos^{\alpha_2} d \left( 0, \cdot \right) f \|_{\infty}?$$
		\end{openquestion}
		To even get the necessary conditions, it would appear that the natural choice of the function $f: \mathbb{S}^n \rightarrow \mathbb{C}$ would be $f \left( x \right) = \chi_{B \left( 0, \lambda \right)} \left( x \right) \sin^{- \alpha_1} d \left( 0, x \right) \cos^{- \alpha_2} d \left( 0, x \right)$. However, with this, the $k$-plane transform is evaluated in terms of Appell's function $F_1$ (see \cite{Olver}), whose behaviour at ${\left( 1, 1 \right)}$ is not well studied. This question, therefore, also warrants a study of the asymptotic behaviour of the Appell's functions near ${\left( 1, 1 \right)}$ and ${\left( - \infty, - \infty \right)}$, which are the two ``end-points" for the variables {of Appell's functions}.
		\begin{openquestion}
			\label{OQ3}
			\normalfont
			What are the necessary and sufficient conditions for the $k$-plane transform to be a bounded operator from $L^p_{\alpha_1, \alpha_2} \left( X \right)$ to $L^q_{\beta_1, \beta_2} \left( {\Xi_k \left( X \right)} \right)$, for a space of constant curvature $X${?}
		\end{openquestion}
		{In this article, we have addressed Open Question \ref{OQ3} for radial functions. However, currently, we are not in a position to give an (affirmative) answer to Open Question \ref{OQ3}.} Such conditions are well-studied in the Euclidean case, some of which we have mentioned in Section \ref{IntroductionSection}. However, obtaining the necessary conditions in the Euclidean case often involves a ``stretching" of a cuboid (or a cylinder), which seems difficult to mimic in non-Euclidean spaces. Therefore, it would be interesting to study various shapes in the hyperbolic space and the sphere, from where one can try and obtain the necessary conditions in solving Open Question \ref{OQ3}. Also, the techniques of Fourier transform that are used in the Euclidean case to obtain $L^p$-$L^q$ estimates do not work well in the non-Euclidean spaces, especially due to the lack of a ``slice theorem".
		
		Apart from this, one may also try to consider the Sobolev boosting estimates for the $k$-plane transform. We could not find much work done in this direction. Particularly, the following question can be asked.
		\begin{openquestion}
			\label{OQ4}
			\normalfont
			{For $p \geq 1$, what } is the best possible regularity {$s > 0$ of the Sobolev space} $W^{s, p} \left( \Xi_k \left( X \right) \right)$ such that the $k$-plane transform {takes $L^p \left( X \right)$ to $W^{s, p} \left( \Xi_k \left( X \right) \right)$ continuously?} This question can also be asked by putting various weights on the measure of $X$ {and $\Xi_k \left( X \right)$}. Here, the Sobolev spaces are considered over the fibres of $\Xi_k \left( X \right)$.
		\end{openquestion}
		
		We end this article with the following question. The motivation to ask this question comes from our study of end-point estimates on Lorentz spaces.
		\begin{openquestion}
			\label{OQ5}
			\normalfont
			What are the best possible values of $p_1, p_2, q_1, q_2 \geq 1$ such that the $k$-plane transform is a bounded operator between the Lorentz spaces $L^{p_1, q_1} \left( X \right)$ and $L^{p_2, q_2} \left( \Xi_k \left( X \right) \right)$? Again, this question can also be asked by adding weights to the measures of $X$ and $\Xi_k \left( X \right)$.
		\end{openquestion}
		The authors of this article could not find literature related to Open Question \ref{OQ5}. This seems to be a new direction of research in the study of $k$-plane transform on Euclidean and non-Euclidean spaces.

	\appendix
	\section{The \texorpdfstring{$k$}{k}-plane transform and its dual of certain functions}
		\label{EvaluationAppendix}
		We now give some explicit computations for evaluating the $k$-plane transform and its dual of certain functions. These computations were used at various places in our work.
		\begin{proposition}
			\label{LPMassBallRn}
			Let $\lambda > 0$ and $B \left( 0, \lambda \right)$ be the ball of radius $\lambda$ centered at $0$ in $\mathbb{R}^n$. Then,
			\begin{equation}
				\label{LPMassBallRnEquation}
				R_k\chi_{B \left( 0, \lambda \right)} = C \left( \lambda^2 - \left( d \left( 0, \xi \right) \right)^2 \right)^{\frac{k}{2}} \chi_{\left( 0, \lambda \right)} \left( d \left( 0, \xi \right) \right),
			\end{equation}
			where, the constant $C > 0$ depends only on $n$ and $k$.
		\end{proposition}
		\begin{proof}
			We have using Equation \eqref{KPlaneTransformRadialRnEquation},
			\begin{align*}
				R_kf \left( \xi \right) &= C \int\limits_{d \left( 0, \xi \right)}^{\infty} \chi_{\left( 0, \lambda \right)} \left( t \right) \left( t^2 - \left( d \left( 0, \xi \right) \right)^2 \right)^{\frac{k}{2} - 1} t \ \mathrm{d}t = C \left( \lambda^2 - \left( d \left( 0, \xi \right) \right)^2 \right)^{\frac{k}{2}} \chi_{\left( 0, \lambda \right)} \left( d \left( 0, \xi \right) \right).
			\end{align*}
		\end{proof}
		We now move to {calculate} the $k$-plane transform of certain ``nice" functions on $\mathbb{H}^n$.
		\begin{proposition}
			\label{KPlaneBallHn}
			Let $\lambda > 0$ and $B \left( 0, \lambda \right)$ be the ball of radius $\lambda$ centered at $0$ in $\mathbb{H}^n$. Then,
			\begin{equation}
				\label{KPlaneballHnEquation}
				\begin{aligned}
					&R_k\chi_{B \left( 0, \lambda \right)} \left( \xi \right) \\
					&= C \tanh^k \lambda \left( 1 - \frac{\sinh^2 d \left( 0, \xi \right)}{\sinh^2 \lambda} \right)^{\frac{k}{2}} {}_2F_1 \left( \frac{k + 1}{2}, \frac{k}{2}; 1 + \frac{k}{2}; \frac{\tanh^2 d \left( 0, \xi \right) \left( \sinh^2 \lambda - \sinh^2 d \left( 0, \xi \right) \right)}{\tanh^2 d \left( 0, \xi \right) \left( \sinh^2 \lambda - \sinh^2 d \left( 0, \xi \right) \right) + \sinh^2 d \left( 0, \xi \right)} \right),
				\end{aligned}
			\end{equation}
			for $d \left( 0, \xi \right) \leq \lambda$. For $d \left( 0, \xi \right) > \lambda$, we have $R_k\chi_{B \left( 0, \lambda \right)} \left( \xi \right) = 0$.
		\end{proposition}
		\begin{proof}
			We have from Equation \eqref{KPlaneTransformRadialHnEquation3},
			\begin{align*}
				R_k\chi_{B \left( 0, \lambda \right)} \left( \xi \right) &= \frac{C}{\cosh^{k - 1} d \left( 0, \xi \right)} \int\limits_{d \left( 0, \xi \right)}^{\lambda} \left[ 1 - \frac{\sinh^2 d \left( 0, \xi \right)}{\sinh^2 t} \right]^{\frac{k}{2} - 1} \sinh^{k - 1} t \ \mathrm{d}t.
			\end{align*}
			It is clear that for $d \left( 0, \xi \right) \geq \lambda$, we have $R_k\chi_{B \left( 0, \lambda \right)} \left( \xi \right) = 0$.
			
			For $d \left( 0, \xi \right) < \lambda$, let us substitute $\sinh^2 t = \left( 1 + y \left( \frac{\sinh^2 \lambda}{\sinh^2 d \left( 0, \xi \right)} - 1 \right) \right) \sinh^2 d \left( 0, \xi \right)$. After simplification, we get
			\begin{align*}
				R_kf \left( \xi \right) &= C \tanh^k d \left( 0, \xi \right) \left( \frac{\sinh^2 \lambda}{\sinh^2 d \left( 0, \xi \right)} - 1 \right)^{\frac{k}{2}} \int\limits_{0}^{1} y^{\frac{k}{2} - 1} \left( 1 + \tanh^2 d \left( 0, \xi \right) \left( \frac{\sinh^2 \lambda}{\sinh^2 d \left( 0, \xi \right)} {- 1} \right) y \right)^{- \frac{1}{2}} \mathrm{d}y \\
				&= C \tanh^k d \left( 0, \xi \right) \left( \frac{\sinh^2 \lambda}{\sinh^2 d \left( 0, \xi \right)} - 1 \right)^{\frac{k}{2}} {}_2F_1 \left( \frac{1}{2}, \frac{k}{2}; 1 + \frac{k}{2}; \tanh^2 d \left( 0, \xi \right) \left( 1 - \frac{\sinh^2 \lambda}{\sinh^2 d \left( 0, \xi \right)} \right) \right) {.}
			\end{align*}
			In the last equality, we have used the integral form of the hypergeometric function given in Equation \eqref{IntegralForm2F1}. The result now follows by using Equation \eqref{Transformation2F1} \eqref{2F1Transformation2}.
		\end{proof}
		\begin{proposition}
			\label{KPlaneBallCoshHn}
			Let $\lambda > 0$ and $B \left( 0, \lambda \right)$ be the ball of radius $\lambda$ centered at $0$ in $\mathbb{H}^n$. Consider the function $f_{\lambda}: \mathbb{H}^n \rightarrow \mathbb{C}$, defined as $f_{\lambda} \left( x \right) = \chi_{B \left( 0, \lambda \right)} \left( x \right) \cosh d \left( 0, x \right)$. Then,
			\begin{equation}
				\label{KPlaneBallCoshHnEquation}
				R_kf_{\lambda} {\left( \xi \right)} = C \frac{\sinh^k \lambda}{\cosh^{k - 1} d \left( 0, \xi \right)} \left( 1 - \frac{\sinh^2 d \left( 0, \xi \right)}{\sinh^2 \lambda} \right)^{\frac{k}{2}} \chi_{\left( 0, \lambda \right)} \left( d \left( 0, \xi \right) \right).
			\end{equation}
		\end{proposition}
		\begin{proof}
			Using Equation \eqref{KPlaneTransformRadialHnEquation}, we have, for $d \left( 0, \xi \right) < \lambda$,
			\begin{align*}
				R_kf_{\lambda} \left( \xi \right) = \frac{C}{\cosh^{k - 1} d \left( 0, \xi \right)} \int\limits_{d \left( 0, \xi \right)}^{\lambda} \cosh t \left( \cosh^2 t - \cosh^2 d \left( 0, \xi \right) \right)^{\frac{k}{2} - 1} \sinh t \ \mathrm{d}t.
			\end{align*}
			By substituting $\cosh^2 t = u$, we get
			\begin{align*}
				R_kf_{\lambda} \left( \xi \right) = \frac{C}{\cosh^{k - 1} d \left( 0, \xi \right)} \int\limits_{\cosh^2 d \left( 0, \xi \right)}^{\cosh^2 \lambda} \left( u - \cosh^2 d \left( 0, \xi \right) \right)^{\frac{k}{2} - 1} \mathrm{d}u = \frac{C}{\cosh^{k - 1} d \left( 0, \xi \right)} \left( \cosh^2 \lambda - \cosh^2 d \left( 0, \xi \right) \right)^{\frac{k}{2}},
			\end{align*}
			which is the same as Equation {\eqref{KPlaneBallCoshHnEquation}}.
		\end{proof}
		\begin{proposition}
			\label{KPlaneBallSinhCoshHn}
			Let $\lambda > 0$ and $B \left( 0, \lambda \right)$ be the ball of radius $\lambda$ centered at $0$ in $\mathbb{H}^n$. For $\alpha \in \mathbb{R}$ and $p \geq 1$, consider the function $f_{\lambda}: \mathbb{H}^n \rightarrow \mathbb{C}$, defined as $f_{\lambda} \left( x \right) = \chi_{B \left( 0, \lambda \right)} \left( x \right) \sinh^{- \frac{\alpha}{p}} d \left( 0, x \right) \cosh d \left( 0, x \right)$. Then, {for $\xi \in \Xi_k \left( \mathbb{H}^n \right)$ with $0 \notin \xi$, we have},
			\begin{equation}
				\label{KPlaneBallSinhCoshHnEquation}
				R_kf_{\lambda} \left( \xi \right) = C \frac{\sinh^{k - \frac{{\alpha}}{p}} \lambda}{\cosh^{k - 1} d \left( 0, \xi \right)} \left( 1 - \frac{\sinh^2 d \left( 0, \xi \right)}{\sinh^2 \lambda} \right)^{\frac{k}{2}} {}_2F_1 \left( \frac{{\alpha}}{2p}, 1; 1 + \frac{k}{2}; 1 - \frac{\sinh^2 d \left( 0, \xi \right)}{\sinh^2 \lambda} \right) \chi_{\left( 0, \lambda \right)} d \left( 0, \xi \right).
			\end{equation}
		\end{proposition}
		\begin{proof}
			Using Equation \eqref{KPlaneTransformRadialHnEquation}, we get for $d \left( 0, \xi \right) < \lambda$,
			\begin{align*}
				R_kf_{\lambda} \left( \xi \right) = \frac{C}{\cosh^{k - 1} d \left( 0, \xi \right)} \int\limits_{d \left( 0, \xi \right)}^{\lambda} \cosh t \sinh^{- \frac{{\alpha}}{p}} t \left( \cosh^2 t - \cosh^2 d \left( 0, \xi \right) \right)^{\frac{k}{2} - 1} \sinh t \ \mathrm{d}t.
			\end{align*}
			By substituting $\cosh^2 t = 1 + \sinh^2 d \left( 0, \xi \right) \left( 1 - z \left( 1 - \frac{\sinh^2 \lambda}{\sinh^2 d \left( 0, \xi \right)} \right) \right)$, and simplifying, we get
			\begin{align*}
				R_kf_{\lambda} \left( \xi \right) = C \frac{\sinh^{k - \frac{{\alpha}}{p}} d \left( 0, \xi \right)}{\cosh^{k - 1} d \left( 0, \xi \right)} \left( \frac{\sinh^2 \lambda}{\sinh^2 d \left( 0, \xi \right)} - 1 \right)^{\frac{k}{2}} \int\limits_{0}^{1} \left( 1 - z \left( 1 - \frac{\sinh^2 \lambda}{\sinh^2 d \left( 0, \xi \right)} \right) \right)^{- \frac{{\alpha}}{2p}} z^{\frac{k}{2} - 1} \ \mathrm{d}z.
			\end{align*}
			Now, using the integral form of the hypergeometric function given in Equation \eqref{IntegralForm2F1}, we get
			\begin{align*}
				R_kf_{\lambda} \left( \xi \right) = C \frac{\sinh^{- \frac{{\alpha}}{p}} d \left( 0, \xi \right)}{\cosh^{k - 1} d \left( 0, \xi \right)} \sinh^{k} \lambda \left( 1 - \frac{\sinh^2 d \left( 0, \xi \right)}{\sinh^2 \lambda} \right)^{\frac{k}{2}} {}_2F_1 \left( \frac{{\alpha}}{2p}, \frac{k}{2}; 1 + \frac{k}{2}; 1 - \frac{\sinh^2 \lambda}{\sinh^2 d \left( 0, \xi \right)} \right).
			\end{align*}
			The result now follows by using Equation \eqref{Transformation2F1}\eqref{2F1Transformation1}.
		\end{proof}
		{\begin{remark}
			\normalfont
			We notice that in Theorem \ref{KPlaneBallSinhCoshHn} we exclude all those totally-geodesic submanifolds of $\mathbb{H}^n$ that pass through the origin. However, the collection of all such submanifolds is a measure zero set in $\Xi_k \left( \mathbb{H}^n \right)$, and hence, Equation \eqref{KPlaneBallSinhCoshHnEquation} makes sense almost everywhere. This is sufficient for us, since we deal with norm estimates, where a set of measure zero causes no change, whatsoever.
		\end{remark}}
		Moving forward, we now consider certain ``nice" functions on $\mathbb{S}^n$, and evaluate their $k$-plane transform. As mentioned earlier, odd functions are in the kernel of the $k$-plane transform. Therefore, in what follows, we consider functions defined on the half sphere $\mathbb{S}^n_+ := \left\lbrace x \in \mathbb{S}^n | \langle x, e_{n + 1} \rangle > 0 \right\rbrace$. These functions can be thought to be extended evenly on $\mathbb{S}^n$.
		\begin{proposition}
			\label{KPlaneBallSnWithoutMass}
			{Let $0 < \lambda < \frac{\pi}{2}$ and $B \left( 0, \lambda \right)$ be the ball of radius $\lambda$ centered at origin in $\mathbb{S}^n$. Consider the function $f_{\lambda}: \mathbb{S}^n_+ \rightarrow \mathbb{C}$, defined as $f_{\lambda} \left( x \right) = \chi_{B \left( 0, \lambda \right)} \left( x \right)$. Then, we have,
			\begin{equation}
				\label{KPlaneBallSnWithoutMassEquation}
				R_kf_{\lambda} \left( \xi \right) = \frac{C}{\cos^{k + 1} d \left( 0, \xi \right)} \sin^k \lambda \left( 1 - \frac{\sin^2 d \left( 0, \xi \right)}{\sin^2 \lambda} \right)^{\frac{k}{2}} {}_2F_1 \left( \frac{1}{2}, \frac{k}{2}; 1 + \frac{k}{2}; 1 - \frac{\cos^2 \lambda}{\cos^2 d \left( 0, \xi \right)} \right) \chi_{\left( 0, \lambda \right)} \left( d \left( 0, \xi \right) \right).
			\end{equation}}
		\end{proposition}
		\begin{proof}
			{Using Equation \eqref{KPlaneTransformRadialSnEquation}, we have for $d \left( 0, \xi \right) < \lambda$,
			$$R_kf_{\lambda} = \frac{C}{\cos^{k - 1} d \left( 0, \xi \right)} \int\limits_{d \left( 0, \xi \right)}^{\lambda} \left( \cos^2 d \left( 0, \xi \right) - \cos^2 t \right)^{\frac{k}{2} - 1} \sin t \ \mathrm{d}t.$$
			Next, we substitute $\cos t = \cos d \left( 0, \xi \right) \sqrt{1 - \left( 1 - \frac{\cos^2 \lambda}{\cos^2 d \left( 0, \xi \right)} \right) v}$ and simply to obtain,
			$$R_kf_{\lambda} \left( \xi \right) = \frac{C}{\cos d \left( 0, \xi \right)} \left( 1 - \frac{\cos^2 \lambda}{\cos^2 d \left( 0, \xi \right)} \right)^{\frac{k}{2}} \int\limits_{0}^{1} v^{\frac{k}{2} - 1} \left( 1 - v \left( 1 - \frac{\cos^2 \lambda}{\cos^2 d \left( 0, \xi \right)} \right) \right)^{- \frac{1}{2}} \mathrm{d}v .$$
			Now, by using the Integral form of the hypergeometric function given in Equation \eqref{IntegralForm2F1}, we have,
			\begin{align*}
				R_kf_{\lambda} \left( \xi \right) &= \frac{C}{\cos d \left( 0, \xi \right)} \left( 1 - \frac{\cos^2 \lambda}{\cos^2 d \left( 0, \xi \right)} \right)^{\frac{k}{2}} {}_2F_1 \left( \frac{1}{2}, \frac{k}{2}; 1 + \frac{k}{2}; 1 - \frac{\cos^2 \lambda}{\cos^2 d \left( 0, \xi \right)} \right) \\
				&= \frac{C}{\cos^{k + 1} d \left( 0, \xi \right)} \sin^k \lambda \left( 1 - \frac{\sin^2 d \left( 0, \xi \right)}{\sin^2 \lambda} \right)^{\frac{k}{2}} {}_2F_1 \left( \frac{1}{2}, \frac{k}{2}; 1 + \frac{k}{2}; 1 - \frac{\cos^2 \lambda}{\cos^2 d \left( 0, \xi \right)} \right).
			\end{align*}
			This completes the proof!}
		\end{proof}
		{We now move a step forward and get expressions for the $k$-plane transform of the characteristic functions of balls centered at the origin of $\mathbb{S}^n$ with certain non-constant densities.}
		\begin{proposition}
			\label{KPlaneBallSn}
			Let $0 < \lambda < \frac{\pi}{2}$ and $B \left( 0, \lambda \right)$ be the ball of radius $\lambda$ centered at the origin in $\mathbb{S}^n$. Consider the function $f_{\lambda}: \mathbb{S}^n_+ \rightarrow \mathbb{C}$, defined as $f_{\lambda} \left( x \right) = \chi_{B \left( 0, \lambda \right)} \left( x \right) \cos d \left( 0, x \right)$. Then, we have,
			\begin{equation}
				\label{KPlaneBallSnEquation}
				R_kf_{\lambda} \left( \xi \right) = \frac{C}{\cos^{k - 1} d \left( 0, \xi \right)} \left( \sin^2 \lambda - \sin^2 d \left( 0, \xi \right) \right)^{\frac{k}{2}} \chi_{\left( 0, \lambda \right)} \left( d \left( 0, \xi \right) \right),
			\end{equation}
			where, the constant $C$ depends only on $n$ and $k$.
		\end{proposition}
		\begin{proof}
			We have using Equation \eqref{KPlaneTransformRadialSnEquation2},
			\begin{align*}
				R_kf_{\lambda} \left( \xi \right) &= {\frac{C}{\cos d \left( 0, \xi \right)}} \int\limits_{d \left( 0, \xi \right)}^{\lambda} \cos t \left[ 1 - \frac{\tan^2 d \left( 0, \xi \right)}{\tan^2 t} \right]^{\frac{k}{2} - 1} \sin^{k - 1} t \ \mathrm{d}t,
			\end{align*}
			provided $d \left( 0, \xi \right) < \lambda$. When $d \left( 0, \xi \right) \geq \lambda$, the {$R_kf_{\lambda} \left( \xi \right) = 0$}. Now, we have, by substituting $\sin t = u$,
			\begin{align*}
				R_kf_{\lambda} \left( \xi \right) &= \frac{C}{\cos d \left( 0, \xi \right)} \int\limits_{\sin d \left( 0, \xi \right)}^{\sin \lambda} u^{k - 1} \left( 1 - \left( \frac{1}{u^2} - 1 \right) \tan^2 d \left( 0, \xi \right) \right)^{\frac{k}{2} - 1} \mathrm{d}u \\
				&= \frac{C}{\cos^{k - 1} d \left( 0, \xi \right)} \int\limits_{\sin d \left( 0, \xi \right)}^{\sin \lambda} \left( u^2 - \sin^2 d \left( 0, \xi \right) \right)^{\frac{k}{2} - 1} u \ \mathrm{d}u \\
				&= \frac{C}{\cos^{k - 1} d \left( 0, \xi \right)} \left( \sin^2 \lambda - \sin^2 d \left( 0, \xi \right) \right)^{\frac{k}{2}}.
			\end{align*}
			This completes the proof!
		\end{proof}
		\begin{proposition}
			\label{KPlaneOriginSn2}
			Let $p \geq 1$, $\alpha \in \mathbb{R}$, $0 < \lambda < \frac{\pi}{2}$ and $B \left( 0, \lambda \right)$ be the ball of radius $\lambda$ centered at the origin in $\mathbb{S}^n$. Consider the function $f_{\lambda}: \mathbb{S}^n_+ \rightarrow \mathbb{C}$, defined as $f_{\lambda} \left( x \right) = \chi_{B \left( 0, \lambda \right)} \left( x \right) \cos d \left( 0, x \right) \sin^{- \frac{\alpha}{p}} d \left( 0, x \right)$. Then, {for $\xi \in \Xi_k \left( \mathbb{S}^n \right)$ with $0 \notin \xi$,} we have,
			\begin{equation}
				\label{KPlaneOriginSn2Equation}
				R_kf_{{\lambda}} \left( \xi \right) = C \frac{\sin^{k - \frac{\alpha}{p}} \lambda}{\cos^{k - 1} d \left( 0, \xi \right)} \left( 1 - \frac{\sin^2 d \left( 0, \xi \right)}{\sin^2 \lambda} \right)^{{\frac{k}{2}}} {}_2F_1 \left( \frac{\alpha}{2p}, 1; 1 + \frac{k}{2}; {\frac{\sin^2 \lambda}{\sin^2 d \left( , \xi \right)}} \left( 1 - \frac{\sin^2 d \left( 0, \xi \right)}{\sin^2 \lambda} \right) \right) \chi_{\left( 0, \lambda \right)} \left( d \left( 0, \xi \right) \right).
			\end{equation}
		\end{proposition}
		\begin{proof}
			From Equation \eqref{KPlaneTransformRadialSnEquation}, we have for $d \left( 0, \xi \right) < \lambda$,
			\begin{align*}
				R_kf \left( \xi \right) &= \frac{C}{\cos^{k - 1} d \left( 0, \xi \right)} \int\limits_{d \left( 0, \xi \right)}^{\lambda} \cos t \sin^{1 - \frac{\alpha}{p}} t \left( \cos^2 d \left( 0, \xi \right) - \cos^2 t \right)^{\frac{k}{2} - 1} \mathrm{d}t.
			\end{align*}
			By substituting $\cos^2 t = \left( 1 - z \left( 1 - \frac{\cos^2 \lambda}{\cos^2 d \left( 0, \xi \right)} \right) \right) \cos^2 d \left( 0, \xi \right)$, and simplifying, we get
			\begin{align*}
				R_kf \left( \xi \right) &= C \cos d \left( 0, \xi \right) \sin^{- \frac{\alpha}{p}} d \left( 0, \xi \right) \left( 1 - \frac{\cos^2 \lambda}{\cos^2 d \left( 0, \xi \right)} \right)^{\frac{k}{2}} \int\limits_{0}^{1} \left[ 1 - {\left( 1 - \frac{\cos^2 \lambda}{\cos^2 d \left( 0, \xi \right)} \right)} \frac{z}{\tan^2 d \left( 0, \xi \right)} \right]^{- \frac{\alpha}{2p}} z^{\frac{k}{2} - 1} \ \mathrm{d}z \\
				&= C \cos d \left( 0, \xi \right) \sin^{- \frac{\alpha}{p}} d \left( 0, \xi \right) \left( 1 - \frac{\cos^2 \lambda}{\cos^2 d \left( 0, \xi \right)} \right)^{\frac{k}{2}} {}_2F_1 \left( \frac{\alpha}{2p}, \frac{k}{2}; 1 + \frac{k}{2}; \cot^2 d \left( 0, \xi \right) {\left( 1 - \frac{\cos^2 \lambda}{\cos^2 d \left( 0, \xi \right)} \right)} \right),
			\end{align*}
			where, we have used Equation \eqref{IntegralForm2F1} in the last step. This is the same as Equation \eqref{KPlaneOriginSn2Equation}.
		\end{proof}
		{\begin{remark}
			\normalfont
			As with the case of Theorem \ref{KPlaneBallSinhCoshHn}, Equation \eqref{KPlaneOriginSn2Equation} makes sense almost everywhere.
		\end{remark}}
		Next, we compute the $k$-plane transform of characteristic functions of the complement (in $\mathbb{S}^n_+$) of a ball. Here as well, we consider these functions with certain densities. We call such functions to be ``\textit{concentrated near the equator}".
		\begin{proposition}
			\label{KPlaneEquatorSnWithoutMass}
			{Let $0 < \lambda < \frac{\pi}{2}$ and $B \left( 0, \lambda \right)$ be a ball of radius $\lambda$ centered at the origin in $\mathbb{S}^n$. Consider the set $E_{\lambda} := \mathbb{S}^n_+ \setminus B \left( 0, \lambda \right)$ and the function $f_{\lambda}: \mathbb{S}^n_+ \rightarrow \mathbb{C}$, defined as $f_{\lambda} \left( x \right) = \chi_{E_{\lambda}} \left( x \right)$. Then,
			\begin{equation}
				\label{KPlaneEquatorSnWithoutMassEquation}
				R_kf_{\lambda} \left( \xi \right) = C \begin{cases}
															1, & d \left( 0, \xi \right) \geq \lambda. \\
															\frac{\cos \lambda}{\cos d \left( 0, \xi \right)} \ {}_2F_1 \left( 1 - \frac{k}{2}, \frac{1}{2}; \frac{3}{2}; \frac{\cos^2 \lambda}{\cos^2 d \left( 0, \xi \right)} \right), & d \left( 0, \xi \right) < \lambda.
														\end{cases}
			\end{equation}}
		\end{proposition}
		\begin{proof}
			{From Equation \eqref{KPlaneTransformRadialSnEquation}, we get
			$$R_kf_{\lambda} \left( \xi \right) = \frac{C}{\cos^{k - 1} d \left( 0, \xi \right)} \int\limits_{\max \left\lbrace d \left( 0, \xi \right), \lambda \right\rbrace}^{\frac{\pi}{2}} \left( \cos^2 d \left( 0, \xi \right) - \cos^2 t \right)^{\frac{k}{2} - 1} \sin t \ \mathrm{d}t.$$
			Now, we substitute $\cos t = \min \left\lbrace \cos d \left( 0, \xi \right), \cos \lambda \right\rbrace \sqrt{y}$ and simplify to obtain,
			$$R_kf_{\lambda} \left( \xi \right) = C \frac{\min \left\lbrace \cos d \left( 0, \xi \right), \cos \lambda \right\rbrace}{\cos d \left( 0, \xi \right)} \int\limits_{0}^{1} y^{- \frac{1}{2}} \left( 1 - y \frac{\min \left\lbrace \cos^2 \lambda, \cos^2 d \left( 0, \xi \right) \right\rbrace}{\cos^2 d \left( 0, \xi \right)} \right)^{\frac{k}{2} - 1} \mathrm{d}y.$$
			From Equation \eqref{IntegralForm2F1}, we now get,
			$$R_kf_{\lambda} \left( \xi \right) = C \frac{\min \left\lbrace \cos d \left( 0, \xi \right), \cos \lambda \right\rbrace}{\cos d \left( 0, \xi \right)} \ {}_2F_1 \left( 1 - \frac{k}{2}, \frac{1}{2}; \frac{3}{2}; \frac{\min \left\lbrace \cos^2 d \left( 0, \xi \right), \cos^2 \lambda \right\rbrace}{\cos^2 d \left( 0, \xi \right)} \right).$$
			This is the same as Equation \eqref{KPlaneEquatorSnWithoutMassEquation}, by a use of Equation \eqref{Behaviour2F1Case1} of Theorem \ref{Behaviour2F1}.}
		\end{proof}
		\begin{proposition}
			\label{KPlaneEquatorSn}
			Let $0 < \lambda < \frac{\pi}{2}$ and $B \left( 0, \lambda \right)$ be the ball of radius $\lambda$ centered at the origin in $\mathbb{S}^n$. Consider the {set $E_{\lambda} := \mathbb{S}^n_+ \setminus B \left( 0, \lambda \right)$, and the} function $f_{\lambda}: \mathbb{S}^n_+ \rightarrow \mathbb{C}$, defined as $f_{\lambda} \left( x \right) = \chi_{{E_{\lambda}}} \left( x \right) \cos d \left( 0, x \right)$. Then, we have,
			\begin{equation}
				\label{KPlaneEquatorSnEquation}
				R_kf_{\lambda} \left( \xi \right) = {C \cos d \left( 0, \xi \right) \begin{cases}
														1, & d \left( 0, \xi \right) > \lambda. \\
														1 - \left( 1 - \frac{\cos^2 \lambda}{\cos^2 d \left( 0, \xi \right)} \right), & d \left( 0, \xi \right) \leq \lambda.
					\end{cases}}
			\end{equation}
		\end{proposition}
		\begin{proof}
			From Equation \eqref{KPlaneTransformRadialSnEquation2}, we have,
			\begin{align*}
				R_kf_{\lambda} \left( \xi \right) &= {\frac{C}{\cos d \left( 0, \xi \right)}} \int\limits_{\max \left\lbrace d \left( 0, \xi \right), \lambda \right\rbrace}^{\frac{\pi}{2}} \cos t \sin^{k - 1} t \left( 1 - \frac{\tan^2 d \left( 0, \xi \right)}{\tan^2 t} \right)^{\frac{k}{2} - 1} \mathrm{d}t.
			\end{align*}
			We now substitute $\sin t = u$ to get
			\begin{align*}
				R_kf_{\lambda} \left( \xi \right) &= \frac{C}{\cos d \left( 0, \xi \right)} \int\limits_{\max \left\lbrace \sin d \left( 0,  \xi \right), \sin \lambda \right\rbrace}^{1} u^{k - 1} \left( 1 - \left( \frac{1}{u^2} - 1 \right) \tan^2 d \left( 0, \xi \right) \right)^{\frac{k}{2} - 1} \mathrm{d}u \\
				&= \frac{C}{\cos^{k - 1} d \left( 0, \xi \right)} \int\limits_{\max \left\lbrace \sin d \left( 0, \xi \right), \sin \lambda \right\rbrace}^{1} \left( u^2 - \sin^2 d \left( 0, \xi \right) \right)^{\frac{k}{2} - 1} u \ \mathrm{d}u \\
				&= \frac{C}{\cos^{k - 1} d \left( 0, \xi \right)} \left[ \left( 1 - \sin^2 d \left( 0, \xi \right) \right)^{\frac{k}{2}} - \left( \left( \max \left\lbrace \sin d \left( 0, \xi \right), \sin \lambda \right\rbrace \right)^2 - \sin^2 d \left( 0, \xi \right) \right)^{\frac{k}{2}} \right].
			\end{align*}
			{This is the same as Equation \eqref{KPlaneEquatorSnEquation}.}
		\end{proof}
		\begin{proposition}
			\label{KPlaneEquatorSn2}
			Let $p \geq 1$, $\alpha < 2p$, $0 < \lambda < \frac{\pi}{2}$ and $B \left( 0, \lambda \right)$ be the ball of radius $\lambda$ centered at the origin in $\mathbb{S}^n$. Consider the {set $E_{\lambda} := \mathbb{S}^n_+ \setminus B \left( 0, \lambda \right)$ and the} function $f_{\lambda}: \mathbb{S}^n_+ \rightarrow \mathbb{C}$, defined as $f_{\lambda} \left( x \right) = \chi_{{E_{\lambda}}} \left( x \right) \cos^{1 - \frac{\alpha}{p}} d \left( 0, x \right)$. Then, we have,
			\begin{equation}
				\label{KPlaneEquatorSn2Equation}
				R_kf \left( \xi \right) = {\frac{C}{\cos d \left( 0, \xi \right)} \begin{cases}
												\cos^{{2} - \frac{\alpha}{p}} d \left( 0, \xi \right), & d \left( 0, \xi \right) > \lambda. \\
												\cos^{2 - \frac{\alpha}{p}} \lambda \ {}_2F_1 \left( 1 - \frac{k}{2}, 1 - \frac{\alpha}{2p}; 2 - \frac{\alpha}{2p}; \frac{\cos^2 \lambda}{\cos^2 d \left( 0, \xi \right)} \right), & d \left( 0, \xi \right) \leq \lambda.			
											\end{cases}}
			\end{equation}
		\end{proposition}
		\begin{proof}
			From Equation \eqref{KPlaneTransformRadialSnEquation}, we have,
			\begin{align*}
				R_kf_{{\lambda}} \left( \xi \right) &= \frac{C}{\cos^{k - 1} d \left( 0, \xi \right)} \int\limits_{\max \left\lbrace d \left( 0, \xi \right), \lambda \right\rbrace}^{\frac{\pi}{2}} \cos^{1 - \frac{\alpha}{p}} t \left( \cos^2 d \left( 0, \xi \right) - \cos^2 t \right)^{\frac{k}{2} - 1} \sin t \ \mathrm{d}t.
			\end{align*}
			By substituting $\cos^2 t = x \min \left\lbrace \cos^2 d \left( 0, \xi \right), \cos^2 \lambda \right\rbrace$ and simplifying, we get
			\begin{align*}
				R_kf_{{\lambda}} \left( \xi \right) &= {\frac{C}{\cos d \left( 0, \xi \right)}} \min \left\lbrace \cos^{2 - \frac{\alpha}{p}} d \left( 0, \xi \right), \cos^{2 - \frac{\alpha}{p}} \lambda \right\rbrace \int\limits_{0}^{1} x^{- \frac{\alpha}{2p}} \left( 1 - x {\frac{\min \left\lbrace \cos^2 d \left( 0, \xi \right), \cos^2 \lambda \right\rbrace}{\cos^2 d \left( 0, \xi \right)}} \right)^{\frac{k}{2} - 1} \mathrm{d}x.
			\end{align*}
			Using Equation \eqref{IntegralForm2F1}, we have,
			\begin{align*}
				R_kf_{{\lambda}} \left( \xi \right) &= {\frac{C}{\cos d \left( 0, \xi \right)}} \min \left\lbrace \cos^{2 - \frac{\alpha}{p}} d \left( 0, \xi \right), \cos^{2 - \frac{\alpha}{p}} \lambda \right\rbrace {}_2F_1 \left( 1 - \frac{k}{2}, 1 - \frac{\alpha}{2p}; 2 - \frac{\alpha}{2p}; {\frac{\min \left\lbrace \cos^2 d \left( 0, \xi \right), \cos^2 \lambda \right\rbrace}{\cos^2 d \left( 0, \xi \right)}} \right) \\
				&= {\frac{C}{\cos d \left( 0, \xi \right)}} \begin{cases}
							\cos^{{2} - \frac{\alpha}{p}} d \left( 0, \xi \right) {}_2F_1 \left( 1 - \frac{k}{2}, 1 - \frac{\alpha}{2p}; 2 - \frac{\alpha}{2p}; \cos^2 d \left( 0, \xi \right) \right), & d \left( 0, \xi \right) > \lambda. \\
							\cos^{2 - \frac{\alpha}{p}} \lambda \ {}_2F_1 \left( 1 - \frac{k}{2}, 1 - \frac{\alpha}{2p}; 2 - \frac{\alpha}{2p}; \cos^2 \lambda \right), & d \left( 0, \xi \right) \leq \lambda.			
						\end{cases}
			\end{align*}
			{In the last equality, we have used Equation \eqref{Behaviour2F1Case1} of Theorem \ref{Behaviour2F1}.}
		\end{proof}
		Moving ahead, we now compute the dual $k$-plane transform of some radial functions on $\Xi_k \left( X \right)$. These results are instrumental in proving the results of Section \ref{ExistenceSection}.
		\begin{proposition}
			\label{DualMixed}
			Let $\varphi: \Xi_k \left( \mathbb{H}^n \right) \rightarrow \mathbb{C}$ be given by $\varphi \left( \xi \right) = \sinh^{\gamma_1} d \left( 0, \xi \right) \cosh^{\gamma_2} d \left( 0, \xi \right)$, with $\gamma_1 > k - n$. Then,
			\begin{equation}
				\label{DualMixedEquation}
				R_k^*\varphi \left( x \right) = \frac{| \mathbb{S}^{k - 1} | | \mathbb{S}^{n - k - 1} |}{2 | \mathbb{S}^{n - 1} |} \sinh^{\gamma_1} d \left( 0, x \right) \cosh^{\gamma_2} d \left( 0, x \right) {}_2F_1 \left( - \frac{\gamma_2}{2}, \frac{k}{2}; \frac{\gamma_1 + n}{2}; \tanh^2 d \left( 0, x \right) \right).
			\end{equation}
		\end{proposition}
		\begin{proof}
			We have from Equation \eqref{DualKPlaneTransformRadialHnEquation},
		\begin{align*}
			R_k^*\varphi \left( x \right) &= \frac{| \mathbb{S}^{k - 1} | | \mathbb{S}^{n - k - 1} |}{| \mathbb{S}^{n - 1} |} \sinh^{2 - n} d \left( 0, x \right) \int\limits_{0}^{\sinh d \left( 0, x \right)} \left( \sinh^2 d \left( 0, x \right) - s^2 \right)^{\frac{k}{2} - 1} s^{\gamma_1 + n - k - 1} \left( 1 + s^2 \right)^{\frac{\gamma_{{2}}}{2}} \mathrm{d}s.
		\end{align*}
		Upon substituting $s = \left( 1 - u \right)^{\frac{1}{2}} \sinh d \left( 0, x \right)$, and simplifying, we get
		\begin{align*}
			R_k^*\varphi \left( x \right) &= \frac{| \mathbb{S}^{k - 1} | | \mathbb{S}^{n - k - 1} |}{2 | \mathbb{S}^{n - 1} |} \sinh^{\gamma_1} d \left( 0, x \right) \cosh^{\gamma_2} d \left( 0, x \right) \int\limits_{0}^{1} u^{\frac{k}{2} - 1} \left( 1 - u \right)^{\frac{\gamma_1 + n - k}{2} - 1} \left( 1 - \tanh^2 d \left( 0, x \right) u \right)^{\frac{\gamma_2}{2}} \mathrm{d}u \\
			&= \frac{| \mathbb{S}^{k - 1} | | \mathbb{S}^{n - k - 1} |}{2 | \mathbb{S}^{n - 1} |} \sinh^{\gamma_1} d \left( 0, x \right) \cosh^{\gamma_2} d \left( 0, x \right) {}_2F_1 \left( - \frac{\gamma_2}{2}, \frac{k}{2}; \frac{\gamma_1 + n}{2}; \tanh^2 d \left( 0, x \right) \right),
		\end{align*}
		where the last equality follows from Equation \eqref{IntegralForm2F1}.
		\end{proof}
		Moving forward, we now give an analogous result for a function defined on $\Xi_k \left( \mathbb{S}^n \right)$.
		\begin{proposition}
			\label{DualMixedSphere}
			Let $\varphi: \Xi_k \left( \mathbb{S}^n \right) \rightarrow \mathbb{C}$ be defined as $\varphi \left( \xi \right) = \sin^{\gamma_1} d \left( 0, \xi \right) \cos^{\gamma_2} d \left( 0, \xi \right)$, for $\gamma_1 > k - n$. Then, we have
			\begin{equation}
				\label{DualMixedSphereEquation}
				R_k^*\varphi \left( x \right) = \frac{| \mathbb{S}^{k - 1} | | \mathbb{S}^{n - k - 1} |}{|\mathbb{S}^{n - 1}|} \sin^{\gamma_1} d \left( 0, x \right) \cos^{\gamma_2} d \left( 0, x \right) {}_2F_1 \left( - \frac{\gamma_2}{2}, \frac{k}{2}; \frac{\gamma_1 + n}{2}; -\tan^2 d \left( 0, x \right) \right).
			\end{equation}
		\end{proposition}
		\begin{proof}
			We have from Equation \eqref{DualKPlaneTransformRadialSnEquation},
			\begin{align*}
				R_k^*\varphi \left( x \right) &= \frac{| \mathbb{S}^{k - 1} | | \mathbb{S}^{n - k -1} |}{|\mathbb{S}^{n - 1}| \sin^{n - 2} d \left( 0, x \right)} \int\limits_{0}^{\sin d \left( 0, x \right)} \left( \sin^2 d \left( 0, x \right) - s^2 \right)^{\frac{k}{2} - 1} s^{n - k - 1} s^{\gamma_1} \left( 1 - s^2 \right)^{\frac{\gamma_2}{2}} \mathrm{d}s.
			\end{align*}
			Substituting $s = u^{\frac{1}{2}} \sin d \left( 0, \xi \right)$ and simplifying, we get,
			\begin{align*}
				R_k^*\varphi \left( x \right) &= \frac{| \mathbb{S}^{k - 1} | | \mathbb{S}^{n - k - 1} |}{|\mathbb{S}^{n - 1}|} \sin^{\gamma_1} d \left( 0, x \right) \int\limits_{0}^{1} \left( 1 - u \right)^{\frac{k}{2} {-1}} u^{\frac{\gamma_1 + n - k}{2} - 1} \left( 1 - \sin^2 d \left( 0, x \right) u \right)^{\frac{\gamma_2}{2}} \mathrm{d}u.
			\end{align*}
			Using the integral representation of the hypergeometric function (Equation \eqref{IntegralForm2F1}), we have,
			\begin{align*}
				R_k^*\varphi \left( x \right) &= \frac{| \mathbb{S}^{k - 1} | | \mathbb{S}^{n - k - 1} |}{|\mathbb{S}^{n - 1}|} \sin^{\gamma_1} d \left( 0, x \right) {}_2F_1 \left( - \frac{\gamma_2}{2}, \frac{\gamma_1 + n - k}{2}; \frac{\gamma_1 + n}{2}; \sin^2 d \left( 0, x \right) \right).
			\end{align*}
			The result now follows from Equation \eqref{Transformation2F1} \eqref{2F1Transformation1}.
		\end{proof}
	\section{Gauss' Hypergeometric Function and its properties}
		\label{HypergoemetricFunctionSubsection}
		In this section, we give the definitions and a few properties of hypergeometric function that are used throughout the article. The material presented here can be found in \cite{Olver}. We denote by $\Gamma \left( \cdot \right)$, the \textit{Gamma function} defined on $\mathbb{C}$, except at the non-positive integers where it has a simple pole. For a description on the Gamma function, we {again} refer the reader to \cite{Olver}. We now give the definition of the Gauss' hypergeometric function.
		\begin{definition}[Gauss' Hypergeometric function]
			\label{HypergeometricFunctionDefinition}
			For $a, b, c \in \mathbb{C}$ and $| z | < 1$, the hypergeometric function is given by the series
			\begin{equation}
				\label{HypergeometricSeries}
				{}_2F_1 \left( a, b; c; z \right) {:}= \sum\limits_{k = 0}^{\infty} \frac{\left( a \right)_k \left( b \right)_k}{\Gamma \left( c + k \right)} \frac{z^k}{k!},
			\end{equation}
			where $\left( x \right)_k$ denotes the Pochhammer's symbol, $\left( x \right)_k {:}= x \left( x + 1 \right) \cdots \left( x + n - 1 \right) = \frac{\Gamma \left( x + n \right)}{\Gamma \left( x \right)}$.
		\end{definition}
		It is well known that the hypergeometric function can be analytically continued to $\mathbb{C} \setminus \left[ 1, \infty \right)$. We use the same symbol ${}_2F_1$ to denote the analytic continuation of Equation \eqref{HypergeometricSeries}. The hypergeometric function can also be written in an integral form.
		\begin{equation}
			\label{IntegralForm2F1}
			{}_2F_1 \left( a, b; c; z \right) = \frac{\Gamma \left( c \right)}{\Gamma \left( b \right) \Gamma \left( c - b \right)} \int\limits_{0}^{1} \frac{t^{b - 1} \left( 1 - t \right)^{c - b - 1}}{\left( 1 - zt \right)^a} \mathrm{d}t,
		\end{equation}
		for $\text{Re } c > \text{Re } b > 0$, and $z \in \mathbb{C} \setminus \left[ 1, \infty \right)$. In the computations involved in this article, we only deal with real parameters and variables. {Moreover, we see that $c > 0$ always holds in our calculations, so that $\Gamma \left( c \right)$ in Equation \eqref{IntegralForm2F1} is always finite.}
		{\begin{remark}
			\label{HypergemetricFunctionBoundedBelow}
			\normalfont
			Before we move ahead and give further properties of the hypergeometric function, let us make the following simple observation. Let $a < 0$, and $c > b > 0$. Then, for any $z < 1$, we have $\left( 1 - zt \right)^{-a} \geq \left( 1 - t \right)^{-a}$, for $0 < t < 1$. Therefore, from Equation \eqref{IntegralForm2F1}, we observe that
			$${}_2F_1 \left( a, b; c; z \right) = \frac{\Gamma \left( c \right)}{\Gamma \left( b \right) \Gamma \left( c - b \right)} \int\limits_{0}^{1} \frac{t^{b - 1} \left( 1 - t \right)^{c - b - 1}}{\left( 1 - zt \right)^{a}} \mathrm{d}t \geq \frac{\Gamma \left( c \right)}{\Gamma \left( b \right) \Gamma \left( c - b \right)} \int\limits_{0}^{1} t^{b - 1} \left( 1 - t \right)^{c - b - a - 1} \mathrm{d}t = \frac{\Gamma \left( c \right) \Gamma \left( c - b - a \right)}{\Gamma \left( c - a \right) \Gamma \left( c - b \right)} > 0.$$
			That is, when $a < 0$ and $c > b > 0$, the hypergeometric function with $z < 1$ is bounded below by a positive constant.
		\end{remark}}
		
		The behaviour of hypergeometric function near $z = 1$ plays an important role in giving necessary conditions for $L^p$-$L^q$ boundedness of the totally-geodesic $k$-plane transform. We state here the three {cases} of {our} interest. 
		\begin{theorem}
			\label{Behaviour2F1}
			{Near the singular point $z = 1$, depending on the value of $c - b - a$, we have the following.}
			\begin{enumerate}
				\item If $\text{Re} \left( c - b - a \right) > 0$, then,
				\begin{equation}
					\label{Behaviour2F1Case1}
					{}_2F_1 \left( a, b; c; 1 \right) = \frac{\Gamma \left( c \right) \Gamma \left( c - b - a \right)}{\Gamma \left( c - a \right) \Gamma \left( c - b \right)}.
				\end{equation}
				\item If $c = a + b$, we have
				\begin{equation}
					\label{Behaviour2F1Case2}
					\lim\limits_{z \rightarrow 1^-} \frac{{}_2F_1 \left( a, b; c; z \right)}{- \ln \left( 1 - z \right)} = \frac{\Gamma \left( a + b \right)}{\Gamma \left( a \right) \Gamma \left( b \right)}.
				\end{equation}
				\item If $\text{Re} \left( c - b - a \right) < 0$, we have
				\begin{equation}
					\label{Behaviour2F1Case3}
					\lim\limits_{z \rightarrow 1^-} \frac{{}_2F_1 \left( a, b; c; z \right)}{\left( 1 - z \right)^{c - b - a}} = \frac{\Gamma \left( c \right) \Gamma \left( a + b - c \right)}{\Gamma \left(  a \right) \Gamma \left( b \right)}.
				\end{equation}
			\end{enumerate}
		\end{theorem}
		A few remarks are in order.
		\begin{remark}~
			\normalfont
			\begin{enumerate}
				\item {Equation \eqref{Behaviour2F1Case1} says that when $\text{Re} \left( c - b - a \right) > 0$, the hypergeometric function behaves like a constant near $z = 1$.}
				\item When $c = a + b$, {we conclude from Equation \eqref{Behaviour2F1Case2} that} the hypergeometric function near $z = 1$ is essentially interchangeable with $-\ln \left( 1 - z \right)$. 
				\item Equation \eqref{Behaviour2F1Case3} gives that when $\text{Re} \left( c - b - a \right) < 0$ and $z$ is close to $1$, we can replace (upto a constant) the hypergeometric function by $\left( 1 - z \right)^{c - b - a}$.
			\end{enumerate}
		\end{remark}
		
		The following transformation formulae for hypergeometric functions are also useful, {and} can be found in \cite{Olver}.
		\begin{equation}
			\label{Transformation2F1}
			{}_2F_1 \left( a, b; c; z \right) {\labelrel{=}{2F1Transformation1}} \left( 1 - z \right)^{-a} {}_2F_1 \left( a, c - b; c; \frac{z}{z - 1} \right) {\labelrel{=}{2F1Transformation2}} \left( 1 - z \right)^{-b} {}_2F_1 \left( c - a, b; c; \frac{z}{z - 1} \right) {\labelrel{=}{2F1Transformation3}} \left( 1 - z \right)^{c - b - a} {}_2F_1 \left( c - a, c - b; c; z \right).
		\end{equation}
		The transformation {formulae} in Equation \eqref{Transformation2F1} holds for all $z \in \mathbb{C} \setminus \left[ 1, \infty \right)$.
		
		We also have for $z \notin \left[ 0, \infty \right)$,
		\begin{equation}
			\label{Transformation2F12}
			\frac{\sin \left( \pi \left( b - a \right) \right)}{\pi} {}_2F_1 \left( a, b; c; z \right) = \frac{\left( -z \right)^{-a}}{\Gamma \left( b \right) \Gamma \left( c - a \right)} {}_2F_1 \left( a, a - c + 1; a - b + 1; \frac{1}{z} \right) - \frac{\left( -z \right)^{-b}}{\Gamma \left( a \right) \Gamma \left( c - b \right)} {}_2F_1 \left( b, b - c  + 1; b - a + 1; \frac{1}{z} \right).
		\end{equation}
	\section{Weighted Boundedness of Riemann-Liouville Fractional Integrals}
		\label{FIBoundednessSection}
		The connection between Radon transforms and fractional integrals has been extensively studied by {Boris} Rubin. To state a few works, we refer the reader to \cite{RubinRFI} and \cite{Rubin} and the references therein. It is observed that due to this connection, mapping properties of Radon transforms (particularly, the $k$-plane transforms) {are} well-understood {with the help of} the corresponding properties of fractional integrals. In what follows, we consider one-dimensional fractional integrals, often called the Riemann-Liouville fractional integrals. These have a close connection with the $k$-plane transform of radial functions, and help in proving {certain} results of Section \ref{RadialFunctionSection}. For a complete review {on} the {topic of one-dimensional fractional integrals}, we refer the reader to \cite{RubinFIP}. Here, we only provide the notation and the results we want in a manner convenient to us. Let $\bar{\mathbb{R}}$ denote the extended real line. For $\varphi: \Omega \subseteq \bar{\mathbb{R}} \rightarrow \mathbb{C}$, a ``nice" function, we {consider} the following Riemann-Liouville fractional integrals of $\varphi$ for $\alpha \in \mathbb{C}$ with $\text{Re } \alpha > 0$.
		\begin{equation}
			\label{FiniteLowerRLFractionalIntegral}
			I_{a+}^{\alpha} \varphi \left( x \right) = \frac{1}{\Gamma \left( \alpha \right)} \int\limits_{a}^{x} \frac{\varphi \left( y \right)}{\left( x - y \right)^{1 - \alpha}} \mathrm{d}y, \text{ for } x > a.
		\end{equation}
		\begin{equation}
			\label{InfiniteUpperRLFractionalIntegral}
			I_-^{\alpha} \varphi \left( x \right) = \frac{1}{\Gamma \left( \alpha \right)} \int\limits_{x}^{\infty} \frac{\varphi \left( y \right)}{\left( x - y \right)^{1 - \alpha}} \mathrm{d}y.
		\end{equation}
		Depending on the set $\Omega$, we introduce the following partitions.
		\begin{equation}
			\label{PartitionOmega}
			\begin{aligned}
				0= a_1 < a_2 < \cdots < a_l = a&, && \text{if } \Omega = \left[ 0, a \right] \text{ for } a < + \infty. \\
				0 = a_1 < a_2 < \cdots < a_l < + \infty&, && \text{if } \Omega = \left[ 0, \infty \right] =: \overline{\mathbb{R}_+}.
			\end{aligned}
		\end{equation}
		{For $\gamma_1, \cdots, \gamma_l, \gamma_{\infty} \in \mathbb{R}$, we consider} the weight function $\rho: \Omega \rightarrow \mathbb{C}$ as
		\begin{equation}
			\label{WeightDomainFI}
			\rho \left( x \right) {:}= \begin{cases}
										\left( 1 + \left| x \right| \right)^{\gamma_{\infty}} \prod\limits_{j = 1}^{l} \left| x - a_j \right|^{\gamma_j}, & \text{when } \left| \Omega \right| = + \infty. \\
										\prod\limits_{j = 1}^{l} \left| x - a_j \right|^{\gamma_j}, & \text{when } \left| \Omega \right| < + \infty.
									\end{cases}
		\end{equation}
		{For $1 < p < +\infty$, we consider $m, q \in \mathbb{R}$ that satisfy the following relations:}
		\begin{equation}
			\label{FractionalIntegralMQ}
			{0 \leq m \leq \alpha, 0 < \alpha < m + \frac{1}{p}, \text{ and } \alpha - m = \frac{1}{p} - \frac{1}{q}.}
		\end{equation}
		Let $\epsilon_0, \cdots, \epsilon_l > 0$ be chosen arbitrarily. We define for $j = 1, \cdots, l$,
		\begin{equation}
			\label{FractionalIntegralDeltaJ}
			\delta_j {:}= \begin{cases}
							\gamma_j, & \text{when } \gamma_j > \alpha - \frac{1}{p}. \\
							\alpha - \frac{1}{p} - m - \epsilon_j, & \text{when } \gamma_j \leq \alpha - \frac{1}{p}.
						\end{cases}
		\end{equation}
		{Moreover, we} define
		$${\delta_{\infty}} {:}= \gamma - m - \sum\limits_{j = 1}^{l} \delta_j.$$
		{Now}, let us define the following weight function.
		\begin{equation}
			\label{WeightCodomainMinus}
			\rho_- \left( x \right) {:}= \left( 1 + x \right)^{\delta_{\infty}} \prod\limits_{j = 1}^{l} \left| x - a_j \right|^{\delta_j}.
		\end{equation}
		{With these notations made, we have the following theorem for when $\Omega = \overline{\mathbb{R}_+}$.}
		\begin{theorem}[\cite{RubinFIP}]
			\label{LPLQBoundednessFI}
			Let $\varphi: \Omega \rightarrow \mathbb{C}$ be such that $\| \rho \varphi \|_{p} < + \infty$. If ${\gamma_{\infty} + \sum\limits_{i = 1}^{l} \gamma_i} > \alpha - \frac{1}{p}$, then, we have,
			\begin{equation}
				\label{BoundednessFIEquation2}
				\| \rho_{-} I_{-}^{\alpha} \varphi \|_q \leq C \| \rho \varphi \|_p.
			\end{equation}
		\end{theorem}
		We also have analogous theorem for the case when $\Omega = \left[ 0, a \right]$, when $a < + \infty$. {We require this in the study of the boundedness properties of $X$-ray transform on the sphere.}
		\begin{theorem}
			\label{LPLQBoundednessFIFiniteOmega}
			Let $\varphi: \Omega = \left[ 0, a \right] \rightarrow \mathbb{C}$ be such that $\| \rho \varphi \|_{p} < + \infty$. Then, for $m$ and $\delta_j$ mentioned {in Equations \eqref{FractionalIntegralMQ} and \eqref{FractionalIntegralDeltaJ}, respectively}, we have
			\begin{equation}
				\label{BoundednessFIFiniteOmegaEquation}
				\| \rho_+ I_{0+}^{\alpha} \varphi \|_q \leq C \| \rho \varphi \|_p.
			\end{equation}
			Here,
			$$\rho_+ \left( x \right) = x^{\gamma_1 - m} \prod\limits_{j = 2}^{l} \left| x - a_j \right|^{\delta_j}.$$
		\end{theorem}
		{In our study, the case $\alpha = \frac{1}{2}$ and $m = 0$ is of particular interest, owing to its connection with the $X$-ray transform of radial functions. We observe that in Theorems \ref{LPLQBoundednessFI} and \ref{LPLQBoundednessFIFiniteOmega}, we have not considered the case when $p = 1$. First, let us observe that when $p = 1$, $\alpha = \frac{1}{2}$ and $m = 0$, we must have, from Equation \eqref{FractionalIntegralMQ}, $q = 2$. In fact, we now see, through a counter-example that such a result concerning weighted $L^1$-$L^2$ boundedness for fractional integrals cannot be expected. The example that we now present comes from a similar example of Hardy and Littlewood (see \cite{HardyLittlewoodFI}) for the unweighted case.}
		
		{We first show that there are functions $\psi$ in $L^1 \left( \left( 0, \infty \right), \rho_1 \right)$ for which $I^{\frac{1}{2}}_+\psi \notin L^2 \left( \left( 0, \infty \right), \rho_2 \right)$, for any choices of ``nice" weights $\rho_1$ and $\rho_2$.}
		\begin{example}
			\label{L1L2NotPossibleIHalfPlus}
			\normalfont
			{Let us consider the case $\Omega = \left[ 0, \infty \right]$, and two partitions
			\begin{equation}
				\label{TwoPartitions1}
				\begin{aligned}
					0 = a_1^{\left( 1 \right)} < a_2^{\left( 1 \right)} < \cdots < a_{l_1}^{\left( 1 \right)} < + \infty, \\
					0 = a_1^{\left( 2 \right)} < a_2^{\left( 2 \right)} < \cdots < a_{l_2}^{\left( 2 \right)} < + \infty.
				\end{aligned}
			\end{equation}
			Corresponding to the partitions mentioned in Equation \eqref{TwoPartitions1} and $\gamma_1^{\left( 1 \right)}, \cdots, \gamma_{l_1}^{\left( 1 \right)}, \gamma_1^{\left( 2 \right)}, \cdots, \gamma_{l_2}^{\left( 2 \right)}, \gamma_{\infty}^{\left( 1 \right)}, \gamma_{\infty}^{\left( 2 \right)} \in \mathbb{R}$, let us consider the weight functions given by,
			\begin{equation}
				\label{TwoWeights1}
				\begin{aligned}
					\rho_1 \left( x \right) &= \left( 1 + x \right)^{\gamma_{\infty}^{\left( 1 \right)}} \prod\limits_{i = 1}^{l_1} \left| x - a_i^{\left( 1 \right)} \right|^{\gamma_i^{\left( 1 \right)}}. \\
					\rho_2 \left( x \right) &= \left( 1 + x \right)^{\gamma_{\infty}^{\left( 2 \right)}} \prod\limits_{i = 1}^{l_2} \left| x - a_i^{\left( 2 \right)} \right|^{\gamma_i^{\left( 2 \right)}}.
				\end{aligned}
			\end{equation}
			Let us consider
			\begin{equation}
				\label{InteriorPoint1}
				a := \begin{cases}
							1, & \text{if } l_1 = 1 \text{ and } l_2 = 1. \\
							a_2^{\left( 1 \right)}, & \text{if } l_1 \geq 2 \text{ and } l_2 = 1. \\
							a_2^{\left( 2 \right)}, & \text{if } l_1 = 1 \text{ and } l_2 \geq 1. \\
							\min \left\lbrace a_2^{\left( 1 \right)}, a_2^{\left( 2 \right)} \right\rbrace, & \text{if } l_1, l_2 \geq 2.
						\end{cases}
			\end{equation}
			Now, we define $\psi: L^1 \left( \left( 0, \infty \right), \rho_1 \right) \rightarrow \mathbb{C}$ as
			\begin{equation}
				\label{CE1FunctionPsi}
				\psi \left( x \right) = \frac{1}{x - \frac{a}{2}} \left( \ln \frac{1}{x - \frac{a}{2}} \right)^{- \gamma} \frac{1}{\rho_1 \left( x \right)} \chi_{\left( \frac{a}{2}, \frac{3a}{4} \right)} \left( x \right),
			\end{equation}
			where, $1 < \gamma \leq \frac{3}{2}$. We now have,
			\begin{align*}
				\| \psi \|_{L^1 \left( \left( 0, \infty \right), \rho_1 \right)} = \int\limits_{\frac{a}{2}}^{\frac{3a}{4}} \frac{1}{x - \frac{a}{2}} \left( \ln \frac{1}{x - \frac{a}{2}} \right)^{- \gamma} \mathrm{d}x.
			\end{align*}
			By substituting $x - \frac{a}{2} = e^{-u}$ and simplifying, we see that
			$$\| \psi \|_{L^1 \left( \left( 0, \infty \right), \rho_1 \right)} = \int\limits_{\ln \frac{4}{a}}^{\infty} u^{- \gamma} \mathrm{d}u < + \infty,$$
			since $\gamma > 1$. On the other hand, we also have,
			\begin{align*}
				\| I^{\frac{1}{2}}_+ \psi \|_{L^2 \left( \left( 0, \infty \right), \rho_2 \right)} &= \left[ \int\limits_{0}^{\infty} \left| \int\limits_{0}^{x} \frac{1}{t - \frac{a}{2}} \left( \ln \frac{1}{t - \frac{a}{2}} \right)^{- \gamma} \frac{1}{\rho_1 \left( t \right)} \chi_{\left( \frac{a}{2}, \frac{3a}{4} \right)} \left( t \right) \frac{1}{\left( x- t \right)^{\frac{1}{2}}} \mathrm{d}t \right|^2 \rho_2 \left( x \right) \mathrm{d}x \right]^{\frac{1}{2}} \\
				&\geq \left[ \int\limits_{\frac{a}{2}}^{\frac{3a}{4}} \left| \int\limits_{\frac{a}{2}}^{x} \frac{1}{t - \frac{a}{2}} \left( \ln \frac{1}{t - \frac{a}{2}} \right)^{- \gamma} \frac{1}{\rho_1 \left( t \right)} \frac{1}{\left( x - t \right)^{\frac{1}{2}}} \mathrm{d}t \right|^2 \rho_2 \left( x \right) \mathrm{d}x \right]^{\frac{1}{2}}
			\end{align*}
			From the definitions of $\rho_1$ and $\rho_2$ it is clear that for $x \in \left( \frac{a}{2}, \frac{3a}{4} \right)$, both $\rho_1 \left( x \right)$ and $\rho_2 \left( x \right)$ is bounded above and below. Also, for $t \geq \frac{a}{2}$, we have $\left( x - t \right)^{- \frac{1}{2}} \geq \left( x - \frac{a}{2} \right)^{\frac{1}{2}}$. Hence, we get
			\begin{align*}
				\| I^{\frac{1}{2}}_+\psi \|_{L^2 \left( \left( 0, \infty \right), \rho_2 \right)} &\geq C \left[ \int\limits_{\frac{a}{2}}^{\frac{3a}{4}} \frac{1}{x - \frac{a}{2}} \left| \int\limits_{\frac{a}{2}}^{x} \frac{1}{t - \frac{a}{2}} \left( \ln \frac{1}{t - \frac{a}{2}} \right)^{- \gamma} \mathrm{d}t \right|^2 \mathrm{d}x \right]^{\frac{1}{2}}.
			\end{align*}
			By substituting $t - \frac{a}{2} = e^{-u}$ in the inner integral and using the fact that $\gamma > 1$, we obtain
			\begin{align*}
				\| I^{\frac{1}{2}}_+\psi \|_{L^2 \left( \left( 0, \infty \right), \rho_2 \right)} &\geq C \left[ \int\limits_{\frac{a}{2}}^{\frac{3a}{4}} \frac{1}{x - \frac{a}{2}} \left| \int\limits_{\ln \frac{1}{x - \frac{a}{2}}}^{\infty} u^{- \gamma} \mathrm{d}u \right|^2 \mathrm{d}x \right]^{\frac{1}{2}} \\
				&= C \left[ \int\limits_{\frac{a}{2}}^{\frac{3a}{4}} \frac{1}{x - \frac{a}{2}} \left( \ln \frac{1}{x - \frac{a}{2}} \right)^{2 \left( 1 - \gamma \right)} \mathrm{d}x \right]^{\frac{1}{2}}.
			\end{align*}
			Again, by substituting $x - \frac{a}{2} = e^{-u}$ and simplifying, we get
			\begin{align*}
				\| I^{\frac{1}{2}}_+\psi \|_{L^2 \left( \left( 0, \infty \right), \rho_2 \right)} &\geq C \left[ \int\limits_{\ln \frac{4}{a}}^{\infty} u^{2 \left( 1 - \gamma \right)} \mathrm{d}u \right]^{\frac{1}{2}} = + \infty,
			\end{align*}
			since $\gamma \leq \frac{3}{2}$.}
		\end{example}
		{Another interesting case that is not handled in Theorems \ref{LPLQBoundednessFI} and \ref{LPLQBoundednessFIFiniteOmega} is when $\alpha = \frac{1}{2}$, $p = 2$ and $m = 0$. Again, if one were to impose the conditions of Equation \eqref{FractionalIntegralMQ}, then one would require $q = \infty$. We now show that we cannot expect a weighted $L^2$-$L^{\infty}$ estimate for $I^{\frac{1}{2}}_+$.}
		\begin{proposition}
			\label{L2LInfinityIHalfPlusNotPossible}
			{Let $\Omega = \left[ 0, \infty \right]$. We consider the partitions mentioned in Equation \eqref{TwoPartitions1} and weight functions $\rho_1$ and $\rho_2$ mentioned in Equation \eqref{TwoWeights1}. Then, {there is no $C > 0$ such that} for all $\varphi \in L^2 \left( \Omega, \rho_1 \right)$,
			\begin{equation}
				\label{L2LInfinityInequalityFINP}
				\| \rho_2 I^{\frac{1}{2}}_- \varphi \|_{L^{\infty} \left( \Omega \right)} \leq C \| \varphi \|_{L^2 \left( \Omega, \rho_1 \right)}.
			\end{equation}}
		\end{proposition}
		\begin{proof}
			{Suppose that Inequality \eqref{L2LInfinityInequalityFINP} holds for all $\varphi \in L^2 \left( \Omega, \rho_1 \right)$. Then, for ``nice" $\psi: \Omega \rightarrow \mathbb{C}$, we have,
			\begin{align*}
				\int\limits_{0}^{\infty} I^{\frac{1}{2}}_- \varphi \left( x \right) \psi \left( x \right) \rho_1 \left( x \right) \mathrm{d}x &= \int\limits_{0}^{\infty} \int\limits_{x}^{\infty} \frac{\varphi \left( t \right)}{\left( t - x \right)^{\frac{1}{2}}} \psi \left( x \right) \rho_1 \left( x \right) \mathrm{d}t \ \mathrm{d}x \\
				&= \int\limits_{0}^{\infty} \varphi \left( t \right) \int\limits_{0}^{t} \frac{\psi \left( x \right) \rho_1 \left( x \right)}{\left( t - x \right)^{\frac{1}{2}}} \mathrm{d}x \ \mathrm{d}t \\
				&= \int\limits_{0}^{\infty} \varphi \left( t \right) \frac{1}{\rho_1 \left( t \right)} I^{\frac{1}{2}}_+ \left( \rho_1 \psi \right) \left( t \right) \rho_1 \left( t \right) \mathrm{d}t.
			\end{align*}
			Therefore, from Inequality \eqref{L2LInfinityInequalityFINP}, we have, for all $\varphi \in L^2 \left( \Omega, \rho_1 \right)$,
			\begin{align*}
				\left| \int\limits_{0}^{\infty} \frac{1}{\rho_1 \left( t \right)} I^{\frac{1}{2}}_+ \left( \rho_1 \psi \right) \left( t \right) \varphi \left( t \right) \rho_1 \left( t \right) \mathrm{d}t \right| &\leq \int\limits_{0}^{\infty} \left| \rho_2 \left( x \right) I^{\frac{1}{2}}_- \varphi \left( x \right) \right| \left| \psi \left( x \right) \right| \frac{\rho_1 \left( x \right)}{\rho_2 \left( x \right)} \mathrm{d}x \\
				&\leq \| \rho_2 I^{\frac{1}{2}}_- \varphi \|_{L^{\infty} \left( \Omega \right)} \| \psi \|_{L^1 \left( \Omega, \frac{\rho_1}{\rho_2} \right)} \\
				&\leq C \| \varphi \|_{L^2 \left( \Omega, \rho_1 \right)} \| \| \psi \|_{L^1 \left( \Omega, \frac{\rho_1}{\rho_2} \right)}.
			\end{align*}
			Therefore, whenever $\psi \in L^1 \left( \Omega, \frac{\rho_1}{\rho_2} \right)$, we have that $\frac{1}{\rho_1} I^{\frac{1}{2}}_+ \left( \rho_1 \psi \right) \in L^2 \left( \Omega, \rho_1 \right)$, and
			$$\Bigg\| I^{\frac{1}{2}}_+ \left( \rho_1 \psi \right) \Bigg\|_{L^2 \left( \Omega, \frac{1}{\rho_1} \right)} = \Bigg\| \frac{1}{\rho_1} I^{\frac{1}{2}}_+ \left( \rho_1 \psi \right) \Bigg\|_{L^2 \left( \Omega, \rho_1 \right)} \leq C \| \psi \|_{L^1 \left( \Omega, \frac{\rho_1}{\rho_2} \right)} = C \| \rho_1 \psi \|_{L^1 \left( \Omega, \frac{1}{\rho_2} \right)}.$$
			The above inequality gives that whenever $\psi: \Omega \rightarrow \mathbb{C}$ is such that $\rho_1 \psi \in L^1 \left( \Omega, \frac{1}{\rho_2} \right)$, we have that $I^{\frac{1}{2}}_+ \left( \rho_1 \psi \right) \in L^2 \left( \Omega, \frac{1}{\rho_1} \right)$. However, this is not possible, as demonstrated by Example \ref{L1L2NotPossibleIHalfPlus}. Therefore, Inequality \eqref{L2LInfinityInequalityFINP} is not true for every $\varphi \in L^2 \left( \Omega, \rho_1 \right)$.}
		\end{proof}
		{We now give similar results for $I^{\frac{1}{2}}_{1-}$ and $I^{\frac{1}{2}}_+$, by considering $\Omega = \left[ 0, 1 \right]$. In what follows, we fix the following partitions of Omega.
		\begin{equation}
			\label{TwoPartitions2}
			\begin{aligned}
				0 = a_1^{\left( 1 \right)} < a_2^{\left( 1 \right)} < \cdots < a_{l_1}^{\left( 1 \right)} = 1. \\
				0 = a_1^{\left( 2 \right)} < a_2^{\left( 2 \right)} < \cdots < a_{l_2}^{\left( 2 \right)} = 1.
			\end{aligned}
		\end{equation}
		Corresponding to the above partition and $\gamma_1^{\left( 1 \right)}, \cdots, \gamma_{l_1}^{\left( 1 \right)}, \gamma_1^{\left( 2 \right)}, \cdots, \gamma_{l_2}^{\left( 2 \right)} \in \mathbb{R}$, we consider the following weight functions.
		\begin{equation}
			\label{TwoWeights2}
			\begin{aligned}
				\rho_1 \left( x \right) &= \prod\limits_{i = 1}^{l_1} \left| x - a_i^{\left( 1 \right)} \right|^{\gamma_i^{\left( 1 \right)}}, \\
				\rho_2 \left( x \right) &= \prod\limits_{i = 1}^{l_2} \left| x - a_i^{\left( 2 \right)} \right|^{\gamma_i^{\left( 2 \right)}}.
			\end{aligned}
		\end{equation}
		With these notations in place, we first show, through a counter-example that the operator $I^{\frac{1}{2}}_{1-}$ does not take functions in $L^1 \left( \Omega, \rho_1 \right)$ to functions in $L^2 \left( \Omega, \rho_2 \right)$.}
		\begin{example}
			\label{L1L2NotPossibleIHalfMinus}
			\normalfont
			{As in Example \ref{L1L2NotPossibleIHalfPlus}, let us consider
			\begin{equation}
				\label{InteriorPointCE2FI}
				a := \min \left\lbrace a_2^{\left( 1 \right)}, a_2^{\left( 2 \right)} \right\rbrace.
			\end{equation}
			Indeed, by the definition of our partition in Equation \eqref{TwoPartitions2}, it is clear that $l_1, l_2 \geq 2$ and we do not have as many cases as in Equation \eqref{InteriorPoint1}. Now, we define $\psi: \Omega \rightarrow \mathbb{C}$ as
			\begin{equation}
				\label{L1L2NotPoissiblePsi}
				\psi \left( x \right) = \frac{1}{\frac{a}{2} - x} \left( \ln \frac{1}{\frac{a}{2} - x} \right)^{- \gamma} \frac{1}{\rho_1 \left( x \right)} \chi_{\left( \frac{a}{4}, \frac{a}{2} \right)} \left( x \right).
			\end{equation}
			Then, following a procedure almost verbatim to that of Example \ref{L1L2NotPossibleIHalfPlus}, it is clear at once that $\psi \in L^1 \left( \Omega, \rho_1 \right)$. On the other hand, we have,
			\begin{align*}
				\| I^{\frac{1}{2}}_{1-}\psi \|_{L^2 \left( \Omega, \rho_2 \right)} &= \left[ \int\limits_{0}^{1} \left| \int\limits_{x}^{1} \frac{1}{\frac{a}{2} - t} \left( \ln \frac{1}{\frac{a}{2} - t} \right)^{- \gamma} \frac{1}{\rho_1 \left( t \right)} \chi_{\left( \frac{a}{4}, \frac{a}{2} \right)} \left( t \right) \frac{1}{\left( t - x \right)^{\frac{1}{2}}} \mathrm{d}t \right|^2 \rho_2 \left( x \right) \mathrm{d}x \right]^{\frac{1}{2}} \\
				&\geq \left[ \int\limits_{\frac{a}{4}}^{\frac{a}{2}} \left| \int\limits_{x}^{\frac{a}{2}} \frac{1}{\frac{a}{2} - t} \left( \ln \frac{1}{\frac{a}{2} - t} \right)^{- \gamma} \frac{1}{\rho_1 \left( t \right)} \frac{1}{\left( t - x \right)^{\frac{1}{2}}} \mathrm{d}t \right|^2 \rho_2 \left( x \right) \mathrm{d}x \right]^{\frac{1}{2}}.
			\end{align*}
			Arguing as in Example \ref{L1L2NotPossibleIHalfPlus}, we observe that for $x \in \left( \frac{a}{4}, \frac{a}{2} \right)$, both $\rho_1 \left( x \right)$ and $\rho_2 \left( x \right)$ are bounded (above and below). Also, for $t \leq \frac{a}{2}$, we have $\left( t - x \right)^{\frac{1}{2}} \geq \left( \frac{a}{2} - x \right)^{\frac{1}{2}}$. Hence, we have, by following the ideas used in Example \ref{L1L2NotPossibleIHalfPlus},
			\begin{align*}
				\| I^{\frac{1}{2}}_{1-}\psi \|_{L^2 \left( \Omega, \rho_2 \right)} &\geq C \left[ \int\limits_{\frac{a}{4}}^{\frac{a}{2}} \frac{1}{\frac{a}{2} - x} \left| \int\limits_{x}^{\frac{a}{2}} \frac{1}{\frac{a}{2} - t} \left( \ln \frac{1}{\frac{a}{2} - t} \right)^{- \gamma} \mathrm{d}t \right|^{2} \mathrm{d}x \right]^{\frac{1}{2}} = + \infty.
			\end{align*}}
		\end{example}
		{Using Example \ref{L1L2NotPossibleIHalfMinus}, we now give an analogue of Proposition \ref{L2LInfinityIHalfPlusNotPossible}.}
		\begin{proposition}
			\label{L2LInfinityIHalfPlusNP}
			{Let us consider $\Omega = \left[ 0, 1 \right]$, the partitions mentioned in Equation \eqref{TwoPartitions2} and the weight functions mentioned in Equation \eqref{TwoWeights2}. Then, the following inequality is not valid for all functions $\varphi \in L^2 \left( \Omega, \rho_1 \right)$.
			\begin{equation}
				\label{L2LInfinityInequalityIHalfPlusNP}
				\| \rho_2 I^{\frac{1}{2}}_+\varphi \|_{L^{\infty} \left( \Omega \right)} \leq C \| \varphi \|_{L^2 \left( \Omega , \rho_1 \right)}.
			\end{equation}}
		\end{proposition}
		\begin{proof}
			{Suppose that Inequality \eqref{L2LInfinityInequalityIHalfPlusNP} holds for every $\varphi \in L^2 \left( \Omega, \rho_1 \right)$. Then, for ``nice" functions $\psi: \Omega \rightarrow \mathbb{C}$, we have,
			\begin{align*}
				\int\limits_{0}^{1} I^{\frac{1}{2}}_+ \varphi \left( x \right) \psi \left( x \right) \rho_1 \left( x \right) \mathrm{d}x &= \int\limits_{0}^{1} \int\limits_{0}^{x} \frac{\varphi \left( t \right)}{\left( x - t \right)^{\frac{1}{2}}} \psi \left( x \right) \rho_1 \left( x \right) \mathrm{d}t \ \mathrm{d}x \\
				&= \int\limits_{0}^{1} \varphi \left( t \right) \int\limits_{t}^{1} \frac{\psi \left( x \right) \rho_1 \left( x \right)}{\left( x - t \right)^{\frac{1}{2}}} \mathrm{d}x \ \mathrm{d}t \\
				&= \int\limits_{0}^{1} \varphi \left( t \right) \frac{1}{\rho_1 \left( t \right)} I^{\frac{1}{2}}_{1-} \left( \rho_1 \psi \right) \left( t \right) \rho_1 \left( t \right) \mathrm{d}t.
			\end{align*}
			Hence, from Inequality \eqref{L2LInfinityInequalityIHalfPlusNP}, we have for all $\varphi \in L^2 \left( \Omega, \rho_1 \right)$,
			\begin{align*}
				\left| \int\limits_{0}^{1} \frac{1}{\rho_1 \left( t \right)} I^{\frac{1}{2}}_{1-} \left( \rho_1 \psi \right) \left( t \right) \varphi \left( t \right) \rho_1 \left( t \right) \mathrm{d}t \right| &\leq \int\limits_{0}^{1} \left| \rho_2 \left( x \right) I^{\frac{1}{2}}_+ \varphi \left( x \right) \right| \left| \psi \left( x \right) \right| \frac{\rho_1 \left( x \right)}{\rho_2 \left( x \right)} \mathrm{d}x \\
				&\leq \| \rho_2 I^{\frac{1}{2}}_+ \varphi \|_{L^{\infty} \left( \Omega \right)} \| \psi \|_{L^1 \left( \Omega, \frac{\rho_1}{\rho_2} \right)} \\
				&\leq C \| \varphi \|_{L^2 \left( \Omega, \rho_1 \right)} \| \psi \|_{L^1 \left( \Omega, \frac{\rho_1}{\rho_2} \right)}.
			\end{align*}
			Therefore, by duality, we conclude that whenever $\psi \in L^1 \left( \Omega , \frac{\rho_1}{\rho_2} \right)$, we have $\frac{1}{\rho_1} I^{\frac{1}{2}}_{1-} \left( \rho_1 \psi \right) \in L^2 \left( \Omega, \rho_1 \right)$, and
			$$\| I^{\frac{1}{2}}_+ \left( \rho_1 \psi \right) \|_{L^2 \left( \Omega, \frac{1}{\rho_1} \right)} = \Bigg\| \frac{1}{\rho_1} I^{\frac{1}{2}}_{1-} \left( \rho_1 \psi \right) \Bigg\|_{L^2 \left( \Omega, \rho_1 \right)} \leq C \| \psi \|_{L^1 \left( \Omega, \frac{\rho_1}{\rho_2} \right)} = C \| \rho_1 \psi \|_{L^1 \left( \Omega, \frac{1}{\rho_2} \right)}.$$
			However, the above inequality is not true, as seen in Example \ref{L1L2NotPossibleIHalfMinus}. Therefore, Inequality \eqref{L2LInfinityInequalityIHalfPlusNP} is not valid for all $\varphi \in L^2 \left( \Omega, \rho_1 \right)$.}
		\end{proof}

	
\end{document}